
\input amstex

\hfuzz 60pt

\expandafter\ifx\csname beta.def\endcsname\relax \else\endinput\fi
\expandafter\edef\csname beta.def\endcsname{%
 \catcode`\noexpand\@=\the\catcode`\@\space}

\let\atbefore @

\catcode`\@=11

\overfullrule\z@
\hsize 6.25truein
\vsize 9.63truein

\let\@ft@\expandafter \let\@tb@f@\atbefore

\newif\ifMag
\ifnum\mag>1000 \Magtrue\fi
=\ifMag cmr8\else cmr9\fi

\newdimen\p@@ \p@@\p@
\def\m@ths@r{\ifnum\mathsurround=\z@\z@\else\maths@r\fi}
\def\maths@r{1.6\p@@} \def\mathsurzero{\def\maths@r{\z@}}

\mathsurround\maths@r
\font\Brm=cmr12 \font\Bbf=cmbx12 \font\Bit=cmti12 \font\ssf=cmss10
\font\Bsl=cmsl10 scaled 1200 \font\Bmmi=cmmi10 scaled 1200
\font\BBf=cmbx12 scaled 1200 \font\BMmi=cmmi10 scaled 1440

\def\atletter{\edef\atrestore{\catcode`\noexpand\@=\the\catcode`\@\space}
 \catcode`\@=11}

\newread\@ux \newwrite\@@x \newwrite\@@cd
\let\@np@@\input
\def\@np@t#1{\openin\@ux#1\relax\ifeof\@ux\else\closein\@ux\relax\@np@@ #1\fi}
\def\input#1 {\openin\@ux#1\relax\ifeof\@ux\wrs@x{No file #1}\else
 \closein\@ux\relax\@np@@ #1\fi}
\def\Input#1 {\relax} 

\def\wr@@x#1{} \def\wrs@x{\immediate\write\sixt@@n}

\def\readldf{\@np@t{\jobname.ldf}}
\def\writeldf{\def\wr@@x{\immediate\write\@@x}
 \def\cl@selbl{\wr@@x{\string\Snodef{\the\Sno}}\wr@@x{\string\endinput}%
 \immediate\closeout\@@x} \immediate\openout\@@x\jobname.ldf}
\let\cl@selbl\relax

\def\tod@y{\ifcase\month\or
 January\or February\or March\or April\or May\or June\or July\or
 August\or September\or October\or November\or December\fi\space\,
\number\day,\space\,\number\year}

\newcount\c@time
\def\h@@r{hh}\def\m@n@te{mm}
\def\wh@tt@me{\c@time\time\divide\c@time 60\edef\h@@r{\number\c@time}%
 \multiply\c@time -60\advance\c@time\time\edef
 \m@n@te{\ifnum\c@time<10 0\fi\number\c@time}}
\def\t@me{\h@@r\/{\rm:}\m@n@te}  \let\whattime\wh@tt@me
\def\today{\tod@y\wr@@x{\string\todaydef{\tod@y}}}
\def\nowtime{\t@me{\let\/\ic@\wr@@x{\string\nowtimedef{\t@me}}}}
\def\todaydef#1{} \def\nowtimedef#1{}

\def\em#1{{\it #1\/}} \def\emph#1{{\sl #1\/}}

\def\itemleft#1{\par\setbox\z@\hbox{\rm #1\enspace}\hangindent\wd\z@
 \hglue-2\parindent\kern\wd\z@\textindent{\rm#1}}
\def\itemflat#1{\par\setbox\z@\hbox{\rm #1\enspace}\hang\ifnum\wd\z@>\parindent
 \noindent\unhbox\z@\ignore\else\textindent{\rm#1}\fi}

\newcount\itemlet
\def\newbi{\itemlet 96} \newbi
\def\bitem{\gad\itemlet \par\hangindent1.5\parindent
 \hglue-.5\parindent\textindent{\rm\rlap{\char\the\itemlet}\hp{b})}}
\def\atem{\newbi\bitem}

\newcount\itemrm

\def\iitem{\gad\itemrm \par\hangindent1.5\parindent
 \hglue-.5\parindent\textindent{\rm\hp{v}\llap{\romannumeral\the\itemrm})}}

\def\center{\par\begingroup\leftskip\z@ plus \hsize \rightskip\leftskip
 \parindent\z@\parfillskip\z@skip \def\\{\unskip\break}}
\def\endcenter{\endgraf\endgroup}

\def\Abstract{\begingroup\narrower\nt{\bf Abstract.}\enspace\ignore}
\def\endAbs{\endgraf\endgroup}

\let\b@gr@@\begingroup \let\B@gr@@\begingroup
\def\b@gr@{\b@gr@@\let\b@gr@@\undefined}
\def\B@gr@{\B@gr@@\let\B@gr@@\undefined}

\def\@fn@xt#1#2#3{\let\@ch@r=#1\def\n@xt{\ifx\t@st@\@ch@r
 \def\n@@xt{#2}\else\def\n@@xt{#3}\fi\n@@xt}\futurelet\t@st@\n@xt}

\def\@fwd@@#1#2#3{\setbox\z@\hbox{#1}\ifdim\wd\z@>\z@#2\else#3\fi}
\def\s@twd@#1#2{\setbox\z@\hbox{#2}#1\wd\z@}

\def\r@st@re#1{\let#1\s@v@} \def\s@v@d@f{\let\s@v@}

\def\p@sk@p#1#2{\par\skip@#2\relax\ifdim\lastskip<\skip@\relax\removelastskip
 \ifnum#1=\z@\else\penalty#1\relax\fi\vskip\skip@
 \else\ifnum#1=\z@\else\penalty#1\relax\fi\fi}
\def\sk@@p#1{\par\skip@#1\relax\ifdim\lastskip<\skip@\relax\removelastskip
 \vskip\skip@\fi}

\newbox\p@b@ld
\def\poorbold#1{\setbox\p@b@ld\hbox{#1}\kern-.01em\copy\p@b@ld\kern-\wd\p@b@ld
 \kern.02em\copy\p@b@ld\kern-\wd\p@b@ld\kern-.012em\raise.02em\box\p@b@ld}

\ifx\plainfootnote\undefined \let\plainfootnote\footnote \fi

\let\s@v@\proclaim \let\proclaim\relax
\def\r@R@fs#1{\let#1\s@R@fs} \let\s@R@fs\Refs \let\Refs\relax
\def\r@endd@#1{\let#1\s@endd@} \let\s@endd@\enddocument
\let\bye\relax

\def\myR@fs{\@fn@xt[\m@R@f@\m@R@fs} \def\m@R@fs{\@fn@xt*\m@r@f@@\m@R@f@@}
\def\m@R@f@@{\m@R@f@[References]} \def\m@r@f@@*{\m@R@f@[]}

\def\Twelvepoint{\twelvepoint \let\Bbf\BBf \let\Bmmi\BMmi
\font\Brm=cmr12 scaled 1200 \font\Bit=cmti12 scaled 1200
\font\ssf=cmss10 scaled 1200 \font\Bsl=cmsl10 scaled 1440
\font\BBf=cmbx12 scaled 1440 \font\BMmi=cmmi10 scaled 1728}

\newif\ifamsppt

\newdimen\b@gsize

\newdimen\r@f@nd \newbox\r@f@b@x \newbox\adjb@x
\newbox\p@nct@ \newbox\k@yb@x \newcount\rcount
\newbox\b@b@x \newbox\p@p@rb@x \newbox\j@@rb@x \newbox\y@@rb@x
\newbox\v@lb@x \newbox\is@b@x \newbox\p@g@b@x \newif\ifp@g@ \newif\ifp@g@s
\newbox\inb@@kb@x \newbox\b@@kb@x \newbox\p@blb@x \newbox\p@bl@db@x
\newbox\ed@b@x \newif\ifed@ \newif\ifed@s \newif\if@fl@b \newif\if@fn@m
\newbox\p@p@nf@b@x \newbox\inf@b@x \newbox\b@@nf@b@x

\newif\ifp@gen@

\@ft@\ifx\csname amsppt.sty\endcsname\relax

\headline={\hfil}
\footline={\ifp@gen@\ifnum\pageno=\z@\else\hfil\foliorm\folio\fi\else
 \ifnum\pageno=\z@\hfil\foliorm\folio\fi\fi\hfil\global\p@gen@true}
\parindent1pc

\font@\tensmc=cmcsc10
\font@\sevenex=cmex7
\font@\sevenit=cmti7
\font@\eightrm=cmr8
\font@\sixrm=cmr6
\font@\eighti=cmmi8 \skewchar\eighti='177
\font@\sixi=cmmi6 \skewchar\sixi='177
\font@\eightsy=cmsy8 \skewchar\eightsy='60
\font@\sixsy=cmsy6 \skewchar\sixsy='60
\font@\eightex=cmex8
\font@\eightbf=cmbx8
\font@\sixbf=cmbx6
\font@\eightit=cmti8
\font@\eightsl=cmsl8
\font@\eightsmc=cmcsc8
\font@\eighttt=cmtt8
\font@\ninerm=cmr9
\font@\ninei=cmmi9 \skewchar\ninei='177
\font@\ninesy=cmsy9 \skewchar\ninesy='60
\font@\nineex=cmex9
\font@\ninebf=cmbx9
\font@\nineit=cmti9
\font@\ninesl=cmsl9
\font@\ninesmc=cmcsc9
\font@\ninemsa=msam9
\font@\ninemsb=msbm9
\font@\nineeufm=eufm9
\font@\eightmsa=msam8
\font@\eightmsb=msbm8
\font@\eighteufm=eufm8
\font@\sixmsa=msam6
\font@\sixmsb=msbm6
\font@\sixeufm=eufm6

\loadmsam\loadmsbm\loadeufm
\input amssym.tex

\def\footnoterule{\kern-3\p@\hrule width5pc\kern 2.6\p@}
\def\m@k@foot#1{\insert\footins
 {\interlinepenalty\interfootnotelinepenalty
 \eightpoint\splittopskip\ht\strutbox\splitmaxdepth\dp\strutbox
 \floatingpenalty\@MM\leftskip\z@\rightskip\z@
 \spaceskip\z@\xspaceskip\z@
 \leavevmode\footstrut\ignore#1\unskip\lower\dp\strutbox
 \vbox to\dp\strutbox{}}}
\def\ftext#1{\m@k@foot{\vsk-.8>\nt #1}}
\def\pr@cl@@m#1{\p@sk@p{-100}\medskipamount\b@gr@\nt\ignore
 \bf #1\unskip.\enspace\sl\ignore}
\outer\def\proclaim{\pr@cl@@m} \s@v@d@f\proclaim \let\proclaim\relax
\def\endproclaim{\endgroup\p@sk@p{55}\medskipamount}
\def\demo#1{\sk@@p\medskipamount\nt{\ignore\it #1\unskip.}\enspace
 \ignore}
\def\enddemo{\sk@@p\medskipamount}

\def\cite#1{{\rm[#1]}} \let\nofrills\relax
 \def\Refs#1#2{\relax}

\def\big@#1#2{{\hbox{$\left#2\vcenter to#1\b@gsize{}%
 \right.\nulldelimiterspace\z@\m@th$}}}
\def\big{\big@\@ne}
\def\Big{\big@{1.5}}
\def\bigg{\big@\tw@}
\def\Bigg{\big@{2.5}}
\normallineskiplimit\p@

\def\tenpoint{\p@@\p@ \normallineskiplimit\p@@
 \mathsurround\m@ths@r \normalbaselineskip12\p@@
 \abovedisplayskip12\p@@ plus3\p@@ minus9\p@@
 \belowdisplayskip\abovedisplayskip
 \abovedisplayshortskip\z@ plus3\p@@
 \belowdisplayshortskip7\p@@ plus3\p@@ minus4\p@@
 \textonlyfont@\rm\tenrm \textonlyfont@\it\tenit
 \textonlyfont@\sl\tensl \textonlyfont@\bf\tenbf
 \textonlyfont@\smc\tensmc \textonlyfont@\tt\tentt
 \ifsyntax@ \def\big##1{{\hbox{$\left##1\right.$}}}%
  \let\Big\big \let\bigg\big \let\Bigg\big
 \else
  \textfont\z@\tenrm \scriptfont\z@\sevenrm \scriptscriptfont\z@\fiverm
  \textfont\@ne\teni \scriptfont\@ne\seveni \scriptscriptfont\@ne\fivei
  \textfont\tw@\tensy \scriptfont\tw@\sevensy \scriptscriptfont\tw@\fivesy
  \textfont\thr@@\tenex \scriptfont\thr@@\sevenex
	\scriptscriptfont\thr@@\sevenex
  \textfont\itfam\tenit \scriptfont\itfam\sevenit
	\scriptscriptfont\itfam\sevenit
  \textfont\bffam\tenbf \scriptfont\bffam\sevenbf
	\scriptscriptfont\bffam\fivebf
  \textfont\msafam\tenmsa \scriptfont\msafam\sevenmsa
	\scriptscriptfont\msafam\fivemsa
  \textfont\msbfam\tenmsb \scriptfont\msbfam\sevenmsb
	\scriptscriptfont\msbfam\fivemsb
  \textfont\eufmfam\teneufm \scriptfont\eufmfam\seveneufm
	\scriptscriptfont\eufmfam\fiveeufm
  \setbox\strutbox\hbox{\vrule height8.5\p@@ depth3.5\p@@ width\z@}%
  \setbox\strutbox@\hbox{\lower.5\normallineskiplimit\vbox{%
	\kern-\normallineskiplimit\copy\strutbox}}%
   \setbox\z@\vbox{\hbox{$($}\kern\z@}\b@gsize1.2\ht\z@
  \fi
  \normalbaselines\rm\dotsspace@1.5mu\ex@.2326ex\jot3\ex@}

\def\eightpoint{\p@@.8\p@ \normallineskiplimit\p@@
 \mathsurround\m@ths@r \normalbaselineskip10\p@
 \abovedisplayskip10\p@ plus2.4\p@ minus7.2\p@
 \belowdisplayskip\abovedisplayskip
 \abovedisplayshortskip\z@ plus3\p@@
 \belowdisplayshortskip7\p@@ plus3\p@@ minus4\p@@
 \textonlyfont@\rm\eightrm \textonlyfont@\it\eightit
 \textonlyfont@\sl\eightsl \textonlyfont@\bf\eightbf
 \textonlyfont@\smc\eightsmc \textonlyfont@\tt\eighttt
 \ifsyntax@\def\big##1{{\hbox{$\left##1\right.$}}}%
  \let\Big\big \let\bigg\big \let\Bigg\big
 \else
  \textfont\z@\eightrm \scriptfont\z@\sixrm \scriptscriptfont\z@\fiverm
  \textfont\@ne\eighti \scriptfont\@ne\sixi \scriptscriptfont\@ne\fivei
  \textfont\tw@\eightsy \scriptfont\tw@\sixsy \scriptscriptfont\tw@\fivesy
  \textfont\thr@@\eightex \scriptfont\thr@@\sevenex
	\scriptscriptfont\thr@@\sevenex
  \textfont\itfam\eightit \scriptfont\itfam\sevenit
	\scriptscriptfont\itfam\sevenit
  \textfont\bffam\eightbf \scriptfont\bffam\sixbf
	\scriptscriptfont\bffam\fivebf
  \textfont\msafam\eightmsa \scriptfont\msafam\sixmsa
	\scriptscriptfont\msafam\fivemsa
  \textfont\msbfam\eightmsb \scriptfont\msbfam\sixmsb
	\scriptscriptfont\msbfam\fivemsb
  \textfont\eufmfam\eighteufm \scriptfont\eufmfam\sixeufm
	\scriptscriptfont\eufmfam\fiveeufm
 \setbox\strutbox\hbox{\vrule height7\p@ depth3\p@ width\z@}%
 \setbox\strutbox@\hbox{\raise.5\normallineskiplimit\vbox{%
   \kern-\normallineskiplimit\copy\strutbox}}%
 \setbox\z@\vbox{\hbox{$($}\kern\z@}\b@gsize1.2\ht\z@
 \fi
 \normalbaselines\eightrm\dotsspace@1.5mu\ex@.2326ex\jot3\ex@}

\def\ninepoint{\p@@.9\p@ \normallineskiplimit\p@@
 \mathsurround\m@ths@r \normalbaselineskip11\p@
 \abovedisplayskip11\p@ plus2.7\p@ minus8.1\p@
 \belowdisplayskip\abovedisplayskip
 \abovedisplayshortskip\z@ plus3\p@@
 \belowdisplayshortskip7\p@@ plus3\p@@ minus4\p@@
 \textonlyfont@\rm\ninerm \textonlyfont@\it\nineit
 \textonlyfont@\sl\ninesl \textonlyfont@\bf\ninebf
 \textonlyfont@\smc\ninesmc \textonlyfont@\tt\ninett
 \ifsyntax@ \def\big##1{{\hbox{$\left##1\right.$}}}%
  \let\Big\big \let\bigg\big \let\Bigg\big
 \else
  \textfont\z@\ninerm \scriptfont\z@\sevenrm \scriptscriptfont\z@\fiverm
  \textfont\@ne\ninei \scriptfont\@ne\seveni \scriptscriptfont\@ne\fivei
  \textfont\tw@\ninesy \scriptfont\tw@\sevensy \scriptscriptfont\tw@\fivesy
  \textfont\thr@@\nineex \scriptfont\thr@@\sevenex
	\scriptscriptfont\thr@@\sevenex
  \textfont\itfam\nineit \scriptfont\itfam\sevenit
	\scriptscriptfont\itfam\sevenit
  \textfont\bffam\ninebf \scriptfont\bffam\sevenbf
	\scriptscriptfont\bffam\fivebf
  \textfont\msafam\ninemsa \scriptfont\msafam\sevenmsa
	\scriptscriptfont\msafam\fivemsa
  \textfont\msbfam\ninemsb \scriptfont\msbfam\sevenmsb
	\scriptscriptfont\msbfam\fivemsb
  \textfont\eufmfam\nineeufm \scriptfont\eufmfam\seveneufm
	\scriptscriptfont\eufmfam\fiveeufm
  \setbox\strutbox\hbox{\vrule height8.5\p@@ depth3.5\p@@ width\z@}%
  \setbox\strutbox@\hbox{\lower.5\normallineskiplimit\vbox{%
	\kern-\normallineskiplimit\copy\strutbox}}%
   \setbox\z@\vbox{\hbox{$($}\kern\z@}\b@gsize1.2\ht\z@
  \fi
  \normalbaselines\rm\dotsspace@1.5mu\ex@.2326ex\jot3\ex@}

\font@\twelverm=cmr10 scaled 1200
\font@\twelveit=cmti10 scaled 1200
\font@\twelvesl=cmsl10 scaled 1200
\font@\twelvebf=cmbx10 scaled 1200
\font@\twelvesmc=cmcsc10 scaled 1200
\font@\twelvett=cmtt10 scaled 1200
\font@\twelvei=cmmi10 scaled 1200 \skewchar\twelvei='177
\font@\twelvesy=cmsy10 scaled 1200 \skewchar\twelvesy='60
\font@\twelveex=cmex10 scaled 1200
\font@\twelvemsa=msam10 scaled 1200
\font@\twelvemsb=msbm10 scaled 1200
\font@\twelveeufm=eufm10 scaled 1200

\def\twelvepoint{\p@@1.2\p@ \normallineskiplimit\p@@
 \mathsurround\m@ths@r \normalbaselineskip12\p@@
 \abovedisplayskip12\p@@ plus3\p@@ minus9\p@@
 \belowdisplayskip\abovedisplayskip
 \abovedisplayshortskip\z@ plus3\p@@
 \belowdisplayshortskip7\p@@ plus3\p@@ minus4\p@@
 \textonlyfont@\rm\twelverm \textonlyfont@\it\twelveit
 \textonlyfont@\sl\twelvesl \textonlyfont@\bf\twelvebf
 \textonlyfont@\smc\twelvesmc \textonlyfont@\tt\twelvett
 \ifsyntax@ \def\big##1{{\hbox{$\left##1\right.$}}}%
  \let\Big\big \let\bigg\big \let\Bigg\big
 \else
  \textfont\z@\twelverm \scriptfont\z@\eightrm \scriptscriptfont\z@\sixrm
  \textfont\@ne\twelvei \scriptfont\@ne\eighti \scriptscriptfont\@ne\sixi
  \textfont\tw@\twelvesy \scriptfont\tw@\eightsy \scriptscriptfont\tw@\sixsy
  \textfont\thr@@\twelveex \scriptfont\thr@@\eightex
	\scriptscriptfont\thr@@\sevenex
  \textfont\itfam\twelveit \scriptfont\itfam\eightit
	\scriptscriptfont\itfam\sevenit
  \textfont\bffam\twelvebf \scriptfont\bffam\eightbf
	\scriptscriptfont\bffam\sixbf
  \textfont\msafam\twelvemsa \scriptfont\msafam\eightmsa
	\scriptscriptfont\msafam\sixmsa
  \textfont\msbfam\twelvemsb \scriptfont\msbfam\eightmsb
	\scriptscriptfont\msbfam\sixmsb
  \textfont\eufmfam\twelveeufm \scriptfont\eufmfam\eighteufm
	\scriptscriptfont\eufmfam\sixeufm
  \setbox\strutbox\hbox{\vrule height8.5\p@@ depth3.5\p@@ width\z@}%
  \setbox\strutbox@\hbox{\lower.5\normallineskiplimit\vbox{%
	\kern-\normallineskiplimit\copy\strutbox}}%
  \setbox\z@\vbox{\hbox{$($}\kern\z@}\b@gsize1.2\ht\z@
  \fi
  \normalbaselines\rm\dotsspace@1.5mu\ex@.2326ex\jot3\ex@}

\font@\twelvetrm=cmr10 at 12truept
\font@\twelvetit=cmti10 at 12truept
\font@\twelvetsl=cmsl10 at 12truept
\font@\twelvetbf=cmbx10 at 12truept
\font@\twelvetsmc=cmcsc10 at 12truept
\font@\twelvettt=cmtt10 at 12truept
\font@\twelveti=cmmi10 at 12truept \skewchar\twelveti='177
\font@\twelvetsy=cmsy10 at 12truept \skewchar\twelvetsy='60
\font@\twelvetex=cmex10 at 12truept
\font@\twelvetmsa=msam10 at 12truept
\font@\twelvetmsb=msbm10 at 12truept
\font@\twelveteufm=eufm10 at 12truept

\def\twelvetruepoint{\p@@1.2truept \normallineskiplimit\p@@
 \mathsurround\m@ths@r \normalbaselineskip12\p@@
 \abovedisplayskip12\p@@ plus3\p@@ minus9\p@@
 \belowdisplayskip\abovedisplayskip
 \abovedisplayshortskip\z@ plus3\p@@
 \belowdisplayshortskip7\p@@ plus3\p@@ minus4\p@@
 \textonlyfont@\rm\twelvetrm \textonlyfont@\it\twelvetit
 \textonlyfont@\sl\twelvetsl \textonlyfont@\bf\twelvetbf
 \textonlyfont@\smc\twelvetsmc \textonlyfont@\tt\twelvettt
 \ifsyntax@ \def\big##1{{\hbox{$\left##1\right.$}}}%
  \let\Big\big \let\bigg\big \let\Bigg\big
 \else
  \textfont\z@\twelvetrm \scriptfont\z@\eightrm \scriptscriptfont\z@\sixrm
  \textfont\@ne\twelveti \scriptfont\@ne\eighti \scriptscriptfont\@ne\sixi
  \textfont\tw@\twelvetsy \scriptfont\tw@\eightsy \scriptscriptfont\tw@\sixsy
  \textfont\thr@@\twelvetex \scriptfont\thr@@\eightex
	\scriptscriptfont\thr@@\sevenex
  \textfont\itfam\twelvetit \scriptfont\itfam\eightit
	\scriptscriptfont\itfam\sevenit
  \textfont\bffam\twelvetbf \scriptfont\bffam\eightbf
	\scriptscriptfont\bffam\sixbf
  \textfont\msafam\twelvetmsa \scriptfont\msafam\eightmsa
	\scriptscriptfont\msafam\sixmsa
  \textfont\msbfam\twelvetmsb \scriptfont\msbfam\eightmsb
	\scriptscriptfont\msbfam\sixmsb
  \textfont\eufmfam\twelveteufm \scriptfont\eufmfam\eighteufm
	\scriptscriptfont\eufmfam\sixeufm
  \setbox\strutbox\hbox{\vrule height8.5\p@@ depth3.5\p@@ width\z@}%
  \setbox\strutbox@\hbox{\lower.5\normallineskiplimit\vbox{%
	\kern-\normallineskiplimit\copy\strutbox}}%
  \setbox\z@\vbox{\hbox{$($}\kern\z@}\b@gsize1.2\ht\z@
  \fi
  \normalbaselines\rm\dotsspace@1.5mu\ex@.2326ex\jot3\ex@}

\font@\elevenrm=cmr10 scaled 1095
\font@\elevenit=cmti10 scaled 1095
\font@\elevensl=cmsl10 scaled 1095
\font@\elevenbf=cmbx10 scaled 1095
\font@\elevensmc=cmcsc10 scaled 1095
\font@\eleventt=cmtt10 scaled 1095
\font@\eleveni=cmmi10 scaled 1095 \skewchar\eleveni='177
\font@\elevensy=cmsy10 scaled 1095 \skewchar\elevensy='60
\font@\elevenex=cmex10 scaled 1095
\font@\elevenmsa=msam10 scaled 1095
\font@\elevenmsb=msbm10 scaled 1095
\font@\eleveneufm=eufm10 scaled 1095

\def\elevenpoint{\p@@1.1\p@ \normallineskiplimit\p@@
 \mathsurround\m@ths@r \normalbaselineskip12\p@@
 \abovedisplayskip12\p@@ plus3\p@@ minus9\p@@
 \belowdisplayskip\abovedisplayskip
 \abovedisplayshortskip\z@ plus3\p@@
 \belowdisplayshortskip7\p@@ plus3\p@@ minus4\p@@
 \textonlyfont@\rm\elevenrm \textonlyfont@\it\elevenit
 \textonlyfont@\sl\elevensl \textonlyfont@\bf\elevenbf
 \textonlyfont@\smc\elevensmc \textonlyfont@\tt\eleventt
 \ifsyntax@ \def\big##1{{\hbox{$\left##1\right.$}}}%
  \let\Big\big \let\bigg\big \let\Bigg\big
 \else
  \textfont\z@\elevenrm \scriptfont\z@\eightrm \scriptscriptfont\z@\sixrm
  \textfont\@ne\eleveni \scriptfont\@ne\eighti \scriptscriptfont\@ne\sixi
  \textfont\tw@\elevensy \scriptfont\tw@\eightsy \scriptscriptfont\tw@\sixsy
  \textfont\thr@@\elevenex \scriptfont\thr@@\eightex
	\scriptscriptfont\thr@@\sevenex
  \textfont\itfam\elevenit \scriptfont\itfam\eightit
	\scriptscriptfont\itfam\sevenit
  \textfont\bffam\elevenbf \scriptfont\bffam\eightbf
	\scriptscriptfont\bffam\sixbf
  \textfont\msafam\elevenmsa \scriptfont\msafam\eightmsa
	\scriptscriptfont\msafam\sixmsa
  \textfont\msbfam\elevenmsb \scriptfont\msbfam\eightmsb
	\scriptscriptfont\msbfam\sixmsb
  \textfont\eufmfam\eleveneufm \scriptfont\eufmfam\eighteufm
	\scriptscriptfont\eufmfam\sixeufm
  \setbox\strutbox\hbox{\vrule height8.5\p@@ depth3.5\p@@ width\z@}%
  \setbox\strutbox@\hbox{\lower.5\normallineskiplimit\vbox{%
	\kern-\normallineskiplimit\copy\strutbox}}%
  \setbox\z@\vbox{\hbox{$($}\kern\z@}\b@gsize1.2\ht\z@
  \fi
  \normalbaselines\rm\dotsspace@1.5mu\ex@.2326ex\jot3\ex@}

\def\m@R@f@[#1]{\mathsurzero{
 \s@ct{}{#1}}\wr@@c{\string\Refcd{#1}{\the\pageno}}\B@gr@
 \frenchspacing\rcount\z@\refkey{[##1]}\refno{[##1]}\widest{AZ}\keyright
 \let\Key\key\let\refin\relax}
\def\widest#1{\s@twd@\r@f@nd{\r@fk@y{#1}\enspace}}
\def\widestno#1{\s@twd@\r@f@nd{\r@fn@{#1}\enspace}}
\def\widestlabel#1{\s@twd@\r@f@nd{#1\enspace}}
\def\refkey{\def\r@fk@y##1} \def\refno{\def\r@fn@##1}
\def\keyright{\def\r@fit@m{\hang\textindent}}
\def\keyflat{\def\r@fit@m##1{\setbox\z@\hbox{\rm ##1\enspace}\hang\noindent
 \ifnum\wd\z@<\parindent\indent\hglue-\wd\z@\fi\unhbox\z@}}

\def\R@fb@x{\global\setbox\r@f@b@x} \def\K@yb@x{\global\setbox\k@yb@x}
\def\ref{\par\b@gr@\rm\R@fb@x\box\voidb@x\K@yb@x\box\voidb@x\@fn@mfalse
 \@fl@bfalse\b@g@nr@f}
\def\c@nc@t#1{\setbox\z@\lastbox
 \setbox\adjb@x\hbox{\unhbox\adjb@x\unhbox\z@\unskip\unskip\unpenalty#1}}
\def\adjust#1{\relax\ifmmode\penalty-\@M\null\hfil$\clubpenalty\z@
 \widowpenalty\z@\interlinepenalty\z@\offinterlineskip\endgraf
 \setbox\z@\lastbox\unskip\unpenalty\c@nc@t{#1}\nt$\hfil\penalty-\@M
 \else\endgraf\c@nc@t{#1}\nt\fi}
\def\adjustnext#1{\P@nct\hbox{#1}\ignore}
\def\adjustend#1{\def\@djp@{#1}\ignore}

\def\cl@s@{\adjust{\@djp@}\endgraf\setbox\z@\lastbox
 \global\setbox\@ne\hbox{\unhbox\adjb@x\ifvoid\z@\else\unhbox\z@\unskip\unskip
 \unpenalty\fi}\egroup\ifnum\c@rr@nt=\k@yb@x\global\fi
 \setbox\c@rr@nt\hbox{\unhbox\@ne\box\p@nct@}\P@nct\null}
\def\@p@n#1{\def\c@rr@nt{#1}\setbox\c@rr@nt\vbox\bgroup\let\@djp@\relax
 \hsize\maxdimen\nt}
\def\b@g@nr@f{\bgroup\@p@n\z@}
\def\key{\cl@s@\ifvoid\k@yb@x\@p@n\k@yb@x\else\@p@n\z@\fi}
\def\label{\cl@s@\ifvoid\k@yb@x\global\@fl@btrue\@p@n\k@yb@x\else\@p@n\z@\fi}
\def\no{\cl@s@\ifvoid\k@yb@x\gad\rcount\global\@fn@mtrue
 \K@yb@x\hbox{\the\rcount}\fi\@p@n\z@}
\def\labelno{\cl@s@\ifvoid\k@yb@x\gad\rcount\@fl@btrue\@p@n\k@yb@x\the\rcount
 \else\@p@n\z@\fi}
\def\by{\cl@s@\@p@n\b@b@x} \def\paper{\cl@s@\@p@n\p@p@rb@x\it\ignore}
\def\jour{\cl@s@\@p@n\j@@rb@x} \def\yr{\cl@s@\@p@n\y@@rb@x}
\def\vol{\cl@s@\@p@n\v@lb@x\bf\ignore} \def\issue{\cl@s@\@p@n\is@b@x}
\def\page{\cl@s@\ifp@g@s\@p@n\z@\else\p@g@true\@p@n\p@g@b@x\fi}
\def\pages{\cl@s@\ifp@g@\@p@n\z@\else\p@g@strue\@p@n\p@g@b@x\fi}
\def\inbook{\cl@s@\@p@n\inb@@kb@x} \def\book{\cl@s@\@p@n\b@@kb@x\it\ignore}
\def\publ{\cl@s@\@p@n\p@blb@x} \def\publaddr{\cl@s@\@p@n\p@bl@db@x}
\def\ed{\cl@s@\ifed@s\@p@n\z@\else\ed@true\@p@n\ed@b@x\fi}
\def\eds{\cl@s@\ifed@\@p@n\z@\else\ed@strue\@p@n\ed@b@x\fi}
\def\info{\cl@s@\@p@n\inf@b@x} \def\paperinfo{\cl@s@\@p@n\p@p@nf@b@x}
\def\bookinfo{\cl@s@\@p@n\b@@nf@b@x} \let\finalinfo\info
\def\P@nct{\global\setbox\p@nct@} \def\nopunct{\P@nct\box\voidb@x}
\def\p@@@t#1#2{\ifvoid\p@nct@\else#1\unhbox\p@nct@#2\fi}
\def\sp@@{\penalty-50 \space\hskip\z@ plus.1em}
\def\c@mm@{\p@@@t,\sp@@} \def\sp@c@{\p@@@t\empty\sp@@} \def\p@@nt{.\kern.3em}
\def\p@tb@x#1#2{\ifvoid#1\else#2\@nb@x#1\fi}
\def\@nb@x#1{\unhbox#1\P@nct\lastbox}
\def\endr@f@{\cl@s@\nopunct
 \R@fb@x\hbox{\unhbox\r@f@b@x \p@tb@x\b@b@x\empty
 \ifvoid\j@@rb@x\ifvoid\inb@@kb@x\ifvoid\p@p@rb@x\ifvoid\b@@kb@x
  \ifvoid\p@p@nf@b@x\ifvoid\b@@nf@b@x
  \p@tb@x\v@lb@x\c@mm@ \ifvoid\y@@rb@x\else\sp@c@(\@nb@x\y@@rb@x)\fi
  \p@tb@x\is@b@x{\c@mm@ no\p@@nt}\p@tb@x\p@g@b@x\c@mm@ \p@tb@x\inf@b@x\c@mm@
  \else\p@tb@x \b@@nf@b@x\c@mm@ \p@tb@x\v@lb@x\c@mm@
  \p@tb@x\is@b@x{\sp@c@ no\p@@nt}%
  \ifvoid\ed@b@x\else\sp@c@(\@nb@x\ed@b@x,\space\ifed@ ed.\else eds.\fi)\fi
  \p@tb@x\p@blb@x\c@mm@ \p@tb@x\p@bl@db@x\c@mm@ \p@tb@x\y@@rb@x\c@mm@
  \p@tb@x\p@g@b@x{\c@mm@\ifp@g@ p\p@@nt\else pp\p@@nt\fi}%
  \p@tb@x\inf@b@x\c@mm@\fi
  \else \p@tb@x\p@p@nf@b@x\c@mm@ \p@tb@x\v@lb@x\c@mm@
  \ifvoid\y@@rb@x\else\sp@c@(\@nb@x\y@@rb@x)\fi
  \p@tb@x\is@b@x{\c@mm@ no\p@@nt}\p@tb@x\p@g@b@x\c@mm@ \p@tb@x\inf@b@x\c@mm@\fi
  \else \p@tb@x\b@@kb@x\c@mm@
  \p@tb@x\b@@nf@b@x\c@mm@ \p@tb@x\p@blb@x\c@mm@
  \p@tb@x\p@bl@db@x\c@mm@ \p@tb@x\y@@rb@x\c@mm@
  \ifvoid\p@g@b@x\else\c@mm@\@nb@x\p@g@b@x p\fi \p@tb@x\inf@b@x\c@mm@ \fi
  \else \c@mm@\@nb@x\p@p@rb@x\ic@\p@tb@x\p@p@nf@b@x\c@mm@
  \p@tb@x\v@lb@x\sp@c@ \ifvoid\y@@rb@x\else\sp@c@(\@nb@x\y@@rb@x)\fi
  \p@tb@x\is@b@x{\c@mm@ no\p@@nt}\p@tb@x\p@g@b@x\c@mm@\p@tb@x\inf@b@x\c@mm@\fi
  \else \p@tb@x\p@p@rb@x\c@mm@\ic@\p@tb@x\p@p@nf@b@x\c@mm@
  \c@mm@\@nb@x\inb@@kb@x \p@tb@x\b@@nf@b@x\c@mm@ \p@tb@x\v@lb@x\sp@c@
  \p@tb@x\is@b@x{\sp@c@ no\p@@nt}%
  \ifvoid\ed@b@x\else\sp@c@(\@nb@x\ed@b@x,\space\ifed@ ed.\else eds.\fi)\fi
  \p@tb@x\p@blb@x\c@mm@ \p@tb@x\p@bl@db@x\c@mm@ \p@tb@x\y@@rb@x\c@mm@
  \p@tb@x\p@g@b@x{\c@mm@\ifp@g@ p\p@@nt\else pp\p@@nt\fi}%
  \p@tb@x\inf@b@x\c@mm@\fi
  \else\p@tb@x\p@p@rb@x\c@mm@\ic@\p@tb@x\p@p@nf@b@x\c@mm@\p@tb@x\j@@rb@x\c@mm@
  \p@tb@x\v@lb@x\sp@c@ \ifvoid\y@@rb@x\else\sp@c@(\@nb@x\y@@rb@x)\fi
  \p@tb@x\is@b@x{\c@mm@ no\p@@nt}\p@tb@x\p@g@b@x\c@mm@ \p@tb@x\inf@b@x\c@mm@
 \fi}}
\def\m@r@f#1#2{\endr@f@\ifvoid\p@nct@\else\R@fb@x\hbox{\unhbox\r@f@b@x
 #1\unhbox\p@nct@\penalty-200\enskip#2}\fi\egroup\b@g@nr@f}
\def\endref{\endr@f@\ifvoid\p@nct@\else\R@fb@x\hbox{\unhbox\r@f@b@x.}\fi
 \parindent\r@f@nd
 \r@fit@m{\ifvoid\k@yb@x\else\if@fn@m\r@fn@{\unhbox\k@yb@x}\else
 \if@fl@b\unhbox\k@yb@x\else\r@fk@y{\unhbox\k@yb@x}\fi\fi\fi}\unhbox\r@f@b@x
 \endgraf\egroup\endgroup}
\def\moreref{\m@r@f;\empty}
\def\transl{\m@r@f;{\unskip\space
 {\sl English translation\ic@}:\penalty-66 \space}}
\def\endRefs{\endgraf\goodbreak\endgroup}

\hyphenation{acad-e-my acad-e-mies af-ter-thought anom-aly anom-alies
an-ti-deriv-a-tive an-tin-o-my an-tin-o-mies apoth-e-o-ses
apoth-e-o-sis ap-pen-dix ar-che-typ-al as-sign-a-ble as-sist-ant-ship
as-ymp-tot-ic asyn-chro-nous at-trib-uted at-trib-ut-able bank-rupt
bank-rupt-cy bi-dif-fer-en-tial blue-print busier busiest
cat-a-stroph-ic cat-a-stroph-i-cally con-gress cross-hatched data-base
de-fin-i-tive de-riv-a-tive dis-trib-ute dri-ver dri-vers eco-nom-ics
econ-o-mist elit-ist equi-vari-ant ex-quis-ite ex-tra-or-di-nary
flow-chart for-mi-da-ble forth-right friv-o-lous ge-o-des-ic
ge-o-det-ic geo-met-ric griev-ance griev-ous griev-ous-ly
hexa-dec-i-mal ho-lo-no-my ho-mo-thetic ideals idio-syn-crasy
in-fin-ite-ly in-fin-i-tes-i-mal ir-rev-o-ca-ble key-stroke
lam-en-ta-ble light-weight mal-a-prop-ism man-u-script mar-gin-al
meta-bol-ic me-tab-o-lism meta-lan-guage me-trop-o-lis
met-ro-pol-i-tan mi-nut-est mol-e-cule mono-chrome mono-pole
mo-nop-oly mono-spline mo-not-o-nous mul-ti-fac-eted mul-ti-plic-able
non-euclid-ean non-iso-mor-phic non-smooth par-a-digm par-a-bol-ic
pa-rab-o-loid pa-ram-e-trize para-mount pen-ta-gon phe-nom-e-non
post-script pre-am-ble pro-ce-dur-al pro-hib-i-tive pro-hib-i-tive-ly
pseu-do-dif-fer-en-tial pseu-do-fi-nite pseu-do-nym qua-drat-ic
quad-ra-ture qua-si-smooth qua-si-sta-tion-ary qua-si-tri-an-gu-lar
quin-tes-sence quin-tes-sen-tial re-arrange-ment rec-tan-gle
ret-ri-bu-tion retro-fit retro-fit-ted right-eous right-eous-ness
ro-bot ro-bot-ics sched-ul-ing se-mes-ter semi-def-i-nite
semi-ho-mo-thet-ic set-up se-vere-ly side-step sov-er-eign spe-cious
spher-oid spher-oid-al star-tling star-tling-ly sta-tis-tics
sto-chas-tic straight-est strange-ness strat-a-gem strong-hold
sum-ma-ble symp-to-matic syn-chro-nous topo-graph-i-cal tra-vers-a-ble
tra-ver-sal tra-ver-sals treach-ery turn-around un-at-tached
un-err-ing-ly white-space wide-spread wing-spread wretch-ed
wretch-ed-ly Brown-ian Eng-lish Euler-ian Feb-ru-ary Gauss-ian
Grothen-dieck Hamil-ton-ian Her-mit-ian Jan-u-ary Japan-ese Kor-te-weg
Le-gendre Lip-schitz Lip-schitz-ian Mar-kov-ian Noe-ther-ian
No-vem-ber Rie-mann-ian Schwarz-schild Sep-tem-ber}

\def\leftheadtext#1{} \def\rightheadtext#1{}

\let\nopagenumber\p@gen@false \let\putpagenumber\p@gen@true
\let\pagefirst\nopagenumber \let\pagenext\putpagenumber

\else

\amsppttrue

\let\twelvepoint\relax \let\Twelvepoint\relax \let\putpagenumber\relax
\let\logo@\relax \let\pagefirst\firstpage@true \let\pagenext\firstpage@false
\def\nopagenumber{\let\f@li@ld\folio\def\folio{\global\let\folio\f@li@ld}}

\def\ftext#1{\footnotetext""{\vsk-.8>\nt #1}}

\def\m@R@f@[#1]{\Refs\nofrills{}\m@th\tenpoint
 {
 \s@ct{}{#1}}\wr@@c{\string\Refcd{#1}{\the\pageno}}
 \def\k@yf@##1{\hss[##1]\enspace} \let\keyformat\k@yf@
 \def\widest##1{\s@twd@\refindentwd{\tenpoint\k@yf@{##1}}}
 \let\Key\key \def\refin{\kern\refindentwd}}
\let\info\finalinfo \r@R@fs\Refs
\def\adjust#1{#1} \let\adjustend\relax
\let\adjustnext\adjust 

\fi

\outer\def\myRefs{\myR@fs} \r@st@re\proclaim
\def\bye{\par\vfill\supereject\cl@selbl\cl@secd\b@e} \r@endd@\b@e
\let\Cite\cite \let\Key\key \def\endpro{\par\endproclaim}
\let\d@c@\document \def\document{\d@c@\tenpoint}
\hyphenation{ortho-gon-al}

\newtoks\@@tp@t \@@tp@t\output
\output=\@ft@{\let\{\noexpand\the\@@tp@t}
\let\{\relax

\newif\ifVersion

\def\s@ct#1#2{\ifVersion
 \skip@\lastskip\ifdim\skip@<1.5\bls\vskip-\skip@\p@n@l{-200}\vsk.5>%
 \p@n@l{-200}\vsk.5>\p@n@l{-200}\vsk.5>\p@n@l{-200}\vsk-1.5>\else
 \p@n@l{-200}\fi\ifdim\skip@<.9\bls\vsk.9>\else
 \ifdim\skip@<1.5\bls\vskip\skip@\fi\fi
 \vtop{\twelvepoint\raggedright\bf\vp1\vsk->\vskip.16ex\s@twd@\parindent{#1}%
 \ifdim\parindent>\z@\adv\parindent.5em\fi\hang\textindent{#1}#2\strut}
 \else
 \p@sk@p{-200}{.8\bls}\vtop{\bf\s@twd@\parindent{#1}%
 \ifdim\parindent>\z@\adv\parindent.5em\fi\hang\textindent{#1}#2\strut}\fi
 \nointerlineskip\nobreak\vtop{\strut}\nobreak\vskip-.6\bls\nobreak}

\def\p@n@l#1{\ifnum#1=\z@\else\penalty#1\relax\fi}

\def\s@bs@ct#1#2{\ifVersion
 \skip@\lastskip\ifdim\skip@<1.5\bls\vskip-\skip@\p@n@l{-200}\vsk.5>%
 \p@n@l{-200}\vsk.5>\p@n@l{-200}\vsk.5>\p@n@l{-200}\vsk-1.5>\else
 \p@n@l{-200}\fi\ifdim\skip@<.9\bls\vsk.9>\else
 \ifdim\skip@<1.5\bls\vskip\skip@\fi\fi
 \vtop{\elevenpoint\raggedright\it\vp1\vsk->\vskip.16ex%
 \s@twd@\parindent{#1}\ifdim\parindent>\z@\adv\parindent.5em\fi
 \hang\textindent{#1}#2\strut}
 \else
 \p@sk@p{-200}{.6\bls}\vtop{\it\s@twd@\parindent{#1}%
 \ifdim\parindent>\z@\adv\parindent.5em\fi\hang\textindent{#1}#2\strut}\fi
 \nointerlineskip\nobreak\vtop{\strut}\nobreak\vskip-.8\bls\nobreak}

\def\gadv{\global\adv} \def\gad#1{\gadv#1\@ne} \def\gadneg#1{\gadv#1-\@ne}

\newcount\t@@n \t@@n=10 \newbox\testbox

\newcount\Sno \newcount\Lno \newcount\Fno

\def\pr@cl#1{\r@st@re\pr@c@\pr@c@{#1}\global\let\pr@c@\relax}

\def\tagg#1{\tag"\rlap{\rm(#1)}\kern.01\p@"}
\def\l@L#1{\l@bel{#1}L} \def\l@F#1{\l@bel{#1}F} \def\<#1>{\l@b@l{#1}F}
\def\Tag#1{\tag{\l@F{#1}}} \def\Tagg#1{\tagg{\l@F{#1}}}
\def\Rem{\demo{\sl Remark}} \def\Ex{\demo{\bf Example}}
\def\Pf#1.{\demo{Proof #1}} \def\epf{\qed\enddemo}
\def\Ap@x{Appendix}
\def\Appendix{\Sno=64 \t@@n\@ne \wr@@c{\string\Appencd}
 \def\sf@rm{\char\the\Sno} \def\sf@rm@{\Ap@x\space\sf@rm} \def\sf@rm@@{\Ap@x}
 \def\s@ct@n##1##2{\s@ct\empty{\setbox\z@\hbox{##1}\ifdim\wd\z@=\z@
 \if##2*\sf@rm@@\else\if##2.\sf@rm@@.\else##2\fi\fi\else
 \if##2*\sf@rm@\else\if##2.\sf@rm@.\else\sf@rm@.\enspace##2\fi\fi\fi}}}
\def\Appcd#1#2#3{\def\Ap@@{\hglue-\l@ftcd\Ap@x}\ifx\@ppl@ne\empty
 \def\l@@b{\@fwd@@{#1}{\space#1}{}}\if*#2\entcd{}{\Ap@@\l@@b}{#3}\else
 \if.#2\entcd{}{\Ap@@\l@@b.}{#3}\else\entcd{}{\Ap@@\l@@b.\enspace#2}{#3}\fi\fi
 \else\def\l@@b{\@fwd@@{#1}{\c@l@b{#1}}{}}\if*#2\entcd{\l@@b}{\Ap@x}{#3}\else
 \if.#2\entcd{\l@@b}{\Ap@x.}{#3}\else\entcd{\l@@b}{#2}{#3}\fi\fi\fi}

\let\s@ct@n\s@ct
\def\s@ct@@[#1]#2{\@ft@\xdef\csname @#1@S@\endcsname{\sf@rm}\wr@@x{}%
 \wr@@x{\string\labeldef{S}\space{\?#1@S?}\space{#1}}%
 {
 \s@ct@n{\sf@rm@}{#2}}\wr@@c{\string\Entcd{\?#1@S?}{#2}{\the\pageno}}}
\def\s@ct@#1{\wr@@x{}{
 \s@ct@n{\sf@rm@}{#1}}\wr@@c{\string\Entcd{\sf@rm}{#1}{\the\pageno}}}
\def\s@ct@e[#1]#2{\@ft@\xdef\csname @#1@S@\endcsname{\sf@rm}\wr@@x{}%
 \wr@@x{\string\labeldef{S}\space{\?#1@S?}\space{#1}}%
 {
 \s@ct@n\empty{#2}}\wr@@c{\string\Entcd{}{#2}{\the\pageno}}}
\def\s@cte#1{\wr@@x{}{
 \s@ct@n\empty{#1}}\wr@@c{\string\Entcd{}{#1}{\the\pageno}}}
\def\theSno#1#2{\dff\?#1@S?{#2}%
 \wr@@x{\string\labeldef{S}\space{#2}\space{#1}}\fi}

\newif\ifd@bn@\d@bn@true
\def\Section{\gad\Sno\ifd@bn@\Fno\z@\Lno\z@\fi\@fn@xt[\s@ct@@\s@ct@}
\def\section{\gad\Sno\ifd@bn@\Fno\z@\Lno\z@\fi\@fn@xt[\s@ct@e\s@cte}
\let\Sect\Section 
\def\subsection{\@fn@xt*\subs@ct@\subs@ct}
\def\subs@ct#1{{
 \s@bs@ct\empty{#1}}\wr@@c{\string\subcd{#1}{\the\pageno}}}
\def\subs@ct@*#1{\vsk->\vsk>{
 \s@bs@ct\empty{#1}}\wr@@c{\string\subcd{#1}{\the\pageno}}}
\let\subsect\subsection \def\Snodef#1{\Sno #1}

\def\l@b@l#1#2{\def\n@@{\csname #2no\endcsname}%
 \if*#1\gad\n@@ \@ft@\xdef\csname @#1@#2@\endcsname{\l@f@rm}\else\def\t@st{#1}%
 \ifx\t@st\empty\gad\n@@ \@ft@\xdef\csname @#1@#2@\endcsname{\l@f@rm}%
 \else\@ft@\ifx\csname @#1@#2@mark\endcsname\relax\gad\n@@
 \@ft@\xdef\csname @#1@#2@\endcsname{\l@f@rm}%
 \@ft@\gdef\csname @#1@#2@mark\endcsname{}%
 \wr@@x{\string\labeldef{#2}\space{\?#1@#2?}\space\ifnum\n@@<10 \space\fi{#1}}%
 \fi\fi\fi}
\def\labeldef#1#2#3{\dff\?#3@#1?{#2}}
\def\Labeldef#1#2#3{\dff\?#3@#1?{#2}\@ft@\gdef\csname @#3@#1@mark\endcsname{}}

\def\l@bel#1#2{\l@b@l{#1}{#2}\?#1@#2?}

\newcount\c@cite
\def\?#1?{\csname @#1@\endcsname}
\def\[{\@fn@xt:\c@t@sect\c@t@}
\def\c@t@#1]{{\c@cite\z@\@fwd@@{\?#1@L?}{\adv\c@cite1}{}%
 \@fwd@@{\?#1@F?}{\adv\c@cite1}{}\@fwd@@{\?#1?}{\adv\c@cite1}{}%
 \relax\ifnum\c@cite=\z@{\bf ???}\wrs@x{No label [#1]}\else
 \ifnum\c@cite=1\let\@@PS\relax\let\@@@\relax\else\let\@@PS\underbar
 \def\@@@{{\rm<}}\fi\@@PS{\?#1?\@@@\?#1@L?\@@@\?#1@F?}\fi}}
\def\(#1){{\rm(\c@t@#1])}}
\def\c@t@s@ct#1{\@fwd@@{\?#1@S?}{\?#1@S?\relax}%
 {{\bf ???}\wrs@x{No section label {#1}}}}
\def\c@t@sect:#1]{\c@t@s@ct{#1}} \let\SNo\c@t@s@ct

\newdimen\l@ftcd \newdimen\r@ghtcd \let\nlc\relax

\def\d@tt@d{\leaders\hbox to 1em{\kern.1em.\hfil}\hfill}
\def\entcd#1#2#3{\item{#1}{#2}\alb\kern.9em\hbox{}\kern-.9em\d@tt@d
 \kern-.36em{#3}\kern-\r@ghtcd\hbox{}\par}
\def\Entcd#1#2#3{\def\l@@b{\@fwd@@{#1}{\c@l@b{#1}}{}}\vsk.2>%
 \entcd{\l@@b}{#2}{#3}}
\def\subcd#1#2{{\adv\leftskip.333em\entcd{}{\it #1}{#2}}}
\def\Refcd#1#2{\def\t@@st{#1}\ifx\t@@st\empty\ifx\r@fl@ne\empty\relax\else
 \R@fcd{\r@fl@ne}{#2}\fi\else\R@fcd{#1}{#2}\fi}
\def\R@fcd#1#2{\sk@@p{.6\bls}\entcd{}{\hglue-\l@ftcd\bf #1}{#2}}
\def\Refline{\def\r@fl@ne} \def\Refempty{\let\r@fl@ne\empty}
\def\Appencd{\par\adv\leftskip-\l@ftcd\adv\rightskip-\r@ghtcd\@ppl@ne
 \adv\leftskip\l@ftcd\adv\rightskip\r@ghtcd\let\Entcd\Appcd}
\def\appline{\def\@ppl@ne} \def\Appempty{\let\@ppl@ne\empty}
\def\Appline#1{\def\@ppl@ne{\s@bs@ct{}{\bf#1}}}
\def\leftcd#1{\adv\leftskip-\l@ftcd\s@twd@\l@ftcd{\c@l@b{#1}\enspace}
 \adv\leftskip\l@ftcd}
\def\rightcd#1{\adv\rightskip-\r@ghtcd\s@twd@\r@ghtcd{#1\enspace}
 \adv\rightskip\r@ghtcd}
\def\C@nt{Contents} \def\Ap@s{Appendices} \def\R@fcs{References}
\def\contents{\@fn@xt*\cont@@\cont@}
\def\cont@{\@fn@xt[\cnt@{\cnt@[\C@nt]}}
\def\cont@@*{\@fn@xt[\cnt@@{\cnt@@[\C@nt]}}
\def\cnt@[#1]{\c@nt@{M}{#1}{44}{\s@bs@ct{}{\bf\Ap@s}}}
\def\cnt@@[#1]{\c@nt@{M}{#1}{44}{}}
\def\endco{\par\penalty-500\vsk>\vskip\z@\endgroup}
\def\readcd{\@np@t{\jobname.cd}}
\def\Cde{\@fn@xt*\Cde@@\Cde@}
\def\Cde@{\@fn@xt[\Cd@{\Cd@[\C@nt]}}
\def\Cde@@*{\@fn@xt[\Cd@@{\Cd@@[\C@nt]}}
\def\Cd@[#1]{\cnt@[#1]\readcd\endco}
\def\Cd@@[#1]{\cnt@@[#1]\readcd\endco}
\def\contlabeldef{\def\c@l@b}

\long\def\c@nt@#1#2#3#4{\s@twd@\l@ftcd{\c@l@b{#1}\enspace}
 \s@twd@\r@ghtcd{#3\enspace}\adv\r@ghtcd1.333em
 \def\@ppl@ne{#4}\def\r@fl@ne{\R@fcs}\s@ct{}{#2}\B@gr@\parindent\z@\let\nlc\nl
 \let\nl\relax\parskip.2\bls\adv\leftskip\l@ftcd\adv\rightskip\r@ghtcd}

\def\writecd{\immediate\openout\@@cd\jobname.cd \def\wr@@c{\write\@@cd}
 \def\cl@secd{\immediate\write\@@cd{\string\endinput}\immediate\closeout\@@cd}
 \def\closecd{\cl@secd\global\let\cl@secd\relax}}
\let\cl@secd\relax \def\wr@@c#1{} \let\closecd\relax

\def\dff{\@ft@\d@f} \def\d@f{\@ft@\def}
\def\edff{\@ft@\ed@f} \def\ed@f{\@ft@\edef}
\def\defi#1#2{\def#1{#2}\wr@@x{\string\def\string#1{#2}}}

\def\qed{\hbox{}\nobreak\hfill\nobreak{\m@th$\,\square$}}
\def\back#1 {\strut\kern-.33em #1\enspace\ignore} 
\def\Text#1{\crcr\noalign{\alb\vsk>\normalbaselines\vsk->\vbox{\nt #1\strut}%
 \nobreak\nointerlineskip\vbox{\strut}\nobreak\vsk->\nobreak}}
\def\inText#1{\crcr\noalign{\penalty\postdisplaypenalty\vskip\belowdisplayskip
 \vbox{\normalbaselines\noindent#1}\penalty\predisplaypenalty
 \vskip\abovedisplayskip}}

\def\hcor#1{\advance\hoffset by #1}
\def\vcor#1{\advance\voffset by #1}
\let\bls\baselineskip \let\ignore\ignorespaces
\ifx\ic@\undefined \let\ic@\/\fi
\def\vsk#1>{\vskip#1\bls} \let\adv\advance
\def\vv#1>{\vadjust{\vsk#1>}\ignore}
\def\vvn#1>{\vadjust{\nobreak\vsk#1>\nobreak}\ignore}
\def\vvv#1>{\vskip\z@\vsk#1>\nt\ignore}
\def\vvgood{\vadjust{\penalty-500}} \def\vgood{\strut\vvn->\vvgood\nl\ignore}
\def\Par{\vsk.5>} \def\setparindent{\edef\Parindent{\the\parindent}}
\def\Type{\vsk.5>\bgroup\parindent\z@\tt\rightskip\z@ plus1em minus1em%
 \spaceskip.3333em \xspaceskip.5em\relax}
\def\endType{\vsk.5>\egroup\nt} 

  \let\dollar\$ \let\ampersand\&
\let\sss\scriptscriptstyle  
\let\vp\vphantom \let\hp\hphantom \let\nt\noindent
\let\cline\centerline \let\lline\leftline \let\rline\rightline
\def\nn#1>{\noalign{\vskip#1\p@@}} \def\NN#1>{\openup#1\p@@}
\def\cnn#1>{\noalign{\vsk#1>}}
 
\let\Lim\lim \def\lim{\Lim\limits} \let\Sum\sum \def\sum{\Sum\limits}
\def\Plus{\bigoplus\limits} 
\let\Prod\prod \def\prod{\Prod\limits} \let\Int\int \def\int{\Int\limits}

\def\tsum{\mathop{\tsize\Sum}\limits} 
\def\tprod{\mathop{\tsize\Prod}\limits}
\def\&{.\kern.1em} \def\>{{\!\;}} \def\]{{\!\!\;}} \def\){\>\]} \def\}{\]\]}
\def\nl{\leavevmode\hfill\break} \def\~{\leavevmode\@fn@xt~\m@n@s\@md@sh}
\def\m@n@s~{\raise.15ex\mbox{-}} \def\@md@sh{\raise.13ex\hbox{--}}
\let\procent\% \def\%#1{\ifmmode\mathop{#1}\limits\else\procent#1\fi}
\let\@ml@t\" \def\"#1{\ifmmode ^{(#1)}\else\@ml@t#1\fi}
\let\@c@t@\' \def\'#1{\ifmmode _{(#1)}\else\@c@t@#1\fi}
\let\colon\: \def\:{^{\vp|}}

\let\texspace\ \def\ {\ifmmode\alb\fi\texspace}

\let\n@wp@ge\newpage \def\newpage{\endgraf\n@wp@ge}
\let\=\m@th \def\mbox#1{\hbox{\m@th$#1$}}
\def\mtext#1{\text{\m@th$#1$}} \def\^#1{\text{\m@th#1}}
\def\Line#1{\kern-.5\hsize\line{\m@th$\dsize#1$}\kern-.5\hsize}
\def\Lline#1{\kern-.5\hsize\lline{\m@th$\dsize#1$}\kern-.5\hsize}
\def\Cline#1{\kern-.5\hsize\cline{\m@th$\dsize#1$}\kern-.5\hsize}
\def\Rline#1{\kern-.5\hsize\rline{\m@th$\dsize#1$}\kern-.5\hsize}

\def\Ll@p#1{\llap{\m@th$#1$}} \def\Rl@p#1{\rlap{\m@th$#1$}}
\def\clap#1{\llap{#1\hss}} \def\Cl@p#1{\llap{\m@th$#1$\hss}}
\def\Llap#1{\mathchoice{\Ll@p{\dsize#1}}{\Ll@p{\tsize#1}}{\Ll@p{\ssize#1}}%
 {\Ll@p{\sss#1}}}
\def\Clap#1{\mathchoice{\Cl@p{\dsize#1}}{\Cl@p{\tsize#1}}{\Cl@p{\ssize#1}}%
 {\Cl@p{\sss#1}}}
\def\Rlap#1{\mathchoice{\Rl@p{\dsize#1}}{\Rl@p{\tsize#1}}{\Rl@p{\ssize#1}}%
 {\Rl@p{\sss#1}}}
 
\def\LRtph#1#2{\setbox\z@\hbox{#1}\dimen\z@\wd\z@\hbox{\hbox to\dimen\z@{#2}}}
\def\LRph#1#2{\LRtph{\m@th$#1$}{\m@th$#2$}}

\def\LLph#1#2{\LRph{#1}{\hss#2}} \def\RRph#1#2{\LRph{#1}{#2\hss}}

\def\Lph#1#2{\mathchoice{\LLph{\dsize#1}{\dsize#2}}{\LLph{\tsize#1}{\tsize#2}}
 {\LLph{\ssize#1}{\ssize#2}}{\LLph{\sss#1}{\sss#2}}}
\def\Rph#1#2{\mathchoice{\RRph{\dsize#1}{\dsize#2}}{\RRph{\tsize#1}{\tsize#2}}
 {\RRph{\ssize#1}{\ssize#2}}{\RRph{\sss#1}{\sss#2}}}
\def\Lto#1{\setbox\z@\mbox{\tsize{#1}}%
 \mathrel{\mathop{\hbox to\wd\z@{\rightarrowfill}}\limits#1}}
\def\Lgets#1{\setbox\z@\mbox{\tsize{#1}}%
 \mathrel{\mathop{\hbox to\wd\z@{\leftarrowfill}}\limits#1}}
 \def\vpp#1{{\vp{\big]}}_{#1}}
\def\lbc{\mathopen{[\![}} \def\rbc{\mathclose{]\!]}}
 
\def\LBc{\mathopen{\Big[\!\}\Big[}} \def\RBc{\mathclose{\Big]\!\}\Big]}}

\let\alb\allowbreak 
\def\ald{\noalign{\alb}} \let\alds\allowdisplaybreaks

\let\o\circ \let\x\times \let\ox\otimes 
\let\sub\subset  \let\tabs\+
\let\le\leqslant \let\ge\geqslant
\let\der\partial \let\8\infty \let\*\star
\let\bra\langle \let\ket\rangle
\let\lto\longrightarrow 
\let\map\mapsto  \let\hto\hookrightarrow
 
 \def\vert{\ |\ } \def\nin{\not\in}
\let\Empty\varnothing 
 
\let\lb\lbrace \let\rb\rbrace
\let\lf\lfloor \let\rf\rfloor
\let\trileft\triangleleft 

\def\lsym#1{#1\alb\ldots\relax#1\alb}
\def\lc{\lsym,}  \def\lx{\lsym\x} \def\lox{\lsym\ox}
\def\llc{\,,\alb\ {\ldots\ ,}\alb\ }
\def\Re{\mathop{\roman{Re}\>}} \def\Im{\mathop{\roman{Im}\>}}
\def\End{\mathop{\roman{End}\>}}

\def\Res{\mathop{\roman{Res}\>}\limits}

\def\id{\roman{id}}  \def\for{\text{for }\,}
\def\1{^{-1}} \def\_#1{_{\Rlap{#1}}}
\def\vst#1{{\lower1.9\p@@\mbox{\bigr|_{\raise.5\p@@\mbox{\ssize#1}}}}}
\def\vrp#1:#2>{{\vrule height#1 depth#2 width\z@}}
\def\vru#1>{\vrp#1:\z@>} \def\vrd#1>{\vrp\z@:#1>}
\def\qqq{\qquad\quad} 
\def\sscr#1{\raise.3ex\mbox{\sss#1}} \def\@@PS{\bold{OOPS!!!}}
\def\ono{\bigl(1+o(1)\bigr)} 

\def\intcl{\mathop
 {\Rlap{\raise.3ex\mbox{\kern.12em\curvearrowleft}}\int}\limits}
\def\intcr{\mathop
 {\Rlap{\raise.3ex\mbox{\kern.24em\curvearrowright}}\int}\limits}

\def\pms{\raise.25ex\mbox{\ssize\pm}\>}
\def\mps{\raise.25ex\mbox{\ssize\mp}\>}
\def\pss{{\sscr+}} \def\mss{{\sscr-}}
\def\pmss{{\sscr\pm}} 

\let\al\alpha
\let\bt\beta
\let\gm\gamma \let\Gm\Gamma \let\Gms\varGamma
\let\dl\delta \let\Dl\Delta 
\let\epe\epsilon \let\eps\varepsilon \let\epsilon\eps

\let\zt\zeta
\let\tht\theta \let\Tht\Theta
\let\thi\vartheta 

\let\ka\kappa
\let\la\lambda \let\La\Lambda

\let\si\sigma 
 \let\Sig\varSigma
\let\pho\phi \let\phi\varphi \let\Pho\varPhi
\let\Pss\varPsi

\let\om\omega \let\Om\Omega 

\def\C{\Bbb C}
\def\R{\Bbb R}
\def\Z{\Bbb Z}

\def\AA{\Bbb A}
\def\BB{\Bbb B}

\def\FF{\Bbb F}
\def\II{\Bbb I}
\def\SS{\Bbb S}

\def\UU{\Bbb U}

\def\YY{\Bbb Y}

\def\Zp{\Z_{\ge 0}} \def\Zpp{\Z_{>0}}
\def\Zn{\Z_{\le 0}} \def\Znn{\Z_{<0}}

\def\difl/{differential} \def\dif/{difference}
\def\cf.{cf.\ \ignore} \def\Cf.{Cf.\ \ignore}
\def\egv/{eigenvector} \def\eva/{eigenvalue} \def\eq/{equation}
\def\lhs/{the left hand side} \def\rhs/{the right hand side}
\def\Lhs/{The left hand side} \def\Rhs/{The right hand side}
\def\gby/{generated by} \def\wrt/{with respect to} \def\st/{such that}
\def\resp/{respectively} \def\off/{offdiagonal} \def\wt/{weight}
\def\pol/{polynomial} \def\rat/{rational} \def\tri/{trigonometric}
\def\fn/{function} \def\var/{variable} \def\raf/{\rat/ \fn/}
\def\inv/{invariant} \def\hol/{holomorphic} \def\hof/{\hol/ \fn/}
\def\mer/{meromorphic} \def\mef/{\mer/ \fn/} \def\mult/{multiplicity}
\def\sym/{symmetric} \def\perm/{permutation} \def\fd/{finite-dimensional}
\def\rep/{representation} \def\irr/{irreducible} \def\irrep/{\irr/ \rep/}
\def\hom/{homomorphism} \def\aut/{automorphism} \def\iso/{isomorphism}
\def\lex/{lexicographical} \def\as/{asymptotic} \def\asex/{\as/ expansion}
\def\ndeg/{nondegenerate} \def\neib/{neighbourhood} \def\deq/{\dif/ \eq/}
\def\hw/{highest \wt/} \def\gv/{generating vector} \def\eqv/{equivalent}
\def\msd/{method of steepest descend} \def\pd/{pairwise distinct}
\def\wlg/{without loss of generality} \def\Wlg/{Without loss of generality}
\def\onedim/{one-dimensional} \def\qcl/{quasiclassical}
\def\hgeom/{hyper\-geometric} \def\hint/{\hgeom/ integral}
\def\hwm/{\hw/ module} \def\emod/{evaluation module} \def\Vmod/{Verma module}
\def\symg/{\sym/ group} \def\sol/{solution} \def\eval/{evaluation}
\def\anf/{analytic \fn/} \def\anco/{analytic continuation}
\def\qg/{quantum group} \def\qaff/{quantum affine algebra}

\def\Rm/{\^{$R$-}matrix} \def\Rms/{\^{$R$-}matrices} \def\YB/{Yang-Baxter \eq/}
\def\Ba/{Bethe ansatz} \def\Bv/{Bethe vector} \def\Bae/{\Ba/ \eq/}
\def\KZv/{Knizh\-nik-Zamo\-lod\-chi\-kov} \def\KZvB/{\KZv/-Bernard}
\def\KZ/{{\sl KZ\/}} \def\qKZ/{{\sl qKZ\/}}
\def\KZB/{{\sl KZB\/}} \def\qKZB/{{\sl qKZB\/}}
\def\qKZo/{\qKZ/ operator} \def\qKZc/{\qKZ/ connection}
\def\KZe/{\KZ/ \eq/} \def\qKZe/{\qKZ/ \eq/} \def\qKZBe/{\qKZB/ \eq/}

\def\h@ph{\discretionary{}{}{-}} \def\$#1$-{\,\^{$#1$}\h@ph}

\def\TFT/{Research Insitute for Theoretical Physics}
\def\HY/{University of Helsinki} \def\AoF/{the Academy of Finland}
\def\CNRS/{Supported in part by MAE\~MICECO\~CNRS Fellowship}
\def\LPT/{Laboratoire de Physique Th\'eorique ENSLAPP}
\def\ENSLyon/{\'Ecole Normale Sup\'erieure de Lyon}
\def\LPTaddr/{46, All\'ee d'Italie, 69364 Lyon Cedex 07, France}
\def\enslapp/{URA 14\~36 du CNRS, associ\'ee \`a l'E.N.S.\ de Lyon,
au LAPP d'Annecy et \`a l'Universit\`e de Savoie}
\def\ensemail/{vtarasov\@ enslapp.ens-lyon.fr}
\def\DMS/{Department of Mathematics, Faculty of Science}
\def\DMO/{\DMS/, Osaka University}
\def\DMOaddr/{Toyonaka, Osaka 560, Japan}
\def\dmoemail/{vt\@ math.sci.osaka-u.ac.jp}
\def\SPb/{St\&Peters\-burg}
\def\home/{\SPb/ Branch of Steklov Mathematical Institute}
\def\homeaddr/{Fontanka 27, \SPb/ \,191011, Russia}
\def\homemail/{vt\@ pdmi.ras.ru}
\def\absence/{On leave of absence from \home/}
\def\UNC/{Department of Mathematics, University of North Carolina}
\def\ChH/{Chapel Hill}
\def\UNCaddr/{\ChH/, NC 27599, USA} \def\avemail/{av\@ math.unc.edu}
\def\grant/{NSF grant DMS\~9501290}	
\def\Grant/{Supported in part by \grant/}

\def\Aomoto/{K\&Aomoto}
\def\Dri/{V\]\&G\&Drin\-feld}
\def\Fadd/{L\&D\&Fad\-deev}
\def\Feld/{G\&Felder}
\def\Fre/{I\&B\&Fren\-kel}
\def\Gustaf/{R\&A\&Gustafson}
\def\Kazh/{D\&Kazhdan} \def\Kir/{A\&N\&Kiril\-lov}
\def\Kor/{V\]\&E\&Kore\-pin}
\def\Lusz/{G\&Lusztig}
\def\MN/{M\&Naza\-rov}
\def\Resh/{N\&Reshe\-ti\-khin} \def\Reshy/{N\&\]Yu\&Reshe\-ti\-khin}
\def\SchV/{V\]\&\]V\]\&Schecht\-man} \def\Sch/{V\]\&Schecht\-man}
\def\Skl/{E\&K\&Sklya\-nin}
\def\Smirn/{F\]\&Smirnov} \def\Smirnov/{F\]\&A\&Smirnov}
\def\Takh/{L\&A\&Takh\-tajan}
\def\VT/{V\]\&Ta\-ra\-sov} \def\VoT/{V\]\&O\&Ta\-ra\-sov}
\def\Varch/{A\&\]Var\-chenko} \def\Varn/{A\&N\&\]Var\-chenko}

\def\AMS/{Amer.\ Math.\ Society}
\def\CMP/{Comm.\ Math.\ Phys.{}}
\def\DMJ/{Duke.\ Math.\ J.{}}
\def\Inv/{Invent.\ Math.{}} 
\def\IMRN/{Int.\ Math.\ Res.\ Notices}
\def\JPA/{J.\ Phys.\ A{}}
\def\JSM/{J.\ Soviet\ Math.{}}
\def\LMP/{Lett.\ Math.\ Phys.{}}
\def\LMJ/{Leningrad Math.\ J.{}}
\def\LpMJ/{\SPb/ Math.\ J.{}}
\def\SIAM/{SIAM J.\ Math.\ Anal.{}}
\def\SMNS/{Selecta Math., New Series}
\def\TMP/{Theor.\ Math.\ Phys.{}}
\def\ZNS/{Zap.\ nauch.\ semin. LOMI}

\def\ASMP/{Advanced Series in Math.\ Phys.{}}

\def\AMSa/{AMS \publaddr Providence}
\def\Birk/{Birkh\"auser}
\def\CUP/{Cambridge University Press} \def\CUPa/{\CUP/ \publaddr Cambridge}
\def\Spri/{Springer-Verlag} \def\Spria/{\Spri/ \publaddr Berlin}
\def\WS/{World Scientific} \def\WSa/{\WS/ \publaddr Singapore}

\newbox\lefthbox \newbox\righthbox

\let\sectsep. \let\labelsep. \let\contsep. \let\labelspace\relax
\let\sectpre\relax \let\contpre\relax
\def\sf@rm{\the\Sno} \def\sf@rm@{\sectpre\sf@rm\sectsep}
\def\c@l@b#1{\contpre#1\contsep}
\def\l@f@rm{\ifd@bn@\sf@rm\labelsep\fi\labelspace\the\n@@}

\def\sectformdef{\def\sf@rm}

\let\DoubleNum\d@bn@true \let\SingleNum\d@bn@false

\def\NoNewNum{\let\writeldf\relax\def\l@b@l##1##2{\if*##1%
 \@ft@\xdef\csname @##1@##2@\endcsname{\mbox{*{*}*}}\fi}}
\def\NoNewTime{\def\todaydef##1{\def\today{##1}}
 \def\nowtimedef##1{\def\nowtime{##1}}}
\def\NoInput{\let\Input\input\let\writeldf\relax}
\def\Fixed{\NoNewTime\NoInput}

\tenpoint

\csname beta.def\endcsname

\Fixed

\expandafter\ifx\csname swan.def\endcsname\relax \else\endinput\fi
\expandafter\edef\csname swan.def\endcsname{%
 \catcode`\noexpand\@=\the\catcode`\@\space}
\catcode`\@=11

\font\Brm=cmr12 scaled 1200
\font\Bit=cmti12 scaled 1200

\ifx\twelvemsa\undefined\font\twelvemsa=msam10 scaled 1200\fi
\ifMag\let\Msam\relax\else\def\Msam{\textfont\msafam\twelvemsa}\fi

\newif\ifPcd \let\noPcd\Pcdfalse \let\makePcd\Pcdtrue
\def\NoPcde{\let\noPcd\relax \let\makePcd\relax}

\def\Th#1{\pr@cl{(\l@F{#1}) Theorem}\ignore}
\def\Lm#1{\pr@cl{(\l@F{#1}) Lemma}\ignore}
\def\Cr#1{\pr@cl{(\l@F{#1}) Corollary}\ignore}
\def\Df#1{\pr@cl{(\l@F{#1}) Definition}\ignore}
\def\Cj#1{\pr@cl{(\l@F{#1}) Conjecture}\ignore}
\def\Prop#1{\pr@cl{(\l@F{#1}) Proposition}\ignore}
\def\Pf#1.{\demo{\let\{\relax Proof #1}\ifPcd\def\t@st@{#1}\ifx\t@st@\empty
 \else\wr@@c{\string\subcd{Proof #1}{\the\pageno}}\fi\fi\ignore}

\let\grt\gg

\def\q{{q\1}} \let\one\id

\def\pb{\bar p}

\def\wb{\bar w}
\def\yb{\bar y}
\def\zb{\bar z}

\def\Cbar{\Rlap{\kern1.2\p@\overline{\vp{\C}\kern5\p@}}\C}
\def\Db{\Rlap{\kern3\p@\overline{\vp{D}\kern4\p@}}D}
\def\Fb{\Rlap{\kern3\p@\overline{\vp{\F}\kern5\p@}}\F}
\def\Hb{\Rlap{\kern2.5\p@\overline{\vp{\H}\kern5\p@}}\H}

\def\Ub{\>\overline{\}\Ud\!}\,}
\def\UUb{\>\overline{\}\UU\}}\>}
\def\Wb{\,\overline{\!W\}}\>}

\def\nb{\bold n}
\def\P{\bold P}
\def\S{\bold S}

\def\F{\Cal F}
\def\H{\Cal H}

\def\Rc{\Cal R}
\def\Zc{\Cal Z}

\def\Ich{\check I}
\def\Iic{\check\Ii}

\def\Bg{\frak B}
\def\eg{\frak e}
\def\Fg{\frak F}

\def\g{\frak g}
\def\ig{\frak i}
\def\jg{\frak j}

\def\lg{\frak l}
\def\mg{\frak m}

\def\A{\roman A}
\def\B{\roman B}
\def\CC{\roman C}
\def\d{\roman d}

\def\Ff{\roman F}
\def\I{\roman I}

\def\Ud{\roman U}

\def\ZZ{\roman Z}

\def\gs{\slanted g}

\def\IIt{\widetilde\II}
\def\Iti{\widetilde I}
\def\Is{\widetilde\I}
  
\def\Lti{\)\widetilde{\]L\)}\]}
\def\Vt{\)\widetilde{\]V\)}\]}
\def\Ti{\widetilde T}

\def\wti{\tilde w}
\def\Pht{\widetilde\Phi}
\def\Pti{\widetilde\Psi}

\def\HOne{H^1(\Omb,\nabla)} \def\Hone{H_1(\Omb,\nabla)}

\def\Aid{\AA_{\>\id}}

\def\Fq{\F_{\!q}} \def\Fqs{\Fq\sing}
\def\Fol{\Fo[\)l\>]} \def\Foll{\Fo[\)\ell\,]}
\def\Fql{\Fq[\)l\>]} \def\Fqll{\Fq[\)\ell\,]} \def\Fqsl{\Fqs[\)l\>]}

\def\Fsi{\Fun_{\vp|\!\}\sss\Sig}}

\let\Fun\Fwh \let\Fo\F

\def\Has{H_{\}\sss A}}

\def\Ii{I^{\sscr\o}}

\def\Omb{\Om^{\sscr\bullet}}
\def\qH{q^{-H}}

\def\Go{\mathchoice{{\%{\vru1.23ex>\smash{G}}^{\>\smash{\sss\o}\}}}}
 {{\%{\vru1.23ex>\smash{G}}^{\>\smash{\sss\o}\}}}}
 {{\%{\vru.92ex>\smash{G}}^{\>\smash{\sss\o}}}}{\@@PS}}
\def\Mo{\mathchoice{{\%{\vru1.19ex>\smash{M}\}}^{\,\smash{\sss\o}\!}}}
 {{\%{\vru1.19ex>\smash{M}\}}^{\,\smash{\sss\o}\!}}}
 {{\%{\vru.893ex>\smash{M}\}}^{\,\smash{\sss\o}\!}}}{\@@PS}}
\def\UUa{\mathchoice{{\%{\vru1.2ex>\smash{\UU}}^{\)\smash{\sss*}\]}}}
 {{\%{\vru1.2ex>\smash{\UU}}^{\)\smash{\sss*}\]}}}
 {{\%{\vru.9ex>\smash{\UU}}^{\)\smash{\sss*}\]}}}{\@@PS}}
\def\UUo{\mathchoice{{\%{\vru1.2ex>\smash{\UU}}^{\)\smash{\sss\o}\]}}}
 {{\%{\vru1.2ex>\smash{\UU}}^{\)\smash{\sss\o}\]}}}
 {{\%{\vru.9ex>\smash{\UU}}^{\)\smash{\sss\o}\]}}}{\@@PS}}
\def\Wo{\mathchoice{{\%{\vru1.2ex>\smash{W}\}}^{\>\smash{\sss\o}\}}}}
 {{\%{\vru1.2ex>\smash{W}\}}^{\>\smash{\sss\o}\}}}}
 {{\%{\vru.9ex>\smash{W}\}}^{\>\smash{\sss\o}\}}}}{\@@PS}}

\def\pii{\pi i}

\def\Cn{\C^{\>n}} \def\Cl{\C^{\,\ell}} \def\Cln{\C^{\,\ell+n}}
\def\Con{\C^{\>1+n}} \def\Cll{\C^{\,\ell-1}}

\def\oxcc[#1]{\mathrel{\ox{\vp|}_{\)\C\)[#1]}}}
 
\def\Pn{\P_{\!n}} \def\Sl{\S^\ell}
 \def\Zpn{\Zp^n}
 \def\Zln{\Zc_\ell^n} \def\Zlm{\Zc_\ell^{n-1}}
\def\Zll{\Zc_{\ell-1}^n}

\def\Ap{A^{\raise.4ex\mbox{\sss\]+}}}
\def\Am{A^{\raise.4ex\mbox{\sss\]-}}}
\def\Apm{A^{\raise.4ex\mbox{\sss\}\pm}}}
\def\Cp{C^{\raise.4ex\mbox{\sss\>+}}}
\def\Cm{C^{\raise.4ex\mbox{\sss\>-}}}
\def\Cpm{C^{\raise.4ex\mbox{\sss\>\pm}}}
\def\CCp{\CC^\pss} \def\CCm{\CC^\mss} \def\CCpm{\CC^\pmss}

\def\DWwi{\det\bigl[I(W_\lg,w_\mg)\bigr]}
\def\DWwo{\det\bigl[I(\Wo_\lg,w_\mg)\bigr]}
\def\DWwj{\det\bigl[I(W_l,w_m)\bigr]_{l,m=1}^n}
\def\DWWJ{\det\bigl[I(W'_l,w'_m)\bigr]_{l,m=1}^n}
\def\DWwp{\det\bigl[I(\Wo_l,w_m)\bigr]_{l,m=1}^{n-1}}
\def\l@min{\limits_{\sss n-1}}
\def\DIo{\det\l@min\bigl[I_k(w_l)\bigr]}
\def\Dwwo{\det\l@min\bigl[\Iti(\wb_l,w_m)\bigr]}
\def\Dgpo{\det\l@min\bigl[I(\gp_l,w_m)\bigr]}
\def\Dgmo{\det\l@min\bigl[I(\gmm_l,w_m)\bigr]}
\def\Dgpmo{\det\l@min\bigl[I(\gpm_l,w_m)\bigr]}
\def\Dgpi{\det\l@min\bigl[I(\gsp_{ln},w_m)\bigr]}
\def\Dgmi{\det\l@min\bigl[I(\gsm_{ln},w_m)\bigr]}
\def\Dgpmi{\det\l@min\bigl[I(\gspm_{ln},w_m)\bigr]}
\def\Det#1{\mathop{\vp e\smash\det}\limits^{\,\sss\"{#1}\!}}

\def\gp{g^{\sss+}} \def\gmm{g^{\sss-}} \def\gpm{g^{\sss\pm}}
 
\def\gsp{\gs^{\sss+}} \def\gsm{\gs^{\sss-}} \def\gspm{\gs^{\sss\pm}}

\def\ZZp{\ZZ^\pss} \def\ZZm{\ZZ^\mss} 
\def\ZZps{\ZZ_{\sss+}} \def\ZZms{\ZZ_{\sss-}}

\def\npZ{\lb\)1\lc\ell\)\rb\cap p\)\Z=\Empty}

\def\6{{\>\overline{\}\8\!}\,}}
\def\l@inf{\lower.21ex\mbox{\ssize\8}} \def\9{_{\kern-.02em\l@inf}}
\def\l@infb{\lower.21ex\mbox{\ssize\6}} \def\0{_{\kern-.02em\l@infb}}

\def\5#1{^{!#1}} \def\7#1{_{!#1}}

\def\Dlq{\Dl{\vp{\big(}}^{\!q}}

\def\smm{\sum_{m=1}} \def\tsmm{\!\tsum_{m=1}}
\def\sun{\smm^n} \def\tsun{\tsmm^n\!}
 
\def\smk{\sum_{k=0}} 
 
\def\sms{\sum_{s=1}} 
\def\smsi{\sms^\8} 

\def\prm{\prod_{m=1}} \def\tprm{\tprod_{m=1}}
\def\pron{\prm^n} \def\tpron{\tprm^n}
\def\prmn{\prm^{n-1}} 
\def\plmn{\prod_{1\le l<m\le n}}
\def\pmln{\prod_{1\le m<l\le n}}
\def\prlm{\prod_{1\le l<m}}	

\def\pros{\prod_{s=0}^{\ell-1}} \def\pral{\prod_{a=1}^\ell}
\def\prab{\prod_{1\le a<b\le\ell}\!\!}

\def\mn{m=1\lc n}  \def\kon{k=0\lc n} \def\lcn{l=1\lc n}
 
 \def\lan{\la_1\lc\la_n} \def\Lan{\La_1\lc\La_n}
\def\phin{\phi_1\lc\phi_n} \def\Dtel{Dt_1\lsym\wedge Dt_\ell}
 \def\yn{y_1\lc y_n} \def\zn{z_1\lc z_n}
\def\zmn{z_1\lc z_m\]+p\lc z_n} \def\znm{z_1\lc p\)z_m\lc z_n}

\def\yz@{\!\pron\! y_m/z_m}
\def\yz{\mathchoice{{\tsize\yz@}}{\yz@}{\yz@}{\yz@}}
\def\zy@{\!\pron\! z_m/y_m} 
\def\zy{\mathchoice{{\tsize\zy@}}{\zy@}{\zy@}{\zy@}}
\def\yzb@{\!\pron\!\yb_m/\zb_m}
\def\yzb{\mathchoice{{\tsize\yzb@}}{\yzb@}{\yzb@}{\yzb@}}
\def\zyb@{\!\pron\!\zb_m/\yb_m} 
\def\zyb{\mathchoice{{\tsize\zyb@}}{\zyb@}{\zyb@}{\zyb@}}
\def\aell{a=1\lc\ell} \def\tell{t_1\lc t_\ell} 
\def\tall{t_1\lc t_{a+1},t_a\lc t_\ell}
\def\sell{s_1\lc s_\ell} \def\uell{u_1\lc u_\ell} \def\xell{x_1\lc x_\ell}
\def\xill{\xi_1\lc\xi_\ell} 
\def\ztn{\zt_1\lc\zt_n} \def\qLan{q^{\La_1}\lc q^{\La_n}}
 \def\lge{\lg^{\)!}}

\def\PV{P_{V_1V_2}\:} \def\RV{R_{V_1V_2}} \def\RVx{\RV\:(x)}
\def\RVz{\RV^q(\zt)} \def\RVo{\RV^q(0)}
\def\VV{V_1\ox V_2} \def\EVV{\End(\VV)}
\def\Vn{V_1\lc V_n}  \def\Vqn{V^q_1\lc V^q_{n{\vp1}}}
\def\Vox{V_1\lox V_n} \def\EV{\End(\Vox)}
\def\Vax{(V_1\lox V_n)^*}
\def\Vqx{V^q_1\lox V^q_{n{\vp1}}}

\def\Voxt{V_{\tau_1}\lox V_{\tau_n}}
\def\Vaxt{(V_{\tau_1}\lox V_{\tau_n})^*}
\def\Vqxt{V^q_{\tau_1}\lox V^q_{\tau_n}}
\def\Vaxtm{(V_{\tau_1}\lox V_{\tau_{m+1}}\ox V_{\tau_m}\lox V_{\tau_n})^*}
\def\Vqxtm{V^q_{\tau_1}\lox V^q_{\tau_{m+1}}\ox V^q_{\tau_m}\lox V^q_{\tau_n}}
\def\Voxzt{V_{\tau_1}(z_{\tau_1})\lox V_{\tau_n}(z_{\tau_n})}
\def\Vaxzt{\bigl(V_{\tau_1}(z_{\tau_1})\lox V_{\tau_n}(z_{\tau_n})\bigr)^*}
\def\Vqxzt{V^q_{\tau_1}(\zt_{\tau_1})\lox V^q_{\tau_n}(\zt_{\tau_n})}
\def\Vqxztt{V^q_{\tau'_1}(\zt_{\tau'_1})\lox V^q_{\tau'_n}(\zt_{\tau'_n})}
\def\elli{_\ell\:}
\def\Vl{(\Vox)\elli} \def\Vql{(\Vqx)\elli} \def\Val{\Vax_\ell}
\def\Vlt{(\Voxt)\elli} \def\Vqlt{(\Vqxt)\elli} \def\Valt{\Vaxt_\ell}
\def\Vaxtt{(V_{\tau'_1}\lox V_{\tau'_n})^*} \def\Valtt{\Vaxtt_\ell}
\def\Vqxtt{V^q_{\tau'_1}\lox V^q_{\tau'_n}} \def\Vqltt{(\Vqxtt)\elli}
 \def\FVal{F\Vax_{\ell-1}}
\def\FValt{F\Vaxt_{\ell-1}}
\def\vn{v_1\lc v_n} \def\vqn{v^q_1\lc v^q_{n{\vp1}}}
\def\vox{v_1\lox v_n} 
\def\voxt{v_{\tau_1}\lox v_{\tau_n}}
\def\vqxt{v^q_{\tau_1}\lox v^q_{\tau_n}}
\def\Fv{F^{\lg_1}v_1\lox F^{\lg_n}v_n} \def\Fva{(\Fv)^*}
\def\Fvq{{F_{\!q}}^{\!\]\lg_1}v^q_1\lox{F_{\!q}}^{\!\]\lg_n}v^q_{n{\vp1}}}
\def\Fvt{F^{\lg_{\tau_1}}v_{\tau_1}\lox F^{\lg_{\tau_n}}v_{\tau_n}}
\def\Fvat{(\Fvt)^*}
\def\Fvqt{{F_{\!q}}^{\!\]\lg_{\tau_1}}v^q_{\tau_1}\lox
 {F_{\!q}}^{\!\]\lg_{\tau_n}}v^q_{\tau_n}}
\def\Fvm{v_1\lox Fv_m\lox v_n}
\def\sing{^{\sscr{\italic{sing}}}}
\def\Vls{(\Vox)_\ell\sing} \def\Vlst{(\Voxt)_\ell\sing}
 
\def\VFV{\Val\big/\FVal} \def\VFVt{\Valt\big/\FValt}
\def\Vqst{(\Vqxt)\sing} \def\Vqlst{\Vqst_\ell}
\def\Vqstt{(\Vqxtt)\sing} \def\Vqlstt{\Vqstt_\ell}

\let\cale\circlearrowleft \let\cari\circlearrowright
\ifMag
\def\DOT{\hbox to 1.5em{\hfill\clap{$\cdot$}\hfill}}
\def\AST{\hbox to 1.5em{\hfill\raise.3ex\clap{$\sss*$}\hfill}}
\def\CDOT{\hbox to 1.5em{\hfill\clap{$\cari$}\clap{$\cdot$}\hfill}}
\def\CAST{\hbox to 1.5em{\hfill\clap{$\cale$}\raise.3ex\clap{$\sss*$}\hfill}}
\def\BLANK{\hbox to 1.5em{\hfill}}
\else
\def\DOT{\hbox to 1.7em{\hfill\raise.07ex\clap{$\cdot$}\hfill}}
\def\AST{\hbox to 1.7em{\hfill\raise.37ex\clap{$\sss*$}\hfill}}
\def\CDOT{\hbox to 1.7em{\hfill\clap{$\cari$}\raise.07ex\clap{$\cdot$}\hfill}}
\def\CAST{\hbox to 1.7em{\hfill\clap{$\cale$}\raise.37ex\clap{$\sss*$}\hfill}}
\def\BLANK{\hbox to 1.7em{\hfill}}
\fi
\def\vertar{\Clap{\left\uparrow\vcenter to 15ex{}\right.\n@space}}
\def\Sinp#1{\sin\bigl(\pi(#1)/p\bigr)} \def\Sinpi#1{\sin\bigl(\pi(#1)\bigr)}
\def\Expp#1{\exp\bigl(\pii\)(#1)/p\bigr)}
\def\Expt#1{\exp\bigl(2\pii\)(#1)/p\bigr)}
\def\Gmpb#1{\Gm\bigl(#1\bigr)} \def\Gmpp#1{\Gmpb{(#1)/p}}

\def\abm{(\al_m-\bt_m)}  
\def\abmn{\mathchoice{\sun\abm}{\!\sun\abm}{\!\sun\!\abm}{\@@PS}}
\def\Lmn{\mathchoice{2\sun\La_m}{2\!\sun\!\La_m}{2\!\sun\!\La_m}{\@@PS}}
\def\lmn{\mathchoice{\sun\la_m}{\!\sun\!\la_m}{\!\sun\!\la_m}{\@@PS}}
\def\tlmn{\tsun\la_m} \def\tLmn{2\tsun\La_m}
\def\M{\Lmn/p} \def\tM{\tLmn/p}

\def\Imu{0<\Im\mu<2\pi} \def\Immu{0\le\Im\mu<2\pi}

\def\dd#1{{\dsize{\der\over\der#1}}} 
\def\dt{\,d^\ell t} 
 \def\zs{z^{\sss\star}}
\def\part#1{\^{\it$#1$-part}}

\def\gl{\frak{gl}_2} \def\gsl{\frak{sl}_2} \def\Usl{U(\gsl)}
\def\Ysl{Y(\gsl)} \def\Uu{U_q(\gsl)} 
\def\Ygl{Y(\gl)} \def\Ugh{U_q(\widetilde{\gl})} \def\Ugg{U_q'(\widetilde{\gl})}

\def\Icy/{\^{$I$-cycle}} \def\Jcy/{\^{$J$-cycle}} \def\azo/{\as/ zone}
\def\regud/{regularized} \def\regul/{regularization}
\def\appr/{approximate} \def\apsol/{\appr/ \sol/}
\def\quea/{quantized universal enveloping algebra} \def\asol/{\as/ \sol/}
\def\conn/{connection coefficient} \def\sconn/{system of \conn/s}
\def\loc/{local system} \def\dloc/{discrete local system}
\def\prim/{primitive factor} \def\hform/{\hgeom/ form}
\def\hpair/{\hgeom/ pairing} \def\hmod/{\hgeom/ module}
\def\hgf/{\hgeom/ space} \def\hgF/{\hgeom/ Fock space}
\def\rhgf/{\rat/ \hgf/} \def\rhgF/{\rat/ \hgF/} \def\rhgm/{\rat/ \hmod/}
\def\thgf/{\tri/ \hgf/} \def\thgF/{\tri/ \hgF/} \def\thgm/{\tri/ \hmod/}
\def\sthgf/{singular \thgf/}
\def\hcg/{\hgeom/ cohomology group} \def\hhg/{\hgeom/ homology group}
\def\trib/{trivial bundle} \def\tvb/{the trivial vector bundle}
\def\wtd/{\wt/ decomposition} \def\ntree/{\^{$n$-}tree}
\def\ivx/{internal vertex} \def\ivces/{internal vertices}
\def\cosub/{coboundary subspace} \def\qlo/{quantum loop algebra}

\def\GM/{Gauss-Manin connection} \def\MB/{Mellin-Barnes}

\def\stype/{\^{$\,\gsl$-}type} \def\gmod/{\^{$\,\g$-}module}
\def\smod/{\^{$\,\gsl$-}module} \def\hwsm/{\hw/ \smod/}
\def\Ymod/{\^{$\,\Ygl$-}module} \def\Ysmod/{\^{$\,\Ysl$-}module}
\def\Umod/{\^{$\,\Uu$-}module} \def\Uhmod/{\^{$\,\Ugg$-}module}
\def\tenco/{tensor coordinates} \def\traf/{transition \fn/}
\def\phf/{phase \fn/} \def\wtf/{\wt/ \fn/} 
\def\rwf/{\rat/ \wtf/} \def\twf/{\tri/ \wtf/}

\def\Vval/{\^{$V\!\!\;$-}valued} \def\SL#1{\^{$\,\Sl$}#1}
\def\p#1{\ifx#1'\def\n@xt##1{\^{$\pb\>$##1}}\else
 \def\n@xt{\^{$p\>$#1}}\fi\n@xt}

\def\ncp/{\ndeg/ \cp/} \def\adm/{admissible} \def\arf/{\adm/ \raf/}
\def\mphf/{modified \phf/} \def\ephf/{extended \phf/}

\setparindent

\newif\ifens
\let\goodbe\relax \let\egood\relax \let\eegood\relax
\let\goodbm\relax \let\mgood\relax 

\ifMag\let\goodbm\goodbreak \let\mgood\vgood 
 \let\goodbreak\relax \let\vgood\relax \let\vvgood\relax \fi

\ifamsppt \let\ensppt\relax \else \def\ensppt{\enstrue
 \let\goodbe\goodbreak \let\egood\vgood \let\eegood\vvgood
 \let\goodbreak\relax \let\vgood\relax \let\vvgood\relax
 \pageno\z@\footline={\ifnum\pageno=0\hfil\else\hfil\foliorm\folio\hfil\fi}}
\fi

\let\fixedpage\newpage

\csname swan.def\endcsname

\setbox\lefthbox\hbox{\eightpoint\sl \VT/ \ and \ \Varch/}
\setbox\righthbox\hbox{\=\eightpoint\sl
$q$-Hypergeometric Functions, Yangians and Quantum Affine Algebras}
\leftheadtext{\copy\lefthbox}
\rightheadtext{\copy\righthbox}

\labeldef{S} {1} {1}

\labeldef{S} {2} {2}
\labeldef{F} {2\labelsep \labelspace 1}  {flat}
\labeldef{F} {2\labelsep \labelspace 2}  {period}
\labeldef{F} {2\labelsep \labelspace 3}  {sl2type}
\labeldef{F} {2\labelsep \labelspace 4}  {Phox}
\labeldef{F} {2\labelsep \labelspace 5}  {Phi}
\labeldef{F} {2\labelsep \labelspace 6}  {Stir}
\labeldef{F} {2\labelsep \labelspace 7}  {Phisym}
\labeldef{F} {2\labelsep \labelspace 8}  {list1}
\labeldef{F} {2\labelsep \labelspace 9}  {act}
\labeldef{F} {2\labelsep \labelspace 10} {dis}
\labeldef{F} {2\labelsep \labelspace 11} {subbundle}
\labeldef{F} {2\labelsep \labelspace 12} {npZ}
\labeldef{F} {2\labelsep \labelspace 13} {Lass}
\labeldef{F} {2\labelsep \labelspace 14} {assum}
\labeldef{F} {2\labelsep \labelspace 15} {kanon}
\labeldef{F} {2\labelsep \labelspace 16} {HFo}
\labeldef{F} {2\labelsep \labelspace 17} {kaon}
\labeldef{F} {2\labelsep \labelspace 18} {Zln}
\labeldef{F} {2\labelsep \labelspace 19} {wlga}
\labeldef{F} {2\labelsep \labelspace 20} {wbasis}
\labeldef{F} {2\labelsep \labelspace 21} {wDw}
\labeldef{F} {2\labelsep \labelspace 22} {HFR}
\labeldef{F} {2\labelsep \labelspace 23} {wtau}
\labeldef{F} {2\labelsep \labelspace 24} {thgf}
\labeldef{F} {2\labelsep \labelspace 25} {triact}
\labeldef{F} {2\labelsep \labelspace 26} {Wlga}
\labeldef{F} {2\labelsep \labelspace 27} {Wlgo}
\labeldef{F} {2\labelsep \labelspace 28} {Wbasis}
\labeldef{F} {2\labelsep \labelspace 29} {Wbasiso}
\labeldef{F} {2\labelsep \labelspace 30} {Wtau}

\labeldef{S} {3} {2a}
\labeldef{F} {3\labelsep \labelspace 1}  {VVF}
\labeldef{F} {3\labelsep \labelspace 2}  {Rdef}
\labeldef{F} {3\labelsep \labelspace 3}  {Rvv}
\labeldef{F} {3\labelsep \labelspace 4}  {RXX}
\labeldef{F} {3\labelsep \labelspace 5}  {Rspec}
\labeldef{F} {3\labelsep \labelspace 6}  {YB}
\labeldef{F} {3\labelsep \labelspace 7}  {Yijrel}
\labeldef{F} {3\labelsep \labelspace 8}  {Yintw}
\labeldef{F} {3\labelsep \labelspace 9}  {Rij}
\labeldef{F} {3\labelsep \labelspace 10} {Kmz}
\labeldef{F} {3\labelsep \labelspace 11} {FR}
\labeldef{F} {3\labelsep \labelspace 13} {Rqdef}
\labeldef{F} {3\labelsep \labelspace 14} {Rqvv}
\labeldef{F} {3\labelsep \labelspace 15} {Rqmore}
\labeldef{F} {3\labelsep \labelspace 16} {Rspeq}
\labeldef{F} {3\labelsep \labelspace 17} {YBq}
\labeldef{F} {3\labelsep \labelspace 18} {Ughr}
\labeldef{F} {3\labelsep \labelspace 19} {Uintw}

\labeldef{S} {4} {2b}
\labeldef{F} {4\labelsep \labelspace 1}  {Tiju}
\labeldef{F} {4\labelsep \labelspace 2}  {YFg}
\labeldef{F} {4\labelsep \labelspace 3}  {HFE}
\labeldef{F} {4\labelsep \labelspace 4}  {FfDf}
\labeldef{F} {4\labelsep \labelspace 5}  {tcfs}
\labeldef{F} {4\labelsep \labelspace 6}  {T12}
\labeldef{F} {4\labelsep \labelspace 7}  {FgVax}
\labeldef{F} {4\labelsep \labelspace 9}  {local}
\labeldef{F} {4\labelsep \labelspace 10} {tcnon}
\labeldef{F} {4\labelsep \labelspace 11} {tcon}
\labeldef{F} {4\labelsep \labelspace 12} {qKZ-GM}
\labeldef{F} {4\labelsep \labelspace 14} {Lijxi}
\labeldef{F} {4\labelsep \labelspace 15} {UFg}
\labeldef{F} {4\labelsep \labelspace 16} {HFEq}
\labeldef{F} {4\labelsep \labelspace 17} {tcfsq}
\labeldef{F} {4\labelsep \labelspace 18} {L21}
\labeldef{F} {4\labelsep \labelspace 19} {FgVqx}
\labeldef{F} {4\labelsep \labelspace 21} {Csing}
\labeldef{F} {4\labelsep \labelspace 22} {localq}
\labeldef{F} {4\labelsep \labelspace 23} {FFo}
\labeldef{F} {4\labelsep \labelspace 24} {FFq}
\labeldef{F} {4\labelsep \labelspace 25} {FgF}
\labeldef{F} {4\labelsep \labelspace 26} {FgFq}

\labeldef{S} {5} {2c}
\labeldef{F} {5\labelsep \labelspace 1}  {IWwt}
\labeldef{F} {5\labelsep \labelspace 2}  {Immu}
\labeldef{F} {5\labelsep \labelspace 3}  {IWw}
\labeldef{F} {5\labelsep \labelspace 5}  {Wwell}
\labeldef{F} {5\labelsep \labelspace 6}  {list2}
\labeldef{F} {5\labelsep \labelspace 7}  {Ianco}
\labeldef{F} {5\labelsep \labelspace 8}  {Ianco0}
\labeldef{F} {5\labelsep \labelspace 9}  {IWw=0}
\labeldef{F} {5\labelsep \labelspace 10} {hpFF}
\labeldef{F} {5\labelsep \labelspace 11} {hpFFo}
\labeldef{F} {5\labelsep \labelspace 12} {1000}
\labeldef{F} {5\labelsep \labelspace 13} {2000}
\labeldef{F} {5\labelsep \labelspace 14} {mu<>0}
\labeldef{F} {5\labelsep \labelspace 15} {mu=0}
\labeldef{F} {5\labelsep \labelspace 16} {Barnes}
\labeldef{F} {5\labelsep \labelspace 17} {BarneS}
\labeldef{F} {5\labelsep \labelspace 18} {sWz}
\labeldef{F} {5\labelsep \labelspace 19} {PsiW}
\labeldef{F} {5\labelsep \labelspace 20} {qKZsol}
\labeldef{F} {5\labelsep \labelspace 21} {3000}
\labeldef{F} {5\labelsep \labelspace 22} {4000}
\labeldef{F} {5\labelsep \labelspace 23} {5000}
\labeldef{F} {5\labelsep \labelspace 24} {6000}

\labeldef{S} {6} {2d}
\labeldef{F} {6\labelsep \labelspace 1}  {syst}
\labeldef{F} {6\labelsep \labelspace 2}  {fzz}
\labeldef{F} {6\labelsep \labelspace 3}  {azone}
\labeldef{F} {6\labelsep \labelspace 4}  {asol}
\labeldef{F} {6\labelsep \labelspace 5}  {RpR}
\labeldef{F} {6\labelsep \labelspace 6}  {Asol}
\labeldef{F} {6\labelsep \labelspace 7}  {AsIWw}

\labeldef{S} {7} {2e}
\labeldef{F} {7\labelsep \labelspace 1}  {AhB}
\labeldef{F} {7\labelsep \labelspace 2}  {Phih}
\labeldef{F} {7\labelsep \labelspace 3}  {Phicl}
\labeldef{F} {7\labelsep \labelspace 4}  {YYtau}
\labeldef{F} {7\labelsep \labelspace 5}  {wticl}
\labeldef{F} {7\labelsep \labelspace 6}  {qcl}
\labeldef{F} {7\labelsep \labelspace 7}  {UUlg}
\labeldef{F} {7\labelsep \labelspace 8}  {qcl0}
\labeldef{F} {7\labelsep \labelspace 9}  {conj}

\labeldef{S} {8} {3}
\labeldef{F} {8\labelsep \labelspace 1}  {eta<>0}
\labeldef{F} {8\labelsep \labelspace 2}  {nabla}
\labeldef{F} {8\labelsep \labelspace 3}  {eta=0}
\labeldef{F} {8\labelsep \labelspace 4}  {cycles}
\labeldef{F} {8\labelsep \labelspace 5}  {cycles0}
\labeldef{F} {8\labelsep \labelspace 6}  {Varch}
\labeldef{F} {8\labelsep \labelspace 7}  {Phi1}
\labeldef{F} {8\labelsep \labelspace 8}  {wma}
\labeldef{F} {8\labelsep \labelspace 9}  {wba}
\labeldef{F} {8\labelsep \labelspace 10} {detM}
\labeldef{F} {8\labelsep \labelspace 11} {D1}
\labeldef{F} {8\labelsep \labelspace 12} {kanon1}
\labeldef{F} {8\labelsep \labelspace 13} {kaon1}
\labeldef{F} {8\labelsep \labelspace 14} {Wma}
\labeldef{F} {8\labelsep \labelspace 15} {Wmo}
\labeldef{F} {8\labelsep \labelspace 16} {Wba}
\labeldef{F} {8\labelsep \labelspace 17} {Wbao}
\labeldef{F} {8\labelsep \labelspace 18} {IWw1}
\labeldef{F} {8\labelsep \labelspace 19} {IWwt1}
\labeldef{F} {8\labelsep \labelspace 20} {ICC}
\labeldef{F} {8\labelsep \labelspace 21} {points}
\labeldef{F} {8\labelsep \labelspace 22} {Ianco1}
\labeldef{F} {8\labelsep \labelspace 23} {Ianco01}
\labeldef{F} {8\labelsep \labelspace 24} {IWws}
\labeldef{F} {8\labelsep \labelspace 25} {IWwd}
\labeldef{F} {8\labelsep \labelspace 27} {sWz1}
\labeldef{F} {8\labelsep \labelspace 28} {qKZsol1}
\labeldef{F} {8\labelsep \labelspace 29} {asol1}
\labeldef{F} {8\labelsep \labelspace 30} {z+p}
\labeldef{F} {8\labelsep \labelspace 31} {PhiwW}
\labeldef{F} {8\labelsep \labelspace 32} {IWwm}
\labeldef{F} {8\labelsep \labelspace 33} {mu<>01}
\labeldef{F} {8\labelsep \labelspace 34} {mu=01}
\labeldef{F} {8\labelsep \labelspace 35} {mu0}
\labeldef{F} {8\labelsep \labelspace 36} {Jj}
\labeldef{F} {8\labelsep \labelspace 37} {IWnwn}
\labeldef{F} {8\labelsep \labelspace 39} {qcl1}
\labeldef{F} {8\labelsep \labelspace 40} {qcl01}
\labeldef{F} {8\labelsep \labelspace 41} {IWmwm}
\labeldef{F} {8\labelsep \labelspace 42} {||<A}

\labeldef{S} {9} {4}
\labeldef{F} {9\labelsep \labelspace 1}  {DetM}
\labeldef{F} {9\labelsep \labelspace 3}  {list3}
\labeldef{F} {9\labelsep \labelspace 4}  {IIx}
\labeldef{F} {9\labelsep \labelspace 5}  {IWwD}
\labeldef{F} {9\labelsep \labelspace 6}  {wWprim}
\labeldef{F} {9\labelsep \labelspace 7}  {ident}
\labeldef{F} {9\labelsep \labelspace 8}  {identi}
\labeldef{F} {9\labelsep \labelspace 9}  {u-eq}
\labeldef{F} {9\labelsep \labelspace 10} {Fufu}
\labeldef{F} {9\labelsep \labelspace 11} {Fnew}
\labeldef{F} {9\labelsep \labelspace 12} {0<<1}
\labeldef{F} {9\labelsep \labelspace 13} {asdet}
\labeldef{F} {9\labelsep \labelspace 14} {Xic}
\labeldef{F} {9\labelsep \labelspace 15} {Wprime}
\labeldef{F} {9\labelsep \labelspace 16} {detN}
\labeldef{F} {9\labelsep \labelspace 17} {Omu}
\labeldef{F} {9\labelsep \labelspace 18} {asIW'w}
\labeldef{F} {9\labelsep \labelspace 19} {shafl}
\labeldef{F} {9\labelsep \labelspace 20} {Wprimd}
\labeldef{F} {9\labelsep \labelspace 21} {I*}
\labeldef{F} {9\labelsep \labelspace 22} {t**}
\labeldef{F} {9\labelsep \labelspace 23} {cone}
\labeldef{F} {9\labelsep \labelspace 25} {asJal}
\labeldef{F} {9\labelsep \labelspace 26} {IWbw}
\labeldef{F} {9\labelsep \labelspace 27} {asIWbw}

\def\ContS{
\Entcd{1}{Introduction}{1}
\Entcd{2}{Discrete flat connections and local systems}{4}
\subcd{Discrete flat connections}{4}
\subcd{Connection coefficients of local systems}{7}
\subcd{The functional space of a rational \stype / local system}{8}
\subcd{Bases in the rational hyper\-geometric space of a fiber}{10}
\subcd{The trigonometric hyper\-geometric space}{11}
\Entcd{3}{\Rms / and the \qKZc /}{13}
\subcd{Highest weight \smod /s}{13}
\subcd{The rational \Rm /}{14}
\subcd{The Yangian $\Ygl $}{15}
\subcd{The rational \qKZc / associated with $\gsl $}{16}
\subcd{The trigonometric \Rm /}{17}
\subcd{The quantum loop algebra $\Ugg $}{18}
\Entcd{4}{Tensor coordinates and module structures on the hyper\-geometric spaces}{20}
\subcd{The rational hyper\-geometric module}{20}
\subcd{Tensor coordinates on the rational hyper\-geometric spaces of fibers}{21}
\subcd{The trigonometric hyper\-geometric module}{23}
\subcd{Tensor coordinates on the trigonometric hyper\-geometric spaces of fibers}{24}
\subcd{Tensor products of the hyper\-geometric modules}{25}
\Entcd{5}{The hyper\-geometric pairing}{27}
\Entcd{6}{Asymptotic solutions to the \qKZe /}{31}
\Entcd{7}{Quasiclassical asymptotics}{35}
\Entcd{8}{The one-dimensional case}{37}
\subcd{One-dimensional discrete cohomologies}{38}
\subcd{One-dimensional discrete homologies}{39}
\subcd{Quasiclassical asymptotics}{46}
\Entcd{9}{The multidimensional case}{47}
\Refcd{References}{63}
}

\def\ContM{
\Entcd{1}{Introduction}{1}
\Entcd{2}{Discrete flat connections and local systems}{6}
\subcd{Discrete flat connections}{6}
\subcd{Connection coefficients of local systems}{9}
\subcd{The functional space of a rational \stype / local system}{11}
\subcd{Bases in the rational hyper\-geometric space of a fiber}{13}
\subcd{The trigonometric hyper\-geometric space}{15}
\Entcd{3}{\Rms / and the \qKZc /}{18}
\subcd{Highest weight \smod /s}{18}
\subcd{The rational \Rm /}{18}
\subcd{The Yangian $\Ygl $}{19}
\subcd{The rational \qKZc / associated with $\gsl $}{21}
\subcd{The trigonometric \Rm /}{22}
\subcd{The quantum loop algebra $\Ugg $}{24}
\Entcd{4}{Tensor coordinates and module structures on the hyper\-geometric spaces}{26}
\subcd{The rational hyper\-geometric module}{26}
\subcd{Tensor coordinates on the rational hyper\-geometric spaces of fibers}{28}
\subcd{The trigonometric hyper\-geometric module}{30}
\subcd{Tensor coordinates on the trigonometric hyper\-geometric spaces of fibers}{32}
\subcd{Tensor products of the hyper\-geometric modules}{33}
\Entcd{5}{The hyper\-geometric pairing}{35}
\Entcd{6}{Asymptotic solutions to the \qKZe /}{41}
\Entcd{7}{Quasiclassical asymptotics}{46}
\Entcd{8}{The one-dimensional case}{49}
\subcd{One-dimensional discrete cohomologies}{50}
\subcd{One-dimensional discrete homologies}{52}
\subcd{Quasiclassical asymptotics}{60}
\Entcd{9}{The multidimensional case}{62}
\Refcd{References}{82}
}

\def\ContSP{
\Entcd{1}{Introduction}{1}
\Entcd{2}{Discrete flat connections and local systems}{4}
\subcd{Discrete flat connections}{4}
\subcd{Connection coefficients of local systems}{7}
\subcd{The functional space of a rational \stype / local system}{8}
\subcd{Bases in the rational hyper\-geometric space of a fiber}{10}
\subcd{The trigonometric hyper\-geometric space}{11}
\Entcd{3}{\Rms / and the \qKZc /}{13}
\subcd{Highest weight \smod /s}{13}
\subcd{The rational \Rm /}{14}
\subcd{The Yangian $\Ygl $}{15}
\subcd{The rational \qKZc / associated with $\gsl $}{16}
\subcd{The trigonometric \Rm /}{17}
\subcd{The quantum loop algebra $\Ugg $}{18}
\Entcd{4}{Tensor coordinates and module structures on the hyper\-geometric spaces}{20}
\subcd{The rational hyper\-geometric module}{20}
\subcd{Tensor coordinates on the rational hyper\-geometric spaces of fibers}{21}
\subcd{The trigonometric hyper\-geometric module}{23}
\subcd{Tensor coordinates on the trigonometric hyper\-geometric spaces of fibers}{24}
\subcd{Tensor products of the hyper\-geometric modules}{25}
\Entcd{5}{The hyper\-geometric pairing}{27}
\Entcd{6}{Asymptotic solutions to the \qKZe /}{31}
\Entcd{7}{Quasiclassical asymptotics}{35}
\Entcd{8}{The one-dimensional case}{37}
\subcd{One-dimensional discrete cohomologies}{38}
\subcd{One-dimensional discrete homologies}{39}
\subcd{Quasiclassical asymptotics}{46}
\Entcd{9}{The multidimensional case}{47}
\subcd{Proof of Lemmas~\[wbasis], \[Wbasis], \[Wbasiso]}{47}
\subcd{Proof of Lemma~\[wDw]}{49}
\subcd{Proof of Theorem~\[Ianco]}{49}
\subcd{Proof of Theorem~\[sWz]}{51}
\subcd{Proof of Theorem~\[AsIWw]}{51}
\subcd{Proof of formula \(BarneS)}{53}
\subcd{Proof of Theorem~\[mu<>0]}{55}
\subcd{Proof of Theorem~\[mu=0]}{56}
\subcd{Proof of Theorem~\[qcl0]}{61}
\Refcd{References}{63}
}

\def\ContMP{
\Entcd{1}{Introduction}{1}
\Entcd{2}{Discrete flat connections and local systems}{6}
\subcd{Discrete flat connections}{6}
\subcd{Connection coefficients of local systems}{9}
\subcd{The functional space of a rational \stype / local system}{11}
\subcd{Bases in the rational hyper\-geometric space of a fiber}{13}
\subcd{The trigonometric hyper\-geometric space}{15}
\Entcd{3}{\Rms / and the \qKZc /}{18}
\subcd{Highest weight \smod /s}{18}
\subcd{The rational \Rm /}{18}
\subcd{The Yangian $\Ygl $}{19}
\subcd{The rational \qKZc / associated with $\gsl $}{21}
\subcd{The trigonometric \Rm /}{22}
\subcd{The quantum loop algebra $\Ugg $}{24}
\Entcd{4}{Tensor coordinates and module structures on the hyper\-geometric spaces}{26}
\subcd{The rational hyper\-geometric module}{26}
\subcd{Tensor coordinates on the rational hyper\-geometric spaces of fibers}{28}
\subcd{The trigonometric hyper\-geometric module}{30}
\subcd{Tensor coordinates on the trigonometric hyper\-geometric spaces of fibers}{32}
\subcd{Tensor products of the hyper\-geometric modules}{33}
\Entcd{5}{The hyper\-geometric pairing}{35}
\Entcd{6}{Asymptotic solutions to the \qKZe /}{41}
\Entcd{7}{Quasiclassical asymptotics}{46}
\Entcd{8}{The one-dimensional case}{49}
\subcd{One-dimensional discrete cohomologies}{50}
\subcd{One-dimensional discrete homologies}{52}
\subcd{Quasiclassical asymptotics}{60}
\Entcd{9}{The multidimensional case}{62}
\subcd{Proof of Lemmas~\[wbasis], \[Wbasis], \[Wbasiso]}{62}
\subcd{Proof of Lemma~\[wDw]}{64}
\subcd{Proof of Theorem~\[Ianco]}{65}
\subcd{Proof of Theorem~\[sWz]}{67}
\subcd{Proof of Theorem~\[AsIWw]}{68}
\subcd{Proof of formula \(BarneS)}{70}
\subcd{Proof of Theorem~\[mu<>0]}{72}
\subcd{Proof of Theorem~\[mu=0]}{73}
\subcd{Proof of Theorem~\[qcl0]}{80}
\Refcd{References}{82}
}

\document

\def\abstext{\Abstract
The \rat/ quantized \KZv/ (\qKZ/) \eq/ associated with the Lie algebra $\gsl$
is a system of linear \deq/s with values in a tensor product of $\gsl\]$
\Vmod/s. We solve the equation in terms of multidimensional \^{$q$-}\hgeom/
\fn/s and define a natural \iso/ of the space of \sol/s and the tensor product
of the corresponding \qg/ $\Uu$ \Vmod/s where the parameter $q$ is related to
the step $p$ of the \qKZe/ via $q=e^{\pi i/p}$.
\par
We construct \asol/s associated with suitable \azo/s and compute the \traf/s
between the \asol/s in terms of the \tri/ \Rms/. This description of the
\traf/s gives a new connection between \rep/ theories of Yangians and \qlo/s
and is analogous to the Kohno-Drinfeld theorem on the monodromy group of
the \difl/ \KZe/.
\par
In order to establish these results we construct a discrete \GM/,
in particular, a suitable discrete \loc/, discrete homology and cohomology
groups with coefficients in this \loc/, and identify an associated \deq/
with the \qKZe/.
\endAbs}


\Sno 0

\ifMag\else\ifamsppt\else\ifx\EnsLapp\relax\enstrue\fi\fi\fi
\ifens
\vbox to\vsize
{\twelvepoint
\vskip.9truein
\vfill
\center
\=
{\bls 20pt
\Bbf
Geometry of {\Bmmi q\)}-Hypergeometric Functions as a Bridge
\\
between Yangians and Quantum Affine Algebras
\par}
\vsk1.5>
\VT/$^{\,\star}$ \ and \ \Varch/$^{\,*}$
\vsk1.5>
{\it
$^\star$Laboratoire de Physique Th\'eorique
{\nineit E\tenit N\it S\Bit L\it A\tenit P\nineit P}$^{\,\diamond}\!$\\
ENS Lyon, \LPTaddr/
\vsk>
$^*$\UNC/\\ \UNCaddr/}
\endcenter
\vfill
\tenpoint
\abstext
\vfill\vfill
{\twelvepoint
\rline{{\ninerm E\tenrm N\rm S\Brm L\rm A\tenrm P\ninerm P}\~~L\~~585/96}
\rline{q-alg/9604011}
\rline{April 1996}
\rline{{\tensl to appear in}}
\rline{\it Inventiones}
\rline{\it Mathematicae}}
\vfill\vfill
\=
\lline{$\]^\diamond\>$\enslapp/}
\vsk.1>
\lline{$\]^\star\>$\absence/}
\lline{\hp{$^\star\>$}\CNRS/}
\vsk.15>
\lline{\twelvepoint\hp{$^\star$}\sl E-mail\/{\rm:} \ensemail/}
\vsk.1>
\lline{$\]^*\>$\Grant/}
\vsk.15>
\lline{\twelvepoint\hp{$^*$}\sl E-mail\/{\rm:} \avemail/}}
\else
\vp1
{\ifMag\vsk-2>\eightpoint\else\vsk->\fi
\vsk-2>
\rline{\sl Inventiones Mathematicae}
\rline{{\bf128} (1997), 501--588}} 
\vsk1.2>
\center
\=
{\bls 16pt
\Bbf
Geometry of {\Bmmi q\)}-Hypergeometric Functions as a Bridge
\\
between Yangians and Quantum Affine Algebras
\par}
\vsk1.5>
\VT/$^{\,\star}$ \ and \ \Varch/$^{\,*}$
\vsk1.5>
{\it
$^\star$\LPT/$^{\,\diamond}\!$, \)\ENSLyon/\\ \LPTaddr/
\vsk.35>
\ifMag
$^*$\UNC/\\ \UNCaddr/
\else
$^*$\UNC/, \UNCaddr/
\fi}
\vsk1.8>
{\sl April \,4, 1996}
\endcenter
\ftext{\=\bls11pt $\]^\diamond\>$\enslapp/\nl
$\]^\star\>$\absence/\vv-.05>\nl \hp{$^\star\>$}\CNRS/\nl
{\tenpoint\hp{$^\star$}\sl E-mail\/{\rm:} \ensemail/\,, \;\homemail/}\nl
$\]^*\>$\Grant/\nl
{\tenpoint\hp{$^*$}\sl E-mail\/{\rm:} \avemail/}}
\fi

\Sect[1]{Introduction}
In this paper we consider the \rat/ quantized \KZv/ (\qKZ/) \eq/ associated
with the Lie algebra $\gsl$ and solve it. The \rat/ \qKZe/ associated with
$\gsl$ is a system of \deq/s for a \fn/ $\Psi(\zn)$ with values in a tensor
product $\Vox$ of \smod/s. The system of \eq/s has the form
$$
\align
\Psi(\zmn)\,=\,
R_{m,m-1}(z_m-z_{m-1}+p)\ldots R_{m,1}(z_m-z_1+p)\>\ka^{-H_m}\;\x &
\\
\nn3>
{}\x\,R_{m,n}(z_m-z_n)\ldots R_{m,m+1}(z_m-z_{m+1})\>\Psi(\zn) &\,,
\endalign
$$
$\mn$, where $p$ and $\ka$ are parameters of the \qKZe/, $H$ is a generator of
the Cartan subalgebra of $\gsl$, $H_m$ is $H$ acting in the \^{$m$-th} factor,
$R_{l,m}(x)$ is the \rat/ \Rm/ $R_{V_lV_m}(x)\in\End(V_l\ox V_m)$ acting in the
\^{$l$-th} and \^{$m$-th} factors of the tensor product of \smod/s. In this
paper we consider only the negative steps $p$. The case of other values of
the step can be treated by \anco/.
\par
The \qKZe/ is an important system of \deq/s. The \qKZe/s had been introduced in
\Cite{FR} as \eq/s for matrix elements of vertex operators of the \qaff/.
An important special case of the \qKZe/ had been introduced earlier in \Cite{S}
as \eq/s for form factors in integrable quantum field theory; relevant \sol/s
for these \eq/s had been given therein. Later, the \qKZe/s were derived as
\eq/s for correlation \fn/s in lattice integrable models, \cf. \Cite{JM} and
references therein.
\par
In the \qcl/ limit the \qKZe/ turns into the \difl/ \KZv/ \eq/ for conformal
blocks of the Wess-Zumino-Witten model of conformal field theory on the sphere.
\par
Asymptotic \sol/s to the \qKZe/ as $p$ tends to zero are closely related to
diagonalization of the transfer-matrix of the corresponding lattice integrable
model by the algebraic \Ba/ method \Cite{TV2}.
\Par
We describe the space of \sol/s to the \qKZe/ in terms of \rep/ theory. Namely,
we consider the \qg/ $\Uu$ with $q=e^{\pii/p}$ and the \Umod/s $\Vqn$ where
$V^q_m$ is the deformation of the \smod/ $V_m$. For every \perm/ $\tau\in\S^n$
we consider the tensor product $\Vqxt$ and establish a natural \iso/ of
the space $\SS$ of \sol/s to the \qKZe/ with values in $\Vox$ and the space
$\Vqxt\ox\FF$, where $\FF$ is the space of \fn/s of $\zn$ which are \p-periodic
\wrt/ each of the \var/s,
$$
C_\tau:\Vqxt\ox\FF\,\to\,\SS.
$$
Notice that if $\Psi(z)$ is a \sol/ to the \qKZe/ and $F(z)$ is a \p-periodic
\fn/, then also $F(z)\Psi(z)$ is a \sol/ to the \qKZe/.
\par
We call the \iso/s $C_\tau$ the tensor coordinates on the space of \sol/s.
The compositions of the isomorphisms are linear maps
$$
C_{\tau,\tau'}(\zn):\Vqxtt\,\to\,\Vqxt
$$
depending on $\zn$ and \p-periodic \wrt/ all \var/s. We call these
compositions the \traf/s. It turns out that the \traf/s are defined in terms of
the \tri/ \Rms/ ${R^q_{V_lV_m}(\zt)\in\End(V^q_l\ox V^q_m)}$ acting in tensor
products of \Umod/s. Namely, for any \perm/ $\tau$ and for any transposition
$(m,m+1)$ the \traf/
$$
C_{\tau,\tau\cdot(m,m+1)}(\zn):\Vqxtm\,\to\,\Vqxt
$$
equals the operator ${P_{V^q_{\tau_{m+1}}\}V^q_{\tau_m}}\:
R^q_{V^q_{\tau_{m+1}}\}V^q_{\tau_m}}\!\!
\bigl(\Expt{z_{\tau_{m+1}}-z_{\tau_m}}\bigr)\)}$
acting in the \^{$m$-th} and \^{$(m+1)$-th} factors, here $P_{V_lV_m}\:$
is the transposition of the tensor factors; \cf. Theorem~\[localq].
\Par
We consider \azo/s $\Re z_{\tau_1}\lsym\ll\Re z_{\tau_n}$ labelled by \perm/s
$\tau\in\S^n$. For every \azo/ we define a basis of \asol/s to the \qKZe/.
We show that for every \perm/ $\tau$ the basis of the corresponding \asol/s is
the image of the standard monomial basis in $\Vqx$ under the map
\ifMag\vv-.3>\fi
$$
\Vqxt\,\to\,\Vqxt\ox 1\,\hookrightarrow\,\Vqxt\ox\FF
\,\ {\overset{C_\tau}\to\lto}\,\ \SS\,,
$$
\cf. Theorem~\[asol]. The last two statements express the \traf/s between
the \asol/s via the \tri/ \Rms/.
\Par
The \rat/ \Rm/ $R_{V_lV_m}\:(x)\in\End(V_l\ox V_m)$ is defined in terms of
the action of the Yangian $\Ygl$ in the tensor product of \smod/s. The Yangian
$\Ygl$ is a Hopf algebra which contains the universal enveloping algebra $\Usl$
as a Hopf subalgebra and has a family of \hom/s $\Ygl\to\Usl$ depending on
a parameter. Therefore, each \smod/ $V_m$ carries a \Ymod/ structure $V_m(x)$
depending on a parameter. For \irr/ \smod/s $V_l,\,V_m$ the Yangian modules
$V_l(x)\ox V_m(y)$ and $V_m(y)\ox V_l(x)$ are \irr/ and isomorphic for generic
$x,y$. The map
$$
P_{V_lV_m}\:R_{V_lV_m}\:(x-y):V_l(x)\ox V_m(y)\,\to\,V_m(y)\ox V_l(x)
$$
is the unique suitably normalized intertwiner \Cite{T}, \Cite{D1}.
\Par
Similarly, the \tri/ \Rm/ $R^q_{V_lV_m}(\zt)\in\End(V^q_l\ox V^q_m)$ is defined
in terms of the action of the \qlo/ $\Ugg$ in the tensor product of \Umod/s.
The \qlo/ $\Ugg$ contains $\Uu$ as a Hopf subalgebra and has a family of \hom/s
$\Ugg\to\Uu$ depending on a parameter. Therefore, each \Umod/ $V_m^q$ has
a \Uhmod/ structure $V^q_m(\zt)$ depending on a parameter.
For \irr/ \Umod/s $V_l,\,V_m$ the \Uhmod/s $V^q_l(\xi)\ox V^q_m(\zt)$ and
$V^q_m(\zt)\ox V^q_l(\xi)$ are \irr/ and isomorphic for generic $\xi,\zt$.
The map
$$
P_{V_lV_m}\:R^q_{V_lV_m}(\xi/\zt):
V^q_l(\xi)\ox V^q_m(\zt)\,\to\,V^q_m(\zt)\ox V^q_l(\xi)
$$
is the unique suitably normalized intertwiner \Cite{T}, \Cite{CP}.
\Par
Our result on the \traf/s between \asol/s together with the indicated
construction of \Rms/ shows that the \qKZe/ establishes a connection between
\rep/ theories of the Yangian $\Ygl$ and the \qlo/ $\Ugg$. Our result is
analogous to the Kohno-Drinfeld theorem on the monodromy group of the \difl/
\KZv/ \eq/
\Cite{K}, \Cite{D2}.
\Par
The \difl/ \KZv/ \eq/ (\KZe/) with values in a tensor product of \smod/s
$V=\Vox$ is a system of \difl/ \eq/s for a \Vval/ \fn/ $\Pss(\zn)$ and has
the form
$$
d\)\Pss\,=\;{1\over p}\,\sum_{l\ne m}
{\Om_{lm}\over z_l-z_m}\,\Pss\>d(z_l-z_m)
$$
where $p$ is a parameter of the \eq/, ${\Om_{lm}\in\End(V_l\ox V_m)}$ is
the Casimir operator. The \KZe/ defines an integrable connection over the
complement in $\Cn$ to the union of the diagonal hyperplanes. The fundamental
group of the complement is the pure braid group $\Pn$. The monodromy group of
the \eq/ is the \rep/ $\Pn\to\End(V)$ defined by \anco/ of \sol/s over loops.
The Kohno-Drinfeld theorem says that this \rep/ is isomorphic to the \Rm/ \rep/
of $\Pn$ in the tensor product of \Umod/s $V^q=\Vqx$, $q=e^{\pii/p}$, where
the \Rm/ \rep/ is defined as follows. Let $R^q_{V_lV_m}\in\End(V^q_l\ox V^q_m)$
be the action of the universal \Rm/ of the \qg/ $\Uu$ in the tensor product of
\Umod/s. Then the \Rm/ \rep/ of $\Pn$ in $V^q$ is defined by elementary
transformations
$$
\Vqxt
\ \Lto{^{\ \ P_{V_{\tau_m}V_{\tau_{m+1}}}\:R^q_{V_{\tau_m}V_{\tau_{m+1}}}\ }}
\ \Vqxtm\,.
$$
\par
The Kohno-Drinfeld theorem establishes a connection between \rep/ theories of
a Lie algebra and the corresponding \qg/, see \Cite{D2}. Using the ideas of
the Kohno-Drinfeld result it was proved in \Cite{KL} that the category of
\rep/s of a \qg/ is \eqv/ to a suitably defined fusion category of \rep/s of
the corresponding affine Lie algebra. Similarly to the Kazhdan-Lusztig theorem
one could expect that our result for the \dif/ \qKZe/ could be a base for
a Kazhdan-Lusztig type result connecting certain categories of \rep/s of
Yangians and \qaff/s, \cf. \Cite{KS}.
\Par
In this paper we consider the \rat/ \qKZe/. There are other types of the
\qKZe/: the trigonometric \qKZe/ \Cite{FR} and the elliptic \qKZBe/ \Cite{F}.
Here \KZB/ stands for \KZv/-Bernard, and the \dif/ \qKZBe/ is a discretization
of the \difl/ \KZB/ \eq/ for conformal blocks on the torus.
\Par
The \tri/ \qKZe/ with values in a tensor product of \Umod/s $V^q=\Vqx$ is
a system of \deq/s for a \^{$V^q$-}valued \fn/ $\Psi(\zn)$ and has the form
$$
\align
\Psi(\znm)\,=\,R^q_{m,m-1}(p\)z_m/z_{m-1})\ldots R^q_{m,1}(p\)z_m/z_1)\>
\ka^{-H_m}\;\x&
\\
\nn3>
{}\x R^q_{m,n}(z_m/z_n)\ldots R^q_{m,m+1}(z_m/z_{m+1})\>\Psi(\zn) &\,,
\endalign
$$
$\mn$, where $p$ and $\ka$ are parameters of the \qKZe/, $q^H$ is a generator
of the Cartan subalgebra of $\Uu$, $H_m$ is $H$ acting in the \^{$m$-th}
factor, $R_{l,m}(x)$ is the \tri/ \Rm/
$R^q_{V^q_lV^q_m}(x)\in\End(V^q_l\ox V^q_m)$ acting in the \^{$l$-th} and
\^{$m$-th} factors of the tensor product of \Umod/s. In the next paper
\Cite{TV3} we will describe for the \tri/ \qKZe/ the analogues of the above
results for the \rat/ \qKZe/. Namely, we will describe the space of \sol/s to
the \tri/ \qKZe/ in terms of modules of the elliptic \qg/ associated with the
Lie algebra $\gsl$ \Cite{F}, \Cite{FV} and will get the \traf/s between \asol/s
in the same way as we did for the \rat/ case. This result for the \tri/ \qKZe/
gives a connection between \rep/ theories of the \qlo/ $\Ugg$ and
the elliptic \qg/ associated with $\gsl$.
\par
In the papers \Cite{FTV1}, \Cite{FTV2} we will describe \sol/s to the elliptic
\dif/ \qKZBe/. The construction of \sol/s for the elliptic \qKZBe/ is similar
to the construction of \sol/s to the \rat/ \qKZe/ described in this paper and
to the \sol/s of the \tri/ \qKZ/ in \Cite{M}, \Cite{R}, \Cite{V3}, \Cite{TV1}.
Nevertheless, we do not know yet how to define \asol/s for the elliptic \qKZBe/
and what could be an elliptic analogue of our result on \traf/s.
\Par
There are three different proofs for the Kohno-Drinfeld theorem. Roughly
speaking, they are analytic \Cite{K}, algebraic \Cite{D2}, and geometric
\Cite{SV2}, \Cite{V2}. In the initial proof \Cite{K}, Kohno expands a monodromy
operator as a series of iterated integrals and studies such expansions.
Drinfeld in \Cite{D2} formalizes algebraic properties of \traf/s between
\asol/s and proves that the monodromy group could be nothing else but
the \Rm/ \rep/.
\par
The leading idea of the geometric proof \Cite{SV2}, \Cite{V2}, \Cite{V4} was
the principle that the monodromy of a \difl/ \eq/ could be computed only if the
\difl/ \eq/ is the \eq/ of the \GM/. The \GM/ is a connection associated with a
locally \trib/ of algebraic manifolds with a \loc/ on the space of the bundle.
One considers the associated holomorphic vector bundle which fiber is the
homology group of the fiber of the initial locally \trib/. Then the vector
bundle has a canonical connection called the \GM/. Having a trivialization of
the vector bundle one realises the connection as a system of \difl/ \eq/s.
Its \sol/s are parametrized by elements of the homology group of the fiber.
The monodromy group of that \difl/ \eq/ is the monodromy group of cycles of
the fiber of the initial locally \trib/ under continuous deformations over
loops in the base. The description of the monodromy group of cycles is
a geometric problem which is easier than studying analytic continuation of
\sol/s of an abstract \difl/ \eq/. In order to apply this idea to the proof of
the Kohno-Drinfeld theorem the \difl/ \KZe/ was solved explicitly in terms of
multidimentional \hint/s and \sol/s were represented as integrals of closed
\difl/ forms over cycles depending on parameters, then the space of cycles was
identified with a tensor product of \Umod/s and the monodromy of cycles was
computed in term of \Rms/.
\par
In this paper, in order to establish a connection between \rep/ theories of
Yangians and \qlo/s we quantize the geometric picture for the \KZe/. First we
solve the \qKZe/ in terms of suitable multidimensional \hint/s of \MB/ type.
We define a discrete analogue of a locally \trib/ and a \loc/ on the space of
bundle. We define a discrete analogue of the \GM/ for the discrete locally
\trib/ with a discrete \loc/ and consider the corresponding \deq/. We identify
that \deq/ with the \dif/ \qKZe/. To realize this idea we introduce a suitable
discrete de~Rham complex and its cohomology group in the spirit of \Cite{A},
then we define the homology group as the dual space to the cohomology group and
construct a family of discrete cycles, elements of the discrete homology group,
using ideas of \Cite{S}. We construct the space of discrete cycles as a certain
space of \fn/s. Having a representative of a discrete cohomology class (a \fn/)
and a discrete cycle (a \fn/ again) we define the pairing (the \hpair/) between
the cohomology class and the cycle as an integral of their product with
a certain fixed ``\hgeom/ \phf/'' over a certain fixed contour of the middle
dimension. We show that there are enough discrete cycles and they form the
space dual to the quotient space of the space of our discrete closed forms
modulo discrete coboundaries. To prove this we compute the determinant of
the period matrix and surprisingly get an explicit formula \(mu<>0) for
the determinant analogous to the determinant formulae for the ``continuous''
\hgeom/ \fn/s \Cite{V1}, \cf. the Loeser determinant formula for the Frobenius
transformation \Cite{L}. The form of our discrete cycles suggests a natural
identification of the space of our discrete cycles with a tensor product of
\Umod/s and this identification allows us to prove the result on \traf/s
between \asol/s.
\Par
As we know the \qKZe/ turns into the \difl/ \KZe/ under the \qcl/ limit.
We show that our discretization of geometry under the \qcl/ limit turns into
the geometry of the \difl/ \KZe/: representatives of our discrete cohomology
classes turn into closed \difl/ forms, our discrete cycles turn into ``honest''
topological cycles.
\Par
Note in conclusion, that our \sol/s to the \qKZe/ in the special case
considered in \Cite{S} are close to the \sol/s constructed therein, but
different. It is also worth mentioning that our description of \traf/s
indicates \qlo/ symmetries in the model of quantum field theory considered
in \Cite{S}.
\Par
The paper is organized as follows. Sections~{\SNo{2}\,--\,\SNo{2e}} contain
constructions and statements. In Section~\SNo{3} we consider the special case
of \onedim/ \hgeom/ \fn/s of the \MB/ type to illustrate ideas and proofs.
Section~\SNo{4} contains proofs in the multidimensional case.
\Par
Parts of this work had been written when the authors visited the University
of Tokyo, the Kyoto University, the University Paris VI, \'Ecole Normale
Sup\'erieure de Lyon, the MSRI at Berkeley. The authors thank those
institutions for hospitality. The authors thank \Feld/, P\]\&Etingof
and E\&Mukhin for valuable discussions.

\Sect[2]{Discrete flat connections and \loc/s}
\subsect*{Discrete flat connections}
Consider a complex vector space $\Cn$ called the \em{base space}. Fix a nonzero
complex number $p$ called the \em{step}. The lattice $\Z^n$ acts on the base
space by translations $z\map z+p\)l$ where $l\in\Z^n\sub\Cn$. Let $\BB$ be
an \inv/ subset of the base space. Say that there is a \em{bundle with
a discrete connection} over $\BB$ if for any $z\in\BB$ there are a vector space
$V(z)$ and linear \iso/s
$$
A_m(\zn):V(\zmn)\,\to\,V(\zn)\,,\qqq\mn\,.
$$
\ifMag\else\newpage\nt\fi
The connection is called \em{flat} (or \em{integrable}) if the \iso/s
$A_1\lc A_n$ commute:
\ifMag
$$
\align
A_l(\zn)\> &A_m(z_1\lc z_l\]+p\lc z_n)\;=
\Tag{flat}
\\
&=\;A_m(\zn)\>A_l(\zmn)\,.
\endalign
$$
\else
$$
A_l(\zn)\>A_m(z_1\lc z_l\]+p\lc z_n)\,=\,A_m(\zn)\>A_l(\zmn)\,.
\Tag{flat}
$$
\fi
Say that a \em{discrete subbundle} in $\BB$ is given if a subspace in every
fiber is distinguished and the family of subspaces is \inv/ \wrt/
the connection.
\par
A section $s:z\map s(z)$ is called \em{periodic} (or \em{horizontal})
if its values are \inv/ \wrt/ the connection:
$$
A_m(\zn)\>s(\zmn)\,=\,s(\zn)\,,\qqq\mn\,.
\Tag{period}
$$
A \fn/ $f(\zn)$ on the base space is called a \em{quasiconstant} if
$$
f(\zmn)\,=\,f(\zn)\,,\qqq\mn\,.
$$
Periodic sections form a module over the ring of quasiconstants.
\par
The \em{dual bundle} with the \em{dual connection} has fiber
$V^*(z)$ and \iso/s
$$
A_m^*(\zn):V^*(\zn)\,\to\,V^*(\zmn)\,.
$$
Let $s_1\lc s_N$ be a basis of sections of the initial bundle.
Then the \iso/s $A_m$ of the connection are given by matrices $\A\"m$:
$$
A_m(\zn)\>s_k(\zmn)\,=\,\sum_{l=1}^N \A\"m_{kl}(\zn)\>s_l(\zn)\,.
$$
For any section ${\psi:z\map\psi(z)}$ of the dual bundle,
denote by ${\Psi:z\map\Psi(z)}$ its coordinate vector,
$\Psi_k(z)=\bra\psi(z),s_k(z)\ket$.
\par
The section $\psi$ is periodic if and only if its coordinate vector satisfies
the system of \deq/s
$$
\Psi(\zmn)\,=\,\A\"m(\zn)\>\Psi(\zn)\,,\qqq\mn\,.
$$
Moreover, all \sol/s to the system have this form. This system of \deq/s is
called the \em{periodic section \eq/}.
\Par
Say that \fn/s $\phin$ in \var/s $\zn$ form a \em{\sconn/} if
$$
\phi_l(\zmn)\>\phi_m(\zn)\,=\,\phi_m(z_1\lc z_l\]+p\lc z_n)\>\phi_l(\zn)
$$
for all $l,m$. These \fn/s define a connection on the trivial complex \onedim/
vector bundle.
\par
There is a simple construction of \conn/s. Fix arbitrary \fn/s $\pho_{lm}$,
$l<m$, in one \var/ and nonzero complex numbers $\ka_m$. Set
$$
\phi_m(\zn)\,=\,
\ka_m\>\Bigl[\,\prod_{1\le l<m}\!\pho_{lm}(z_l-z_m-p)\>\Bigr]\1\>
\prod_{m<l\le n}\>\pho_{ml}(z_m-z_l)\,.
$$
The \sconn/ of this form is called \em{decomposable}, the \fn/s $\pho_{lm}$
are called \em{\prim/s} and $\ka_m$ are called \em{scaling parameters}.
\Par
A \fn/ $\Phi(\zn)$ is called a \em{\phf/} of a \sconn/ if
$$
\Phi(\zmn)\,=\,\phi_m(\zn)\>\Phi(\zn)\,,\qqq\mn\,.
$$
Similarly, a \fn/ $\Pho(x)$ is called a \em{\phf/} of a \fn/ $\pho(x)$ in
one \var/ if $\Pho(x+p)=\pho(x)\>\Pho(x)$. Note that the \phf/s are not unique.
\par
If the \conn/s are decomposable, if $\Pho_{lm}$ are \phf/s of \prim/s,
and if $K_m$ are \phf/s of scaling parameters, then
$$
\Phi(\zn)\,=\,\pron K_m(z_m)\,\prod_{l<m}\Pho_{lm}(z_l-z_m)
$$
is a \phf/ of the system of \conn/s.
\Par
For any \fn/ $f(\zn)$ define new \fn/s $Q_1f\lc Q_nf$ and $D_1f\lc D_nf$ by
$$
\gather
(Q_mf)(\zn)\,=\,\phi_m(\zn)\>f(\zmn)\,,
\\
\Text{and}
D_mf\,=\,Q_mf-f\,.
\endgather
$$
The \fn/s $D_1f\lc D_nf$ are the \em{discrete partial derivatives} of the \fn/
$f$. We have $D_lD_mf=D_mD_lf$.
\par
Let $\Ff$ be a vector space of \fn/s of $\zn$ \st/ the operators $Q_1\lc Q_n$
induce linear \iso/s of $\Ff$:
$$
Q_m:\Ff\,\to\,\Ff\,.
$$
Say that the space $\Ff$ and the \conn/s $\phin$ form a \onedim/ \em{\dloc/}
on $\Cn\!$. $\Ff$ is called the \em{\fn/al space} of the \loc/.
\Par
Define the \em{de~Rham complex} $\bigl(\Omb(\Ff),D\bigr)$ of the \loc/
in a standard way. Namely, set
$$
\Om^a\,=\,\bigl\lb\>\om=\tsum_{k_1\lc k_a}f_{k_1\lc k_a}\>
Dz_{k_1}\lsym\wedge Dz_{k_a}\>\bigr\rb
$$
where $Dz_1\lc Dz_n$ are formal symbols and the coefficients $f_{k_1\lc k_a}$
belong to $\Ff$. Define the \difl/ of a \fn/ by $Df=\sun D_mf\>Dz_m$,
and the \difl/ of a form by
$$
D\)\om\;=
\tsum_{k_1\lc k_a}Df_{k_1\lc k_a}\wedge Dz_{k_1}\lsym\wedge Dz_{k_a}\,.
$$
The cohomology groups $H^1\lc H^n$ of this complex are called the
\em{cohomology groups of $\Cn$ with coefficients in the \dloc/}. In particular,
the top cohomology group is $H^n=\Ff/D\Ff$ where $D\Ff=\sun D_m\Ff$. The dual
spaces $H_a=(H^a)^*$ are called the \em{homology groups}.
\Par
There is a geometric construction of bundles with discrete flat connections.
This is a discrete version of the \GM/ construction.
\par
Let $\pi:\Cln\to\Cn$ be an affine projection onto the base with fiber $\Cl\!$.
$\Cln$ will be called the \em{total space}. Let $\zn$ be coordinates on the
base, $\tell$ coordinates on the fiber, so that $\tell$, $\zn$ are coordinates
on the total space. When it is convenient, we will denote the coordinates $\zn$
by $t_{\ell+1}\lc t_{\ell+n}$.
\par
Let $\Ff$, $\phi_1\lc\phi_{\ell+n}$ be a \loc/ on $\Cln\!$. For a point
$z\in\Cn$ define a \loc/ $\Ff(z)$, $\phi_a(\cdot;z)$, $\aell$, on the fiber
over $z$. Set
$$
\Ff(z)\,=\,\lb\>f\vst{\pi\1(z)}\vert f\in\Ff\>\rb\qquad
\text{and}\qquad
\phi_a(\cdot;z)\,=\,\phi_a\vst{\pi\1(z)}\,.
$$
The de~Rham complex, cohomology and homology groups of the fiber are denoted
by $\bigl(\Omb(z),D(z)\bigr)$, $H^a(z)$ and $H_a(z)$, \resp/.
\par
There is a natural \hom/ of the de~Rham complexes
$$
(\Omb(\Cln\!,\Ff),D)\,\to\,(\Omb(z),D(z))\,,\qqq \om\map\om|_{\pi\1(z)}\,,
$$
where the restriction of a form is defined in a standard way: all symbols
$Dz_1\lc Dz_n$ are replaced by zero and all coefficients of the remaining
monomials $Dt_{k_1}\lsym\wedge Dt_{k_a}$ are restricted to the fiber.
\par
For a fixed $a$ the vector spaces $H^a(z)$ form a bundle with a discrete flat
connection. The linear maps
$$
A_m(\zn):H^a(\zmn)\,\to\,H^a(\zn)
$$
are defined as follows. Define $Q_m:\Om^a(\Cln,\Ff)\to\Om^a(\Cln,\Ff)$ by
$$
\om\;\map\tsum_{k_1\lc k_a}Q_mf_{k_1\lc k_a}\>
Dz_{k_1}\lsym\wedge Dz_{k_a}\,.
$$
\vgood
Then $Q_m$ induces a homomorphism of the de~Rham complexes
$$
\bigl(\Omb(\zmn),D(\zmn)\bigr)\,\to\,\bigl(\Omb(\zn),D(\zn)\bigr)\,.
$$
We set $A_m(\zn)$ to be equal to the induced map of the cohomology spaces.
This connection is called the \em{discrete \GM/}.
\par
The \GM/ on the cohomological bundle induces the dual flat connection
on the homological bundle:
$$
A^*_m(\zn):H_a(\zn)\,\to\,H_a(\zmn)\,.
$$
\par
In this paper we study the \GM/ for a class of \dloc/s.

\subsect{Connection coefficients of \loc/s}
There are three important classes of \loc/s: \rat/, \tri/ and elliptic.
\par
Consider a \loc/ with decomposable \conn/s and \prim/s of the form
$$
\pho_{ab}(x)\,=\;{\tau(x+\al_{ab})\over \tau(x+\bt_{ab})}
$$
where $\tau(x)$ is a \fn/ of one \var/ and $\al_{ab},\bt_{ab}$ are suitable
complex numbers. A \loc/ is called \em{\rat/}, \em{\tri/} or
\em{elliptic} if
$$
\tau(x)=x\,,\qqq\tau(x)=\sin(\gm x)\,,\qqq\tau(x)=\tht(\gm x)\,,
$$
\resp/. Here $\tht(x)$ is a theta-\fn/ and $\gm$ is a nonzero complex
number. Note that $\tau(x)=\gm x$ for all $\gm\ne 0$ gives the same \prim/s.
\par
Say that a decomposable \sconn/ on the total space is of the \em{\stype/} if
the constants $\al_{ab}, \bt_{ab}$, and the scaling parameters
$\ka_1\lc\ka_{\ell+n}$ have the following form:
$$
\alignat2
\al_{ab}\, &{}=\,-\bt_{ab}\,=\,-h\qqq &&\for\quad a<b\le\ell\,,
\Tag{sl2type}
\\
\al_{ab}\, &{}=\,-\bt_{ab}\,=\,\La_{b-\ell}\qqq &&\for\quad a\le\ell<b\,,
\\
\al_{ab}\, &{}=\,-\bt_{ab}\,=\,0\qqq &&\for\quad \ell<a<b\,,
\\
\ka_a\, &{}=\,\ka\qqq &&\for\quad a\le\ell\,,
\\
\ka_a\, &{}=\,1\qqq &&\for\quad \ell<a\,.
\endalignat
$$
Such a \sconn/ depends on $n+2$ complex numbers $\Lan$, $\ka,h$.
\Par
In this paper we study \rat/ systems of the \stype/, for the \tri/ case see
\Cite{TV3} and for the elliptic case see \Cite{FTV1}, \Cite{FTV2}.
\Par
The \prim/s of a \rat/ \stype/ \loc/ have the form
$$
\NN3>
\alignat2
\pho_{ab}(x)\,&{}=\;{x-h\over x+h}\qqq &&\for\quad a<b\le\ell\,,
\\
\pho_{ab}(x)\,&{}=\;{x+\La_{b-\ell}\over x-\La_{b-\ell}}\qqq &&
\for\quad a\le\ell<b\,,
\\
\pho_{ab}(x)\,&{}=\;1\qqq &&\for\quad \ell<a<b\,.
\endalignat
$$
Rescaling $\Lan$ and $x$ we can set $h=1$, so we assume that the \prim/s of
a \rat/ \stype/ \loc/ have the form
$$
\NN3>
\alignat2
\pho_{ab}(x)\,&{}=\;{x-1\over x+1}\qqq &&\for\quad a<b\le\ell\,,
\\
\pho_{ab}(x)\,&{}=\;{x+\La_{b-\ell}\over x-\La_{b-\ell}}\qqq &&
\for\quad a\le\ell<b\,,
\\
\pho_{ab}(x)\,&{}=\;1\qqq &&\for\quad \ell<a<b\,.
\endalignat
$$
The \conn/s of a \rat/ \stype/ \loc/ have the form
$$
\NN4>
\gather
\phi_a(t,z)\,=\,\ka\pron{t_a-z_m+\La_m\over t_a-z_m-\La_m}\,
\prod_{a<b\le\ell}{t_a-t_b-1\over t_a-t_b+1}\;
\prod_{1\le b<a}{t_a-t_b-1+p\over t_a-t_b+1+p}\;,
\ifMag\else\qqq\aell\,,\fi
\\
\ifMag
\nn-2>
\rline{$\aell$,}
\\
\fi
\phi_{\ell+m}(t,z)\,=\,
\pral\,{t_a-z_m-\La_m-p\over t_a-z_m+\La_m-p}\;,\qqq\mn\,.
\endgather
$$
A \phf/ of a \prim/ $(x+\al)/(x-\al)$ has the form
$$
\Pho(x;\al)\,=\,{\Gmpp{x+\al}\over\Gmpp{x-\al}}
\Tag{Phox}
$$
and, therefore, a \phf/ of the \sconn/ is given by
\ifMag
$$
\align
\Phi( &\tell,\zn)\;=
\Tagg{Phi}
\\
&{}=\,\exp\bigl(\mu\}\tsum_{a=1}^\ell\}t_a/p\bigr)
\pron\,\pral\>\Pho(t_a-z_m;\La_m)\prab\Pho(t_a-t_b;-1)
\endalign
$$
\else
$$
\kern1.5em
\Phi(\tell,\zn)\,=\,\exp\bigl(\mu\}\tsum_{a=1}^\ell\}t_a/p\bigr)
\pron\,\pral\>\Pho(t_a-z_m;\La_m)\prab\Pho(t_a-t_b;-1)
\kern-1.5em
\Tagg{Phi}
$$
\fi
where parameters $\ka$ and $\mu$ are connected by the \eq/ $\ka=e^\mu$.
\Par
The Stirling formula gives the following \as/s for the \phf/ \(Phox) of
the \prim/ as $x\to\8$:
$$
\Pho(x;\al)\,=\,(x/p)^{2\al/p}\>\ono\,,\qqq |\arg(x/p)|<\pi\,.
\Tag{Stir}
$$
This formula defines \as/s at infinity of the \phf/ of the \sconn/.
\par
The \phf/ \(Phox) of the \prim/ has a symmetry property
$$
\Pho(-x;\al)\,=\,\Pho(x;\al)\;{(x+\al)\>\Sinp{x+\al}\over(x-\al)\>\Sinp{x-\al}}
$$
which leads to a symmetry property
$$
\align
& \Phi(\tall,\zn)\;=
\Tag{Phisym}
\\
\nn6>
&\!{}=\,\Phi(\tell,\zn)\;{(t_a-t_{a+1}-1)\>\Sinp{t_a-t_{a+1}-1}\over
(t_a-t_{a+1}+1)\>\Sinp{t_a-t_{a+1}+1}}
\endalign
$$
of the \phf/ of the \sconn/. This property later motivates definitions \(act)
and \(triact) of certain actions of the \symg/.

\subsect{The \fn/al space of a \rat/ \{\stype/ \loc/}
\par
Define the \fn/al space $\Fun$ of a \rat/ \stype/ \loc/ as the space of \raf/s
on the total space with at most simple poles at the following hyperplanes:
$$
\alignat2
& t_a=z_m-\La_m+(s+1)p\,, &\qqq &t_a=z_m+\La_m-sp\,,
\Tag{list1}
\\
\nn2>
& t_a=t_b-1-(s+1)p\,, && t_a=t_b+1+sp\,,
\endalignat
$$
$1\le b<a\le\ell$, $\mn$, $s\in\Zp$. It is easy to check that the \fn/al space
\vv.1>
is \inv/ \wrt/ all operators $Q_1^{\pm1}\lc Q_n^{\pm1}\!$.
\par
Define an action of the \symg/ $\Sl$ on the \fn/al space:
$$
\si:\Fun\,\to\,\Fun,\qqq f\map[\)f\)]_\si\,,\qqq\si\in\Sl,
\Tag{act}
$$
by the following action of simple transpositions:
$$
[\)f\)]_{(a,a+1)}(\tell,\zn)\,=\,f(\tall,\zn)\;
{t_a-t_{a+1}-1\over t_a-t_{a+1}+1}\;,
$$
$\aell-1$. The operators $Q_1\lc Q_{\ell+n}$ and $D_1\lc D_{\ell+n}$ commute
with the action of the \symg/.
\par
We extend the \SL-action to the de~Rham complex assuming that it respects
the exterior product and
$$
\si:Dt_a\,\map\,Dt_{\si_a}\,,\qqq\si:Dz_m\,\map\,Dz_m\,,\qqq\si\in\Sl.
$$
The same formulae define an action of the \symg/ on the de~Rham complex of
a fiber. The \hom/ of the restriction of the de~Rham complex of the total space
to the de~Rham complex of a fiber commutes with the action of the \symg/.
The action of the \symg/ induces an action of the \symg/ on the homology and
cohomology groups. The \GM/ commutes with this action.
\par
If a \symg/ acts on a vector space $V\!$, we denote by $V_{\]\sss\Sig}$
the subspace of \inv/ vectors and by $V_{\}\sss A}$ the subspace of skew-\inv/
vectors.
\par
In this paper we are interested in the skew-\inv/ part $\Has^\ell(z)$ of
the top cohomology group of a fiber. This subspace is \gby/ forms
${f\)\Dtel}$ where $f$ runs through the space $\Fsi(z)$ of \inv/ \fn/s.
\par
Introduce an important \em{\rhgf/} $\Fo\sub\Fsi$ as the subspace of \fn/s of
the form
$$
P(\tell,\zn)\,\pron\,\pral\,{1\over t_a-z_m-\La_m\!}\;
\prab\ {t_a-t_b\over t_a-t_b+1}
$$
where $P$ is a \pol/ with complex coefficients which is \sym/ in \var/s $\tell$
and has degree less than $n$ in each of the \var/s $\tell$. The restriction of
the \rhgf/ to a fiber defines the \em{\rhgf/} $\Fo(z)\sub\Fsi(z)$ of the fiber
which is a complex \fd/ vector space. A form ${f\)\Dtel}$ with the coefficient
in the \rhgf/ of a fiber is called a \em{\hform/}.
\par
The subspace $\H(z)\sub\Has^\ell(z)$ of the top cohomology group of a fiber
\gby/ the \hform/s is called the \em{\hcg/}.
\par
The union of the hyperplanes
$$
z_l+\La_l-z_m+\La_m\,=\, r+p\)s\,,\qqq r=0\lc\ell-1\,,\qquad s\in\Z\,,
\Tag{dis}
$$
$l,\mn$, $l\ne m$, in the base space $\Cn$ is called the \em{discriminant}.
The complement to the discriminant will be denoted by $\BB$.
\Th{subbundle}
\back\Cite{V3},\;\Cite{TV1}
The family of subspaces ${\lb\>\H(z)\>\rb}_{z\in\BB}$ is \inv/ \wrt/ the \GM/
and, therefore, defines a discrete subbundle.
\endpro
\nt
This subbundle will be called the \em{\hgeom/ subbundle}.
\Par
Later on we often make the following assumptions. We assume that the step $p$
is real negative and \st/
$$
\npZ\,,
\Tag{npZ}
$$
the \wt/s $\Lan$ are \st/
$$
2\La_m-s\,\nin\,p\)\Z\,,\qqq \mn\,,\qquad s=1-\ell\lc\ell-1\,,
\Tag{Lass}
$$
and the coordinates $\zn$ obey the condition
$$
z_l\pm\La_l-z_m\pm\La_m-s\,\nin\,p\)\Z\,,\qqq l,\mn\,,\qquad l\ne m\,,
\Tag{assum}
$$
for any $s=1-\ell\lc\ell-1$ and for an arbitrary combination of signs.
\Th{kanon}
Let $\ka\ne 1$. Let $p<0$. Let \(npZ)\;--\;\(assum) hold. Then
$$
\dim\H(z)\,=\,\dim\Fo(z)\,=\,{n+\ell-1\choose n-1}\,.
$$
\vvv-1.5>
\endpro
\nt
This means that
\vv-.5>
$$
\H(z)\,\simeq\,\Fo(z)\,.
\Tag{HFo}
$$
\vvv->
\goodbreak
\Th{kaon}
Let $\ka=1$. Let $p<0$. Let \(npZ)\;--\;\(assum) hold.
If $\tLmn-s\nin{p\)\Znn}$ for all $s=\ell-1\lc 2\ell-2$,
then $\dsize\dim\H(z)\,=\,{n+\ell-2\choose n-2}$.
\endpro
\nt
Theorems~\[kanon] and \[kaon] are proved in Section~\SNo{4}.
\Par
Theorem~\[kanon] means that if the scaling parameter $\ka$ is not equal to $1$,
then every nonzero \hform/ defines a nonzero cohomology class. On the contrary,
if $\ka=1$, then by Theorem~\[kaon] there are exact \hgeom/ forms. We describe
them in Lemma~\[wDw].

\subsect{Bases in the \rhgf/ of a fiber}
The \fd/ \rhgf/ $\Fo(z)$ of a fiber has $n!$ remarkable bases. These bases
will allow us to identify geometry of an \stype/ \loc/ with \rep/ theory.
The bases are labelled by elements of the \symg/ $\S^n\!$. First we define
the basis corresponding to the unit element of the \symg/.
\par
Let
$$
\gather
\\
\nn-24>
\Zln\,=\,\lb\>\lg\in\Zpn\vert\tsun\lg_m=\ell\>\rb\,.
\Tag{Zln}
\\
\nn-8>
\endgather
$$
Set $\lg^m=\sum_{k=1}^m\lg_k$. In particular, $\lg^0=0$, $\lg^n=\ell$.
For any $\lg\in\Zln$ define a \raf/ $w_\lg\in\Fo$ as follows:
\ifMag
$$
\align
\kern1.5em
w_\lg( & \tell,\zn)\;=
\Tagg{wlga}
\\
\nn2>
&{}=\sum_{\si\in\Sl}\,\Bigl[\,\pron\>{1\over\lg_m!}\>
\prod_{a\in\Gm_{\Rph lm}}\Bigl(\,{1\over t_a-z_m-\La_m\!}
\ \prlm\,{t_a-z_l+\La_l\over t_a-z_l-\La_l}\,\Bigr)\>\Bigr]_\si
\kern-1.5em
\endalign
$$
\else
$$
\kern1.5em
w_\lg(\tell,\zn)\,=
\sum_{\si\in\Sl}\,\Bigl[\,\pron\>{1\over\lg_m!}\>
\prod_{a\in\Gm_{\Rph lm}}\Bigl(\,{1\over t_a-z_m-\La_m\!}
\ \prlm\,{t_a-z_l+\La_l\over t_a-z_l-\La_l}\,\Bigr)\>\Bigr]_\si
\kern-1.5em
\Tagg{wlga}
$$
\fi
where $\Gm_m=\lb 1+\lg^{m-1}\,\lc\lg^m\rb$, $\mn$. The \fn/s $w_\lg$ are called
the \em{\rwf/s}.
\Ex
For $\ell=1$ the \rwf/s have the form
$$
w_{\eg(m)}(t_1,\zn)\,=\,
{1\over t_1-z_m-\La_m\!}\ \prlm\,{t_1-z_l+\La_l\over t_1-z_l-\La_l}
$$
where $\eg(m)=(0\lc\%1_{\sss\^{$m$-th}}\lc 0)$, $\mn$.
\enddemo
\Ex
For $n=1$ the \fn/ $w\'\ell$ has the form
$$
w\'\ell(\tell,z_1)\,=\,
\pral\,{1\over t_a-z_1-\La_1\}}\;\prab\ {t_a-t_b\over t_a-t_b+1}\;.
$$
\enddemo
\Ex
For $\ell=2$ and $n=2$ the \fn/s $w_\lg$ have the form
$$
\NN6>
\align
w_{(2,0)}(t_1,t_2,z_1,z_2)\, &{}=\;
{1\over(t_1-z_1-\La_1)\>(t_2-z_1-\La_1)}\,{t_1-t_2\over t_1-t_2+1}\;,
\\
\ald
w_{(1,1)}(t_1,t_2,z_1,z_2)\, &{}=\;{1\over(t_1-z_1-\La_1)\>(t_2-z_2-\La_2)}\,
{t_2-z_1+\La_1\over t_2-z_1-\La_1}\;+{}
\\
\nn-2>
&{}\)+\;{1\over(t_2-z_1-\La_1)\>(t_1-z_2-\La_2)}\,
{t_1-z_1+\La_1\over t_1-z_1-\La_1}\,{t_1-t_2-1\over t_1-t_2+1}\;,
\\
\ald
w_{(0,2)}(t_1,t_2,z_1,z_2)\, &{}=\;{1\over(t_1-z_2-\La_2)\>(t_2-z_2-\La_2)}\,
{(t_1-z_1+\La_1)\>(t_2-z_1+\La_1)\over(t_1-z_1-\La_1)\>(t_2-z_1-\La_1)}\,
{t_1-t_2\over t_1-t_2+1}\;.
\endalign
$$
\enddemo
\Lm{wbasis}
The \fn/s $w_\lg$, $\lg\in\Zln$, restricted to a fiber over $z$ form a basis
in the \rhgf/ $\Fo(z)$ of the fiber provided that for any $s=0\lc\ell-1$,
$$
z_l-\La_l-z_m-\La_m+s\,\ne\,0\,,\qqq 1\le l<m\le n\,.
$$
\endpro
\nt
Lemma~\[wbasis] is proved in Section~\SNo{4}.
\Lm{wDw}
Let $\ka=1$. Then for any $\lg\in\Zll$ the following relation holds:
$$
\sun(\lg_m+1)\>(2\La_m-\lg_m)\>w_{\lg+\eg(m)}\,=\,
\sum_{a=1}^\ell\,D_a\bigl[\)w_\lg(t_2\lc t_\ell)\)\bigr]_{(1,a)}\,,
$$
where $(1,a)\in\Sl$ are transpositions. Moreover, if $\Rc(z)$ is the subspace
in $\Fo(z)$ \gby/ the elements in \lhs/ of the relations, then
$$
\dim\Fo(z)/\Rc(z)\,=\,{n+\ell-2\choose n-2}
$$
provided that \ ${z_l-\La_l-z_m-\La_m+s\ne 0}$, \ $1\le l\le m\le n$, \ for
any $s=0\lc\ell-1$.
\endpro
\nt
The subspace $\Rc(z)\sub\Fo(z)$ is called the \em{\cosub/}.
\par
The relations \(wDw) induce relations
$$
\sun\lf(\lg_m+1)\>(2\La_m-\lg_m)\>w_{\lg+\eg(m)}\)\Dtel\rf\,=\,0\,,\qqq
\lg\in\Zll\,,
$$
in the cohomology group $H^\ell(z)$, where $\lf\al\rf$ denotes
the cohomological class of a form $\al$. For $\ka=1$ under assumptions
of Theorem~\[kaon] we have
$$
\H(z)\,\simeq\,\Fo(z)/\Rc(z)\,.
\Tag{HFR}
$$
\goodbe
\par
For any \perm/ ${\tau\in\S^n}$ we define a basis
${\lb\)w^\tau_\lg\)\rb_{\lg\)\in\Zln}\:}$ in the \rhgf/ of a fiber
by formulae similar to \(wlga). Namely, set
$$
\ifMag\kern2em\else\kern1em\fi
w^\tau_\lg(\tell,\zn;\Lan)\,=\,w_{{\vp(}^\tau\!\lg}^{\vp1}
(\tell,z_{\tau_1}\lc z_{\tau_n};\La_{\tau_1}\lc\La_{\tau_n})
\ifMag\kern-2em\else\kern-1em\fi
\Tag{wtau}
$$
where $^\tau\}\lg=(\lg_{\tau_1}\lc\lg_{\tau_n})$.
\goodbm
\Ex
For ${\ell=1}$ and a \perm/ ${\tau=(n,n-1\lc 1)}$ the \fn/s have the form
$$
w^\tau_{\eg(m)}(t_1,\zn)\,=\,
{1\over t_1-z_m-\La_m\!}\ \prod_{m<l\le n}{t_1-z_l+\La_l\over t_1-z_l-\La_l}\;.
$$
\enddemo

\subsect{The \thgf/}
In our study of the \GM/ an important role is played by the following
\em{\thgf/}. The \thgf/ is a \tri/ counterpart of the \rhgf/
introduced above.
\Par
The \thgf/ $\Fq$ is the space of \fn/s of \var/s $\tell,\zn$ which have
the form
\ifMag
$$
\gather
{\align
\kern1.5em
P( & \xill,\ztn)\;\x
\Tagg{thgf}
\\
\nn2>
&{}\x \pron\,\pral\,{\Expp{z_m-t_a}\over\Sinp{t_a-z_m-\La_m}}\;
\prab\ {\Sinp{t_a-t_b}\over\Sinp{t_a-t_b+1}}
\kern-1.5em
\endalign}
\\
\Text{where}
\xi_a=\exp(2\pii t_a/p)\,,\qqq \zt_m=\exp(2\pii z_m/p)\,,
\endgather
$$
\else
$$
\gather
\kern1.5em
P(\xill,\ztn)\,\pron\,\pral\,{\Expp{z_m-t_a}\over\Sinp{t_a-z_m-\La_m}}\;
\prab\ {\Sinp{t_a-t_b}\over\Sinp{t_a-t_b+1}}
\kern-1.5em
\Tag{thgf}
\\
\Text{where}
\xi_a=\exp(2\pii t_a/p)\,,\qqq \zt_m=\exp(2\pii z_m/p)\,,
\endgather
$$
\fi
and $P$ is a \pol/ with complex coefficients which is \sym/ in \var/s $\xill$
and has degree less than $n$ in each of the \var/s $\xill$.
\par
Introduce the \em{\sthgf/} $\Fqs\sub\Fq$ as the space of \fn/s of the form
\(thgf) \st/ the polynomial $P$ is divisible by the product
$\xi_1\ldots\xi_\ell$.
\par
The restriction of the \thgf/s to a fiber defines the \em{\thgf/s}
$\Fqs(z)\sub\Fq(z)$ of the fiber. The \thgf/ $\Fq(z)$ is a complex \fd/ vector
space of the same dimension as the \rhgf/ of the fiber.
\par
The \thgf/s of fibers over $z$ and $z'$ are naturally identified if the points
$z$ and $z'$ lie in the same orbit of the \^{$\Z^n\]$-}action on the base
space, since all elements of the \thgf/ are \p-periodic \fn/s.
\Par
Introduce a new action of the \symg/ $\Sl$ on \fn/s,
$$
f\,\map\,\lbc\)f\)\rbc\vpp\si\,,\qqq\si\in\Sl,
\Tag{triact}
$$
by the following action of simple transpositions:
\ifMag
$$
\align
\lbc\)f\)\rbc_{(a,a+1)} & (\tell,\zn)\;=
\\
&{}\!=\,f(\tall,\zn)\;{\Sinp{t_a-t_{a+1}-1}\over\Sinp{t_a-t_{a+1}+1}}\;,
\endalign
$$
\else
$$
\lbc\)f\)\rbc_{(a,a+1)}(\tell,\zn)\,=\,
f(\tall,\zn)\;{\Sinp{t_a-t_{a+1}-1}\over\Sinp{t_a-t_{a+1}+1}}\;,
$$
\fi
$\aell-1$. The \thgf/ is \inv/ \wrt/ this action. The action commutes with
the restriction of \fn/s to a fiber.
\Par
The \thgf/ of a fiber has $n!$ remarkable bases. The bases are labelled by
elements of the \symg/ $\S^n\!$. First we define the basis corresponding to
the unit element of the \symg/.
For any $\lg\in\Zln$ define a \fn/ $W_\lg\in\Fq$ as follows:
$$
\align
\kern1.5em
\ifMag\else W_\lg( \fi
&
\ifMag\kern1.7em W_\lg(\fi
\tell,\zn)\,=\>\pron\,\prod_{s=1}^{\lg_m}\,{\sin(\pi/p)\over\sin(\pi s/p)}\ \x
\Tagg{Wlga}
\\
\nn2>
&{}\x\sum_{\si\in\Sl}\,\LBc\,\pron\,
\prod_{a\in\Gm_{\Rph lm}}\Bigl(\,{\Expp{z_m-t_a}\over\Sinp{t_a-z_m-\La_m}}
\ \prlm\,{\Sinp{t_a-z_l+\La_l}\over\Sinp{t_a-z_l-\La_l}}\,\Bigr)\>\RBc_\si
\kern-1.5em
\endalign
$$
where $\Gm_m=\lb 1+\lg^{m-1}\,\lc\lg^m\rb$, $\mn$. Also for any $\lg\in\Zlm\!$
define a \fn/ $\Wo_\lg\in\Fqs$ as follows:
\ifMag
$$
\NN4>
\align
\kern2em
\Wo_\lg &(\tell,\zn)\;=
\Tagg{Wlgo}
\\
&{}\x\prmn\,\prod_{s=1}^{\lg_m}\,{\sin(\pi/p)\over\sin(\pi s/p)}\,
\Sinp{z_m-\La_m-z_{m+1}-\La_{m+1}+s-1}\>\x{}
\\
&{}\x\sum_{\si\in\Sl}\,\LBc\,\prmn\,\prod_{a\in\Gm_{\Rph lm}}
\Bigl(\,{1\over\Sinp{t_a-z_m-\La_m}\>\Sinp{t_a-z_{m+1}-\La_{m+1}}}\,\x{}
\kern-2em
\\
&\kern4em{}\x\,\prlm\,
{\Sinp{t_a-z_l+\La_l}\over\Sinp{t_a-z_l-\La_l}}\,\Bigr)\>\RBc_\si.
\endalign
$$
\else
$$
\NN4>
\gather
\Rline{\Wo_\lg(\tell,\zn)\,=\>\prmn\,\prod_{s=1}^{\lg_m}\,
{\sin(\pi/p)\over\sin(\pi s/p)}\,\Sinp{z_m-\La_m-z_{m+1}-\La_{m+1}+s-1}
\>\x{}\,\;\hp.}
\Tag{Wlgo}
\\
\Rline{{}\x\sum_{\si\in\Sl}\,\LBc\,\prmn\,\prod_{a\in\Gm_{\Rph lm}}
\Bigl(\,{1\over\Sinp{t_a-z_m-\La_m}\>\Sinp{t_a-z_{m+1}-\La_{m+1}}}
\,\x{}\,\;\hp.}
\\
\Rline{{}\x\,\prlm\,{\Sinp{t_a-z_l+\La_l}\over\Sinp{t_a-z_l-\La_l}}
\,\Bigr)\>\RBc_\si\,.\}}
\endgather
$$
\fi
The \fn/s $W_\lg$ and $\Wo_\lg$ are called the \em{\twf/s}.
\Lm{Wbasis}
The \fn/s $W_\lg$, $\lg\in\Zln$, restricted to a fiber over $z$ form a basis
in the \thgf/ $\Fq(z)$ of the fiber, provided that for any $s=0\lc\ell-1$,
$$
z_l-\La_l-z_m-\La_m+s\,\nin\,p\)\Z\,,\qqq 1\le l<m\le n\,.
$$
\endpro
\Lm{Wbasiso}
The \fn/s $\Wo_\mg$, $\mg\in\Zlm\!$, restricted to a fiber over $z$ form
a basis in the \sthgf/ $\Fqs(z)$ of the fiber, provided that for any
$s=0\lc\ell-1$,
$$
z_l-\La_l-z_m-\La_m+s\,\nin\,p\)\Z\,,\qqq 1\le l<m\le n\,.
$$
\endpro
\nt
Lemmas~\[Wbasis], \[Wbasiso] are proved in Section~\SNo{4}.
\Ex
For $\ell=1$ the \twf/s have the form
$$
W_{\eg(m)}(t_1,\zn)\,=\,{\Expp{z_m-t_1}\over\Sinp{t_1-z_m-\La_m}}
\ \prlm\,{\Sinp{t_1-z_l+\La_l}\over\Sinp{t_1-z_l-\La_l}}\;.
$$
The \sthgf/ ${\Fqs(z)\sub\Fq(z)}$ has dimension ${(n-1)}$ and is \gby/
the \fn/s
$$
\Wo_{\eg(m)}\,=\,W_{\eg(m)}\exp(-\pii\La_m/p)-W_{\eg(m+1)}\exp(\pii\La_{m+1}/p)
\,,\qqq \mn-1\,.
$$
\enddemo
\Ex
For $n=1$ the \fn/ $W\'\ell$ has the form
$$
W\'\ell(\tell,z_1)\,=\,\pral\,{\Expp{z_1-t_a}\over\Sinp{t_a-z_1-\La_1}}\;
\prab\ {\Sinp{t_a-t_b}\over\Sinp{t_a-t_b+1}}\;.
$$
\enddemo
\Ex
For $\ell=2$ and $n=2$ the \fn/s $W_\lg$ have the form
$$
\NN6>
\alignat2
W_{(2,0)}(t_1,t_2,z_1,z_2)\, &{}=\;
{\Expp{2z_1-t_1-t_2}\over\Sinp{t_1-z_1-\La_1}\>\Sinp{t_2-z_1-\La_1}}\>
{\Sinp{t_1-t_2}\over\Sinp{t_1-t_2+1}}\;,
\\
\ald
W_{(1,1)}(t_1,t_2,z_1,z_2)\, &{}=\;
{\Expp{z_1+z_2-t_1-t_2}\over\Sinp{t_1-z_1-\La_1}\>\Sinp{t_2-z_2-\La_2}}\>
{\Sinp{t_2-z_1+\La_1}\over\Sinp{t_2-z_1-\La_1}} &&{}\;+{}
\\
\nn-2>
&& \Llap{{}+\;
{\Expp{z_1+z_2-t_1-t_2}\over\Sinp{t_2-z_1-\La_1}\>\Sinp{t_1-z_2-\La_2}}\>
{\Sinp{t_1-z_1+\La_1}\over\Sinp{t_1-z_1-\La_1}}} &{}\,\,\x{}
\\
&& \Llap{{}\x\;{\Sinp{t_1-t_2-1}\over\Sinp{t_1-t_2+1}}} &\;,
\\
W_{(0,2)}(t_1,t_2,z_1,z_2)\, &{}=\;
{\Expp{2z_2-t_1-t_2}\over\Sinp{t_1-z_2-\La_2}\>\Sinp{t_2-z_2-\La_2}}\,\,\x{}&&
\\
\nn-2>
&{}\>\x\;{\Sinp{t_1-z_1+\La_1}\>\Sinp{t_2-z_1+\La_1}\over
\Sinp{t_1-z_1-\La_1}\>\Sinp{t_2-z_1-\La_1}}\>
{\Sinp{t_1-t_2}\over\Sinp{t_1-t_2+1}}\;, &&
\endalignat
$$
The \sthgf/ $\Fqs(z)\sub\Fq(z)$ is \onedim/ and is \gby/ the \fn/
$$
\Wo_{(2)}\,=\,W_{(2,0)}\Expp{1-2\La_1}-W_{(1,1)}\Expp{\La_2-\La_1}+
W_{(0,2)}\Expp{2\La_2-1}\,.
$$
\enddemo
For any \perm/ ${\tau\in\S^n}$ we define a basis
${\lb\)W^\tau_\lg\)\rb_{\lg\)\in\Zln}\:}$ in the \thgf/ of a fiber
by formulae similar to \(Wlga). Namely, set
$$
\ifMag\kern2em\else\kern1em\fi
W^\tau_\lg(\tell,\zn;\Lan)\,=\,W_{{\vp(}^\tau\!\lg}^{\vp1}
(\tell,z_{\tau_1}\lc z_{\tau_n};\La_{\tau_1}\lc\La_{\tau_n})
\ifMag\kern-2em\else\kern-1em\fi
\Tagg{Wtau}
$$
where $^\tau\}\lg=(\lg_{\tau_1}\lc\lg_{\tau_n})$.
\Ex
For ${\ell=1}$ and a \perm/ ${\tau=(n,n-1\lc 1)}$ the \fn/s have the form
$$
W^\tau_{\eg(m)}(t_1,\zn)\,=\,{\Expp{z_m-t_1}\over\Sinp{t_1-z_m-\La_m}}
\ \prod_{m<l\le n}\>{\Sinp{t_1-z_l+\La_l}\over\Sinp{t_1-z_l-\La_l}}\;.
$$
\enddemo
\goodbm

\Sect[2a]{\{\Rms/ and the \{\qKZc/}
\subsect*{Highest \wt/ \{\smod/s}
Let $E,F,H$ be generators of the Lie algebra $\gsl$, $[H,E]=E$, $[H,F]=-F$,
$[E,F]=2H$.
\Par
For an \smod/ $V$ let $V=\Plus_\la V_\la$ be its \wtd/.
Let $V^*=\Plus_\la V^*_\la$ be its restricted dual.
Define a structure of an \smod/ on $V^*\!$ by
$$
\bra E\phi,x\ket=\bra\phi,Fx\ket\,,\qqq
\bra F\phi,x\ket=\bra\phi,Ex\ket\,,\qqq\bra H\phi,x\ket=\bra\phi,Hx\ket\,.
$$
This \smod/ structure on $V^*\!$ will be called the \em{dual} module structure.
\goodbreak
\par
Let $\Vn$ be \smod/s with \hw/s $\Lan$, \resp/. We have the \wtd/s
\ifMag\vv-.2>\else\vv->\fi
$$
\Vox\,=\,\Plus_{\ell=0}^\8\,\Vl\qquad\text{and}\qquad
\Vax\,=\,\Plus_{\ell=0}^\8\,\Val
$$
\ifMag\vvv-.5>\else\vvv->\fi
where $()\elli$ denotes the eigenspace of $H$ with \eva/
$\sun\La_m-\ell$.
\par
Let $\FVal\sub\Val$ be the image of the operator $F$. Let $\Vls\sub\Vox$
be the kernel of the operator $E$. There is a natural pairing
$$
\Vls\ox\VFV\,\to\,\C\,.
\Tag{VVF}
$$
Let $\Vn$ be \Vmod/s, then this pairing is nondegenerate provided
$$
\pron\,\pros\,(2\La_m-s)\ne 0\,.
$$
\vsk->
\vsk->

\subsect{The \rat/ \{\Rm/}
Let $V_1,V_2$ be \Vmod/s for $\gsl$ with \hw/s $\La_1,\La_2$ and \gv/s
$v_1,v_2$, \resp/. Consider an \^{$\EVV$-}valued \mef/ $\RVx$ with
the following properties:
$$
\gather
[\RVx,F\ox\one+\one\ox F]\,=\,0\,,
\Tag{Rdef}
\\
\nn4>
\RVx\>(H\ox F-F\ox H+x\)F\ox\one)\,=\,(F\ox H-H\ox F+x\)F\ox\one)\>\RVx\,,
\endgather
$$
in $\EVV$ and
$$
\RVx\>v_1\ox v_2\,=\,v_1\ox v_2\,.
\Tag{Rvv}
$$
Such a \fn/ $\RVx$ exists and is uniquely determined. $\RVx$ is called
the \em{$\gsl\!$ \rat/ \Rm/} for the tensor product $\VV$.
\Par
It turns out that $\RVx$ commutes with the standard diagonal action of $\gsl$
in $\VV$:
$$
[\RVx,X\ox\one+\one\ox X]\,=\,0\,,\qqq X\in\gsl\,.
\Tag{RXX}
$$
In particular, $\RVx$ respects the \wtd/ of $\VV$. $\RVx$ also satisfies
the following relation
$$
\RVx\>(E\ox H-H\ox E+x\)E\ox\one)\,=\,(H\ox E-E\ox H+x\)E\ox\one)\>\RVx\,.
$$
\par
The \rat/ \Rm/ $\RVx$ satisfies the symmetry relation
$$
\align
\PV\>\RVx\,&{}=\,R_{V_2V_1}\:(x)\>\PV
\\
\nn9>
\Text{where $\PV:\VV\to V_2\ox V_1$ is the \perm/ map:
$\PV(v\ox v')=v'\ox v$, and the inversion relation}
\RVx\,&{}=\,\RV\1(-x)\,.
\endalign
$$
The following \as/s holds as $x\to\8$:
$$
\RVx\,=\,\one\ox\one +
x\1(2\La_1\La_2\,\one\ox\one-2H\ox H-E\ox F-F\ox E) + O(x^{-2})\,.
$$
\par
Let $\VV=\Plus_{l=0}^\8 V\"l\}$ be the decomposition of the \smod/ $\VV$ into
the direct sum of \irr/s, where the \irr/ module $V\"l\!$ is \gby/ a singular
\vvgood
vector of \wt/ $\La_1+\La_2-l$. Let $\Pi\"l$ be the projector onto $V\"l\!$
along the other summands. Then we have
$$
\RVx\,=\,\sum_{l=0}^\8\,\Pi\"l\cdot\prod_{s=0}^{l-1}\,
{x+\La_1+\La_2-s\over x-\La_1-\La_2+s}\;.
\Tag{Rspec}
$$
\par
Let $V_1,\,V_2,\,V_3$ be \Vmod/s. The corresponding \Rms/ satisfy the \YB/:
$$
\RV\:(x-y)\>R_{V_1V_3}\:(x)\>R_{V_2V_3}\:(y)\,=
\,R_{V_2V_3}\:(y)\>R_{V_1V_3}\:(x)\>\RV\:(x-y)\,.
\Tag{YB}
$$
All of the properties of $\RVx$ given above are well known (\cf. \Cite{KRS},
\Cite{FTT}, \Cite{T}\)).

\subsect{The Yangian $\{\Ygl$}
The \rat/ \Rm/ is connected with an action of the Yangian $\Ygl$ in a tensor
product of \smod/s. The Yangian $\Ygl$ is a remarkable Hopf algebra which
contains $\Usl$ as a Hopf subalgebra. We recall the necessary facts about
$\Ygl$ in this section.
\Par
The \em{Yangian} $\Ygl$ is a unital associative algebra with an infinite set of
generators $T_{ij}\"s\!$, $i,j=1,2$, $s=1,2,\ldots$, subject to the relations
$$
[\)T_{ij}\"r,\)T_{kl}\"{s+1}\)]\,-\,[\)T_{ij}\"{r+1}\!,\)T_{kl}\"s\)]\,=\,
T_{kj}\"r\)T_{il}\"s -\,T_{kj}\"s\)T_{il}\"r\,,
\Tag{Yijrel}
$$
$i,j,k,l=1,2$, $r,s=1,2\ldots$.
Here $T_{ij}\"0=\>\dl_{ij}$ and $\dl_{ij}$ is the Kronecker symbol.
\Par
The Yangian $\Ygl$ is a Hopf algebra with a coproduct $\Dl:\Ygl\to\Ygl\ox\Ygl$:
$$
\Dl:T_{ij}\"s\;\map\,\sum_{k=1}^2\,\sum_{r=0}^s\, T_{ik}\"r\ox\>T_{kj}\"{s-r}.
$$
There is an important one-parametric family of \aut/s $\rho_x\::\Ygl\to\Ygl$:
$$
\rho_x\::T_{ij}\"s\;\map\,\sum_{r=1}^s\,{s-1\choose r-1}\,x^{s-r}\>T_{ij}\"r.
$$
The Yangian $\Ygl$ contains $\Usl$ as a Hopf subalgebra;
the embedding is given by
$$
E\>\map\>T_{21}\"1\,,\qqq F\>\map\>T_{12}\"1\,,\qqq
H\>\map\>\bigl(T_{11}\"1\!-T_{22}\"1\bigr)/2\,.
$$
There is also an \em{evaluation \hom/} $\epe:\Ygl\to\Usl$:
$$
\alignat2
& \epe\):\)T_{11}\"s\>\map\>H\dl_{1s}\,,\qqq &&
\epe\):\)T_{12}\"s\>\map\>F\dl_{1s}\,,
\\
\nn3>
& \epe\):\)T_{21}\"s\>\map\>E\dl_{1s}\,,\qqq &&
\epe\):\)T_{22}\"s\>\map\>-H\dl_{1s}\,,
\endalignat
$$
$s=1,2,\ldots$. Both the \aut/s $\rho_x\:$ and $\epe$ restricted
to the subalgebra $\Usl$ are the identity maps.
\Par
Introduce the generating series $T_{ij}(u)=\dl_{ij}+\smsi T_{ij}\"su^{-s}\!$.
In terms of these series the coproduct, the \aut/s $\rho_x$ and the evaluation
\hom/ look like
$$
\gather
\Dl:T_{ij}(u)\,\map\,\tsum_k\) T_{ik}(u)\ox T_{kj}(u)\,,
\\
\nn2>
\rho_x\::T(u)\,\map\,T(u-x)\,,
\\
\nn6>
\alignedat2
& \epe\):\)T_{11}(u)\>\map\>Hu\1,\qqq && \epe\):\)T_{12}(u)\>\map\>Fu\1,
\\
\nn3>
& \epe\):\)T_{21}(u)\>\map\>Eu\1,\qqq && \epe\):\)T_{22}(u)\>\map\>-Hu\1.
\endalignedat
\endgather
$$
Let $e_{ij}$, $i,j=1,2$, be the ${2{\x}2}$ matrix with the only nonzero
entry $1$ at the intersection of the $i$-th row and $j$-th column. Set
\ifMag\else\vv-.3>\fi
$$
R(x)\,=\>\sum_{i,j=1}^2(x\>e_{ii}\ox e_{jj}\,+\,e_{ij}\ox e_{ji})\,.
$$
Then relations \(Yijrel) in the Yangian $\Ygl$ have the form
$$
R(x-y)\>T\'1(x)\>T\'2(y)\,=\,T\'2(y)\>T\'1(x)\>R(x-y)\,,
$$
where $T\'1(u)=\sum_{ij} e_{ij}\ox 1\ox T_{ij}(u)$ and
$T\'2(u)=\sum_{ij}1\ox e_{ij}\ox T_{ij}(u)$.
\Par
For any \smod/ $V$ denote by $V(x)$ the \Ymod/ which is obtained
from the module $V$ via the \hom/ ${\epe\o\rho_x\:}$.
The module $V(x)$ is called the \em{\emod/}.
\par
Let $V_1,V_2$ be \Vmod/s for $\gsl$ with \gv/s $v_1,v_2$, \resp/. For generic
complex numbers $x,y$ the \Ymod/s $V_1(x)\ox V_2(y)$ and $V_2(y)\ox V_1(x)$
are isomorphic and the \rat/ \Rm/ $\PV\RV\:(x-y)$ intertwines them \Cite{T},
\Cite{D1}. The vectors $v_1\ox v_2$ and $v_2\ox v_1$ are respective \gv/s
of the \Ymod/s $V_1(x)\ox V_2(y)$ and $V_2(y)\ox V_1(x)$. The \rat/ \Rm/
$\RV\:(x-y)$ can be defined as the unique element of $\EVV$ with property
\(Rvv) and \st/
$$
\PV\RV\:(x-y)\>:\>V_1(x)\ox V_2(y)\,\to\,V_2(y)\ox V_1(x)
\Tag{Yintw}
$$
is an \iso/ of the \Ymod/s.
\Par
For a \Ymod/ $V$ let $V=\Plus_\la V_\la$ be its \wtd/ as an \smod/.
Let $V^*=\Plus_\la V^*_\la$ be its restricted dual.
Define a structure of a \Ymod/ on $V^*\!$ by
\ifMag\vv->\fi
$$
\alignat2
& \bra T_{11}(u)\phi,x\ket=\bra\phi,T_{11}(u)x\ket\,,\qqq &&
\bra T_{12}(u)\phi,x\ket=\bra\phi,T_{21}(u)x\ket\,,
\\
\nn4>
& \bra T_{21}(u)\phi,x\ket=\bra\phi,T_{12}(u)x\ket\,,\qqq &&
\bra T_{22}(u)\phi,x\ket=\bra\phi,T_{22}(u)x\ket\,.
\endalignat
$$
This \Ymod/ structure on $V^*\!$ will be called the \em{dual} module structure.

\subsect{The \rat/ \{\qKZc/ associated with $\{\gsl$}
Let $\Vn$ be \smod/s. The \qKZc/ is a discrete connection on the \trib/
over $\Cn$ with fiber $\Vox$. We define it below.
\Par
Let $\Vn$ be \Vmod/s with \hw/s $\Lan$, \resp/. Let $R_{V_iV_j}\:(x)$
be the \rat/ \Rms/. Let $R_{ij}(x)\in\EV$ be defined in a standard way:
$$
R_{ij}(x)\,=\>\sum\one\lox\%{r(x)}_{\^{$i$-th}\,}\lox\%{r'(x)}_{\^{$j$-th}\,}
\lox\one
\Tag{Rij}
$$
provided that $R_{V_iV_j}\:(x)=\sum r(x)\ox r'(x)\in\End(V_i\ox V_j)$.
For any $X\in\gsl$ set
$$
X_m\,=\,\one\lox\%{X}_{\sss\^{$m$-th}}\lox\one\,.
$$
Let $p,\ka$ be complex numbers. For any $\mn$ set
$$
\align
\ifMag\kern2.5em\else\kern1em\fi
K_m(\zn)\,=\, R_{m,m-1}(z_m-z_{m-1}+p)\ldots R_{m,1}(z_m-z_1+p)\>
\ka^{\La_m-H_m}\;\x&
\Tagg{Kmz}
\\
\nn3>
{}\x R_{m,n}(z_m-z_n)\ldots R_{m,m+1}(z_m-z_{m+1}) &\,.
\ifMag\kern-2.5em\else\kern-1em\fi
\endalign
$$
\Th{FR}
\back\Cite{FR}
The linear maps $K_m(z)$ obey the flatness conditions
$$
K_l(\zmn)\>K_m(\zn)\,=\,K_m(z_1\lc z_l\]+p\lc z_n)\>K_l(\zn)\,,
$$
$l,\mn$.
\endpro
The maps $K_1(z)\lc K_n(z)$ define a flat connection on the \trib/
over $\Cn$ with fiber $\Vox$. This connection is called the \em{\qKZc/}.
\goodbreak
\Par
By \(RXX) the operators $K_m(z)$ commute with the diagonal action of $H$
in $\Vox$:
$$
[K_m(\zn)\>,\)H]=0\,,\qqq\mn\,,
$$
and, therefore, preserve the \wtd/ of $\Vox$. Hence, the \qKZc/ induces
the dual flat connection on the \trib/ over $\Cn$ with fiber $\Vax$.
This connection will be called the \em{dual \qKZc/}.
\Par
Let $\BB\sub\Cn$ be the complement to the discriminant \(dis).
\Lm*
For any $z\in\BB$ the linear maps $K^*_1(z)\lc K^*_n(z)$ define \iso/s
of $\Val$.
\endpro
\nt
This statement follows from formulae \(Rspec) and \(Kmz).
\Par
If $\ka=1$, then the dual \qKZc/ commutes with the diagonal action of $\gsl$
in $\Vax$:
$$
[K^*_m(\zn)\>,\)X]=0\,,\qqq X\in\gsl\,,\qquad\mn\,,
$$
and, therefore, admits a trivial discrete subbundle with fiber $\FVal$,
moreover, it induces a flat connection on the \trib/ with fiber
$\Val\big/\FVal$.
\Par
Let $\Vn$ be \smod/s. The \em{\qKZe/} for a \^{$\Vox$-}valued \fn/
$\Psi(\zn)$ is the following system of \eq/s
$$
\Psi(\zmn)\,=\,K_m(\zn)\>\Psi(\zn)\,, \qqq \mn\,.
$$
The \qKZe/ is a remarkable \deq/, see \Cite{S}, \Cite{FR}, \Cite{JM},
\Cite{Lu}.

\subsect{The \tri/ \{\Rm/}
Let $q$ be a nonzero complex number which is not a root of unity.
Let $E_q,F_q,q^{\pm H}$ be generators of $\Uu$:
\vv-1.5>
$$
\NN3>
\gather
q^H\qH=\>\qH q^H=\>1\,,
\\
q^H\]E_q=q\)E_q\>q^H,\qqq q^H\]F_q=\q F_q\>q^H,
\\
[E_q,F_q]\;=\,{q^{2H}-q^{-2H}\over q-\q}\;.
\endgather
$$
A coproduct $\Dlq:\Uu\,\to\,\Uu\ox\Uu$ is given by
$$
\NN2>
\align
\Dlq(q^H\])=q^H\]\ox q^H, &\qqq \Dlq(\qH\])=\qH\]\ox\qH,
\\
\Dlq(E_q)=E_q\ox q^H\]+\qH\]\ox E_q\,, &\qqq
\Dlq(F_q)=F_q\ox q^H\]+\qH\]\ox F_q\,.
\endalign
$$
The coproduct defines a \Umod/ structure on a tensor product of \Umod/s.
\par
Let $V_1,V_2$ be \Vmod/s for $\Uu$ with \hw/s $q^{\La_1}\!,\,q^{\La_2}$ and
\gv/s $v_1,v_2$, \resp/. Consider an \^{$\EVV$-}valued \mef/ $\RVz$ with
the following properties:
$$
\align
\kern1em
\RVz\>(F_q\ox q^H\]+\>\qH\]\ox F_q)\, {}& =\,
(F_q\ox\qH\]+\>q^H\]\ox F_q)\>\RVz\,,
\kern-1em
\Tagg{Rqdef}
\\
\nn4>
\kern1em
\RVz\>(F_q\ox\qH\]+\>\zt\)q^H\]\ox F_q)\, {}& =
\,(F_q\ox q^H\]+\>\zt\)\qH\]\ox F_q)\>\RVz
\kern-1em
\endalign
$$
in $\EVV$ and
$$
\RVz\>v_1\ox v_2=v_1\ox v_2\,.
\Tag{Rqvv}
$$
Such a \fn/ $\RVz$ exists and is uniquely determined. $\RVz$ is called
the \em{$\gsl\!$ \tri/ \Rm/} for the tensor product $\VV$.
\Par
The \tri/ \Rm/ $\RVz$ also satifies the following relations
$$
\NN4>
\align
\kern1em
\RVz\>(E_q\ox q^H\]+\>\qH\]\ox E_q)\, {}& =\,
(E_q\ox\qH\]+\>q^H\]\ox E_q)\>\RVz
\kern-1em
\Tagg{Rqmore}
\\
\kern1em
\RVz\>(\zt\)E_q\ox\qH\]+\>q^H\]\ox E_q)\, {}& =
\,(\zt\)E_q\ox q^H\]+\>\qH\]\ox E_q)\>\RVz\,,
\kern-1em
\\
\RVz\,q^H\]\ox q^H {}& =\,q^H\]\ox q^H\RVz\,.
\endalign
$$
In particular, $\RVz$ respects the \wtd/ of $\VV$.
\par
$\RVz$ satisfies the inversion relation
$$
\PV\>\RVz\,=\,\bigl(R_{V_2V_1}^q(\zt\1)\bigr)\1\>\PV
$$
where $\PV:\VV\to V_2\ox V_1$ is the \perm/ map.
\par
Let $\VV=\Plus_{l=0}^\8 V\"l$ be the decomposition of the \Umod/ $\VV$ into
the direct sum of \irr/s, where the \irr/ module $V\"l\!$ is \gby/ a singular
vector of \wt/ $q^{\La_1+\La_2-l}\!$. Let $\Pi\"l$ be the projector onto
$V\"l\!$ along the other summands. Then we have
$$
\gather
\RVz\,=\,\RVo\,\sum_{l=0}^\8\,\Pi\"l\cdot\prod_{s=0}^{l-1}\,
{1-\zt q^{2s-2\La_1-2\La_2}\over1-\zt q^{2\La_1+2\La_2-2s}}\;.
\Tag{Rspeq}
\\
\Text{where}
\nn-4>
\RVo\,=\,q^{2\La_1\La_2-2H\ox H}\smk^\8(q^2-1)^{2k}
\tprod_{s=1}^k(1-q^{2s})\1\>(q^H\]F_q\ox\qH\]E_q)^k\,.
\endgather
$$
\goodbe
\par
Let $V_1,\,V_2,\,V_3$ be \Vmod/s. The corresponding \Rms/ satisfy the \YB/:
$$
\RV^q(\xi/\zt)\>R_{V_1V_3}^q(\xi)\>R_{V_2V_3}^q(\zt)\,=
\,R_{V_2V_3}^q(\zt)\>R_{V_1V_3}^q(\xi)\>\RV^q(\xi/\zt)\,.
\Tag{YBq}
$$
All of the properties of $\RVz$ given above are well known (\cf. \Cite{T},
\Cite{D1}, \Cite{J}, \Cite{CP}\)).
\goodbm
\Par
Similar to the \rat/ case one can define the \qKZc/ associated with
the \tri/ \Rm/ (\cf. \Cite{FR}). We study this \tri/ \qKZc/ in \Cite{TV3}.

\subsect{The \qlo/ $\{\Ugg$}
The \tri/ \Rm/ is connected with an action of the \qlo/ $\Ugg$ in a tensor
product of \Umod/s. The \qlo/ $\Ugg$ is a Hopf algebra which contains $\Uu$
as a Hopf subalgebra. We give the necessary facts about $\Ugg$ in this section.
\Par
Let $q$ be a complex number, $q\ne\pm1$. The \em{\qlo/} $\Ugh$ is a unital
associative algebra with generators $L_{ij}\"{+0},\ L_{ji}\"{-0}\!$,
$1\le j\le i\le 2$, and $L_{ij}\"s\!$, $i,j=1,2$, $s=\pm 1,\pm 2,\ldots$,
subject to relations \(Ughr)\,\Cite{RS}, \Cite{DF}.
\par
Let $e_{ij}$, $i,j=1,2$, be the ${2{\x}2}$ matrix with the only nonzero entry
$1$ at the intersection of the $i$-th row and $j$-th column. Set
$$
\align
R(\xi)\,=\,{} &(\xi q-\q)\>(e_{11}\ox e_{11}+e_{22}\ox e_{22})\;+
\\
\nn3>
{} +\,{} &(\xi-1)\>(e_{11}\ox e_{22}+e_{22}\ox e_{11})\,+\,
\xi(q-\q)\>e_{12}\ox e_{21}\,+\,(q-\q)\>e_{21}\ox e_{12}\,.
\endalign
$$
Introduce the generating series
$L^{\pm}_{ij}(u)=L_{ij}\"{\pm0}+\smsi L_{ij}\"{\pm s}u^{\pm s}\!$.
The relations in $\Ugh$ have the form
\vv-.7>
$$
\NN4>
\gather
L_{ii}\"{+0}L_{ii}\"{-0}\,=\,1\,,\qquad L_{ii}\"{-0}L_{ii}\"{+0}\,=\,1\,,
\Rlap{\qqq i=1,2\,,}
\Tag{Ughr}
\\
{\align
& R(\xi/\zt)\>L^+\'1(\xi)\>L^+\'2(\zt)\,=\,
L^+\'2(\zt)\>L^+\'1(\xi)\>R(\xi/\zt)\,,
\\
& R(\xi/\zt)\>L^+\'1(\xi)\>L^-\'2(\zt)\,=\,
L^-\'2(\zt)\>L^+\'1(\xi)\>R(\xi/\zt)\,,
\\
& R(\xi/\zt)\>L^-\'1(\xi)\>L^-\'2(\zt)\,=\,
L^-\'2(\zt)\>L^-\'1(\xi)\>R(\xi/\zt)\,,
\endalign}
\endgather
$$
where $L^\nu\'1(\xi)=\sum_{ij} e_{ij}\ox 1\ox L_{ij}^\nu(\xi)$ and
$L^\nu\'2(\xi)=\sum_{ij}1\ox e_{ij}\ox L_{ij}^\nu(\xi)$,\quad $\nu=\pm$.
\goodbreak
\par
Elements $L_{11}\"{+0}L_{22}\"{+0}$, $L_{22}\"{+0}L_{11}\"{+0}$,
$L_{11}\"{-0}L_{22}\"{-0}$, $L_{22}\"{-0}L_{11}\"{-0}$ are central in $\Ugh$.
Impose the following relations:
$$
L_{11}\"{+0}L_{22}\"{+0}\>=\>1\,,\qquad
L_{22}\"{+0}L_{11}\"{+0}\>=\>1\,,\qquad
L_{11}\"{-0}L_{22}\"{-0}\>=\>1\,,\qquad
L_{22}\"{-0}L_{11}\"{-0}\>=\>1\,,
$$
in addition to relations \(Ughr). Denote the corresponding quotient algebra
by $\Ugg$.
\Par
The \qlo/ $\Ugg$ is a Hopf algebra with a coproduct $\Dlq:\Ugg\to\Ugg\ox\Ugg$:
$$
\Dlq:L^\nu_{ij}(\xi)\,\map\,\tsum_k L^\nu_{kj}(\xi)\ox L^\nu_{ik}(\xi)\,,\qqq
\nu=\pm\,.
$$
\Rem
Notice that we take the coproduct $\Dlq$ for the \qlo/ $\Ugg$ which is
in a sense opposite to the coproduct $\Dl$ taken for the Yangian $\Ygl$
(\cf. Theorems~\[FgF], \[FgFq]\)).
\enddemo
There is an important one-parametric family of \aut/s $\rho^q_\zt:\Ugg\to\Ugg$:
$$
\gather
\rho^q_\zt:L_{ij}^\nu(\xi)\,\map\,L_{ij}^\nu(\xi/\zt)\,,\qqq \nu=\pm\,,
\\
\nn-2>
\Text{that is}
\nn-2>
\rho^q_\zt:L_{ij}\"{\pm0}\,\map\,L_{ij}\"{\pm0}\qquad\text{and}\qquad
\rho^q_\zt:L_{ij}\"{s}\map\zt^{-s}L_{ij}\"{s}\,,\qquad s\in\Z_{\ne0}\,.
\endgather
$$
\egood
The \qlo/ $\Ugg$ contains $\Uu$ as a Hopf subalgebra;
the embedding is given by
$$
E_q\>\map\>-L_{21}\"{+0}\big/(q-\q)\,,\qqq
F_q\>\map\>L_{12}\"{-0}\big/(q-\q)\,,\qqq q^H\>\map\>L_{11}\"{-0}\,.
$$
There is also an \em{evaluation \hom/} $\epe^q:\Ugg\to\Uu$:
\vv-.5>
$$
\NN4>
\alignat2
& \epe^q\):\)L^+_{11}(\xi)\,\map\,\qH-q^H\xi\,, &&
\epe^q\):\)L^+_{12}(\xi)\,\map\,-F_q\)(q-\q)\>\xi\,,
\\
& \epe^q\):\)L^+_{21}(\xi)\,\map\,-E_q\)(q-\q)\,, &&
\epe^q\):\)L^+_{22}(\xi)\,\map\,q^H-\qH\xi\,,
\\
\ald
\nn4>
& \epe^q\):\)L^-_{11}(\xi)\,\map\,q^H-\qH\xi\1\,, &&
\epe^q\):\)L^-_{12}(\xi)\,\map\,F_q\)(q-\q),
\\
& \epe^q\):\)L^-_{21}(\xi)\,\map\,E_q\)(q-\q)\>\xi\1\,,\qqq
\ifMag\quad\else\qqq\fi
&& \epe^q\):\)L^-_{22}(\xi)\,\map\,\qH-q^H\xi\1,
\endalignat
$$
\vvv-.5>
that is
\ifMag\vv-.5>\else\vv-.8>\fi
$$
\NN4>
\alignat4
& \epe^q\):\)L_{11}\"{+0}\,\map\,\qH,\qquad &&
\epe^q\):\)L_{11}\"1\,\map\,-q^H, &&
\qqq\Rlap{\epe^q\):\)L_{12}\"1\,\map\,-F_q\)(q-\q)\,,} &&
\\
&\qqq \Rlap{\epe^q\):\)L_{21}\"{+0}\,\map\,-E_q\)(q-\q)\,,} &&&&
\epe^q\):\)L_{22}\"{+0}\,\map\,q^H, && \epe^q\):\)L_{22}\"1\,\map\,-\qH,
\\
\ald
\nn4>
& \epe^q\):\)L_{11}\"{-0}\,\map\,q^H, &&
\epe^q\):\)L_{11}\"{-1}\,\map\,-\qH, \ifMag\ \else\qqq\fi &&
\qqq\Rlap{\epe^q\):\)L_{12}\"{-0}\,\map\,F_q\)(q-\q)\,,} &&
\\
&\qqq \Rlap{\epe^q\):\)L_{21}\"{-1}\,\map\,E_q\)(q-\q)\,,} &&&&
\epe^q\):\)L_{22}\"{-0}\,\map\,\qH,\qquad &&
\epe^q\):\)L_{22}\"{+1}\,\map\,-q^H
\endalignat
$$
and ${\epe^q\):\)L_{ij}\"s\,\map\,0}$ for all other generators $L_{ij}\"s$.
\par
Both the \aut/s $\rho^q_\xi$ and $\epe^q$ restricted to the subalgebra $\Uu$
are the identity maps.
\par
For any \Umod/ $V$ denote by $V(\xi)$ the \Uhmod/ which is obtained from
the module $V$ via the \hom/ ${\epe^q\o\rho^q_\xi}$. The module $V(\xi)$
is called the \em{\emod/}.
\par
Let $V_1,V_2$ be \Vmod/s for $\Uu$ with \gv/s $v_1,v_2$, \resp/. For generic
complex numbers $\xi,\zt$ the \Uhmod/s $V_1(\xi)\ox V_2(\zt)$ and
$V_2(\zt)\ox V_1(\xi)$ are isomorphic and the \tri/ \Rm/ $\PV\RV^q(\xi/\zt)$
intertwines them \Cite{T}, \Cite{CP}. The vectors $v_1\ox v_2$ and
$v_2\ox v_1$ are respective \gv/s of the \Uhmod/s $V_1(\xi)\ox V_2(\zt)$ and
$V_2(\zt)\ox V_1(\xi)$. The \tri/ \Rm/ $\RV^q(\xi/\zt)$ can be defined as
the unique element of $\EVV$ with property \(Rqvv) and \st/
$$
\PV\RV^q(\xi/\zt)\>:\>V_1(\xi)\ox V_2(\zt)\,\to\,V_2(\zt)\ox V_1(\xi)
\Tag{Uintw}
$$
is an \iso/ of the \Uhmod/s.

\Sect[2b]{Tensor coordinates and module structures on the \hgeom/ spaces}
In this section we identify the \GM/ and the \qKZc/. In addition we also
describe a structure of a \Ymod/ on the \rhgf/ and a structure of a \Uhmod/
on the \thgf/, \resp/.

\subsect{The \rhgm/}
Let ${\Fol}$ be the \rhgf/ defined for the projection ${\C^{\,l+n}\!\to\Cn\!}$.
In particular, $\Fo[0]=\C$ and, in our previous notations, we have $\Foll=\Fo$.
Consider the direct sum
$$
\Fg\,=\,\Plus_{l\ge 0}\)\Fol
$$
which will be called the \em{\rhgF/}.
\Par
Let $T_{ij}(u)$, $i,j=1,2$, be the generating series for the Yangian $\Ygl$
introduced in Section~\SNo{2a}. Set
$$
\Ti_{ij}(u)\,=\,T_{ij}(u)\,\pron\,{u-z_m\over u-z_m-\La_m\!}\;,\qqq i,j=1,2\,,
$$
where the \raf/ in \rhs/ is understood as its Laurent series expansion
at $u=\8$. It is clear that the coefficients of the series $\Ti_{ij}(u)$
generate $\Ygl$. Introduce an action of the coefficients of the series
$\Ti_{ij}(u)$ in the space $\Fg$. Namely, for any $f\in\Fol$ set:
\ifMag
$$
\NN4>
\gather
\Rline{(\Ti_{11}(u)f)(t_1\lc t_l)\,=\,
f(t_1\lc t_l)\,\pron\,{u-z_m+\La_m\over u-z_m-\La_m}\;
\prod_{a=1}^l\,{u-t_a-1\over u-t_a}\ \,+\hp{_{(a,l)}}}
\Tagg{Tiju}
\\
\Rline{{}+\,\prod_{a=1}^l\,{u-t_a-1\over u-t_a}\;
\sum_{a=1}^l\,\Bigl[\,{f(t_1\lc t_{l-1},u)\over u-t_l-1}\;
\pron\,{t_l-z_m+\La_m\over t_l-z_m-\La_m}\,\Bigr]_{(a,l)}\,,\!}
\\
\ald
\nn6>
{\align
(\Ti_{22}(u)f)(t_1\lc t_l)\, &{}=\,f(t_1\lc t_l)\,
\prod_{a=1}^l\,{u-t_a+1\over u-t_a}\ \,-
\\
& \>{}-\,\prod_{a=1}^l\,{u-t_a+1\over u-t_a}\;\sum_{a=1}^l\,
\Bigl[\,{f(u,t_2\lc t_l)\over u-t_1+1}\,\Bigr]_{(1,a)}\,,
\endalign}
\\
\ald
\nn6>
\Lline{\quad(\Ti_{12}(u)f)(t_1\lc t_{l+1})\,=\,
\sum_{a=1}^{l+1}\,\Bigl[\,{f(t_2\lc t_{l+1})\over u-t_1}\ \,\x}
\\
\Rline{{}\x\,\Bigl(\;\pron\,{t_1-z_m+\La_m\over t_1-z_m-\La_m}\;
\prod_{b=2}^{l+1}\>\>{u-t_b+1\over u-t_b}\,
{t_1-t_b-1\over t_1-t_b+1}\ \;-\hp{_{(1,a)}\,}\hp{_{(1,a)}\,}}
\\
\nn2>
\Rline{{}-\,\pron\,{u-z_m+\La_m\over u-z_m-\La_m}\,
\prod_{b=2}^{l+1}\,{u-t_b-1\over u-t_b}\;\Bigr)\>\Bigr]_{(1,a)}\,-
\hp{_{(1,a)}}}
\\
\ifMag\ald\fi
\Rline{-\,\prod_{a=1}^{l+1}\,{u-t_a+1\over u-t_a}
\sum_{\tsize{a,b=1\atop a\ne b}}^{l+1}\,
\Bigl[\,{f(u,t_2\lc t_l)\over(u-t_1+1)\>(u-t_{l+1}+1)}\;
\pron\,{t_{l+1}-z_m+\La_m\over t_{l+1}-z_m-\La_m}\,\Bigr]_{\si^{ab}}\,,\!}
\\
\ald
\nn4>
(\Ti_{21}(u)f)(t_1\lc t_{l-1})\,=\,f(t_1\lc t_{l-1},u)\,
\prod_{a=1}^{l-1}\,{u-t_a-1\over u-t_a}\;,\qqq l>0\,,
\endgather
$$
\else
$$
\NN4>
\gather
{\align
(\Ti_{11}(u)f)(t_1\lc t_l)\, &{}=\,f(t_1\lc t_l)\,
\pron\,{u-z_m+\La_m\over u-z_m-\La_m}\;\prod_{a=1}^l\,{u-t_a-1\over u-t_a}\ +
\Tagg{Tiju}
\\
&\>{}+\,\prod_{a=1}^l\,{u-t_a-1\over u-t_a}\;
\sum_{a=1}^l\,\Bigl[\,{f(t_1\lc t_{l-1},u)\over u-t_l-1}\;
\pron\,{t_l-z_m+\La_m\over t_l-z_m-\La_m}\,\Bigr]_{(a,l)}\,,
\endalign}
\\
\ald
\nn6>
(\Ti_{22}(u)f)(t_1\lc t_l)\,=\,f(t_1\lc t_l)\,
\prod_{a=1}^l\,{u-t_a+1\over u-t_a}\ -
\,\prod_{a=1}^l\,{u-t_a+1\over u-t_a}\;\sum_{a=1}^l\,
\Bigl[\,{f(u,t_2\lc t_l)\over u-t_1+1}\,\Bigr]_{(1,a)}\,,
\\
\ald
\nn6>
{\align
(\Ti_{12}(u)f)(t_1\lc t_{l+1})\, &{}=\,\sum_{a=1}^{l+1}\,\Bigl[\,
{f(t_2\lc t_{l+1})\over u-t_1}\ \,\x
\\
&
\aligned
\hp{\sum}\>{}\x\,\Bigl(\;\pron\,{t_1-z_m+\La_m\over t_1-z_m-\La_m}\;
\prod_{b=2}^{l+1}\>\>{u-t_b+1\over u-t_b}\,
{t_1-t_b-1\over t_1-t_b+1}\ & \;-
\\
\nn2>
{}-\,\pron\,{u-z_m+\La_m\over u-z_m-\La_m}\,
\prod_{b=2}^{l+1}\,{u-t_b-1\over u-t_b}\;\Bigr)\> &\Bigr]_{(1,a)}\,-
\endaligned
\\
&\>-\,\prod_{a=1}^{l+1}\,{u-t_a+1\over u-t_a}
\sum_{\tsize{a,b=1\atop a\ne b}}^{l+1}\,
\Bigl[\,{f(u,t_2\lc t_l)\over(u-t_1+1)\>(u-t_{l+1}+1)}\;
\pron\,{t_{l+1}-z_m+\La_m\over t_{l+1}-z_m-\La_m}\,\Bigr]_{\si^{ab}}\,,
\endalign}
\\
\ald
\nn4>
(\Ti_{21}(u)f)(t_1\lc t_{l-1})\,=\,f(t_1\lc t_{l-1},u)\,
\prod_{a=1}^{l-1}\,{u-t_a-1\over u-t_a}\;,\qqq l>0\,,
\endgather
$$
\fi
and $\Ti_{21}(u)f=0$ for $f\in\Fo[0]$. Here $(1,a)$, $(a,l)$ are
transpositions and $\si^{ab}\in\S^{l+1}$ is the following \perm/
$$
\si^{ab}\!:i\>\map\)i\quad\for i=2\lc l\,,\qqq
\si^{ab}\!:1\>\map\)a\,,\qquad\si^{ab}\!:l+1\>\map\)b\,.
$$
\Rhs/s of formulae \(Tiju) are \raf/s in $u$, and the precise meaning of each
of the formulae is that \lhs/ equals the Laurent series expansion of
the respective right hand side at $u=\8$.
\Lm{YFg}
Formulae \(Tiju) define a \Ymod/ structure in the \rhgF/ $\Fg$.
\endpro
\nt
The proof is given by direct verification.
\Par
Let ${\Fg(z)=\Plus_{l\ge 0}\Fol(z)}$ be the \rhgF/ of a fiber.
The \Ymod/ structure in $\Fg$ clearly induces a \Ymod/ structure in $\Fg(z)$.
This module will be called the \em{\rhgm/}.
\Par
For the action of the generators of the subalgebra $\Usl$ formulae \(Tiju)
simplify and for ${f\in\Fol}$ look as follows:
\vv-.1>
$$
\gather
(H\]f)(t_1\lc t_l)\,=\,\bigl(\tsun\La_m-l\bigr)\>f(t_1\lc t_l)\,,
\Tag{HFE}
\\
\nn3>
(F\]f)(t_1\lc t_{l+1})\,=\,\sum_{a=1}^{l+1}\,\Bigl[\,f(t_2\lc t_{l+1})\,
\Bigl(\;\pron\,{t_1-z_m+\La_m\over t_1-z_m-\La_m}\,
\prod_{b=2}^{l+1}\,{t_1-t_b-1\over t_1-t_b+1}\,-1\Bigr)\,\Bigr]_{(1,a)}\,,
\\
\nn8>
(Ef)(t_1\lc t_{l-1})\,=\,\bigl(\)t_l\>f(t_1\lc t_l)\bigl)\vst{t_l=\8}\,,
\Rlap{\qqq l>0\,,}
\endgather
$$
and $Ef=0$ for $f\in\Fo[0]$. Here $(1,a)\in\S^{l+1}$ are transpositions.
\Rem
It is worth to mention that for any \fn/ $w_\lg$ we have
$$
Fw_\lg\,=\,\sun(\lg_m+1)\>(2\La_m-\lg_m)\>w_{\lg+\eg(m)}\,.
$$
\cf. \(wDw). Hence, ${\Rc[\)l\>](z)\sub F\bigl(\Fo[\)l-1\)](z)\bigr)}$,
where ${\Rc[\)l\>](z)}$ is the \cosub/.
\enddemo
\Rem
Let $\ka=1$. Then for any \fn/ ${f\in\Fo[\)\ell-1\)]}$ we have
\ifMag\else\vv-.5>\fi
$$
(F\]f)(\tell)\,=\,
\sum_{a=1}^\ell\,D_a\bigl[\)f(t_2\lc t_\ell)\)\bigr]_{(1,a)}\,.
\Tag{FfDf}
$$
\enddemo

\subsect{Tensor coordinates on the \rhgf/s of fibers}
Let $\Vn$ be $\gsl$ Verma modules with \hw/s $\Lan$ and \gv/s $\vn$, \resp/.
Consider the \wt/ subspace $\Vl$ with a basis given by monomials $\Fv$,
$\lg\in\Zln$. The dual space $\Val$ has the dual basis denoted by $\Fva\!$,
$\lg\in\Zln$.
\par
For any $z\in\Cn$ and for any $\tau\in\S^n$ denote by $B_\tau(z)$
the following \hom/:
$$
\gather
B_\tau(z):\Valt\,\to\,\Fo(z)\,,
\\
\nn3>
B_\tau(z):\Fvat\,\map\,w^\tau_\lg(t,z)\,,
\endgather
$$
where $\Fo(z)$ is the \rhgf/ of the fiber (\cf. \(wlga), \(wtau)\)).
The \hom/s $B_\tau(z)$ are called the \em{\tenco/} on the \rhgf/ of a fiber.
The composition maps
$$
B_{\tau,\tau'}(z):\Valtt\,\to\,\Valt\,,\qqq
B_{\tau,\tau'}(z)\,=\,B\1_\tau(z)\o B_{\tau'}(z)\,,
$$
are called the \em{\traf/s}, \cf. \Cite{V3}.
\Lm{tcfs}
Let $z_l+\La_l-z_m+\La_m\nin\lb\)0\lc\ell-1\)\rb$ for any $l\ne m$, $l,\mn$.
Then for any \perm/ $\tau$ the linear map $B_\tau(z):\Valt\to\Fo(z)$ is \ndeg/.
\endpro
\nt
The statement follows from Lemma~\[wbasis].
\Par
Consider a tensor product $\Voxzt$ of \emod/s over $\Ygl$ coinciding with
$\Voxt$ as an \smod/.
\Lm{T12}
For any $\phi\in\Valt$ we have
\ifMag
$$
\align
\bra\phi,T_{12}(t_1) & \ldots T_{12}(t_\ell)\>\vox\ket\;=
\\
\nn2>
{}={}\,\, &\!\bigl(B_\tau(z)\)\phi\bigr)(\tell)\,
\pral\,\pron\,{t_a-z_m-\La_m\over t_a-z_m}\!\prab\){t_a-t_b+1\over t_a-t_b}\;.
\endalign
$$
\else
$$
\bra\phi,T_{12}(t_1)\ldots T_{12}(t_\ell)\>\voxt\ket\,=\,
\bigl(B_\tau(z)\)\phi\bigr)(\tell)\,\pral\,\pron\,{t_a-z_m-\La_m\over t_a-z_m}
\!\prab\){t_a-t_b+1\over t_a-t_b}\;.
$$
\fi
\endpro
\nt
The \raf/ in \rhs/ above is understood as its power series expansion
in $t_1\1\lc t_\ell\1\!$. It is easy to see that the \raf/ is regular
at the hyperplanes $t_a=t_b$, $1\le a<b\le\ell$, so the power series
is well defined without additional prescriptions.
\Th{FgVax}
For any \perm/ $\tau\in\S^n$ the map
$$
B_{\tau}(z)\):\)\Vaxzt\,\to\,\Fg(z)
$$
is an intertwiner of \Ymod/s.
\endpro
\Cr*
Let $z_l+\La_l-z_m+\La_m\nin\Z$ for any $l\ne m$, $l,\mn$.
Then for any \perm/ $\tau\in\S^n$ the map ${B_{\tau}(z):\Vaxzt\,\to\,\Fg(z)}$
is an \iso/ of \Ymod/s.
\endpro
\nt
The statement follows from Theorem~\[FgVax] and Lemma~\[tcfs].
\Th{local}
\back\Cite{V3}
For any ${\tau\in\S^n}$ and any transposition $(m,m+1)$, $\mn-1$, the \traf/
$$
B_{\tau,\tau\cdot(m,m+1)}(z):\Vaxtm_\ell\,\to\,\Valt
$$
equals the operator
${\bigl(P_{V_{\tau_m}\}V_{\tau_{m+1}}}\:
R_{V_{\tau_m}\}V_{\tau_{m+1}}}\:\!(z_{\tau_m}\!-z_{\tau_{m+1}})\bigr)^*}$
acting in the \^{$m$-th} and \^{$(m+1)$-th} factors.
\endpro
\nt
The theorem follows from Lemma~\[T12] and formula \(Yintw).
\Par
Each $B_\tau(z)$ induces a linear map $\Valt\,\to\,\H(z)$
which also will be denoted by $B_\tau(z)$.
\par
\Th{tcnon}
Let $\ka\ne 1$. Let $p<0$. Let \(npZ)\;--\;\(assum) hold.
Then for any $\tau\in\S^n$ the map $B_\tau(z):\Valt\to\H(z)$ is an \iso/.
\endpro
\nt
This statement follows from Theorem~\[kanon] and Lemma~\[tcfs].
\Par
It is easy to see that for any $\tau\in\S^n$ the image of $\FValt$ under
the map $B_\tau(z)$ coincides with the \cosub/ $\Rc(z)\sub\Fo(z)$.
\Th{tcon}
Let $\ka=1$. Let $p<0$. Let \(npZ)\;--\;\(assum) hold.
If ${\tLmn-s\nin p\)\Znn}$ for all $s=\ell-1\lc 2\ell-2$,
then for any $\tau\in\S^n$ the map $B_\tau(z)$ induces an \iso/
$$
\VFVt\)\to\,\H(z)\,.
$$
\endpro
\nt
The statement follows from Theorem~\[kaon] and Lemmas~\[wbasis], \[wDw].
\Par
Taking into account formula \(VVF) we get an \iso/
$$
\bigl(\Vlst\bigr)^*\,\to\,\H(z)\,.
$$
\Th{qKZ-GM}
\back\Cite{V3},\;\Cite{TV1}
For any $\mn$, the following diagram is commutative:
$$
\kern1.4em
\CD
\Valt @>{\tsize\ K^*_m(\zn)\ }>> \Valt
\\
@V{\tsize{\vrp24pt:18pt>}B_\tau(\zmn)\ }VV
@VV{\tsize\ B_\tau(\zn{\vrp24pt:18pt>})}V
\\
\H(\zmn) @>>{\tsize\ A_m(\zn)\ }> \H(\zn)
\endCD
\kern-1.4em
$$
Here $A_m(z)$ are the operators of the \GM/, $K^*_m(z)$ are the operators dual
to $K_m(z)$, and $K_m(z)$ are the operators of the \qKZc/ in $\Vlt$ defined by
\(Kmz).
\endpro
\goodbreak
\Cr*
The construction above identifies the \qKZc/ and the \GM/ restricted to
the \hgeom/ subbundle.
\endpro

\subsect{The \thgm/}
Let ${\Fql}$ be the \thgf/ defined for the projection ${\C^{\,l+n}\!\to\Cn\!}$.
In particular, $\Fq[0]=\C$ and, in our previous notations, we have $\Fqll=\Fq$.
Consider the direct sum
$$
\Fg_q\,=\,\Plus_{l\ge 0}\)\Fql
$$
which will be called the \em{\thgF/}.
\Par
Let $q=\exp(\pii/p)$. Let $L^{\pm}_{jk}(u)$, $j,k=1,2$, be the generating
series for the \qlo/ $\Ugg$ introduced in Section~\SNo{2a}. Set
$$
\Lti^\pm_{jk}(\xi)\,=\,L^\pm_{jk}(\xi)\,
\pron\,{\pm i\exp\bigl(\pm\pii(z_m-u)/p\bigr)\over2\Sinp{u-z_m-\La_m}}\;,\qqq
\Rlap{j,k=1,2\,,}
$$
where ${\xi=\exp(2\pii\)u/p)}$. The products in \rhs/ are \raf/s in $\xi$.
The precise meaning is that ${\Lti^-_{jk}(\xi)}$ equals the Laurent series
expansion of the corresponding right hand side at $\xi=\8$, and
$\Lti^+_{jk}(\xi)$ equals the Taylor series expansion of the corresponding
right hand side at $\xi=0$. It is clear that the coefficients of the series
$\Lti^\pm_{jk}(\xi)$ generate $\Ugg$. Introduce an action of the coefficients
of the series $\Lti^\pm_{jk}(\xi)$ in the space $\Fg_q$. Namely, for any
$f\in\Fql$ set:
\ifMag
$$
\NN4>
\gather
\Lline{\kern5em (\Lti^\pm_{11}(\xi)f)(t_1\lc t_l)\;=}
\Tagg{Lijxi}
\\
\Rline{
\aligned
&{}=\,f(t_1\lc t_l)\,\pron\,{\Sinp{u-z_m+\La_m}\over\Sinp{u-z_m-\La_m}}\;
\prod_{a=1}^l\,{\Sinp{u-t_a-1}\over\Sinp{u-t_a}}\ \,+
\\
&{}\)+\,\sin(\pi/p)\,\prod_{a=1}^l\,{\Sinp{u-t_a-1}\over\Sinp{u-t_a}}\ \,\x
\\
&\>\x\,\sum_{a=1}^l\,\LBc\,f(t_1\lc t_{l-1},u)\;
{\Expp{u-t_l}\over\Sinp{u-t_l-1}}\;
\pron\,{\Sinp{t_l-z_m+\La_m}\over\Sinp{t_l-z_m-\La_m}}\,\RBc_{(a,l)}\,,
\endaligned}
\kern-1em
\\
\ald
\nn6>
\Lline{(\Lti^\pm_{22}(\xi)f)(t_1\lc t_l)\,=\,f(t_1\lc t_l)\,
\prod_{a=1}^l\,{\Sinp{u-t_a+1}\over\Sinp{u-t_a}}\ \,-}
\\
\Rline{{}-\,\sin(\pi/p)\,
\prod_{a=1}^l\,{\Sinp{u-t_a+1}\over\Sinp{u-t_a}}\;\sum_{a=1}^l\,
\LBc\,f(u,t_2\lc t_l)\;{\Expp{u-t_1}\over\Sinp{u-t_1+1}}\,\RBc_{(1,a)}\,,}
\kern-1em
\\
\noalign{\fixedpage}
\ald
\nn6>
\Lline{(\Lti^\pm_{12}(\xi)f)(t_1\lc t_{l+1})\,=\,
\sin(\pi/p)\,\sum_{a=1}^{l+1}\,\LBc\,f(t_2\lc t_{l+1})\;
{\Expp{u-t_1}\over\Sinp{u-t_1}}\ \,\x}
\\
\Rline{
\aligned
{}\x\,\Bigl(\;\pron\,{\Sinp{t_1-z_m+\La_m}\over\Sinp{t_1-z_m-\La_m}}\;
\prod_{b=2}^{l+1}\>\>{\Sinp{u-t_b+1}\over\Sinp{u-t_b}}\,
{\Sinp{t_1-t_b-1}\over\Sinp{t_1-t_b+1}}\ \;- &
\\
\nn2>
\Llap{{}-\,\pron\,{\Sinp{u-z_m+\La_m}\over\Sinp{u-z_m-\La_m}}\,
\prod_{b=2}^{l+1}\,{\Sinp{u-t_b-1}\over\Sinp{u-t_b}}\;\Bigr)\>
\RBc_{(1,a)}\,-} &
\endaligned}
\\
\ald
\Rline{
\aligned
&{}\)-\,\sin^2(\pi/p)\,\prod_{a=1}^{l+1}\,{\Sinp{u-t_a+1}\over\Sinp{u-t_a}}\,
\sum_{\tsize{a,b=1\atop a\ne b}}^{l+1}\,\LBc\,f(u,t_2\lc t_l)\;\x
\\
&\>{}\x\;{\Expp{2u-t_1-t_{l+1}}\over\Sinp{u-t_1+1}\>\Sinp{u-t_{l+1}+1}}\;
\pron\,{\Sinp{t_{l+1}-z_m+\La_m}\over\Sinp{t_{l+1}-z_m-\La_m}}\,
\RBc_{\si^{ab}}\,,\!\}
\endaligned}
\\
\ald
\nn4>
(\Lti^\pm_{21}(\xi)f)(t_1\lc t_{l-1})\,=\,f(t_1\lc t_{l-1},u)\,
\prod_{a=1}^{l-1}\,{\Sinp{u-t_a-1}\over\Sinp{u-t_a}}\;,\qqq l>0\,,
\endgather
$$
\else
$$
\NN4>
\gather
\Rline{(\Lti^\pm_{11}(\xi)f)(t_1\lc t_l)\,=\,
f(t_1\lc t_l)\,\pron\,{\Sinp{u-z_m+\La_m}\over\Sinp{u-z_m-\La_m}}\;
\prod_{a=1}^l\,{\Sinp{u-t_a-1}\over\Sinp{u-t_a}}\ \,+\!}
\Tagg{Lijxi}
\\
\Rline{
\aligned
&{}+\,\sin(\pi/p)\,\prod_{a=1}^l\,{\Sinp{u-t_a-1}\over\Sinp{u-t_a}}\ \,\x
\\
&{}\)\x\,\sum_{a=1}^l\,\LBc\,f(t_1\lc t_{l-1},u)\;
{\Expp{u-t_l}\over\Sinp{u-t_l-1}}\;
\pron\,{\Sinp{t_l-z_m+\La_m}\over\Sinp{t_l-z_m-\La_m}}\,\RBc_{(a,l)}\,,
\endaligned}
\kern-1em
\\
\ald
\nn6>
\Lline{(\Lti^\pm_{22}(\xi)f)(t_1\lc t_l)\,=\,f(t_1\lc t_l)\,
\prod_{a=1}^l\,{\Sinp{u-t_a+1}\over\Sinp{u-t_a}}\ \,-}
\\
\Rline{{}-\,\sin(\pi/p)\,
\prod_{a=1}^l\,{\Sinp{u-t_a+1}\over\Sinp{u-t_a}}\;\sum_{a=1}^l\,\LBc\,
f(u,t_2\lc t_l)\;{\Expp{u-t_1}\over\Sinp{u-t_1+1}}\,\RBc_{(1,a)}\,,}
\kern-1em
\\
\ald
\nn6>
\Lline{
(\Lti^\pm_{12}(\xi)f)(t_1\lc t_{l+1})\,=\,\sin(\pi/p)\,\sum_{a=1}^{l+1}\,
\LBc\,f(t_2\lc t_{l+1})\;{\Expp{u-t_1}\over\Sinp{u-t_1}}\ \,\x}
\\
\Lline{
\alignedat2
&\!\}{}\x\,\Bigl(\;\pron\,{\Sinp{t_1-z_m+\La_m}\over\Sinp{t_1-z_m-\La_m}}\;
\prod_{b=2}^{l+1}\>\>{\Sinp{u-t_b+1}\over\Sinp{u-t_b}}\,
{\Sinp{t_1-t_b-1}\over\Sinp{t_1-t_b+1}}\ \;- &&
\\
\nn2>
&& \Llap{{}-\,\pron\,{\Sinp{u-z_m+\La_m}\over\Sinp{u-z_m-\La_m}}\,
\prod_{b=2}^{l+1}\,{\Sinp{u-t_b-1}\over\Sinp{u-t_b}}\;\Bigr)\>
\RBc_{(1,a)}\,-\qqq} &
\\
&\!{}-\,\sin^2(\pi/p)\,
\prod_{a=1}^{l+1}\,{\Sinp{u-t_a+1}\over\Sinp{u-t_a}}\ \;\x
\\
&\}{}\x\,\sum_{\tsize{a,b=1\atop a\ne b}}^{l+1}\,\LBc\,f(u,t_2\lc t_l)\;
{\Expp{2u-t_1-t_{l+1}}\over\Sinp{u-t_1+1}\>\Sinp{u-t_{l+1}+1}}\;
\pron\,{\Sinp{t_{l+1}-z_m+\La_m}\over\Sinp{t_{l+1}-z_m-\La_m}}\,
\RBc_{\si^{ab}}\,, &&
\endalignedat}
\\
\ald
\nn4>
(\Lti^\pm_{21}(\xi)f)(t_1\lc t_{l-1})\,=\,f(t_1\lc t_{l-1},u)\,
\prod_{a=1}^{l-1}\,{\Sinp{u-t_a-1}\over\Sinp{u-t_a}}\;,\qqq l>0\,,
\endgather
$$
\fi
and $\Lti^\pm_{21}(u)f=0$ for $f\in\Fq[0]$. Here ${\xi=\exp(2\pii\)u/p)}$,
$(1,a)$, $(a,l)$ are transpositions and $\si^{ab}\in\S^{l+1}$
is the following \perm/
$$
\si^{ab}\!:i\>\map\)i\quad\for i=2\lc l\,,\qqq
\si^{ab}\!:1\>\map\)a\,,\qquad\si^{ab}\!:l+1\>\map\)b\,.
$$
\Rhs/s of formulae \(Lijxi) are \raf/s in $\xi$, and the precise meaning of
each of the formulae is that $\Lti^-_{jk}(\xi)$ equals the Laurent series
expansion of the corresponding right hand side at $\xi=\8$, and
$\Lti^+_{jk}(\xi)$ equals the Taylor series expansion of the corresponding
right hand side at $\xi=0$.
\Lm{UFg}
Formulae \(Lijxi) define a \Uhmod/ structure in the \thgF/ $\Fg_q$.
\endpro
\nt
The proof is given by direct verification.
\Par
Let $\Fg_q(z)=\Plus_{l\ge 0}\Fql(z)$ be the \thgF/ of a fiber. The \Uhmod/
structure in $\Fg_q$ clearly induces a \Uhmod/ structure in $\Fg_q(z)$.
This module will be called the \em{\thgm/}.
\Par
For the action of the generators of the subalgebra $\Uu$ formulae \(Lijxi)
simplify and for ${f\in\Fql}$ look as follows:
\ifMag\vadjust{\fixedpage}\else\vvn->\fi
$$
\gather
(q^{\pm H}\]f)(t_1\lc t_l)\,=\,q^{\pm(\!\sun\La_m-l)}\>f(t_1\lc t_l)\,,
\Tag{HFEq}
\\
\nn3>
\ifMag
\alignedat2
(F_qf) &(t_1\lc t_{l+1})\,=\,\exp\bigl(-\pii\bigl(l+\tsun\La_m\bigr)/p\bigr)\>
\sum_{a=1}^{l+1}\,\Bigl[\,f(t_2\lc t_{l+1})\;\x &&
\\
\nn6>
\x\,&\Bigl(\exp(2\pii\>l/p)\,
\pron\,{\Sinp{t_1-z_m+\La_m}\over\Sinp{t_1-z_m-\La_m}}\,\,
\prod_{b=2}^{l+1}\,{\Sinp{t_1-t_b-1}\over\Sinp{t_1-t_b+1}}\ \;- &&
\\
\nn4>
&& \Llap{-\;\exp\bigl(2\pii\tsun\La_m/p\bigr)\Bigr)\,\Bigr]_{(1,a)}\!} & \,\,,
\endalignedat
\else
{\align
& \Lline{(F_qf)(t_1\lc t_{l+1})\,=\,
\exp\bigl(-\pii\bigl(l+\tsun\La_m\bigr)/p\bigr)\>
\sum_{a=1}^{l+1}\,\Bigl[\,f(t_2\lc t_{l+1})\;\x}
\\
\nn4>
& \Rline{\x\,\Bigl(\exp(2\pii\>l/p)\,
\pron\,{\Sinp{t_1-z_m+\La_m}\over\Sinp{t_1-z_m-\La_m}}\,\,
\prod_{b=2}^{l+1}\,{\Sinp{t_1-t_b-1}\over\Sinp{t_1-t_b+1}}\,-\)
\exp\bigl(2\pii\tsun\La_m/p\bigr)\Bigr)\,\Bigr]_{(1,a)}\,,\!}
\endalign}
\fi
\\
\nn8>
\ifMag
{\align
(E_qf)(t_1\lc t_{l-1})\,=\,-\bigl(2i\sin(\pi/p)\bigr)\1\>
\exp\bigl(\pii\bigl(l-1+\tsun\La_m\bigr)/p\bigr)\;\x &
\\
\nn4>
\x\;f(\tell)\vst{\exp(2\pii t_\ell/p)=0} &\,,\qqq\Rlap{l>0\,,}
\endalign}
\else
\Line{(E_qf)(t_1\lc t_{l-1})\,=\,-\bigl(2i\sin(\pi/p)\bigr)\1\>
\exp\bigl(\pii\bigl(l-1+\tsun\La_m\bigr/p\bigr)\>
f(\tell)\vst{\exp(2\pii t_\ell/p)=0}\,,\hfill l>0\,,\!}
\fi
\endgather
$$
and $E_qf=0$ for $f\in\Fq[0]$. Here $(1,a)\in\S^{l+1}$ are transpositions.

\subsect{Tensor coordinates on the \tri/ \hgeom/ spaces of fibers}
Let $q=\exp(\pii/p)$. Let $\Vqn$ be $\Uu$ Verma modules with \hw/s $\qLan$ and
\gv/s $\vqn$, \resp/. Consider a \wt/ subspace $\Vql$ with a basis given by
monomials $\Fvq$, $\lg\in\Zln$. For any $z\in\BB$ and for any $\tau\in\S^n$
denote by $C_\tau(z)$ the following \hom/:
$$
\gather
{\align
C_\tau(z):\Vqlt\, &{}\to\,\Fq(z)\,,
\\
\nn3>
C_\tau(z):\Fvqt\, &{}\map\,c_\lg\:\>W^\tau_\lg(t,z)\,,
\endalign}
\\
\nn6>
\Text{where $\Fq(z)$ is the \thgf/ of the fiber and}
\nn4>
c_\lg\:\,=\pron\,\prod_{s=0}^{\lg_m-1}
{\Sinp{s+1}\)\Sinp{2\La_m-s}\over\sin(\pi/p)}\;,
\endgather
$$
(\cf. \(wlga), \(wtau)\)). The \hom/s $C_\tau(z)$ are called the \em{\tenco/}
on the \thgf/ of a fiber. The composition maps
$$
C_{\tau,\tau'}(z):\Vqltt\,\to\,\Vqlt\,,\qqq
C_{\tau,\tau'}(z)\,=\,C\1_\tau(z)\o C_{\tau'}(z)\,,
$$
are called the \em{\traf/s}, \cf. \Cite{V3}.
\Lm{tcfsq}
Let ${z_l+\La_l-z_m+\La_m-s\,\nin\,p\)\Z}$ for any $s=0\lc\ell-1$, and for any
$l,\mn$. Then for any \perm/ $\tau$ the linear map
$C_\tau(z):\Vqlt\to\Fq(z)$ is \ndeg/.
\endpro
\nt
The statement follows from Lemma~\[Wbasis].
\Par
Let ${\zt_m=\exp(2\pii\)z_m/p)}$. Consider a tensor product $\Vqxzt$ of \emod/s
over $\Ugg$ coinciding with $\Vqxt$ as a \Umod/.
\Lm{L21}
For any $v\in\Vqlt$ we have
$$
\align
& L^\pm_{21}(\xi_1)\ldots L^\pm_{21}(\xi_\ell)\>v\,=\,
\bigl(C_\tau(z)\)v\bigr)(\tell)\;\x
\\
\nn4>
&\quad\x\,\pral\,\pron\,
{2\Sinp{u-z_m-\La_m}\over\pm i\exp\bigl(\pm\pii(z_m-u)/p\bigr)}\,
\prab\,{\Sinp{t_a-t_b+1}\over\Sinp{t_a-t_b}}\ \vqxt\;,
\endalign
$$
where ${\xi_a=\exp(2\pii\>t_a/p)}$, $\aell$.
\endpro
\nt
It is easy to see that \rhs/ above is a \pol/ in $\xi_1\lc\xi_\ell$ for
the case of the upper signs, and is a \pol/ in $\xi_1\1\lc\xi_\ell\1\!$ for
the case of the lower signs, so the formula makes sense without additional
prescriptions.
\Th{FgVqx}
For any \perm/ $\tau\in\S^n$ the map
$$
C_{\tau}(z)\):\)\Vqxzt\,\to\,\Fg_q(z)
$$
is an intertwiner of \Uhmod/s.
\endpro
\Cr*
Let ${z_l+\La_l-z_m+\La_m-s\,\nin\,p\)\Z}$ for any $s\in\Zp$,
and for any $l,\mn$. Then for any \perm/ $\tau\in\S^n$ the map
\ifMag\nl\fi
${C_{\tau}(z):\Vqxzt\,\to\,\Fg_q(z)}$ is an \iso/ of \Uhmod/s.
\endpro
\nt
The statement follows from Theorem~\[FgVqx] and Lemma~\[tcfsq].
\Cr{Csing}
For any $\tau\in\S^n$ the \hom/ $C_\tau(z)$ maps $\Vqlst\!$
into the \sthgf/ $\Fqs(z)$ of a fiber. The map
$$
C_{\tau}(z)\):\)\Vqlst\to\,\Fqs(z)
$$
is an \iso/ provided that ${z_l+\La_l-z_m+\La_m-s\,\nin\,p\)\Z}$
for any $s=0\lc\ell-1$, and for any $l,\mn$.
\endpro
\nt
The statement follows the last formula in \(HFEq).
\Th{localq}
\back\Cite{V3}
For any $\tau\in\S^n$ and any transposition $(m,m+1)$, $\mn-1$, the \traf/
$$
C_{\tau,\tau\cdot(m,m+1)}(z):\Vqxtm\,\to\,\Vqxt
$$
equals the operator
${P_{V^q_{\tau_{m+1}}\}V^q_{\tau_m}}\:R^q_{V^q_{\tau_{m+1}}\}V^q_{\tau_m}}\!\!
\bigl(\Expt{z_{\tau_{m+1}}-z_{\tau_m}}\bigr)\)}$
acting in the \^{$m$-th} and \^{$(m+1)$-th} factors.
\endpro
\nt
The theorem follows from Lemma~\[L21] and formula \(Uintw).

\subsect{Tensor products of the \hmod/s}
Let ${\Fo[z_1\lc z_m;\La_1\lc\La_m;l\,]}$ and
${\Fq[z_1\lc z_m;\La_1\lc\La_m;l\,]}$ be \resp/ the \rat/ and the \thgf/s
defined for the projection ${\C^{\,l+m}\!\to\C^{\>m}\!}$. In particular,
in our previous notations we have
$$
\Fo=\Fo[\zn;\Lan;\ell\,]\qquad\text{and}\qquad\Fq=\Fq[\zn;\Lan;\ell\,]\,.
$$
There are maps
$$
\NN6>
\align
\chi:\)\Fo[z_1\lc z_k;\La_1\lc\La_k;j\)] &{} \ox
\Fo[z_{k+1}\lc z_{k+m};\La_{k+1}\lc\La_{k+m};l\,]\,\to{}
\\
& {}\quad\ \;\to\,\Fo[z_1\lc z_{k+m};\La_1\lc\La_{k+m};j+l\,]
\\
\nn-12>
\Text{and}
\ifMag\else\nn-10>\fi
\chi_q:\]\Fq[z_1\lc z_k;\La_1\lc\La_k;j\)] &{} \ox
\Fq[z_{k+1}\lc z_{k+m};\La_{k+1}\lc\La_{k+m};l\,]\,\to{}
\\
& {}\quad\ \;\to\,\Fq[z_1\lc z_{k+m};\La_1\lc\La_{k+m};j+l\,]
\endalign
$$
which are \resp/ defined by ${\chi:f\ox g\map f\*g}$ and
${\chi_q\]:f\ox g\map f*g}$ where
\ifMag
$$
\align
(f\*g) &(t_1\lc t_{j+l})\;=
\\
&\!{}=\;{1\over j!\,l\)!}\sum_{\si\in\S^{j+l}}\,
\Bigl[\>f(t_1\lc t_j)\>g(t_{j+1}\lc t_{j+l})\,\prod_{i=1}^k\,
\prod_{a=1}^l\,{t_{a+j}-z_i+\La_i\over t_{a+j}-z_i-\La_i}\,\Bigr]_\si
\\
\Text{and}
(f*g) &(t_1\lc t_{j+l})\;=
\\
&\!{}=\;{1\over j!\,l\)!}\sum_{\si\in\S^{j+l}}\,
\LBc\>f(t_1\lc t_j)\>g(t_{j+1}\lc t_{j+l})\,\prod_{i=1}^k\,\prod_{a=1}^l\,
{\Sinp{t_{a+j}-z_i+\La_i}\over\Sinp{t_{a+j}-z_i-\La_i}}\,\RBc_\si.
\endalign
$$
\else
$$
\gather
(f\*g)(t_1\lc t_{j+l})\,=\,{1\over j!\,l\)!}\sum_{\si\in\S^{j+l}}\,\Bigl[\>
f(t_1\lc t_j)\>g(t_{j+1}\lc t_{j+l})\,\prod_{i=1}^k\,
\prod_{a=1}^l\,{t_{a+j}-z_i+\La_i\over t_{a+j}-z_i-\La_i}\,\Bigr]_\si,
\\
\Text{and}
(f*g)(t_1\lc t_{j+l})\,=\,{1\over j!\,l\)!}\sum_{\si\in\S^{j+l}}\,\LBc\>
f(t_1\lc t_j)\>g(t_{j+1}\lc t_{j+l})\,\prod_{i=1}^k\,\prod_{a=1}^l\,
{\Sinp{t_{a+j}-z_i+\La_i}\over\Sinp{t_{a+j}-z_i-\La_i}}\,\RBc_\si.
\endgather
$$
\fi
We have the next lemmas.
\Lm{FFo}
Assume that ${(z_i-\La_i-z_{k+j}-\La_{k+j}+s)\ne 0}$ for any $i=1\lc k$,
$j=1\lc m$, $s=0\lc l-1$. Then the map
$$
\alignat2
\chi:\)\Plus_{i+j=l} &
\Fo[z_1\lc z_k;\La_1\lc\La_k;i\)]\)\bigl((z_1\lc z_k)\bigr)\ox{} &&
\\
\nn-2>
&\!{}\ox\Fo[z_{k+1}\lc z_{k+m};\La_{k+1}\lc\La_{k+m};j\)]\)
\bigl((z_{k+1}\lc z_{k+m})\bigr)\,\to{} &&
\\
\nn6>
&&\Llap{{}\to\,\Fo[z_1\lc z_{k+m};\La_1\lc\La_{k+m};l\,]\)
\bigl((z_1\lc z_{k+m})\bigr)}\! &
\endalignat
$$
defined by linearity is an \iso/ of the \rhgf/s of fibers.
\endpro
\Lm{FFq}
Assume that ${(z_i-\La_i-z_{k+j}-\La_{k+j}+s)\nin p\)\Z}$
for any $i=1\lc k$, $j=1\lc m$, $s=0\lc l-1$. Then the map
$$
\alignat2
\chi_q\]:\)\Plus_{i+j=l} &
\Fq[z_1\lc z_k;\La_1\lc\La_k;i\)]\)\bigl((z_1\lc z_k)\bigr)\ox{} &&
\\
\nn-2>
&\!{}\ox\Fq[z_{k+1}\lc z_{k+m};\La_{k+1}\lc\La_{k+m};j\)]\)
\bigl((z_{k+1}\lc z_{k+m})\bigr)\,\to{} &&
\\
\nn6>
&&\Llap{{}\to\,\Fq[z_1\lc z_{k+m};\La_1\lc\La_{k+m};l\,]\)
\bigl((z_1\lc z_{k+m})\bigr)}\! &
\endalignat
$$
defined by linearity is an \iso/ of the \thgf/s of fibers.
\endpro
Let
\vvn->
$$
\align
\Fg[z_1\lc z_m;\La_1\lc\La_m]\, & {}=\,
\Plus_{l=0}^\8\)\Fo[z_1\lc z_m;\La_1\lc\La_m;l\,]
\\
\nn-3>
\Text{and}
\nn-3>
\Fg_q[z_1\lc z_m;\La_1\lc\La_m]\, & {}=\,
\Plus_{l=0}^\8\)\Fq[z_1\lc z_m;\La_1\lc\La_m;l\,]
\endalign
$$
be the \rat/ and the \thgF/s, \resp/.
Extend the maps $\chi,\chi_q$ to the \resp/ maps
$$
\NN4>
\align
\chi:{} & \Fg[z_1\lc z_k;\La_1\lc\La_k]\)\bigl((z_1\lc z_k)\bigr)\;\ox
\\
&\,\){}\ox\,\Fg[z_{k+1}\lc z_{k+m};\La_{k+1}\lc\La_{k+m}]\)
\bigl((z_{k+1}\lc z_{k+m})\bigr)\,\to{}
\\
&\kern6.6em {}\to\,\Fg[z_1\lc z_{k+m};\La_1\lc\La_{k+m}]\)
\bigl((z_1\lc z_{k+m})\bigr)\,,
\\
\nn4>
\chi_q\]:\){} &
\Fg_q[z_1\lc z_k;\La_1\lc\La_k]\)\bigl((z_1\lc z_k)\bigr)\;\ox
\\
&\ \;{}\ox\,\Fg_q[z_{k+1}\lc z_{k+m};\La_{k+1}\lc\La_{k+m}]\)
\bigl((z_{k+1}\lc z_{k+m})\bigr)\,\to{}
\\
&\kern7.45em {}\to\,\Fg_q[z_1\lc z_{k+m};\La_1\lc\La_{k+m}]\)
\bigl((z_1\lc z_{k+m})\bigr)\,.
\endalign
$$
\Th{FgF}
The map
$$
\NN4>
\align
\chi\o\}P\):\)\Fg[z_{k+1}\lc z_{k+m}; &
\La_{k+1}\lc\La_{k+m}]\)\bigl((z_{k+1}\lc z_{k+m})\bigr)\,\ox{}
\\
 & {}\!\]\ox\,\Fg[z_1\lc z_k;\La_1\lc\La_k]\)\bigl((z_1\lc z_k)\bigr)\,\to{}
\\
&\hp{\ox\Fg\!}{}\to\,
\Fg[z_1\lc z_{k+m};\La_1\lc\La_{k+m}]\)\bigl((z_1\lc z_{k+m})\bigr)
\endalign
$$
is an intertwiner of \Ymod/s. Here $P$ is the \perm/ map.
The map ${\chi\)\o P}$ is an \iso/ provided that
${(z_i-\La_i-z_{k+j}-\La_{k+j})\nin\Zn}$ for any $i=1\lc k$, $j=1\lc m$.
\endpro
\Th{FgFq}
The map
$$
\NN4>
\align
\chi_q\]:\){} &
\Fg_q[z_1\lc z_k;\La_1\lc\La_k]\)\bigl((z_1\lc z_k)\bigr)\;\ox
\\
&\ \;{}\ox\,\Fg_q[z_{k+1}\lc z_{k+m};\La_{k+1}\lc\La_{k+m}]\)
\bigl((z_{k+1}\lc z_{k+m})\bigr)\,\to{}
\\
&\kern7.45em {}\to\,\Fg_q[z_1\lc z_{k+m};\La_1\lc\La_{k+m}]\)
\bigl((z_1\lc z_{k+m})\bigr)\,.
\endalign
$$
is an intertwiner of \Uhmod/s. The map $\chi_q$ is an \iso/ provided that
\ifMag\else\nl\fi
${(z_i-\La_i-z_{k+j}-\La_{k+j}+s)\nin p\)\Z}$ for any $i=1\lc k$, $j=1\lc m$,
$s\in\Zp$.
\endpro
\goodbe
It is clear that for any \fn/s $f,g,h$ we have $(f\*g)\*h\,=\,f\*(g\*h)$ and
for any \fn/s $f,g,h$ we have $(f*g)*h\,=\,f*(g*h)$. Lemmas~\[FFo], \[FFq] and
Theorems~\[FgF], \[FgFq] can be extended naturally to an arbitrary number of
factors.

\Sect[2c]{The \hpair/}
In this section we define the main object of this paper, the \hpair/.
We define a pairing between the \rat/ and the \tri/ \hgf/s of a fiber.
For any \fn/s $w\in\Fo(z)$ and $W\in\Fq(z)$ we define the \em{\hint/} by
$$
I(W,w)\,=\int_{\){\IIt^\ell}\]}\Phi(t)\>w(t)\>W(t)\>\dt
\Tag{IWwt}
$$
where $\Phi$ is the \phf/ \(Phi) and $\IIt^\ell$ is a suitable deformation
of the imaginary subspace
$$
\II^\ell=\lb\)t\in\Cl\vert\Re t_1=0\llc\Re t_\ell=0\)\rb\,.
$$
\par
We always assume that the step $p$ is real and negative.
The case of arbitrary step can be treated by \anco/.
\par
The \phf/ $\Phi$ has a factor $\exp\bigl(\mu\}\tsum_{a=1}^\ell\}t_a/p\bigr)$
where the parameter $\mu$ is connected with the parameter $\ka$
in the definition of the \conn/s by $\ka=e^\mu$.
We choose the parameter $\mu$ so that it satifies
\ifMag\else\vvn-.5>\fi
$$
\Immu\,.
\Tag{Immu}
$$
\par
We define the \hint/ as follows. First we assume that the real parts
of the \wt/s $\Lan$ are large negative and set
$$
I(W,w)\,=\int_{\){\II^\ell}\]}\Phi(t)\>w(t)\>W(t)\>\dt\,.
\Tag{IWw}
$$
\Lm*
Let $\Imu$. Let the real parts of the \wt/s $\Lan$ be large negative. Then the
\hint/ $I(W,w)$ is well defined for any \fn/s $w\in\Fo(z)$ and $W\in\Fq(z)$.
\endpro
\Pf.
It follows from \(Phi), \(Stir) and \(thgf) that the integrand of the \hint/
decays exponentially as $t$ goes to infinity.
\epf
Let $\Fqs(z)\sub\Fq(z)$ be the \sthgf/.
\Lm{Wwell}
Let $\Im\mu=0$. Let the real parts of the \wt/s $\Lan$ be large negative.
Then the \hint/ $I(W,w)$ is well defined for any \fn/s $w\in\Fo(z)$ and
$W\in\Fqs(z)$.
\endpro
\nt
The proof is similar to the previous lemma.
\Par
The \hint/ for generic $\Lan$, $\zn$ and arbitrary negative $p$ is defined
by \anco/ \wrt/ $\Lan$, $\zn$ and $p$. This \anco/ makes sense since the
integrand is analytic in $\Lan$, $\zn$ and $p$, \cf. \(Phi), \(wlga), \(Wlga).
More precisely, first we define the \hint/ for basic \fn/s $w_\lg$, $W_\mg$
and then extend the definition by linearity to arbitrary \fn/s $w\in\Fo(z)$,
$W\in\Fq(z)$. The result of \anco/ can be represented as an integral of the
integrand over a suitably deformed imaginary subspace. Namely, the poles of
the integrand of the \hint/ $I(W_\lg,w_\mg)$ are located at the hyperplanes
$$
t_a=z_m\pm(\La_m+sp)\,,\qqq t_a=t_b\pm(1-sp)\,,
\Tag{list2}
$$
$1\le b<a\le\ell$, $\mn$, $s\in\Zp$. We deform $\Lan$, $\zn$ and $p$ in such
a way that the topology of the complement in $\Cl$ to the union of hyperplanes
\(list2) does not change. We deform accordingly the imaginary subspace
$\II^\ell$ so that it does not intersect the hyperplanes \(list2) at every
moment of the deformation. The deformed imaginary subspace is denoted by
$\IIt^\ell$. Then the \anco/ of the integral \(IWw) is given by formula
\(IWwt).
\Th{Ianco}
Let $\Imu$. Then for any $\lg,\mg\in\Zln$ the \hint/ $I(W_\lg,w_\mg)$ can be
analytically continued as a holomorphic univalued \fn/ of complex \var/s $p$,
$\Lan$, $\zn$ to the region:
\ifMag\else\vv-.8>\fi
$$
\gather
p<0\,,\qqq \npZ\,,
\\
\nn2>
2\La_m-s\,\nin\, p\)\Z\,,\qqq \mn\,,\qquad s=1-\ell\lc\ell-1\,,
\\
\nn3>
z_l\pm\La_l-z_m\pm\La_m-s\,\nin\, p\)\Z\,,\qqq l,\mn\,,\qquad l\ne m\,,
\endgather
$$
for an arbitrary combination of signs {\rm(}\cf. \(assum)\){\rm)}.
\endpro
\Th{Ianco0}
Let $\Im\mu=0$. Then for any $\lg\in\Zlm\!$, $\mg\in\Zln$ the \hint/
$I(\Wo_\lg,w_\mg)$ can be analytically continued as a holomorphic univalued
\fn/ of complex \var/s $p$, $\Lan$, $\zn$ to the region:
\ifMag\else\vv-.5>\fi
$$
\gather
p<0\,,\qqq \npZ\,,
\\
\nn2>
2\La_m-s\,\nin\, p\)\Z\,,\qqq \mn\,,\qquad s=1-\ell\lc\ell-1\,,
\\
\nn3>
z_l\pm\La_l-z_m\pm\La_m-s\,\nin\, p\)\Z\,,\qqq l,\mn\,,\qquad l\ne m\,,
\endgather
$$
for an arbitrary combination of signs {\rm(}\cf. \(assum)\){\rm)}.
\endpro
\nt
The theorems are proved in Section~\SNo{4}.
\Par
Let $\Rc(z)\sub\Fo(z)$ be the \cosub/.
\Lm{IWw=0}
Let $\mu=0$. Let $p<0$. Let \(npZ)\;--\;\(assum) hold.
Then the \hint/ $I(W,w)$ equals zero for any $w\in\Rc(z)$ and $W\in\Fqs(z)$.
\endpro
\nt
The lemma is proved in Section~\SNo{4}.
\Par
The \hint/ defines the \em{\hpair/}
$$
\gather
I:\Fq(z)\ox\Fo(z)\,\to\,\C
\Tag{hpFF}
\\
\Text{for $\Imu$, and}
\Ii\!:\Fqs(z)\ox\Fo(z)/\Rc(z)\,\to\,\C
\Tag{hpFFo}
\endgather
$$
for $\mu=0$. According to \(HFo) and \(HFR) this can be written as
$$
\gather
I:\Fq(z)\ox\H(z)\,\to\,\C
\Tag{1000}
\\
\Text{and}
\Ii\!:\Fqs(z)\ox\H(z)\,\to\,\C\, ,
\Tag{2000}
\endgather
$$
\resp/.
\goodbreak
\Th{mu<>0}
Let $\Imu$. Let $p<0$. Let \(npZ)\;--\;\(assum) hold.
Then the \hpair/ $I:\Fq(z)\ox\Fo(z)\to\C$ is \ndeg/. Moreover
\ifMag
$$
\NN4>
\align
\DWwi_{\lg,\mg\in\Zln}\) &{}=\,(2i)^{\ell\tsize{n+\ell-1\choose n-1}}\,
\ell\)!^{\tsize{n+\ell-1\choose n-1}}\,\x
\\
&\kern-2.6em{}\x\,(e^\mu-1)^{-\M\,\cdot{\tsize{n+\ell-1\choose n}}
+2n/p\,\cdot\]{\tsize{n+\ell-1\choose n+1}}}\,\x
\\
\nn-1>
&\kern-2.6em{}\x\,\exp\Bigl(\mu\tsun z_m/p\cdot\!{n+\ell-1\choose n}\Bigr)\;\x
\\
\nn2>
&\kern-2.6em{}\x\,\exp\Bigl((\mu+\pii)\>\Bigl(\,\tsun\La_m/p\cdot\!
{n+\ell-1\choose n}\)-\>n/p\cdot\!{n+\ell-1\choose n+1}\Bigr)\Bigr)\;\x
\\
&\kern-2.6em{}\x\,\pros\,\Bigl[\>\Gmpb{-(s+1)/p}^n\>\Gmpb{-1/p}^{-n}\,
\pron\Gmpp{2\La_m-s}\;\x
\\
\nn2>
&\kern-2.6em{}\x\,
\plmn\,{\Gmpp{z_l+\La_l-z_m+\La_m-s}\over\Gmpp{z_l-\La_l-z_m-\La_m+s}}\,
\Bigr]^{\tsize{n+\ell-s-2\choose n-1}}.
\endalign
$$
\else
$$
\NN4>
\align
\DWwi_{\lg,\mg\in\Zln}\) &{}=\,(2i)^{\ell\tsize{n+\ell-1\choose n-1}}\,
\ell\)!^{\tsize{n+\ell-1\choose n-1}}\>
(e^\mu-1)^{-\M\cdot{\tsize{n+\ell-1\choose n}}
+2n/p\cdot\]{\tsize{n+\ell-1\choose n+1}}}\,\x
\\
&\>{}\x\,\exp\Bigl(\mu\tsun z_m/p\cdot\!{n+\ell-1\choose n}\Bigr)\;\x
\\
\nn2>
&\>{}\x\,\exp\Bigl((\mu+\pii)\>\Bigl(\,\tsun\La_m/p\cdot\!
{n+\ell-1\choose n}\)-\>n/p\cdot\!{n+\ell-1\choose n+1}\Bigr)\Bigr)\;\x
\\
&\>{}\x\,\pros\,\Bigl[\>\Gmpb{-(s+1)/p}^n\>\Gmpb{-1/p}^{-n}\,
\pron\Gmpp{2\La_m-s}\;\x
\\
&\>{}\x\,\plmn\,{\Gmpp{z_l+\La_l-z_m+\La_m-s}\over\Gmpp{z_l-\La_l-z_m-\La_m+s}}
\,\Bigr]^{\tsize{n+\ell-s-2\choose n-1}}.
\endalign
$$
\fi
Here $0\le\arg(e^\mu-1)<2\pi$.
\endpro
\Th{mu=0}
Let $\mu=0$. Let $p<0$. Let \(npZ)\;--\;\(assum) hold.
If $\tLmn-s\nin p\)\Znn$ for all $s=\ell-1\lc 2\ell-2$, then the \hpair/
${\Ii\!:\Fqs(z)}{}\ox\Fo(z)/\Rc(z)\to\C$ is \ndeg/. Moreover
\ifMag
$$
\NN4>
\align
\kern1em
&\kern-1.1em
\DWwo_{\lg,\mg\in\Zlm}\)=\,(2i)^{\ell\tsize{n+\ell-2\choose n-2}}\,
\ell\)!^{\tsize{n+\ell-2\choose n-2}}\;\x
\\
&\>{}\x\,\pros\,\Bigl[\>\Gmpb{-(s+1)/p}^{n-1}\>\Gmpb{-1/p}^{1-n}\;\x
\\
&\>{}\x\,\Gmpb{1+\tM+(s+2-2\ell)/p}\1\,\Gmpb{1+(2\La_n-s)/p}\;\x
\\
&\>{}\x\,\prmn\Gmpp{2\La_m-s}\;
\plmn\,{\Gmpp{z_l+\La_l-z_m+\La_m-s}\over\Gmpp{z_l-\La_l-z_m-\La_m+s}}\,
\Bigr]^{\tsize{n+\ell-s-3\choose n-2}}.
\kern-1em
\endalign
$$
\else
$$
\NN4>
\align
&\kern-1.1em
\DWwo_{\lg,\mg\in\Zlm}\)=\,(2i)^{\ell\tsize{n+\ell-2\choose n-2}}\,
\ell\)!^{\tsize{n+\ell-2\choose n-2}}\;\x
\\
&\>{}\x\,\pros\,\Bigl[\>\Gmpb{-(s+1)/p}^{n-1}\>\Gmpb{-1/p}^{1-n}\,
\Gmpb{1+\tM+(s+2-2\ell)/p}\1\,\Gmpb{1+(2\La_n-s)/p}\;\x
\\
&\>{}\x\,\prmn\Gmpp{2\La_m-s}\;
\plmn\,{\Gmpp{z_l+\La_l-z_m+\La_m-s}\over\Gmpp{z_l-\La_l-z_m-\La_m+s}}\,
\Bigr]^{\tsize{n+\ell-s-3\choose n-2}}.
\endalign
$$
\fi
Here we identify $\mg\in\Zlm\!$ with $(\mg,0)\in\Zln$.
\endpro
Theorems~\[mu<>0] and \[mu=0] are proved in section \SNo{4}.
\Ex
Theorem~\[mu<>0] for $n=1$, $\ell=1$ and
Theorem~\[mu=0] for $n=2$, $\ell=1$ give
$$
\align
\int_{\!-i\8}^{\,i\8}\Gm(a+s)\>\Gm(a-s)\>u^{2s}\>ds\, &{}=\,
2\pii\>\Gm(2a)\>(u+u\1)^{-2a}\,,
\Tagg{Barnes}
\\
\int_{\!-i\8}^{\,i\8}\Gm(a+s)\>\Gm(b+s)\>\Gm(c-s)\>\Gm(d-s)\>ds\, &{}=
\,2\pii\;{\Gm(a+c)\>\Gm(a+d)\>\Gm(b+c)\>\Gm(b+d)\over\Gm(a+b+c+d)}\,,
\endalign
$$
which are formulae for the Barnes integrals \Cite{WW}.
\par
For arbitrary $\ell$, Theorem~\[mu<>0] for $n=1$ and Theorem~\[mu=0] for $n=2$
give the following Mellin-Barnes integrals, which are generalizations of
the famous Selberg integral:
$$
\align
\ifMag\kern1.9em\else\kern1em\fi
\int_{-i\8}^{i\8}\!\cdots\!\int_{-i\8}^{i\8}\; &
\prod_{k=1}^\ell\,\Bigl(u^{2s_k}\>\Gm(a+s_k)\>\Gm(a-s_k)\,
\prod_{j=1}^{k-1}\,{\Gm(s_j-s_k+x)\Gm(s_k-s_j+x)\over\Gm(s_j-s_k)\Gm(s_k-s_j)}
\>\Bigr)\,d^\ell s\;=
\ifMag\kern-1.9em\else\kern-1em\fi
\Tagg{BarneS}
\\
=\>&
\;(2\pii)^\ell\>(u+u\1)^{-\ell(2a+(\ell-1)x)}\,
\prod_{k=1}^\ell\,{\Gm(1+kx)\over\Gm(1+x)}\,\Gmpb{2a+(k-1)x}\,,
\endalign
$$
\vv-.3>
\ifMag
$$
\align
&\aligned
\int_{-i\8}^{i\8}\!\cdots\!\int_{-i\8}^{i\8}\;
\prod_{k=1}^\ell\,\Bigl(\Gm(a+s_k)\>\Gm(b+s_k)\> & \Gm(c-s_k)\>\Gm(d-s_k)\;\x
\\
\nn-3>
{}\x{}&
\prod_{j=1}^{k-1}\,{\Gm(s_j-s_k+x)\Gm(s_k-s_j+x)\over\Gm(s_j-s_k)\Gm(s_k-s_j)}
\>\Bigr)\,d^\ell s\;=
\endaligned
\\
\nn-8>
&=\,(2\pii)^\ell\>\prod_{k=1}^\ell\,\Bigl({\Gm(1+kx)\over\Gm(1+x)}\;\x
\\
\nn6>
&\>\x\,{\Gmpb{a+c+(k-1)x}\>\Gmpb{a+d+(k-1)x}\>
\Gmpb{b+c+(k-1)x}\>\Gmpb{b+d+(k-1)x}\over\Gmpb{a+b+c+d+(2\ell-k-1)x}}\>
\Bigr)\,,
\endalign
$$
\else
$$
\align
\int_{-i\8}^{i\8}\!\cdots\!\int_{-i\8}^{i\8}\;
\prod_{k=1}^\ell\,\Bigl(\Gm(a+s_k)\>\Gm(b+s_k)\>\Gm(c-s_k)\>\Gm(d-s_k)\,
\prod_{j=1}^{k-1}\,{\Gm(s_j-s_k+x)\Gm(s_k-s_j+x)\over\Gm(s_j-s_k)\Gm(s_k-s_j)}
\>\Bigr)\,d^\ell s\,=&
\\
\aligned
=\>&\;(2\pii)^\ell\>\prod_{k=1}^\ell\,\Bigl({\Gm(1+kx)\over\Gm(1+x)}\;\x
\\
\nn6>
&\kern2em \x\,{\Gmpb{a+c+(k-1)x}\>\Gmpb{a+d+(k-1)x}\>
\Gmpb{b+c+(k-1)x}\>\Gmpb{b+d+(k-1)x}\over\Gmpb{a+b+c+d+(2\ell-k-1)x}}\>
\Bigr)\!\Rlap{\,,}
\endaligned &
\endalign
$$
\fi
where $\Re a,b,c,d,u,x >0$.
\enddemo
\Rem
After this paper was written we found out that the second formula in \(BarneS)
had appeared in \Cite{G}. In Section~\SNo{4} we give a proof of the first
formula in \(BarneS) and use the formula to prove Theorems~\[mu<>0], \[mu=0].
\enddemo
\Rem
We also obtain determinant formulae similar to \(mu<>0) and \(mu=0) for the
\hpair/ in the \tri/ case \Cite{TV3}. Under the same specialization as above,
those formulae give multidimensional generalizations of the Askey-Roy formula
\Cite{GR, (4.11.2)}, and, on the other hand, can be viewed as a generalization
of the famous \^{$q$-}Selberg integral, \cf. \Cite{Ka}, \Cite{AK}.
\enddemo
\Rem
It is plausible that the assumptions on $p$, $\Lan$, $\zn$ of Theorems~\[mu<>0]
and \[mu=0] as well as of Theorems~\[kanon], \[kaon], \[tcnon], \[tcon],
\[IWw=0], \[asol], \[Asol], \[AsIWw] could be replaced by
the following weaker assumptions: the step $p$ is \st/
${\lb\)2\lc\ell\)\rb\cap p\)\Zpp=\Empty}$, the \wt/s $\Lan$ are \st/
$$
2\La_m-s\,\nin\, p\)\Z\,,\qqq \mn\,,\qquad s=0\lc\ell-1\,.
$$
and the coordinates $\zn$ obey the condition
$$
z_l+\La_l-z_m+\La_m-s\,\nin\, p\)\Z\,,\qqq l,\mn\,,\qquad l\ne m\,,
$$
for any $s=0\lc\ell-1$, so that $z\in\BB$.
\enddemo
Let $W$ be any element of the \thgf/ $\Fq$. The restriction of the \fn/ $W$
to a fiber defines an element $W|_z\in\Fq(z)$ of the \thgf/ of the fiber.
The \hpair/ allows us to consider the element $W|_z\in\Fq(z)$ as an element
$s_W\:(z)$ of the space $\H^*(z)$ dual to the \hcg/ $\H(z)$. This construction
defines a section of the bundle over $\Cn$ with fiber $\H^*(z)$.
\par
There is a simple but important statement.
\Th{sWz}
Let either {$\Imu$ and $W\in\Fq$} or {$\mu=0$ and $W\in\Fqs$}. Let $p<0$.
Let $\Lan$ obey \(Lass). Then the section $s_W\:$ is a periodic section
\wrt/ the \GM/.
\endpro
\nt
The theorem is proved in Section~\SNo{4}.
\Par
The section $s_W\:$ and the \tenco/ $B_\tau$ induce a section
$$
\Psi_W:z\,\map\,B_\tau^*\cdot W|_z\in\Vlt
\Tag{PsiW}
$$
of the \trib/ with fiber $\Vlt$. Theorems~\[sWz] and \[qKZ-GM] imply
\Cr{qKZsol}
The section $\Psi_W$ is a \sol/ to the \qKZe/.
\endpro
The \tenco/ $B_\tau(z)$, $C_{\tau'}(z)$ induce the \em{\hpair/}
$$
\gather
I_{\tau,\tau'}(z):\Vqltt\>\ox\>\Valt\,\to\,\C
\Tag{3000}
\\
\Text{if $\Imu$ and}
\nn2>
\Ii_{\tau,\tau'}(z):\Vqlstt\>\ox\>\VFVt\,\to\,\C
\Tag{4000}
\endgather
$$
if $\mu=0$, which also can be considered as maps
$$
\align
\Ich_{\tau,\tau'}(z):\Vqltt\, &{}\to\,\Vlt
\Tag{5000}
\\
\Text{and}
\Iic_{\tau,\tau'}(z):\Vqlstt\) &{}\to\Vlst.
\Tag{6000}
\endalign
$$
If $v\in\Vqltt$, then the \hpair/ defines a section
$$
\Psi_v:z\,\map\,\Ich_{\tau,\tau'}(z)\cdot v\in\Vlt
$$
and, if $v\in\Vqlstt$, then the \hpair/ defines a section
$$
\Psi_v:z\,\map\,\Iic_{\tau,\tau'}(z)\cdot v\in\Vlst
$$
\Cr*
Let $\Imu$ and, therefore, $\ka\ne 1$. Then for any $v\in\Vqltt$ the section
$\Psi_v$ is a \sol/ to the \qKZe/ with values in $\Vlt$. Under conditions of
Theorem~\[mu<>0] all \sol/s are constructed in this way.
\endpro
Therefore, for $\ka\ne 1$ we constructed the \hgeom/ maps
$$
\Ich_{\tau,\tau'}(z):\Vqxztt\,\to\,\Voxzt
$$
from \qlo/ modules to Yangian modules. Here $\zt_m=\exp(2\pii z_m/p)$, $\mn$.
The maps have the following properties:
$$
\align
& \Ich_{\tau\cdot(m,m+1),\tau'}(z)\,=\,P_{V_{\tau_m}\}V_{\tau_{m+1}}}\:
R_{V_{\tau_m}\}V_{\tau_{m+1}}}\:\!(z_{\tau_m}\!-z_{\tau_{m+1}})\bigr)\,
\Ich_{\tau,\tau'}(z)\,,
\\
\nn6>
& \Ich_{\tau,\tau'\cdot(m,m+1)}(z)\,=\,\Ich_{\tau,\tau'}(z)\,
P_{V^q_{\tau_{m+1}}\}V^q_{\tau_m}}\:R^q_{V^q_{\tau_{m+1}}\}V^q_{\tau_m}}\!\!
\bigl(\Expt{z_{\tau_{m+1}}-z_{\tau_m}}\bigr)\,,
\\
\inText{and as \fn/s of $z$ they satisfy the \qKZe/s:}
& \Ich_{\tau,\tau'}(z_1\lc z_{\tau_m}\!+p\lc z_n)\,=\,
K_m(z_{\tau_1}\lc z_{\tau_n})\>\Ich_{\tau,\tau'}(\zn)\,.
\endalign
$$
\Cr*
Let $\mu=0$ and, therefore, $\ka=1$. Then for any $v\in\Vqlstt$ the section
$\Psi_v$ is a \sol/ to the \qKZe/ with values in $\Vlst\!$. Under conditions
of Theorem~\[mu=0] all \sol/s are constructed in this way.
\endpro
\Rem
Let $\Vox$ be a tensor product of $\gsl\!$ \Vmod/s, $\Vt_1\lox\Vt_n$ the tensor
product of the corresponding \irr/ \smod/s, and ${S:\Vox\to\Vt_1\lox\Vt_n}$ the
natural projection. If $\Psi(z)$ is a \sol/ to the \qKZe/ with values in $\Vox$
then ${S\)\Psi(z)}$ is a \sol/ to the \qKZe/ with values in $\Vt_1\lox\Vt_n$.
\par
This observation shows that the previous constructions give all \sol/ to
the \qKZe/ with values in $(\Vt_1\lox\Vt_n)\elli$ if $\ell\le\dim\Vt_m$
for all $\mn$. Moreover, the space of \sol/s to the \qKZe/ with values
in $(\Vt_1\lox\Vt_n)\elli$ in this case is identified with the space
$(\Vt^q_1\lox\Vt^q_n)\elli\otimes\FF$ where $\Vt^q_1\lox\Vt^q_n$ is the tensor
product of the corresponding \irr/ \Umod/s, and $\FF$ is the space of \fn/s
of $\zn$ which are \p-periodic \wrt/ each of the \var/s.
\par
In a separate paper we shall explain how the construction of this paper gives
all \sol/s to the \qKZe/ with values in a tensor product of \irr/ \smod/s.
\enddemo

\Sect[2d]{Asymptotic \sol/s to the \{\qKZe/}
One of the most important characteristics of a \difl/ \eq/ is the monodromy
group of its \sol/s. For the \difl/ \KZ/ \eq/ with values in a tensor product
of \rep/s of a simple Lie algebra its monodromy group is described in terms of
the corresponding \qg/. This fact establishes a remarkable connection between
\rep/ theories of simple Lie algebras and their quantum groups, see \Cite{K},
\Cite{D2}, \Cite{KL}, \Cite{SV2}, \Cite{V2}, \Cite{V4}.
\par
The substitution of the monodromy group for \deq/s is the set of \traf/s
between \asol/s. For a \deq/ one defines suitable \azo/s in the domain of
the definition of the \eq/ and then an \asol/ for every zone. Thus, for
every pair of \azo/s one gets a \traf/ between the corresponding \asol/s.
\par
In this section we describe \azo/s, \asol/s, and their \traf/s for the \qKZe/
with values in a tensor product of \smod/s when the parameter $\ka$ is
different from $1$. A remarkable fact is that the \traf/s are described in
terms of the \tri/ \Rms/ acting in the tensor product of the corresponding
\Umod/s. This fact establishes a correspondence between \rep/ theories of
Yangians and \qlo/s, since the \qKZe/ is defined in terms of the \rat/ \Rm/
action in the tensor product of \smod/s (and, therefore, in terms of
the Yangian action), and the \tri/ \Rm/ action in the tensor product of
\Umod/s is defined in terms of the action of the \qlo/.
\Par
Let $V$ be a vector space of dimension $N$ for some $N$.
Consider an integrable system of \deq/s for a \Vval/ \fn/ $\Psi(\zn)$:
$$
\Psi(\zmn)\,=\,\A_m(\zn)\>\Psi(\zn)\,,\qqq\mn\,.
\Tag{syst}
$$
Let $\AA$ be a domain in $\Cn\!$. Say that a basis $\Psi_1\lc \Psi_N$ of
\sol/s to system \(syst) form an \em{\asol/} in the domain if
$$
\Psi_j(z)\,=\,\exp\bigl(\tsun a_{jm}z_m/p\bigr)
\pmln\!(z_l-z_m)^{b_{jlm}}\>\bigl(v_j+o(1)\bigr)\,,
\Tag{fzz}
$$
where $a_{jm}$ and $b_{jlm}$ are suitable numbers, $v_1\lc v_N$ are vectors
which form a basis in $V$, and $o(1)$ tends to $0$ as $z$ tends to infinity
in $\AA$. We will call the domain an \em{\azo/}.
\Par
Consider the \qKZe/ with parameter $\ka\ne 1$ and values in $\Vl$.
We describe its \asol/s in suitable \azo/s.
\par
For every \perm/ $\tau\in\S^n$ we consider an \azo/ $\AA_\tau$ in $\Cn\!$
given by
$$
\AA_\tau\,=\,\lb\)z\in\Cn\vert\Re z_{\tau_{1}}\lsym\ll\Re z_{\tau_{n}}\)\rb\,.
\Tag{azone}
$$
Say that $z\to\8$ in $\AA_\tau$ if $\Re(z_{\tau_m}-z_{\tau_{m+1}})\to-\8$
for all $\mn-1$.
\Par
Recall that for every \perm/ $\tau\in\S^n$ we constructed a basis
$W_\lg^\tau$, $\lg\in\Zln$, in the \thgf/. This basis defines a basis
$\Psi_{W_\lg^\tau}$, $\lg\in\Zln$, of \sol/s to the \qKZe/, \cf. \(PsiW).
\Th{asol}
Let $p<0$. Assume that the \wt/s $\Lan$ obey condition \(Lass).
Let $\Imu$ and, therefore, $\ka\ne 1$. Then for any \perm/ $\tau\in\S^n$
the basis $\Psi_{W_\lg^\tau}$, $\lg\in\Zln$, is an \asol/ in the \azo/
$\AA_\tau$. Namely,
$$
\align
\Psi_{W_\lg^\tau}(z)\,=\,\Tht_{\lg}\>\exp\bigl(\mu\tsun\lg_m\)z_m/p\bigr)
\plmn\!\bigl((z_{\tau_l}-z_{\tau_m})/p\bigr)^{2(\lg_{\tau_l}\)\La_{\tau_m}+
\lg_{\tau_m}\)\La_{\tau_l}-\lg_{\tau_l}\)\lg_{\tau_m})/p}\;\x &
\\
\x\;\bigl(\Fv+o(1) &\bigr)
\endalign
$$
as $z\to\8$ in $\AA_\tau$ \st/ at any moment assumption \(assum) holds.
Here $|\arg\bigl((z_{\tau_l}-z_{\tau_m})/p\bigr)|<\pi$ and $\Tht_\lg$ is
a constant independent of the \perm/ $\tau$ and given by
$$
\align
\kern.5em
\Tht_\lg\, &{}=\,(2i)^\ell\,\ell\)!\,\Gm(-1/p)^{-\ell}\,
\pron\)\Bigl[\)(e^\mu-1)^{(\lg_m(\lg_m-1)-2\lg_m\La_m)/p}\;\x
\\
&{}\>\x\,\)\exp\bigl((\mu+\pii)(\lg_m\La_m-\lg_m(\lg_m-1)/2)/p\bigr)\,
\prod_{s=0}^{\lg_m-1}\Gmpp{2\La_m-s}\>\Gmpb{-(s+1)/p}\)\Bigr]\,,
\kern-.5em
\endalign
$$
where $0\le\arg(e^\mu-1)<2\pi$.
\endpro
\par\nt
The theorem is proved in Section~\SNo{4}.
\Rem
The \qKZo/s $K_m(z)$ have the following \as/s in the \azo/ $\AA_\tau$,
$$
K_m(z)\,=\,\ka^{\La_m-H_m}\>\ono\,,\qqq\mn\,.
$$
The vectors $\Fv$ form an eigenbasis of the operator $\ka^{\La_m-H_m}$ with
\eva/s $\ka^{\lg_m}$.
\enddemo
\Rem
The \qKZe/ and the basis of \sol/s $\Psi_{W_\lg^\tau}$, $\lg\in\Zln$,
depend \mer/ally on parameters $\mu$, $\Lan$. The \as/s of the basis
$\Psi_{W_\lg^\tau}$, $\lg\in\Zln$, determine the basis uniquely. Namely,
if a basis of \sol/s \mer/ally depends on the parameters $\mu$, $\Lan$ and has
\as/s in $\AA_\tau$ described in Theorem~\[asol], then such a basis coincides
with the basis $\Psi_{W_\lg^\tau}\!$. In fact, elements of any such basis are
linear combinations of the \fn/s $\Psi_{W_\lg^\tau}\!$ with coefficients
\mer/ally depending on $\mu$, $\Lan$ and \p-periodic in $\zn$. To preserve
the \as/s one can add to an element $\Psi_{W_\lg^\tau}\!$ any other \fn/s
$\Psi_{W_{\lg'}^\tau}\!$ having smaller \as/s. If $\mu<0$, then one can add
only the \fn/s $\Psi_{W_{\lg'}^\tau}\!$ with $\lg'$ \lex/ly greater than
$\lg$, and if $\mu>0$, then one can add only the \fn/s $\Psi_{W_{\lg'}^\tau}\!$
with $\lg'$ \lex/ly smaller than $\lg$. Since the coefficients of added terms
are \mer/ they have to be zero.
\enddemo
\Ex
Theorem~\[asol] allows us to write a \tri/ \Rm/ as an infinite product of
\rat/ \Rms/. Namely, consider the \qKZe/ with values in the tensor product of
two $\gsl$ \Vmod/s $\VV$. Then there are two \azo/s $\Re z_1\ll\Re z_2$ and
$\Re z_1\grt\Re z_2$. Our result on the \traf/ from the first \azo/
to the second one is the following statement.
\par
For any $\gsl$ \Vmod/ $V$ let $V^q\!$ be $\Uu$ \Vmod/ corresponding to $V$.
Let $\La$ be the \hw/ of module $V$ and let $v,v^q$ be the respective \gv/s
of modules $V,V^q\!$. Define a map $G:V\to V^q$:
\ifMag\else\vv-.5>\fi$$
G\>:\>F^l\]v\,\map\,F_{\!q}^lv^q\,
\prod_{s=0}^{l-1}\>\Gmpb{1+(s-2\La)/p}\>\Gmpb{1+(s+1)/p}\,.
$$
Let $p,\mu$ be complex numbers \st/ $p<0$ and $\Imu$. Let $q=e^{\pii/p}\!$. Set
$$
\RV\:(x;\mu,p)\,=\,\exp\bigl(\mu x(\one\ox H)/p\bigr)\>\RVx\>
\exp\bigl(-\mu x(\one\ox H)/p\bigr)\,.
$$
and ${\,J(s,\mu)\,=\,
(G\ox G)\>\bigl(-is\)(e^{\mu/2}\}-e^{-\mu/2})\bigr)^{2H\ox H/p}\!}$,
\,where ${|\arg(-i(e^{\mu/2}\}-e^{-\mu/2}))|<\pi/2}$. Then
$$
\kern.6em
\lim_{s\to\8}\Bigl(J(s,\mu)\>\Bigl(\prod_{r=-s}^s\!\RV\:(x+rp;\mu,p)\Bigr)\,
J(s,\mu{\vp{\big)})\]}\1\>\Bigr)\,=\,
R^q_{V^q_1V^q_2}\bigl(\exp(-2\pii\)x/p)\bigr)\,.
\kern-.6em
\Tag{RpR}
$$
Here the factors of the product are ordered in such a way that $r$ grows from
right to left.
\par
Notice that the minus sign in the argument of the \Rm/ in \rhs/ of formula
\(RpR) above reflects the fact that we use the coproducts $\Dl$ and $\Dlq$
for the Yangian $\Ygl$ and the \qlo/ $\Ugg$ which are in a sense opposite to
each other.
\par
The restriction of formula \(RpR) to the \wt/ subspace $(\VV)_1\:$ of \wt/
$\La_1+\La_2-1$ can be transformed to the infinite product formula for
$2{\)\x}2$ matrices (\cf. \Cite{RF}\)), which looks as follows.
\par
\ifMag\else\vsk.1>\fi
Let $a,b,c,d,\thi$ be complex numbers, $\Re\thi>0$.
Set $\la=\sqrt{\vp{1^1}a^2-bc\)}$,
$$
h\,=\,\pmatrix 1 & 0\\ 0 & -1 \endpmatrix\,,\kern4em
A(u)\,=\;{1\over d+u}\,\pmatrix a+u & b\\ c & a-u\endpmatrix\,.
$$
and ${A(u;\thi)=\thi^{uh}\>A(u)\>\thi^{-uh}\!}$.
Assume that $-bc\)\ne s\)(s+2a)$ for any $s\in\Z$. Then
$$
\lim_{s\to\8}\Bigl(\)s^{-ah}\>h^s\>\Bigl(\prod_{r=-s}^s\!A(u+r;\thi)\Bigr)\,
h^s\)s^{ah}\>\Bigr)\,=\,A^q(u)
$$
where the factors of the product are ordered in such a way that $r$ grows
from right to left and
$$
A^q(u)\,=\;{1\over\Sinpi{d+u}}\, \pmatrix \Sinpi{a+u} &
\;{\dsize{\pi\)b\>(\thi+\thi\1)^{2a}\over\Gm(1+a+\la)\>\Gm(1+a-\la)}} \\
{\dsize{\pi\)c\>(\thi+\thi\1)^{-2a}\over\Gm(1+\la-a)\>\Gm(1-\la-a)}}\; &
\Sinpi{a-u} \endpmatrix\,.
$$
\enddemo
Theorem~\[asol] admits the following generalization. Fix a nonnegative integer
$k$ not greater than $n$. Let $n_0\lc n_k$ be nonnegative integers \st/
$$
\gather
0=n_0<n_1\lsym< n_k=n\,.
\\
\nn-6>
\Text{Set}
\ifMag
{\align
\Fo^i[\)l\>]\, &{}=\,
\Fo[z_{n_{i-1}+1}\lc z_{n_i};\La_{n_{i-1}+1}\lc\La_{n_i};l\,]
\\
\nn-6>
\Text{and}
\nn-6>
\Fq^i[\)l\>]\, &{}=\,
\Fq[z_{n_{i-1}+1}\lc z_{n_i};\La_{n_{i-1}+1}\lc\La_{n_i};l\,]
\endalign}
\else
\nn6>
\Line{\Fo^i[\)l\>]\,=\,
\Fo[z_{n_{i-1}+1}\lc z_{n_i};\La_{n_{i-1}+1}\lc\La_{n_i};l\,]\hfil
\text{and}\hfil \Fq^i[\)l\>]\,=\,
\Fq[z_{n_{i-1}+1}\lc z_{n_i};\La_{n_{i-1}+1}\lc\La_{n_i};l\,]\,.}
\fi
\endgather
$$
Then for any nonnegative integers $\ell_1\lc\ell_k$
\st/ $\sum_{i=1}^k\ell_i=\ell$ we have embeddings
$$
\Fo^1[\)\ell_1\)]\lox\Fo^k[\)\ell_k\)]\,\hto\,\Fo\qqq\text{and}\qqq
\Fq^1[\)\ell_1\)]\lox\Fq^k[\)\ell_k\)]\,\hto\,\Fq
$$
\wrt/ the tensor products introduced in Section~\SNo{2b}.
We consider an \azo/ in $\Cn$ given by
$$
\align
\AA^\nb\,=\,\lb\) z\in\Cn\vert &\Re z_{m_1}\lsym\ll\Re z_{m_k}\,,
\\
\nn3>
&\ \text{for all $\,m_1\lc m_k\,$ \st/ $\,n_{i-1}<m_i\le n_i\,$,
$\,i=1\lc k\)\rb\,$.}
\endalign
$$
We say that $z\to\8$ in $\AA^\nb$ if $\Re(z_l-z_m)\to-\8$ for all $l,m$ \st/
$n_{i-1}<l\le n_i<m\le n_{i+1}$ for some $i=1\lc k-1$, and $z_l-z_m$ remains
bounded for all $l,m$ \st/ $n_{i-1}<l,m\le n_i$ for some $i=1\lc k$
\Par
For any ${W\in\Fq^i[\)l\>]}$ let $\Psi_W(z_{n_{i-1}+1}\lc z_{n_i})$ be
the \sol/ to the \qKZe/ with values in $(V_{n_{i-1}+1}\lox V_{n_i})_l\:$
corresponding to $W$ (\cf. \(PsiW)\)).
\Th{Asol}
Let $p<0$. Assume that the \wt/s $\Lan$ obey condition \(Lass). Let $\Imu$
and, therefore, $\ka\ne 1$. Let $\ell_1\lc\ell_k$ be nonnegative integers
\st/ $\ell_1\lsym+\ell_k=\ell$. Let ${W_i\in\Fq^i[\)\ell_i\)]}$, $i=1\lc k$.
Let $W=W_1\lsym*W_k$. Then the \sol/ $\Psi_W(\zn)$ to the \qKZe/ with values
in $\Vl$ has the following \as/s as $z\to\8$ in $\AA^\nb$ \st/ at any moment
assumption \(assum) holds\/{\rm:}
$$
\align
\Psi_W(\zn)\,=\;{\ell\>!\over\ell_1!\ldots\ell_k!}\,
\prod_{1\le i<j\le k}\!\bigl((z_{n_i}-z_{n_j})/p\bigr)
^{2(\ell_i\!\!\!\sum_{m\in\Gms_j}\!\!\!\La_m\>+\,
\ell_j\!\!\!\sum_{m\in\Gms_i}\!\!\!\La_m\>-\,\ell_i\ell_j)/p} \;\x &
\\
\nn6>
\x\;\bigl(\Psi_{W_1}(z_1\lc z_{n_1})\lox\Psi_{W_k}(z_{n_{k-1}+1}\lc z_n)
+ o(1)\bigr)\! &\,.
\endalign
$$
Here $\Gms_i=\lb\)n_{i-1}+1\lc n_i\)\rb$ and
$|\arg\bigl((z_l-z_m)/p\bigr)|<\pi$.
\endpro
Theorem~\[asol] for $\tau=\id$ follows from Theorem~\[Asol] for $k=n$ so that
$n_j=j$, $j=0\lc n$, and the first formula in \(BarneS). Theorem~\[asol] for
a general \perm/ $\tau$ reduces to the same theorem for $\tau=\id$.
\par
Theorem~\[Asol] follows from the next statement on \as/s of the \hpair/.
\Th{AsIWw}
{\bls 1.1\bls
Let $p<0$. Assume that the \wt/s $\Lan$ obey condition \(Lass). Let $\Imu$ and,
therefore, $\ka\ne 1$. Let $\ell_1\lc\ell_k$ and $\ell'_1\lc\ell'_k$ be
nonnegative integers \st/ $\ell_1\lsym+\ell_k=\ell$ and
$\ell'_1\lsym+\ell'_k=\ell$. Let ${w_i\in\Fo^i[\)\ell_i\)]}$ and
${W_i\in\Fq^i[\)\ell'_i\)]}$, $i=1\lc k$. Let $w=w_1\lsym\*w_k$ and
$W=W_1\lsym*W_k$. Then the \hint/ $I(W,w)$ has the following \as/s as $z\to\8$
in $\AA^\nb$ so that at any moment assumption \(assum) holds\/{\rm:}\par\nt}
$$
\align
I(W,w)\,=\;{\ell\>!\over\ell_1!\ldots\ell_k!}\,
\prod_{1\le i<j\le k}\!\bigl((z_{n_i}-z_{n_j})/p\bigr)
^{2(\ell_i\!\!\!\sum_{m\in\Gms_j}\!\!\!\La_m\>+\,
\ell_j\!\!\!\sum_{m\in\Gms_i}\!\!\!\La_m\>-\,\ell_i\ell_j)/p} \;\x &
\\
\nn-6>
\x\;\bigl(\)\tprod_{i=1}^k\dl_{\ell_i\ell'_i}\>I(W_i,w_i)+o(1)\bigr)\! &\,.
\endalign
$$
Here ${\Gms_i=\lb\)n_{i-1}+1\lc n_i\)\rb}$, $|\arg\bigl((z_l-z_m)/p\bigr)|<\pi$
and $\dl_{lm}$ is the Kronecker symbol.
\endpro
\Rem
In a separate paper we will describe \azo/s and \asol/s for the \qKZe/, if the
parameter $\ka$ of the \eq/ equals $1$. In this case the \azo/s are essentially
the same as the \azo/s for the \KZ/ \difl/ \eq/ and the \asol/s are similar,
\cf. \Cite{V4}. If $\ka=1$, then the \azo/s of the \qKZe/ are labelled by
\perm/s in $\S^n$ and suitable planar trees $T$. For every \perm/ $\tau$ and
a tree $T$ we define an \azo/ and a basis $\Bg_{T,\tau}$ in the space of
singular vectors $\Vls\!$, a basis of ``iterated singular vectors'', see
\Cite{V4}. For every \perm/ $\tau$ and a tree $T$ we also define a basis
$\Wo_{T,\tau}$ in the \sthgf/. This basis defines a basis of \sol/s to
the \qKZe/ with values in $\Vls\!$. This basis gives an \asol/ to the \qKZe/
in the \azo/ corresponding the \perm/ and the tree. Moreover, the leading terms
of \as/s in this case are proportional to elements of the basis $\Bg_{T,\tau}$
and the coefficients of proportionality are products of powers of linear \fn/s
like in \(fzz) with no exponential factors unlike in the case of $\ka\ne 1$.
\par
If $\ka=1$ then the \qKZo/s $K_m(z)$ have the following \as/s
$$
K_m(z)\,=\,1+o(1)_m\,,\qqq\mn\,,
$$
as all differences $z_i-z_j$ tend to infinity. In every \azo/ the leading terms
of $o(1)_m$ form a system of commuting operators, see (2.2.3) in \Cite{V4}.
The vectors of the basis $\Bg_{T,\tau}$ form an eigenbasis of those commuting
operators.
\par
As an illustrating example consider the \eq/ $f(z+p)=(1+a/z)f(z)$.
The \eq/ has a \sol/ $\Gmpp{z+a}\big/\Gmpb{z/p}$ with \as/s $(z/p)^{a/p}$
as $z$ tends to infinity.
\enddemo

\Sect[2e]{Quasiclassical \as/s}
Consider a system of \deq/s
$$
\Psi(\zmn)\,=\,\A\"m(\zn;h)\,\Psi(\zn)\,,\qqq\mn\,,
$$
depending on a parameter $h$ and assume that
$$
\A\"m(z_1/h\lc z_n/h;h)\,=\,1+h\)\B\"m(\zn)+o(h)
\Tag{AhB}
$$
as $h\to 0$. Introduce new coordinates $y_m=hz_m$, $\mn$, and a new \fn/
$$
\Pti(\yn)\,=\,\Psi(y_1/h\lc y_n/h)\,.
$$
Then the system of \deq/s takes the form
$$
\Pti(y_1\lc y_m+hp\lc y_n)\,=\,\bigl(1+h\)\B\"m(\yn)+o(h)\bigr)\>\Pti(\yn)\,,
$$
$\mn$, and turns into a system of \difl/ \eq/s
$$
p\>\dd{y_m}\Pti(\yn)\,=\,\B\"m(\yn)\>\Pti(\yn)\,,\qqq\mn\,,
$$
as $h$ tends to zero. We call this system of \difl/ \eq/s the \em{\qcl/ \as/s}
of the initial system of \deq/s.
\Par
Consider the \qKZe/ with values in $\Vl$ and parameter $\ka=e^{h\eta}$ where
$\eta$ is a given number and $h$ is an additional parameter. Then the \qKZe/
has property \(AhB) and its \qcl/ \as/s is the \KZ/ \difl/ \eq/
$$
p\>\dd{y_m}\Pti(\yn)\,=\,\eta\)H_m\>\Pti(\yn)\,+
\sum_{\tsize{l=1\atop l\ne m}}^n\,{\Om_{lm}\over y_m-y_l}\;\Pti(\yn)\,,
$$
$\mn$, where $\Om_{lm}=2\La_l\La_m-2H_lH_m-E_lF_m-F_lE_m$.
\par
In the previous sections we constructed \sol/s to the \qKZe/. The \sol/s were
labelled by elements of a suitable subspace of a tensor product of \Umod/s.
We show that these \sol/s have \qcl/ \as/s and turn into the \hgeom/ \sol/s
to the \KZ/ \difl/ \eq/ which are described in \Cite{SV1}. To show this fact
we study \qcl/ \as/s of the \hpair/.
\Par
Let $h$ be a real positive number. Assume that $\Im\eta\ge 0$.
We connect the parameter $\mu$ in the \phf/ \(Phi) with the parameter $\eta$
by an \eq/ $\mu=h\eta$.
\par
The case $\Im\eta<0$ can be treated similarly. The parameters $\mu$ and $\eta$
have to be connected by an \eq/ $\mu=2\pii-h\eta$, if $\Im\eta<0$.
\Par
The \as/s \(Stir) of the \phf/ of a \prim/ gives the following \as/s
for the \phf/ \(Phi) as $h\to +0$:
$$
\gather
\Phi(u/h,y/h)\,=\,h^{\ell(\ell-1-\Lmn)/p}\,\Pht(u,y)\>\ono\,,
\Tag{Phih}
\\
\lline{where}
\endgather
$$
\vvn-2>
\ifMag
$$
\align
\kern2em
\Pht( &\uell,\yn)\,={}
\Tagg{Phicl}
\\
{}&=\,\exp\bigl(\eta\}\tsum_{a=1}^\ell\}u_a/p\bigr)
\pron\,\pral\>\bigl((u_a-y_m)/p\bigr)^{2\La_m/p}\!\!\!
\prab\bigl((u_a-u_b/p\bigr)^{-2/p}.
\kern-2em
\endalign
$$
\else
$$
\Rline{\Pht(\uell,\yn)\,=\,\exp\bigl(\eta\}\tsum_{a=1}^\ell\}u_a/p\bigr)
\pron\,\pral\>\bigl((u_a-y_m)/p\bigr)^{2\La_m/p}\!\!\!
\prab\bigl((u_a-u_b)/p\bigr)^{-2/p}.}
\Tagg{Phicl}
$$
\fi
Here we fix a branch of the \fn/ $(x/p)^\al$ by $|\arg(x/p)|<\pi$.
\Par
Consider a domain $\YY$ given by
$$
\YY\,=\,\lb\)y\in\Cn\vert\Im y_1\lsym<\Im y_n\)\rb\,.
\Tag{YYtau}
$$
For every $y\in\YY$ and each $\mn$ we consider an imaginary interval
$$
\gather
\Ud_m\,=\,\lb\)x\in\C\vert\Re u=0\,,\ \,
\Im y_{m-1}\le\Im x\le\Im y_m\)\rb\,,\qqq y_0=-i\8\,,
\\
\nn6>
\Text{and a chain}
\nn-6>
\Ub_m\,=\,
\sum_{l=1}^m\exp\bigl(4\pii\!\tsum_{1\le k<l}\!\La_k/p\bigr)\>\Ud_l\,.
\endgather
$$
\vvv-.3>
For any $\lg\in\Zln$ we define a chain $\UUb_\lg$ in the imaginary subspace
in $\Cl$ by
$$
\UUb_\lg\,=\,\underbrace{\Ub_1\lx\Ub_1}_{\lg_1}\lx
\underbrace{\Ub_n\lx\Ub_n}_{\lg_n}\,.
$$
For any $\lg\in\Zln$ we also define a \raf/ $\wti_\lg(u,y)$ by
$$
\wti_\lg(\uell,\yn)\,=\sum_{\si\in\Sl}\,\pron\>(\lg_m!)\1\>
\prod_{a\in\Gm_{\Rph lm}}(u_{\si_a}-y_m)\1
\Tag{wticl}
$$
where $\Gm_m=\lb 1+\lg^{m-1}\,\lc\lg^m\rb$, $\mn$.
\Th{qcl}
Let $p<0$. Let $\Re\La_m<0$ and let $\Re y_m=0$ for all $\mn$. Let $\mu=h\eta$,
$\Im\eta>0$. Then for any $\lg,\mg\in\Zln$ the \hint/ $I(W_\lg,w_\mg)$ has
the following \as/s as $h\to+0$ and $y\in\YY${\rm:}
$$
\align
I(W_\lg,w_\mg)\, &{}=\,(-2i)^\ell\,\ell\)!\,h^{\ell(\ell-1-\Lmn)/p}\,
\pron\,\prod_{s=1}^{\lg_m}\,{\sin(\pi/p)\over\sin(\pi s/p)}\ \x
\\
\nn4>
&\ifMag\kern-2.5em\fi
{}\>\x\,
\exp\bigl(\pii\tsun\La_m(\lg^{m-1}+\lg^m-2\ell)/p\bigr)\,
\int_{\UUb_\lg}\Pht(u,y)\>\wti_\mg(u,y)\>d^\ell u\,\ono\,.
\endalign
$$
\endpro
\Rem
Recall that the \hint/ $I(W_\lg,w_\mg)$ is defined by \(IWw), the \fn/s
$W_\lg$ and $w_\mg$ are given by formulae \(Wlga) and \(wlga), \resp/,
and we replace in these formulae $\zn$ by $y_1/h\lc y_n/h$.
\enddemo
For any $\lg\in\Zlm$ consider a domain $\UU_\lg$ in the imaginary subspace
in $\Cl$ defined by
$$
\align
\kern1.5em
\UU_\lg\,=\,\lb\,& u\in\Cl\,\vert\,\Re u_a=0\,,\quad\aell\,,
\Tagg{UUlg}
\\
\nn3>
& \Im y_m\le\Im u_{1+\lg^{m-1}}\lsym\le\Im u_{\lg^m}\le\Im y_{m+1}\,,\quad
\mn-1\)\rb\,.
\kern-1.5em
\endalign
$$
\Th{qcl0}
Let $p<0$. Let $\Re\La_m<0$ and let $\Re y_m=0$ for all $\mn$. Let $\mu=h\eta$,
$\Im\eta=0$. Then for any $\lg\in\Zlm$ and any $\mg\in\Zln$ the \hint/
$I(\Wo_\lg,w_\mg)$ has the following \as/s as $h\to+0$ and $y\in\YY${\rm:}
$$
\align
I(\Wo_\lg,w_\mg)\, &{}=\,(2i)^\ell\,\ell\)!\,h^{\ell(\ell-1-\Lmn)/p}\,\x
\\
\nn4>
&\ifMag\kern-2.5em\fi
{}\>\x\,
\exp\bigl(2\pii\tsun\La_m(\ell-\lg^{m-1})/p\bigr)\,
\int_{\UU_\lg}\Pht(u,y)\>\wti_\mg(u,y)\>d^\ell u\,\ono\,.
\endalign
$$
\endpro
\Rem
Recall that the \hint/ $I(\Wo_\lg,w_\mg)$ is defined by \(IWw), the \fn/s
$\Wo_\lg$ and $w_\mg$ are given by formulae \(Wlgo) and \(wlga), \resp/,
and we replace in those formulae $\zn$ by $y_1/h\lc y_n/h$.
\enddemo
Theorems~\[qcl] and \[qcl0] essentially follow from formulae \(Wlga), \(Wlgo)
and \(Phih).
\Cj{conj}
The claims of Theorems~\[qcl] and \[qcl0] remain valid for any $\Lan$ which
obey condition \(Lass) if other assumptions of the theorems hold and
the integrals in \rhs/s of \(qcl), \(qcl0) are defined by \anco/.
\endpro
\Rem
If $\eta=0$, that is $\ka=1$, then the limiting \phf/ \(Phicl) has no
exponential factor and is a product of powers of linear \fn/s. In particular,
if the numbers $\La_m/p$ and $2/p$ are all \rat/, then the limiting integral
is an integral of an algebraic \fn/. From this point of view our initial
\hint/s are a deformation of periods of algebraic \difl/ forms,
and the subject of our study is a \p-deformation of algebraic geometry.
\enddemo

\Sect[3]{The \onedim/ case}
In this section we consider in details the \onedim/ case $\ell=1$.
So we consider the affine projection $\pi:\Con\}\to\Cn$ and a discrete
\rat/ \stype/ \loc/ on $\Con$ and study its de~Rham complex. Our main goal
of doing this is methodological. Since this case is technically simpler
than the general case, the ideas of the proofs become more clear and visual.
The case $\ell=1$ can be viewed as a \p-deformation of the following example.
\par
Let $\zn$ be \pd/ points in $\C$. Let $\Fun$ be the space of \raf/s in $t$
which are regular in $\C\setminus\lb\zn\rb$. Consider the holomorphic
de~Rham complex $\Omb$ on $\C\setminus\lb\zn\rb$ with coefficients
in $\Fun$ associated with the \difl/ ${\nabla=\d+\om\wedge{\cdot}}$,
$\om=\eta dt+\sun\la_m\om_m$, where ${\om_m=dt/(t-z_m)}$.
\Th{eta<>0}
Let ${\eta\ne 0}$. Then for generic $\lan$ the forms $\om_1\lc\om_n$ form
a basis in $\HOne$.
\endproclaim
For $\eta=0$ the \difl/ of $1$ gives a relation in $\HOne$
$$
\sun\la_m\om_m\,\sim\,0\,.
\Tag{nabla}
$$
\Th{eta=0}
Let ${\eta=0}$. Then for generic $\lan$ the forms $\om_1\lc\om_n$ span $\HOne$.
Moreover, relation \(nabla) is the only independent relation between them.
\endproclaim
Let ${\zn\in i\)\R}$, $\Im z_1\lsym<\Im z_n$, $z_0=-i\8$, $z_{n+1}=+i\8$.
Consider the following intervals:
$$
I_k\,=\,\lb\)t\in\C\vert\Re t=0\,,\ \,\Im z_k\le\Im t\le\Im z_{k+1}\)\rb\,,\qqq
\kon\,.
$$
Set
$$
I_k(\om)\,=\,\int_{I_k}\exp(\eta t)\>\tpron(t-z_m)^{\la_m}\>\om
$$
(the integral must be appropriately \regud/). Here we assume that
$0\le\arg(t-z_m)<2\pi$, thus fixing a branch of the integrand.
The intervals $I_k$ become linear \fn/als on the space of \difl/ forms.
For a \fn/ $f$ we have
$$
I_k(\nabla f)\,=\,0\,,\qqq\kon\,.
$$
This means that the linear \fn/als on \difl/ forms defined by intervals $I_k$
can be considered as elements of the space $\Hone$ of linear \fn/als
on $\HOne$.
\Th{cycles}
Let $\lan$ be generic. Then
\vvn.1>
\atem For any $\eta$, $\Im\eta>0$, the intervals $I_1\lc I_n$ form a basis
in $\Hone$.
\bitem For any $\eta$, $\Im\eta<0$, the intervals $I_0\lc I_{n-1}$ form
a basis in $\Hone$.
\endpro
\Th{cycles0}
Let $\eta=0$. Let $\lan$ be generic. Then the intervals $I_1\lc I_{n-1}$ form
a basis in $\Hone$.
\endpro
\Rem
Theorems~\[cycles] and \[cycles0] follow from elementary topological
considerations. Theorem~\[cycles0] can be also deduced from the following
formula \Cite{V1}\>:
\ifMag
$$
\align
&\det\Bigl[\int_{z_k}^{z_{k+1}}\!{\,\la_l\over t-z_l}\,
\tpron(t-z_m)^{\la_m}\>dt\>\Bigr]_{k,l=1}^{n-1}\;=
\Tag{Varch}
\\
&\hp{\det}=\;\Gm(1+\tlmn)\1\pron\Gm(1+\la_m)\,\prod_{l\ne m}(z_l-z_m)^{\la_m}.
\endalign
$$
\else
$$
\kern-5em
\rline{$\dsize
\det\Bigl[\int_{z_k}^{z_{k+1}}\!{\,\la_l\over t-z_l}\,
\tpron(t-z_m)^{\la_m}\>dt\>\Bigr]_{k,l=1}^{n-1}=
\,\Gm(1+\tlmn)\1\pron\Gm(1+\la_m)\,\prod_{l\ne m}(z_l-z_m)^{\la_m}.$}
\kern-5em
\Tag{Varch}
$$
\fi
\enddemo

\subsect{One-dimensional discrete cohomologies}
Consider the affine projection $\pi:\Con\to\Cn$ and a discrete \rat/ \stype/
\loc/ on $\Con$. In this case the \conn/s are equal to
$$
\gather
\phi_1(t,z)\,=\,\ka\pron\,{t-z_m+\La_m\over t-z_m-\La_m}\;,
\\
\nn6>
\phi_{m+1}(t,z)\,=\;{t-z_m-\La_m-p\over t-z_m+\La_m-p}\;,
\endgather
$$
$\mn$, and the \phf/ takes the form
$$
\Phi(t)=\exp(\mu t/p)\pron
\ {\Gm((t-z_m+\La_m)/p)\over\Gm((t-z_m-\La_m)/p)}\;.
\Tag{Phi1}
$$
The \fn/al space $\Fun$ is the space of \raf/s in $t$ and $\zn$ with at most
simple poles at the following hyperplanes
$$
t=z_m-\La_m+(s+1)p\,, \qqq t=z_m+\La_m-sp\,,
$$
$\mn$, $s\in\Zp$. The \rhgf/ $\Fo\sub\Fun$ is the subspace consisting of \fn/s
of the form
$$
P(t,\zn)\,\pron\,{1\over t-z_m-\La_m\!}
$$
where $P$ is a \pol/ of degree less than $n$ in the \var/ $t$.
The discriminant $\BB\sub\Cn$ is the union of the hyperplanes
$$
z_l-z_m+\La_l+\La_m\,=\,p\)s\,,\qqq s\in\Z\,,
$$
$l,\mn$, $l\ne m$, in the base space $\Cn$.
\Par
To simplify notations in this section we write $w_m(t,z)$ instead of
$w_{\eg(m)}(t,z)$. Recall that
$$
w_m(t,\zn)\,=\,{1\over t-z_m-\La_m\!}\ \prlm\,{t-z_l+\La_l\over t-z_l-\La_l}\,,
\qqq\mn\,.
\Tag{wma}
$$
\Lm{wba}
\back{{(\rm\cf. \(wbasis)\))}}
For any $z\in\BB$ the \fn/s $w_1\lc w_n$ restricted to a fiber over $z$
form a basis in the \rhgf/ $\Fo(z)$ of the fiber.
\endpro
\Pf.
Consider \fn/s
$$
g_m(t,z)\,=\,t^{m-1}\,\pron\,{1\over t-z_m-\La_m\!}\;,\qqq\mn\,.
$$
Their restrictions to a fiber over $z$ form a basis of the space $\Fo(z)$.
Define a matrix $M(z)$ by
$$
w_l(t,z)=\sun M_{lm}(z)\>g_m(t,z)\,,\qqq\lcn\,.
$$
The lemma follows from the formula
$$
\det M\,=\plmn\!(z_l-\La_l-z_m-\La_m)\,.
\Tag{detM}
$$
The last formula is similar to the Vandermonde determinant formula.
\epf
The \cosub/ $\Rc(z)$ is \onedim/ and is spanned by $\sun\La_m w_m$.
Relation \(wDw) has the form
\vv->
$$
D(z)\){\cdot}1\,=\,2\sun\La_m w_m\>dt\,,
\Tag{D1}
$$
where $D(z)$ is the \difl/ of the de~Rham complex of a fiber over $z$.
Consider the de~Rham complex of the fiber,
$$
0\,\to\,\Om^0(z)\,\to\,\Om^1(z)\,\to 0\,.
$$
Let $\H(z)\sub H^1(z)$ be the image of the \rhgf/ of the fiber.
\Th{kanon1}
Let $\ell=1$. Let $\ka\ne 1$. Assume that $p<0$ and $2\La_m\nin p\)\Z$ for
any $\mn$. Let $z\in\BB$. Then $\dim\H(z)=n$, that is $\H(z)\simeq\Fo(z)$.
\endpro
\Th{kaon1}
Let $\ell=1$. Let $\ka=1$. Assume that $p<0$ and $2\La_m\nin p\)\Z$
for any $\mn$. Let $z\in\BB$. If $\tLmn\nin p\)\Znn$, then $\dim\H(z)=n-1$,
that is $\H(z)\simeq\Fo(z)/\Rc(z)$.
\endpro
Theorems~\[kanon1] and \[kaon1] can be proved by rather straightforward
calculations. Nevertheless, we will give further another proof which can
be naturally extended to the general case.
\Rem
Assume that the weights $\Lan$ are \st/ $2\La_m\nin p\)\Zp$ for any $\mn$.
Let $z\in\BB$. Then it is easy to check the following.
\vsk.2>
\atem If $\ka\ne 1$, then we have $\H(z)=H^1(z)$ and $\dim\H(z)=n$.
\bitem If $\ka=1$ and $2\!\tsun\!\La_m\nin p\)\Znn$, then also $\H(z)=H^1(z)$,
but $\dim\H(z)=n-1$.
\par\nt
Otherwise, we have $\dim H^1(z)/\H(z)=1$ and $\dim\H(z)$ can be $n-2$ or $n-1$.
\enddemo

\subsect{One-dimensional discrete homologies}
The \thgf/ $\Fq$ is the space of \fn/s of the form
$$
\gather
P(\xi,\ztn)\,\pron\,
{\Expp{z_m-t}\over\Sinp{t-z_m-\La_m}}
\\
\Text{where}
\xi=\exp(2\pii t/p)\,,\qqq \zt_m=\exp(2\pii z_m/p)\,,
\endgather
$$
and $P$ is a \pol/ of degree less than $n$ in the \var/ $\xi$.
\par
We write $W_m(t,z)$ instead of $W_{\eg(m)}(t,z)$ and $\Wo_m(t,z)$
instead of $\Wo_{\eg(m)}(t,z)$. Recall that
$$
\gather
\ifMag\kern-6em\rline\bgroup$\dsize\else\kern.9em\fi
W_m(t,\zn)\,=\,{\Expp{z_m-t}\over\Sinp{t-z_m-\La_m}}
\ \prlm\,{\Sinp{t-z_m+\La_m}\over\Sinp{t-z_m-\La_m}}\,,
\ifMag$\egroup\kern-6em\else\kern-.9em\fi
\Tag{Wma}
\\
\Text{$\mn$, and}
\nn4>
\ifMag\kern-6em\rline\bgroup$\dsize\else\kern1em\fi
\Wo_m\,=\,W_m\exp(-\pii\La_m/p)-W_{m+1}\exp(\pii\La_{m+1}/p)
\,,\qqq \mn-1\,.
\ifMag$\egroup\kern-6em\else\kern-1em\fi
\Tag{Wmo}
\endgather
$$
\Lm{Wba}
\back{{(\rm\cf. \(Wbasis)\))}}
For any $z\in\BB$ the \fn/s $W_1\lc W_n$ restricted to a fiber over $z$
form a basis in the \thgf/ $\Fq(z)$ of the fiber.
\endpro
\Pf.
Consider \fn/s
$$
G_m(t,z)\,=\,\exp\bigl(2\pii\)(m-1)\)t/p\bigr)\,
\pron\,2{\exp(-\pii\)t/p)\over\Sinp{t-z_m-\La_m}}\;,
\qqq\mn\,.
$$
The restrictions of these \fn/s to a fiber over $z$ form a basis of
the space $\Fq(z)$. Define a matrix $M^q(z)$ by
$$
W_l(t,z)=\sun M^q_{lm}(z)\>G_m(t,z)\,,\qqq\lcn\,.
$$
The lemma follows from the formula
$$
\det M^q\,=\,(2i)^{n(1-n)/2}\>\exp\bigl(\pii\tsun z_m/p\bigr)
\plmn\!\Sinp{z_l-\La_l-z_m-\La_m}\,,
$$
(\cf. \(detM)\)).
\epf
\Lm{Wbao}
\back{{(\rm\cf. \(Wbasiso)\))}}
For any $z\in\BB$ the \fn/s $\Wo_1\lc\Wo_{n-1}$ restricted to a fiber
over $z$ form a basis in the \sthgf/ $\Fqs(z)$ of the fiber.
\endpro
\nt
The proof is similar to the proof of Lemma~\[Wba].
\Par
Let $\II$ be the imaginary axis in the space $\C$ with coordinate $t$
oriented from $-i\8$ to $+i\8$. Recall that the \hint/ $I(W,w)$ for
\fn/s $w\in\Fo(z)$, $W\in\Fq(z)$ is defined as the \anco/ of the integral
$$
I(W,w)\,=\,\int_{\II\,}\Phi(t)\>w(t)\>W(t)\>dt
\Tag{IWw1}
$$
\wrt/ $\Lan$ and $\zn$, starting from large real negative $\Lan$ and
imaginary $\zn$. The \anco/ can be written as an integral over a deformed
imaginary space
$$
I(W,w)\,=\,\int_{\IIt\,}\Phi(t)\>w(t)\>W(t)\>dt\,.
\Tag{IWwt1}
$$
The deformation of the imaginary space is not unique. Below we describe
an example of the deformed imaginary axis $\IIt$ which is involved in
the integral \(IWwt1).
\par
The deformed imaginary axis $\IIt$ is a sum of three terms:
$$
\IIt\,=\,\Is+\CCp+\CCm\,,
\Tag{ICC}
$$
which are defined below. First we assume that all the points
$$
z_m\pm(\La_m+sp)\,,\qqq\mn\,,\qquad s\in\Zp\,,
\Tag{points}
$$
are not imaginary. In this case we set $\Is=\II$.
To define the terms $\CCpm\!$ consider the following sets:
$$
\NN2>
\align
\ZZp\, &{}=\,\lb\)z_m+\La_m+ps\vert
\Re(z_m+\La_m+ps)>0\,,\ \,\mn\,,\ \,s\in\Zp\)\rb\,,
\\
\ZZm\, &{}=\,\lb\)z_m-\La_m-ps\vert
\Re(z_m-\La_m-ps)<0\,,\ \,\mn\,,\ \,s\in\Zp\)\rb\,,
\\
\ZZps\, &{}=\,\lb\)z_m-\La_m-ps\vert
\Re(z_m-\La_m-ps)>0\,,\ \,\mn\,,\ \,s\in\Zp\)\rb\,,
\\
\ZZms\, &{}=\,\lb\)z_m+\La_m+ps\vert
\Re(z_m+\La_m+ps)<0\,,\ \,\mn\,,\ \,s\in\Zp\)\rb\,.
\endalign
$$
We define $\CCp\!$ to be the sum of small circles with centers at the points
of $\ZZp$ oriented anticlockwise. Similarly, $\CCm\!$ is the union of small
circles with centers at the points of $\ZZm$ oriented clockwise.
We assume that the circles are so small that there are no points of the sets
$\ZZps$, $\ZZms$ inside them and they do not intersect the imaginary axis.
\par
{\bls 1.1\bls
If some of the points \(points) are imaginary, then we take $\Is$ to be
an appropriate deformation of the imaginary axis. Namely,
if $\Re(z_m+\La_m+ps)=0$, then we replace the small interval $\Re t=0$,
$|\Im(t-z_m-\La_m-ps)|\le\eps$, by a small semicircle
$|\>t-z_m-\La_m-ps|=\eps$, $\Re(t-z_m-\La_m-ps)\ge 0$. Similarly,
if $\Re(z_m-\La_m-ps)=0$, then we replace the small interval $\Re t=0$,
$\Im(t-z_m+\La_m+ps)|\le\eps$, by a small semicircle
$|\>t-z_m+\La_m+ps|=\eps$, $\Re(t-z_m+\La_m+ps)\le 0$.
The terms $\CCpm\!$ remain the same.}
\goodbm
\Ex
Let $n=1$. In this case the deformed imaginary axis $\IIt$ looks like
$$
{\Msam
\AST\BLANK\AST\CDOT\AST\CDOT\AST\CDOT\AST\vertar
\DOT\CAST\DOT\CAST\DOT\CAST\DOT\BLANK\DOT}
$$
where asterisks and dots stand for points $z_1+\La_1+ps$ and $z_1-\La_1-ps$,
$s\in\Zp$, \resp/.
\enddemo
\Lm{Ianco1}
Let $\Imu$. Then for any $l,\mn$ the \hint/ $I(W_l,w_m)$ can be analytically
continued as a univalued holomorphic \fn/ of complex \var/s $p$, $\Lan$, $\zn$
to the region
$$
p<0\,,\qqq z\in\BB\,,\qqq 2\La_m\nin p\)\Zn\,,\qquad\mn\,.
$$
\endpro
\Pf.
The only thing to be shown is convergence of the integral in \rhs/ of formula
\(IWwt) for \fn/s $W=W_l$, $w=w_m$. The convergence is clear since
$$
\Phi(t)=t^{\M}\exp(\mu t/p)\ono\,,\qqq t\to\pm i\8\,,
$$
and therefore, under the assumptions of the lemma the integrand decays
exponentially as $t$ goes to infinity.
\epf
\Lm{Ianco01}
Let $\Im\mu=0$. Then for any $\lcn-1$, $\mn$ the \hint/ $I(\Wo_l,w_m)$ can be
analytically continued as a univalued holomorphic \fn/ of complex \var/s $p$,
$\Lan$, $\zn$ to the region
$$
p<0\,,\qqq z\in\BB\,,\qqq 2\La_m\nin p\)\Zn\,,\qquad\mn\,.
$$
\endpro
\nt
The proof is similar to the proof of the previous lemma.
\Par
In what follows we need to consider the \hint/ $I(W,w)$ for \fn/s $w$ from
the \fn/al space $\Fun(z)$ of a fiber. The definition is similar to
the definition of the \hint/ for $w\in\Fo(z)$. Below we describe explicitly
the \anco/ of the \hint/ $I(W,w)$ for any \fn/ $w\in\Fun(z)$ as an integral
over a suitable deformation of the imaginary line.
\par
For any integer $s$ let $\IIt[s]$ be the deformation of the imaginary axis
which is defined similarly to $\IIt$ but the parameters $\Lan$ are replaced
by $\La_1+ps\lc\La_n+ps$, \resp/. In particular, $\IIt[0]=\IIt$.
\par
For a \fn/ $w\in\Fun(z)$ we have
$$
I(W,w)\,=\int_{\IIt[s]}\Phi(t)\>w(t)\>W(t)\>dt
\Tag{IWws}
$$
where the integer $s$ is chosen so that the integrand has no poles at
the points $z_m\pm(\La_m+pr)$ for $r<s$, $r\in\Z$. Under this assumption
\rhs/ of \(IWws) does not depend on $s$.
\Par
\goodbm
Let $D\Fun(z)=\lb\)Dw\vert w\in\Fun(z)\)\rb$.
\Lm{IWwd}
Let either {$\Imu$ and $W\in\Fq(z)$} or {$\mu=0$ and $W\in\Fqs(z)$}.
Assume that $p<0$ and $2\La_m\nin p\)\Z$ for any $\mn$. Let $z\in\BB$. Then
\atem The \hint/ $I(W,w)$ is well defined for any \fn/ $w\in\Fun(z)$.
\bitem The \hint/ $I(W,w)$ equals zero for any \fn/ $w\in D\Fun(z)$.
\endpro
\Pf.
The proof of claim a) is similar to the proof of Lemma~\[Ianco1]. Claim b)
follows from the next observation. Let $\IIt_p[s]$ be the contour obtained from
$\IIt[s]$ by the translation $t\map t+p$. Then for a given \fn/ $w\in\Fun(z)$
and a large negative $s$ the contour $\IIt[s]$ and $\IIt_p[s]$ are homologous
in the complement of the set of poles of the \fn/ $\Phi(t)\>w(t)\>W(t)$.
\epf
\Lm*
Let $\mu=0$. Assume that $p<0$ and $2\La_m\nin p\)\Zn$ for any $\mn$.
Then the \hint/ $I(W,w)$ equals zero for any $w\in\Rc(z)$ and $W\in\Fqs(z)$.
\endpro
\Pf.
The lemma follows from formula \(D1) and Lemma~\[IWwd].
\epf
The \hint/ defines linear \fn/als $I(W,\cdot)$ on the \fn/al space of a fiber.
Lemma~\[IWwd] means that these linear \fn/s can be considered as elements of
the homology group $H_1(z)$, the dual space to the cohomology group of
the de~Rham complex of the discrete \loc/ of the fiber.
\par
Let $W$ be any element of the \thgf/ $\Fq$. Let $W|_z\in\Fq(z)$ be its
restriciton to a fiber. Consider an element $s_W\:(z)=I(W|_z,\cdot)$ of
the homology group $H_1(z)$.
\Th{sWz1}
\back{{(\rm\cf. \(sWz)\))}}
Let $\ell=1$. Let either {$\Imu$ and $W\in\Fq$} or {$\mu=0$ and $W\in\Fqs$}.
Assume that $p<0$ and $2\La_m\nin p\)\Z$ for any $\mn$. Then the section
$s_W\:$ is a periodic section \wrt/ the \GM/.
\endpro
\Pf.
Let the contour $\IIt_m[s]$ be defined similar to $\IIt[s]$ but the parameter
$z_m$ is replaced by $z_m-p$. The statement of the theorem means that for any
\fn/ $w\in\Fo(z)$ and each $\mn$ we have the equality
$$
I(W,w)\,=\int_{\IIt_m[s]}\Phi(t)\>w(t)\>W(t)\>dt\,,
$$
where $s$ is a sufficiently large negative integer. The last equality holds
since the integrand $\Phi(t)\>w(t)\>W(t)$ has no poles separating the contours
$\IIt_m[s]$ and $\IIt[s]$.
\epf
Consider a section $\Psi_W$ of the \trib/ over $\Cn$ with fiber $\Vl$:
$$
\Psi_W(z)\,=\,\sun I(W|_z, w_m|_z)\,\Fvm\,.
$$
\Cr{qKZsol1}
\back{{\rm(\cf. \(qKZsol)\))}}
The section $\Psi_W$ is a \sol/ to the \qKZe/.
\endpro
Our further strategy is as follows. First we show that if $\Imu$, then
the basis of sections $\Psi_{W_m}$, $\mn$, is an \asol/ to the \qKZe/,
(\cf. Theorem~\[asol1]). Using this fact we prove that the \hpair/
${I:\Fq(z)\ox\Fo(z)\to\C}$ is \ndeg/ if $\Imu$ (\cf. Theorem~\[mu<>01]).
Studying the \as/ behaviour of the \hint/ as $\mu$ tends to zero we will show
that for $\mu=0$ the \hpair/ ${\Ii\!\]:\Fqs(z)\ox\Fo(z)/\Rc(z)\to\C}$ is \ndeg/
(\cf. Theorem~\[mu=01]). At the end of the section we will describe the \qcl/
\as/s of the \hint/ for $\ell=1$ (\cf. Theorems~\[qcl1], \[qcl01]).
\Par
For every \perm/ $\tau\in\S^n$, consider the \azo/ in $\Cn\!$
given by
$$
\AA_\tau\,=\,\lb\)z\in\Cn\vert\Re z_{\tau_{1}}\lsym\ll\Re z_{\tau_{n}}\)\rb\,,
$$
and say that $z\to\8$ in $\AA_\tau$ if $\Re(z_{\tau_m}-z_{\tau_{m+1}})\to-\8$
for all $\mn-1$.
\Th{asol1}
\back{{\rm(\cf. \(asol)\))}}
Let $\ell=1$. Let $\Imu$ and, therefore, $\ka\ne 1$. Assume that $p<0$ and
$2\La_m\nin p\)\Z$ for any $\mn$. Then for any \perm/ $\tau\in\S^n$ the basis
$\Psi_{W_m^\tau}$, $\mn$, is an \asol/ in the \azo/ $\AA_\tau$. Namely,
$$
\align
\Psi_{W_m^\tau}(z)\,=\,\Tht_m\>\exp(\mu z_m/p)
\prod_{1\le l<\tau\1_m}\!\bigl((z_{\tau_l}-z_m)/p\bigr)^{2\La_{\tau_l}}\!
\prod_{\tau\1_m<l\le n}\!\bigl((z_m-z_{\tau_l})/p\bigr)^{2\La_{\tau_l}}\;\x &
\\
\nn4>
\x\;\bigl(\Fvm+o(1)\bigr) &\,.
\endalign
$$
as $z\to\8$ in $\AA_\tau$ so that $z\in\BB$ at any moment.
Here $|\arg\bigl((z_k-z_l)/p\bigr)|<\pi$ and $\Tht_m$ is a constant
independent of the \perm/ $\tau$ and given by
$$
\Tht_m\,=\,2i\>(e^\mu-1)^{-2\La_m/p}\>\exp\bigl((\mu+\pii)\La_m/p\bigr)\,
\Gm(2\La_m/p)\,,
$$
where $0\le\arg(e^\mu-1)<2\pi$.
\endpro
\Pf.
To simplify notations we will give a proof only for $\tau=\id$.
A simple but important fact is that for any $W\in\Fq$
$$
\Psi_W(z_1+p\lc z_n+p)\,=\,\ka\>\Psi_W(\zn)\,.
\Tag{z+p}
$$
It allows us to fix freely the real part of one of the coordinates $\zn$.
\par
Consider the \hint/ $I(W_m,w_m)$. The corresponding integrand
$\Phi(t)\>w_m(t)\>W_m(t)$ can be rewritten as follows:
$$
\NN6>
\align
\kern1.5em
\Phi(t)\>w_m(t)\>W_m(t)\, &{}=\,
(-\pi p)\1\>\exp\bigl((\mu-\pii)t/p+\pii z_m/p\bigr)\;\x
\Tagg{PhiwW}
\\
&{}\>\x\,\Gmpp{t-z_m+\La_m}\>\Gmpp{z_m+\La_m-t}\;\x
\\
&{}\>\x\,
\prlm{\Gmpp{z_l+\La_l-t}\over\Gmpp{z_l-\La_l-t}}\>
\prod_{m<l\le n}{\Gmpp{t-z_l+\La_l}\over\Gmpp{t-z_l-\La_l}}\;.
\kern-1.5em
\endalign
$$
This \fn/ has no poles at points $z_l-\La_l-sp$, $s\in\Z$, for $l<m$ and has
no poles at points $z_l+\La_l+sp$, $s\in\Z$, for $l>m$. Moreover, due to
\(z+p), \wlg/, we can assume that $z$ tends to infinity in $\Aid$ so that
$\Re z_l\to-\8$ for $l<m$, $\Re z_m$ remain finite, and $\Re z_l\to+\8$
for $l>m$. Under this assumption the integrand has no poles at the points
$z_l+\La_l+sp$, $s\in\Z$, for $l<m$ in the halfplane $\Re t>0$ and has no poles
at the points $z_l-\La_l-sp$, $s\in\Z$, for $l>m$ in the halfplane $\Re t<0$.
Therefore, we can ``straighten'' the contour and write
$$
I(W_m,w_m)\,=\int_{\IIt_m}\Phi(t)\>w_m(t)\>W_m(t)\>dt
\Tag{IWwm}
$$
where the contour $\IIt_m$ is the contour defined above for \anco/ of
the integral
$$
\int_{\II\,}
\exp\bigl((\mu-\pii)t/p\bigr)\>\Gmpp{t-z_m+\La_m}\>\Gmpp{z_m+\La_m-t}\;dt\,.
$$
\par
The remaining part of the calculation is a standard exercise.
We replace the integrand in the integral \(IWwm) by its asymptotics
as $z\to\8$ in $\Aid$ and obtain
$$
\align
I(W_m,w_m)\,=\,
(-\pi p)\1\>\exp(\mu z_m/p)\prlm\!\bigl((z_l-z_m)/p\bigr)^{2\La_l/p}\!
\prod_{m<l\le n}\!\bigl((z_m-z_l)/p\bigr)^{2\La_l/p}\;\x &
\\
\nn6>
\x\int_{\,\IIt_m}\exp\bigl((\mu-\pii)(t-z_m)/p\bigr)\>
\Gmpp{t-z_m+\La_m}\>\Gmpp{z_m-t+\La_m}\>dt\>\ono &\,.
\endalign
$$
The last integral reduces to the Barnes integral \(Barnes) and is calculated
explicitly. Finally, we have
$$
\align
I(W_m,w_m)\,=\,2i\>(e^\mu-1)^{-2\La_m/p}\>\exp\bigl((\mu+\pii)\La_m/p\bigr)\,
\Gm(2\La_m/p)\>\exp(\mu z_m/p)\;\x
\\
\nn6>
\x\prlm\!\bigl((z_l-z_m)/p\bigr)^{2\La_l/p}\!
\prod_{m<l\le n}\!\bigl((z_m-z_l)/p\bigr)^{2\La_l/p}\,\ono & \,,
\endalign
$$
as $z\to\8$ in $\Aid$. Here $0\le\arg(e^\mu-1)<2\pi$.
\par
\goodbreak
\goodbe
The \hint/ $I(W_m,w_l)$ for $l\ne m$ can be treated similarly to the \hint/
$I(W_m,w_m)$ considered above. The final answer is
$$
I(W_m,w_l)\,=\,I(W_m,w_m)\>o(1)\,,
$$
which completes the proof of Theorem~\[asol1].
\epf
\Th{mu<>01}
\back{{\rm(\cf.\(mu<>0)\))}}
Let $\ell=1$. Let $\Imu$. Assume that $p<0$ and $2\La_m\nin p\)\Zn$ for any
$\mn$. Let $z\in\BB$. Then the \hpair/ $I:\Fq(z)\ox\Fo(z)\to\C$ is \ndeg/.
Moreover,
$$
\NN4>
\align
\DWwj\>=\,(2i)^n\>(e^\mu-1)^{-\M}
\exp\bigl((\mu+\pii)\tsun\La_m/p+\mu\tsun z_m/p\bigr)\;\x &
\\
\x\;\pron\Gm(2\La_m/p)\>
\plmn\,{\Gmpp{z_l+\La_l-z_m+\La_m}\over\Gmpp{z_l-\La_l-z_m-\La_m}}\> &\,.
\endalign
$$
Here $0\le\arg(e^\mu-1)<2\pi$.
\endpro
\Pf.
Denote by $F(z)$ the determinant $\DWwj$ and by $G(z)$ \rhs/ of the formula
above. Since for every $\lcn$ the section $\Psi_{W_l}$ is
a \^{$(\Vox)_1\:$-}valued \sol/ to the \qKZe/, $F(z)$ solves the next system
of \deq/s
$$
F(\zmn)\,=\,\Det1 K_m(\zn)\>F(\zn)\,.
$$
Here $\Det1 K_m(z)$ stands for the determinant of the operator
$K_m(z)$ \(Kmz) acting in the weight subspace $(\Vox)_1\:$. Using formula
\(Rspec) we see that
$$
\Det1 K_m(\zn)\,=\,
\ka\>\prlm{z_m+\La_m-z_l+\La_l+p\over z_m-\La_m-z_l-\La_l+p}\>
\prod_{m<l\le n}{z_m+\La_m-z_l+\La_l\over z_m-\La_m-z_l-\La_l}\;.
$$
Therefore, the ratio $F(z)/G(z)$ is a \p-periodic \fn/ of each
of the \var/s $\zn$:
$$
{F\over G}\)(\zmn)\,=\,{F\over G}\)(\zn)\;.
$$
Theorem~\[asol1] implies that the ratio $F(z)/G(z)$ tends to $1$ as $z$ tends
to infinity in the \azo/ $\Aid$. Hence, this ratio equals $1$ identically,
which completes the proof.
\epf
\Th{mu=01}
Let $\ell=1$. Let $\mu=0$. Assume that $p<0$ and $2\La_m\nin p\)\Z$
for any $\mn$. Let $z\in\BB$. If ${\tLmn\nin p\)\Znn}$, then the \hpair/
${\Ii\!:\Fqs(z)}{}\ox\Fo(z)/\Rc(z)\to\C$ is \ndeg/. Moreover,
$$
\NN4>
\align
\DWwp\>=\,(2i)^{n-1}\,\Gmpb{1+\tM}\1\,\Gm(1+2\La_n/p)\;\x &
\\
\x\;\prmn\Gm(2\La_m/p)\;
\plmn\,{\Gmpp{z_l+\La_l-z_m+\La_m}\over\Gmpp{z_l-\La_l-z_m-\La_m}}\> &\,.
\endalign
$$
\endpro
\Pf.
Since both sides of the formula above are analytic \fn/s of $\Lan$, it
suffices to prove the formula under the assumption
$$
0\,<\,\tM\,<\,1\,.
$$
To prove the theorem we first assume that $\mu\ne 0$ and study the \as/s
of the determinant
\ifMag\else\nl\fi
$\DWwj$ as $\mu$ tends to zero. We will show that
$$
\NN3>
\align
\quad\DWwj\>=\,(ip/\La_n)\>\exp\bigl(\pii\tsun\La_m/p\bigr)\>\mu^{-\M}\;\x &
\Tagg{mu0}
\\
\x\;\Gmpb{1+\tM}\1\,\DWwp\ono &
\endalign
$$
as $\mu\to 0$, $0<\arg\mu<\pi$. Due to formula \(mu<>01) the last formula
will imply the required formula for $\DWwo$.
\par
First we change bases in the \rat/ and \tri/ \hgeom/ spaces of a fiber. We set
$$
\alignat3
W'_m &{}=\Wo_m\,,\qquad &&\mn-1\,,\qquad\qquad & W'_n &{}=W_n\,,
\\
\nn-2>
\Text{and}
\nn-8>
w'_m &{}=w_m\,, &&\mn-1\,, & w'_n &{}=\tsun\La_m w_m\,.
\endalignat
$$
\vvv-.7>
We have
\vvn-.3>
$$
\DWWJ\>=\,\La_n\>\exp\bigl(-\pii\>\tsmm^{n-1}\La_m)\>\DWwj\,.
\Tag{Jj}
$$
As $\mu$ tends to zero, the entries $I(W'_l,w'_m)$, $l,\mn-1$, have finite
limits $I(\Wo_l,w_m)$, \resp/. Similarly, the entries $I(W'_l,w'_n)$, $\lcn-1$,
tend to zero since ${D(z)\){\cdot}1}=2w'_ndt$ at $\mu=0$ and, therefore,
$I(\Wo_l,w_n)=0$ at $\mu=0$. More precisely, we have $I(W'_n,w'_m)=O(\mu)$ as
$\mu\to 0$. The behaviour of the entries $I(W'_n,w'_m)$, $\mn$, is described
in the next lemma.
\Lm{IWnwn}
Let $0<\M<1$. Let $\mu\to 0$, $0<\arg\mu<\pi$. Then
\vv-.5>
$$
I(W_n,w_m)\,=\,2i\>\exp(\pii\La_n/p)\,\mu^{-\M}\,\Gmpb{\tM}\,\ono.
$$
\vvv-.3>
\endpro
\Pf.
As $t\to-i\8$, the integrand of the \hint/ $I(W_n,w_m)$ has
the following \as/s:
\vv-.5>
\ifMag
$$
\align
\Phi(t)\>w_m(t)\>W_n(t)\,=\,
(-2i/p)\>\exp\bigl(\mu\)t/p-\pii\La_n/p-2\pii\tsmm^{n-1}\La_m/p\bigr)\>
(t/p)^{-1+\M}\,\x &
\\
\x\;\ono &\,.
\endalign
$$
\else
$$
\Phi(t)\>w_m(t)\>W_n(t)\,=\,
(-2i/p)\>\exp\bigl(\mu\)t/p-\pii\La_n/p-2\pii\tsmm^{n-1}\La_m/p\bigr)\>
(t/p)^{-1+\M}\ono\,.
$$
\fi
Denote by $F(t)$ \lhs/ of the equality above and by $G(t)$ \rhs/ without
the factor $1+o(1)$.
\par
Let $s$ be a positive number \st/ $s>\max\lb\>|z_1|\>\lc\>|z_n|\>\rb$.
Let $\IIt_s$ be the part of the deformed imaginary axis $\IIt$ in
the halfplane $\Im t>-s$. We have
$$
\gather
I(W_n,w_m)\,=\,\Bigl(\int_{-i\8}^{\;0} - \int_{-is}^{\;0}\,\Bigr)\>G(t)\>dt+
\int_{-i\8}^{-is}\bigl(F(t)-G(t)\bigr)\>dt + \int_{\IIt_s}F(t)\>dt\,.
\endgather
$$
The first integral in \rhs/ above can be calculated explicitly since
$$
\int_{-i\8}^{\;0}\exp(\mu t/p)\>(t/p)^{-1+\M}\>dt/p\,=\,
-\exp\bigl(2\pii\tsun\La_m\bigr)\,\mu^{-\M}\,\Gmpb{\tM}\,,
$$
and the three other integrals have finite limits as $\mu\to 0$.
The lemma is proved.
\epf
\Cr*
Let $0<\M<1$. Let $\mu\to 0$, $0<\arg\mu<\pi$. Then
$$
I(W'_n,w'_n)\,=\,2i\>\exp(\pii\La_n/p)\,\mu^{-\M}\,\Gmpb{1+\tM}\,\ono.
$$
\endpro
Finally, we have
$$
\DWWJ\,=\,\DWwp\>I(W'_n,w'_n)\,\ono\,.
$$
Using formula \(Jj) and Corollary~\[*] we get formula \(mu0).
Theorem~\[mu=01] is proved.
\epf
\Pf of Theorems~\[kanon1] and \[kaon1].
Theorem~\[kanon1] and \[kaon1] follow from Theorem~\[mu<>01] and \[mu=01],
\resp/, and Lemma~\[IWwd].
\epf

\subsect{Quasiclassical \as/s}
Recall that to study the \qcl/ \as/s of the \hint/ we introduced new parameters
$h$ and $\eta=\mu/h$, and new coordinates $u=ht$ and $y_m=hz_m$, $\mn$.
The \qcl/ \as/s of a \hint/ is the \as/s of the integral as $h\to 0$ while
the coordinates $\yn$ and the parameter $\eta$ remain fixed.
\Par
For each $\mn$ we defined an imaginary interval
$$
\gather
\Ud_m\,=\,\lb\)u\in\C\vert\Re u=0\,,\ \,
\Im y_{m-1}\le\Im u\le\Im y_m\)\rb\,,\qqq y_0=-i\8\,,
\\
\nn4>
\Text{a chain}
\nn-4>
\Ub_m\,=\,
\sum_{l=1}^m \exp\bigl(4\pii\!\tsum_{1\le k<l}\!\La_k/p\bigr)\>\Ud_l\,.
\\
\Text{and a \raf/}
\wti_m(u,\yn)\,=\;{1\over u-y_m}\;.
\\
\Text{Set}
\nn-4>
\Pht(u,\yn),=\,\exp(\eta u/p)\pron\>\bigl((u-y_m)/p\bigr)^{2\La_m/p}
\endgather
$$
\vvv-.5>
where $|\arg\bigl((u-y_m)/p\bigr)|<\pi$.
\Th{qcl1}
\back{{\rm(\cf. \(qcl)\))}}
Let $\ell=1$. Let $p<0$. Let $\Re\La_m<0$ and let $\Re y_m=0$ for all $\mn$.
Let $\mu=h\eta$, $\Im\eta>0$. Then for any $l,\mn$ the \hint/ $I(W_l,w_m)$
has the following \as/s as $h\to+0$ and $y\in\YY${\rm:}
\ifMag
$$
\align
& I(W_l,w_m)\;=
\\
\nn2>
&\,=\,-2i\>h^{-2\!\sum_{k=1}^n\La_k/p}\,
\exp\bigl(-\pii\La_l/p-2\pii\!\tsum_{1\le k<l}\!\La_k/p\bigr)\,
\int_{\Ub_l}\Pht(u,y)\>\wti_m(u,y)\>du\,\ono\,,
\endalign
$$
\else
$$
I(W_l,w_m)\,=\,-2i\>h^{-2\!\sum_{k=1}^n\La_k/p}\,
\exp\bigl(-\pii\La_l/p-2\pii\!\tsum_{1\le k<l}\!\La_k/p\bigr)\,
\int_{\Ub_l}\Pht(u,y)\>\wti_m(u,y)\>du\,\ono\,.
$$
\fi
\vvv-.5>
\endpro
\Rem
Recall that the \hint/ $I(W_l,w_m)$ is defined by \(IWw) where $\ell=1$,
the \fn/s $W_l$ and $w_m$ are given by formulae \(Wma) and \(wma), \resp/,
and we replace in these formulae $\zn$ by $y_1/h\lc y_n/h$.
\enddemo
\Th{qcl01}
\back{{\rm(\cf. \(qcl0), \(conj)\))}}
Let $\ell=1$. Let $p<0$. Let $\Re\La_m<0$ and let $\Re y_m=0$ for all $\mn$.
Let $\mu=h\eta$, $\Im\eta=0$. Then for any $\lcn-1$ and any $\mn$ the \hint/
$I(\Wo_l,w_m)$ has the following \as/s as $h\to+0$ and $y\in\YY${\rm:}
$$
I(\Wo_l,w_m)\,=\,2i\>h^{-2\!\sum_{k=1}^n\!\La_k/p}\,
\exp\bigl(2\pii\!\tsum_{1\le k\le l}\!\La_k/p\bigr)\,
\int_{\Ud_{l+1}}\Pht(u,y)\>\wti_m(u,y)\>du\,\ono\,.
$$
\vvv->
\endpro
\Rem
Recall that the \hint/ $I(\Wo_l,w_m)$ is defined by \(IWw), the \fn/s
$\Wo_l$ and $w_m$ are given by formulae \(Wmo) and \(wma), \resp/,
and we replace in these formulae $\zn$ by $y_1/h\lc y_n/h$.
\enddemo
\Rem
The claims of Theorems~\[qcl1] and \[qcl01] remain valid for any $\Lan$ \st/
$\La_m\nin p\)\Zn$ for all $\mn$, if the other assumptions of the theorems
hold and the integrals in \rhs/s of \(qcl1), \(qcl01) are regularized in
the standard way. We omit the proof since it is not essential for our purpose
in this paper.
\enddemo
The idea of the proofs of Theorems~\[qcl1], \[qcl01] is simple.
After a suitable renormalization, the \qcl/ \as/s of the \fn/ $W_l$ is given
by a linear combination of the characteristic \fn/s of the intervals
$\Ud_1\lc\Ud_l$ with the coefficients defined by the chain $\Ub_l$.
Similarly, after a suitable renormalization, the \qcl/ \as/s of the \fn/
$\Wo_l$ is given by the characteristic \fn/ of the interval $\Ud_{l+1}$.
Therefore, modulo renormalization factors the \qcl/ \as/s of the \hint/s
$I(W_l,w_m)$, $I(\Wo_l,w_m)$ are given by integrals of products of powers of
linear \fn/s over the chain $\Ub_l$ or over the interval $\Ud_{l+1}$, \resp/.
\Pf of Theorem~\[qcl1].
To simplify notations we will give a proof only for $l=m$.
Consider the \hint/ $I(W_m,w_m)$. It is given by
$$
I(W_m,w_m)\,=\,\int_{\II\,}\Phi(t)\>w_m(t)\>W_m(t)\>dt\,.
\Tag{IWmwm}
$$
Let $h$ be a positive number. The factors of the integrand above have
the following \qcl/ \as/s as ${h\to+0}$ while the parameter $\eta=\mu/h$,
the \var/ $u=ht$ and the coordinates $y_m=hz_m$, $\mn$, remain fixed:
\vv-.5>
$$
\align
\Phi(t,\zn)\, &
{}=\,\exp(\mu t/p)\pron\>\bigl((t-z_m)/p\bigr)^{2\La_m/p}\>\ono\,,
\\
\nn4>
w_m(t,\zn)\, &{}=\;{1\over t-z_m}\,\ono\,,
\\
\nn12>
W_m(t,\zn)\, &{}=\,2i\>\exp\bigl(2\pii(z_m-t)/p+\pii\La_m/p+
2\pii\!\tsum_{1\le k<m}\!\La_k/p\bigr)\>\ono
\\
\Text{if $\Im z_m<\Im t$ and}
\nn6>
W_m(t,\zn)\, &{}=\,-2i\>\exp\bigl(-\pii\La_m/p
+2\pii\!\tsum_{1\le k<l}\!\La_k/p-2\pii\!\tsum_{l\le k<m}\!\La_k/p\bigr)\>\ono
\endalign
$$
if $\Im z_l<\Im t<\Im z_{l+1}$, $l=0\lc m-1$. Here $z_0=-i\8$.
\par
To compute the \qcl/ \as/s of the \hint/ $I(W_m,w_m)$ we replace
the integrand in \rhs/ of formula \(IWmwm) by its \qcl/ \as/s,
and after simple transformations we obtain that
\ifMag
$$
\align
& I(W_l,w_m)\;=
\\
\nn8>
&\,=\,-2i\>h^{-2\!\sum_{k=1}^n\La_k/p}\,
\exp\bigl(-\pii\La_l/p-2\pii\!\tsum_{1\le k<l}\!\La_k/p\bigr)\,
\int_{\Ub_l}\Pht(u,y)\>\wti_m(u,y)\>du\,\ono\,.
\endalign
$$
\else
$$
I(W_m,w_m)\,=\,-2i\>h^{-2\!\sum_{k=1}^n\La_k/p}\,
\exp\bigl(-\pii\La_m/p-2\pii\!\tsum_{1\le k<m}\!\La_k/p\bigr)\,
\int_{\Ub_m}\Pht(u,y)\>\wti_m(u,y)\>du\,\ono\,.
$$
\fi
This step can be justified in a standard way using the next lemma.
\Lm{||<A}
Let ${\Re\al>0}$. Then there is a constant $A$ \st/ for any real $s$
the following estimates hold:
\ifMag\else\vv->\fi
$$
\gather
|\)(\al^2+s^2)^{1-\al}\>\exp(-\pi\)|s|\>)\>\Gm(\al+is)\>\Gm(\al-is)\)|\,<\,A\,,
\\
\nn6>
|\)(\al^2+s^2)^{-\al}\>\Gm(is+\al)/\Gm(is-\al)\)|\,<\,A\,.
\endgather
$$
\vvv-.5>
\endpro
\Pf.
The required formula follows from the next specialization
of the Stirling formula
$$
\bigl|\,\log\Gm(x)-(x-1/2)\>\log x + x -\log\sqrt{2\pi}\,\bigr|\,
<\,{K\over\Re x}\:,\qqq\Re x>0\,,
$$
where $K$ is some constant \Cite{WW}.
\epf
\nt
Theorem~\[qcl1] is proved.
\epf
The proof of Theorem~\[qcl01] is similar to the proof of Theorem~\[qcl1].

\makePcd
\Sect[4]{The multidimensional case}
This section contains proofs of the statements formulated in
Sections~{\SNo{2}\,--\,\SNo{2e}}. We start from the lemmas which describe
bases in the \rat/ and \thgf/s of a fiber.
\Pf of Lemmas\{~\{\[wbasis], \{\[Wbasis], \{\[Wbasiso].
First we have to show that \fn/s $w_\lg$, $W_\lg$ and $\Wo_\lg$ lie in
the \rat/, \tri/ and \sthgf/s, \resp/. The arguments in all the cases
are similar, so we will consider only the \rat/ case.
\par
It is clear from definiton \(wlga) that the \fn/ $w_\lg(t,z)$ has the form
\ifMag\else\vv-.5>\fi
$$
Q(\tell,\zn)\,\pron\,\pral\,{1\over t_a-z_m-\La_m\!}\;\prab\ {1\over t_a-t_b+1}
$$
where $Q$ is a \pol/ which has degree less than $n+\ell-1$ in each of
the \var/s $\tell$. Furthermore, by construction the \fn/ $w_\lg$ as a \fn/ of
$\tell$ is invariant \wrt/ the action \(act) of the \symg/ $\Sl\!$, which means
that the \pol/ $Q$ is skew\sym/ \wrt/ the \var/s $\tell$. Hence, the \pol/ $Q$
is divisible by $\!\!\dsize\prab(t_a-t_b)$ and the ratio is a \pol/ which is
\sym/ in \var/s $\tell$ and has degree less than $n$ in each of the \var/s
$\tell$; that is the \fn/ $w_\lg$ is in the \rhgf/.
\Par
For any $\lg\in\Zln$ (\cf. \(Zln)\)), let $P_\lg(\uell)$ be the following
\sym/ \pol/
$$
P_\lg(\uell)\,=\;{1\over\lg_1!\ldots\lg_n!}\,\sum_{\si\in\Sl}\,
\pron\,\prod_{a\in\Gm_{\Rph lm}} u_{\si_a}^{m-1}\,.
$$
Here $\Gm_m=\lb 1+\lg^{m-1}\,\lc\lg^m\rb$, $\mn$. Consider the following \fn/s
$$
g_\lg(t,z)\,=\,P_\lg\)(\tell)\,\pron\,\pral\,{1\over t_a-z_m-\La_m\!}\;
\prab\ {t_a-t_b\over t_a-t_b+1}\;,
$$
\ifMag
$$
\align
G_\lg( & t,z)\,=\,P_\lg\)(\xill)\;\x
\\
\nn2>
&{}\x \pron\,\pral\,{\exp(-\pii t_a/p)\over\Sinp{t_a-z_m-\La_m}}\;
\prab\ {\Sinp{t_a-t_b}\over\Sinp{t_a-t_b+1}}
\endalign
$$
\else
$$
G_\lg(t,z)\,=\,
P_\lg\)(\xill)\,\pron\,\pral\,{\exp(-\pii t_a/p)\over\Sinp{t_a-z_m-\La_m}}\;
\prab\ {\Sinp{t_a-t_b}\over\Sinp{t_a-t_b+1}}
$$
\fi
\ifMag
$$
\align
\Go_\lg( & t,z)\,=\,P_\lg\)(\xill)\;\x
\\
\nn2>
&{}\x \pron\,\pral\,{\exp(\pii t_a/p)\over\Sinp{t_a-z_m-\La_m}}\;
\prab\ {\Sinp{t_a-t_b}\over\Sinp{t_a-t_b+1}}
\endalign
$$
\else
$$
\Go_\lg(t,z)\,=\,
P_\lg\)(\xill)\,\pron\,\pral\,{\exp(\pii t_a/p)\over\Sinp{t_a-z_m-\La_m}}\;
\prab\ {\Sinp{t_a-t_b}\over\Sinp{t_a-t_b+1}}
$$
\fi
where $\xi_a=\exp(2\pii t_a/p)$, $\aell$. Restrictions of the \fn/s
$g_\lg(t,z)$, $\lg\in\Zln$, to a fiber over $z$ form a basis of the \rhgf/
of the fiber. Restrictions of the \fn/s $G_\lg(t,z)$, $\lg\in\Zln$,
(resp.\ $\Go_\lg(t,z)$, $\lg\in\Zlm$) to a fiber over $z$ form a basis
of the \tri/ (resp.\ the singular \tri/) \hgf/ of the fiber.
\par
Define matrices $M(z)$, $M^q(z)$ and $\Mo(z)$ by
$$
\alignat2
w_\lg(t,z)\, &{}=\>\sum_{\mg\in\Zln}M_{\lg\mg}(z)\>g_\mg(t,z)\,, &&
\lg\in\Zln,
\\
W_\lg(t,z)\, &{}=\>\sum_{\mg\in\Zln}M^q_{\lg\mg}(z)\>G_\mg(t,z)\,, &&
\lg\in\Zln,
\\
\Wo_\lg(t,z)\, &{}=\>\sum_{\mg\in\Rph{\Zln}{\Zlm}}
\Mo_{\lg\mg}(z)\>\Go_\mg(t,z)\,, & \qqq &\lg\in\Zlm,
\endalignat
$$
\Lm{DetM}
\ifMag\vvn->\else\vvn-2>\fi
$$
\gather
\det M\,=\,\pros\,\plmn\!\bigl(z_l-\La_l-z_m-\La_m+s\bigr)_{\vp|}^
{\tsize{n+\ell-s-2\choose n-1}},
\\
\nn12>
{\align
\det M^q\,= {}&\,(2i)^{n(1-n)/2\cdot\!\!\tsize{n+\ell-1\choose n}}\>
\exp\Bigl(\pii\tsun z_m/p\cdot\!{n+\ell-1\choose n}\Bigr)\;\x
\\
\nn5>
&{}\x\,\pros\,\plmn\!\!\Sinp{z_l-\La_l-z_m-\La_m+s}_{\vp|}^
{\tsize{n+\ell-s-2\choose n-1}},
\endalign}
\\
\ald
\nn12>
\det\Mo\,=\,(2i)^{(1-n)(n-2)/2\cdot\!\!\tsize{n+\ell-2\choose n-1}}\,
\pros\,\plmn\!\!\Sinp{z_l-\La_l-z_m-\La_m+s}_{\vp|}^
{\tsize{n+\ell-s-3\choose n-2}}.
\endgather
$$
\endpro
\Pf.
The first and the second formulae are \eqv/ to Lemmas 5.2 and 2.2 in \Cite{T},
\resp/. The third formula can be reduced to the second one by a suitable
change of \var/s.
\epf
Lemmas~\[wbasis], \[Wbasis], \[Wbasiso] clearly follow from Lemma~\[DetM].
\epf
\Pf of Lemma\{~\{\[wDw].
\Rhs/ of formula \(wDw) can be rewritten as
$$
\align
\kern2em
\pron\,\prod_{s=0}^{\lg_m-1}\>(2\La_m-s)/p\,\sum_{\si\in\Sl}\,\Bigl[\>
\Bigl(\pron\>{t_1-z_m+\La_m\over t_1-z_m-\La_m}\,
\prod_{a=2}^\ell\>{t_1-t_a-1\over t_1-t_a+1}\;-\,1\Bigr)\;\x &
\Tagg*
\\
\nn4>
\x\;\prod_{a\in\Gm'_{\Rph lm}}\Bigl(\,{1\over t_a-z_m-\La_m\!}
\ \prlm\,{t_a-z_l+\La_l\over t_a-z_l-\La_l}\,\Bigr)\>\Bigr]_\si\!\! &
\kern-2em
\endalign
$$
where $\Gm'_m=\lb 2+\lg^{m-1}\,\lc\lg^m,1+\lg^m\rb$, $\mn$. Set
\ifMag
$$
\align
f_m(t,z)\,=\,\Bigl(\,{t_1-z_m+\La_m\over t_1-z_m-\La_m}\,
\prod_{a\in\Gm'_{\Rph lm}}\>{t_1-t_a-1\over t_1-t_a+1}\;-\,1\Bigr)\;\x &
\\
\nn6>
\x\;\prlm\Bigl[\,{t_1-z_l+\La_l\over t_1-z_l-\La_l}\,
\prod_{a\in\Gm'_l}\>{t_1-t_a-1\over t_1-t_a+1}\,\Bigr] &\,,
\endalign
$$
\else
$$
f_m(t,z)\,=\,\Bigl(\,{t_1-z_m+\La_m\over t_1-z_m-\La_m}\,
\prod_{a\in\Gm'_{\Rph lm}}\>{t_1-t_a-1\over t_1-t_a+1}\;-\,1\Bigr)
\prlm\Bigl[\,{t_1-z_l+\La_l\over t_1-z_l-\La_l}\,
\prod_{a\in\Gm'_l}\>{t_1-t_a-1\over t_1-t_a+1}\,\Bigr]\,,
$$
\fi
\vvv-.5>
so that
\vv-.5>
$$
\sun f_m(t,z)\,=\,\pron\>{t_1-z_m+\La_m\over t_1-z_m-\La_m}\,
\prod_{a=2}^\ell\>{t_1-t_a-1\over t_1-t_a+1}\;-\,1\,,
$$
and expression \(*) equals
\ifMag
$$
\align
\sum_{k=1}^n\>\Bigr\lb\] &\pron\,\prod_{s=0}^{\lg_m-1}\>(2\La_m-s)/p\;\x
\\
\nn2>
&{}\x\sum_{\si\in\Sl}\,\Bigl[\>f_k(t,z)\,
\pron\,\prod_{a\in\Gm'_{\Rph lm}}\Bigl(\,{1\over t_a-z_m-\La_m\!}
\ \prlm\,{t_a-z_l+\La_l\over t_a-z_l-\La_l}\,\Bigr)\>\Bigr]_\si\)\Bigr\rb\,.
\endalign
$$
\else
$$
\sum_{k=1}^n\>\Bigr\lb\]\pron\,\prod_{s=0}^{\lg_m-1}\>(2\La_m-s)/p\,
\sum_{\si\in\Sl}\,\Bigl[\>f_k(t,z)\,
\pron\,\prod_{a\in\Gm'_{\Rph lm}}\Bigl(\,{1\over t_a-z_m-\La_m\!}
\ \prlm\,{t_a-z_l+\La_l\over t_a-z_l-\La_l}\,\Bigr)\>\Bigr]_\si\)\Bigr\rb\,.
$$
\fi
Lemma~\[wDw] now follows from the formulae
\ifMag
$$
\align
(& \lg_k+1)\>(2\La_k-\lg_k)\>w_{\lg+\eg(k)}\;=
\\
\nn4>
&{}=\sum_{\si\in\Sl}\,\Bigl[\>f_k(t,z)\,\pron\>{1\over\lg_m!}\>
\prod_{a\in\Gm'_{\Rph lm}}\Bigl(\,{1\over t_a-z_m-\La_m\!}
\ \prlm\,{t_a-z_l+\La_l\over t_a-z_l-\La_l}\,\Bigr)\>\Bigr]_\si
\endalign
$$
\else
$$
(\lg_k+1)\>(2\La_k-\lg_k)\>w_{\lg+\eg(k)}\,=\>
\sum_{\si\in\Sl}\,\Bigl[\>f_k(t,z)\,\pron\>{1\over\lg_m!}\>
\prod_{a\in\Gm'_{\Rph lm}}\Bigl(\,{1\over t_a-z_m-\La_m\!}
\ \prlm\,{t_a-z_l+\La_l\over t_a-z_l-\La_l}\,\Bigr)\>\Bigr]_\si
$$
\fi
(\cf. \(wlga)\)).
\epf
Lemmas~\[T12], \[L21] follow from formulae (A.3), (A.5) in [IK], \resp/,
and the definitions of the \emod/s, by induction \wrt/ the number of factors
of the tensor products.
\par
Theorems~\[FgVax], \[FgVqx] follow from formulae (A.5)\,--\,(A.8) in [Ko]
and Lemmas~\[T12], \[L21], \resp/.
\Par
Lemmas~\[FFo], \[FFq] follow from formulae \(wlga), \(Wlga) for the \rat/
and \twf/s, \resp/, and Lemma~\[DetM].
\par
The claims of Theorems~\[FgF], \[FgFq] that the maps ${\chi\)\o P}$ and
$\chi_q$ are intertwiners can be verified direcly from formulae \(Tiju) and
\(Lijxi), \resp/, though the calculations are cumbersome. These claims also
follow from Theorems~\[FgVax], \[FgVqx], \resp/. The claims that the maps
${\chi\)\o P}$ and $\chi_q$ are \iso/s follow from Lemmas~\[FFo], \[FFq],
\resp/.
\Par
The proofs of Theorems~\[Ianco] and \[Ianco0] are based on the following
simple lemma.
\Lm{list3}
Consider a configuration of hyperplanes in $\Cl$
$$
t_a=z_m\pm\La_m+sp\,,\qqq t_a=t_b\pm 1+sp\,,
$$
$1\le b<a\le\ell$, $\mn$, $s\in\Z$. The dimensions of all edges of
the configuration do not depend on $p$, $\Lan$, $\zn$ provided that
assumptions \(npZ)\,--\,\(assum) hold.
\endpro
\Pf.
The initial configuration of hyperplanes induces a configuration of hyperplanes
in any edge of the initial configuration. The dimensions of all edges of
the initial configuration remain the same if and only if all the induced
configurations do not have nonstandard coinciding hyperplanes.
This is obviously true if assumptions \(npZ)\,--\,\(assum) hold.
\epf
This lemma implies that the topology of the complement of configuration
\(list3) of hyperplanes in $\Cl$ remains the same for all $p$, $\Lan$ and
$\zn$ satisfying conditions \(npZ)\,--\,\(assum).
\Pf of Theorem\{~\{\[Ianco].
The theorem is proved by induction \wrt/ the number of integration \var/s
in the \hint/.
\par
Recall that the \hint/ $I(W_\lg,w_\mg)$ is defined by
$$
I(W_\lg,w_\mg)\,=\int_{\){\II^\ell}\]}\Phi(t)\>w_\mg(t)\>W_\lg(t)\>\dt\,,
$$
if $\Lan$ are large negative. (\cf. \(IWw)\)). We can replace the imaginary
subspace $\II^\ell$ in the last formula by any subspace of the form
$$
\II_x\,=\,\lb\)t\in\Cl\vert\Re t_a=x_a\,,\quad\aell\)\rb\,,
\Tag{IIx}
$$
where $\xell$ are small \pd/ real numbers without changing the integral.
\par
For the \anco/ we move $\Lan$, $\zn$ and $p$ and preserve the integration
contour $\II_x$ as long as it does not touch the hyperplanes of configuration
\(list2). If a hyperplane $\Pi$ of configuration \(list2) goes through
the integration contour $\II_x$, the integration contour should be deformed
to avoid the intersection. Deforming the integration contour we add a tube
over the intersection of $\II_x$ and the hyperplane $\Pi$. The result of
the deformation is the sum $\II_x+\II'_x\x C_\Pi\:$, where $\II'_x\sub\Pi$ is
a suitable subspace of real dimension $(\ell-1)$, and $C_\Pi\:$ is a small
circle around the hyperplane $\Pi$. For example, if $\Pi$ is given by an \eq/
$t_\ell-z_n-\La_n-ps=0$, then $\Pi$ has coordinates $t_1\lc t_{\ell-1}$,
in these coordinates
$$
\gather
\II'_x\,=\,\lb\)t\in\Pi\vert\Re t_a=x_a\,,\quad\aell-1\)\rb\,,
\\
\Text{and the circle $C_\Pi\:$ is given by}
C_\Pi\:\,=\,\lb\)t_\ell\in\C\vert |\>t_\ell-z_n-\La_n-ps|=\rho\)\rb\,,
\endgather
$$
$\rho$ is a small positive number. The \anco/ of the initial \hint/
$I(W_\lg,w_\mg)$ equals the sum of two integrals
$$
\int_{\){\II_x}\]}\Phi(t)\>w_\mg(t)\>W_\lg(t)\>\dt\ \;+
\ \int_{\){\II_x'}\]}\Res_{t\in\Pi}
\bigl(\Phi(t)\>w_\mg(t)\>W_\lg(t)\bigr)\>d^{\ell-1}t\,,
$$
and the second integral is of the same type as the first one but of a smaller
dimension. Therefore, under the \anco/ the passage of a hyperplane of
the configuration through the integration contour results in the appearence
of a new \hint/ with a smaller number of integrations. This reason shows that
the \hint/ can be analytically continued to the region described in
Theorem~\[Ianco].
\par
Now we show the univaluedness of the \hint/ by induction on the number of
integration variables. Denote the domain described in Theorem~\[Ianco] by $U$.
Consider its fundamental group $\pi_1(U, z^*)$. The generators of the
fundamental group can be chosen in a special form. Namely, for any hyperplane
$\Pi$ lying at the boundary of the domain $U$ choose a curve $\al_\Pi\:$ in $U$
from the base point $\zs$ to a generic point $z_\Pi\:$ of the hyperplane $\Pi$
and fix a loop $\gm_\Pi\:$ in $U$ which goes from $\zs$ to $\Pi$ along
$\al_\Pi\:$, then turns around $\Pi$ along a small circle $\bt_\Pi\:$ and
returns back to $\zs$ along the same curve $\alpha_\Pi\:$.
The loops $\gamma_\Pi\:$ generate the fundamental group.
\par
Let us show that the \hint/ $I(W_\lg,w_\mg)$ has the trivial monodromy under
the \anco/ along the curve $\gamma_\Pi\:$. In fact, under the \anco/ from the
base point $\zs$ to the hyperplane $\Pi$ along the curve $\al_\Pi\:$ we create
smaller dimensional integrals each time one of the hyperplanes of singularities
hits the integration contour. Under the \anco/ of the integral along the circle
$\bt_\Pi\:$ the hyperplanes of singularities do not touch the integration
contour if the point $z_\Pi\:$ is generic. Now under the \anco/ along the curve
$\al_\Pi\:$ from $\Pi$ to $\zs$ we create again smaller dimensional integrals
each time one of the hyperplanes of singularities hits the integration
contour. But the corresponding integrals created on the way to $\Pi$ and on the
way from $\Pi$ come with the opposite signs. Moreover, they are equal according
to the induction assumptions. Hence, the monodromy of the integral along
the loop $\gm_\Pi\:$ is trivial. Theorem~\[Ianco] is proved.
\epf
The proof of Theorem~\[Ianco0] is similar to the proof of Theorem~\[Ianco].
\Par
As in the case $\ell=1$ we extend the notion of the \hint/ $I(W,w)$ and
consider the \hint/ for any \fn/ $w$ in the \fn/al space $\Fun(z)$ of a fiber.
Namely, let $w(t,z)\in\Fun(z)$ be a \fn/ of the form
\ifMag
$$
\NN4>
\alignat2
&\Rlap{P(\tell,\zn,\Lan)\;\x} &&
\\
&\ \ {}\x{}\,& \prod_{s=0}^r\,\Bigl[\,
\pron\,\pral\,{1\over(t_a-z_m-\La_m+sp)\>(t_a-z_m+\La_m-(s+1)p)}\;\,\x &
\\
&& \x\prab\ {1\over(t_a-t_b+1+sp)\>(t_a-t_b-1-(s+1)p)}\,\Bigr]\> &
\endalignat
$$
\else
$$
\align
P(\tell,\zn,\Lan)\,\prod_{s=0}^r\,\Bigl[\,
\pron\,\pral\,{1\over(t_a-z_m-\La_m+sp)\>(t_a-z_m+\La_m-(s+1)p)}\;\,\x &
\\
\nn4>
\x\prab\ {1\over(t_a-t_b+1+sp)\>(t_a-t_b-1-(s+1)p)}\,\Bigr]\> &
\endalign
$$
\fi
where $P$ is a \pol/. If the real parts of the \wt/s $\Lan$ are large negative
and $p$ is small, then we define the \hint/ $I(W,w)$ by formula \(IWw).
The \hint/ is well defined if either $\Imu$ and $W\in\Fq(z)$ or $\Im\mu=0$
and $W\in\Fqs(z)$, since the integrand vanishes exponentially at infinity.
For generic $\Lan$, $\zn$ and $p$ we define \resp/ the \hint/ $I(W_\lg,w)$ or
the \hint/ $I(\Wo_\lg,w)$ by \anco/ \wrt/ $\Lan$, $\zn$ and $p$. Similar to the
proof of Theorem~\[Ianco] one can show that these \hint/s can be analytically
continued as holomorphic univalued \fn/s of complex \var/s $p$, $\Lan$, $\zn$
to the region described in Theorem~\[Ianco]. For arbitrary \fn/s $w\in\Fun(z)$,
$W\in\Fq(z)$ we define the \hint/ by linearity.
\Par
Let $D\Fun(z)=\lb\)Dw\vert w\in\Fun(z)\)\rb$.
\Lm{IWwD}
Let $p<0$. Let \(npZ)\;--\;\(assum) hold.
Let either {$\Imu$ and $W\in\Fq(z)$} or {$\mu=0$ and $W\in\Fqs(z)$}. Then
\atem The \hint/ $I(W,w)$ is well defined for any \fn/ $w\in\Fun(z)$.
\bitem The \hint/ $I(W,w)$ equals zero for any \fn/ $w\in D\Fun(z)$.
\endpro
\Pf.
Claim a) holds by the definition of the \hint/ $I(W,w)$. Claim b) is clear,
if the real parts of the \wt/s $\Lan$ are large negative and $p$ is small.
Then the \anco/ of the integral gives claim b) for generic $p$, $\Lan$, $\zn$.
\epf
Lemma~\[IWw=0] follows from Lemmas~\[wDw] and \[IWwD].
\Par
The \hint/ defines linear \fn/als $I(W,\cdot)$ on the \fn/al space of a fiber.
Lemma~\[IWwD] means that these linear \fn/als can be considered as elements of
the top homology group $H_\ell(z)$, the dual space to the top cohomology group
of the de~Rham complex of the discrete \loc/ of the fiber.
\Pf of Theorem\{~\{\[sWz].
The section $s_W\:$ is defined by
$$
s_W\:(z)\,=\,I(W|_z,\cdot)
$$
where $W|_z$ denotes the restriction of the \fn/ $W(t,z)$ to the fiber over
$z$. The theorem is a direct corollary of the periodicity of the \fn/ $W$
\wrt/ each of the \var/s $\zn$:
$$
W(t,\zmn)\,=\,W(t,\zn)\,,\qqq\mn\,,
$$
\cf. the case $\ell=1$ in Section~\SNo{3}.
\epf
Our further strategy is the same as in the case $\ell=1$. First we prove
Theorem~\[AsIWw] which implies Theorem~\[asol]. Using Theorem~\[asol] we prove
that the \hpair/ ${I:\Fq(z)\ox\Fo(z)\to\C}$ is \ndeg/ if $\Imu$
(\cf. Theorem~\[mu<>0]). Studying the \as/ behaviour of the \hint/ as $\mu$
tends to zero we show that for $\mu=0$ the \hpair/
${\Ii\!\]:\Fqs(z)\ox\Fo(z)/\Rc(z)\to\C}$ is \ndeg/ (\cf. Theorem~\[mu=0]).
\par
Theorems~\[kanon] and \[kaon] will follow from Theorems~\[mu<>0] and \[mu=0],
\resp/, and Lemma~\[IWwD].
\Pf of Theorem\{~\{\[AsIWw].
To simplify notations we will give a proof only for the case $k=n$, so that
$n_m=m$, $\mn$. The general case is similar.
\par
Let ${w\"m\'l\in\Fo[z_m;\La_m;l\,]}$ and ${W\"m\'l\in\Fq[z_m;\La_m;l\,]}$
be the following \fn/s:
$$
\gather
w\"m\'l(t_1\lc t_l,z_m)\,=\;{1\over l\)!}\,\sum_{\si\in\S^l}\,\Bigl[\,
\prod_{a=1}^l\,{1\over t_a-z_m-\La_m\!}\,\Bigr]_\si,\qqq
\Tag{wWprim}
\\
\nn8>
W\"m\'l(t_1\lc t_l,z_m)\,=\,
\prod_{s=1}^{l}\,{\sin(\pi/p)\over\sin(\pi s/p)}\,\sum_{\si\in\S^l}\,\LBc\,
\prod_{a=1}^l\,{\Expp{z_m-t_a}\over\Sinp{t_a-z_m-\La_m}}\,\RBc_\si,
\endgather
$$
(\cf. \(wlga), \(Wlga)\)). We have the equalities
$$
w_\lg\,=\,w\"1\'{\lg_1}\lsym\* w\"n\'{\lg_n}\qquad\text{and}\qquad
W_\lg\,=\,W\"1\'{\lg_1}\lsym* W\"n\'{\lg_n}
$$
Therefore, we have to study the \as/s of the \hint/s $I(W_\lg,w_\mg)$.
\par
Consider the \hint/ $I(W_\lg,w_\mg)$. Due to property \(Phisym) all the terms
in the definition \(Wlga) of the \fn/ $W_\lg$ give the same contribution to
the integral. So we can replace the integrand $\Phi(t)\>w_\mg(t)\>W_\lg(t)$
by the following integrand
$$
\align
\ifMag
\aligned
& F(t)\;=
\\
&{}=\,\pi^{-\ell}\,\ell\)!\,w_\mg(t)\>\exp\bigl(\pii\tsun\lg_mz_m/p\bigr)
\!\prab\;{\Gmpp{t_a-t_b-1}\over\Gmpp{t_a-t_b+1}}\;
\pron\,\Bigl[\,\prod_{s=1}^{\lg_m}\,{\sin(\pi/p)\over\sin(\pi s/p)}\ \x
\endaligned
\else
F(t)\,=\,\pi^{-\ell}\,\ell\)!\,w_\mg(t)\>\exp\bigl(\pii\tsun\lg_mz_m/p\bigr)
\!\prab\;{\Gmpp{t_a-t_b-1}\over\Gmpp{t_a-t_b+1}}\;
\pron\,\Bigl[\,\prod_{s=1}^{\lg_m}\,{\sin(\pi/p)\over\sin(\pi s/p)}\ \x
\fi &
\\
\nn5>
\x\,\prod_{a\in\Gm_{\Rph lm}}\exp\bigl((\mu-\pii)\)t_a/p\bigr)\>
\Gmpp{t_a-z_m+\La_m}\>\Gmpb{1-(t_a-z_m-\La_m)/p}\;\x &
\\
\nn2>
\Llap{\x\,\prlm{\Gmpb{1-(t_a-z_l-\La_l)/p}\over\Gmpb{1-(t_a-z_l+\La_l)/p}}\>
\prod_{m<l\le n}{\Gmpp{t_a-z_l+\La_l}\over\Gmpp{t_a-z_l-\La_l}}\,\Bigr]\,,\!}&
\endalign
$$
where $\Gm_m=\lb 1+\lg^{m-1}\,\lc\lg^m\rb$, $\mn$.
\par
Assume that the real parts of the \wt/s $\Lan$ are negative. If all $\zn$
are imaginary, then we have
$$
I(W_\lg,w_\mg)\,=\,\int_{\II^\ell}\>F(t)\>\dt\,.
$$
The \anco/ of $I(W_\lg,w_\mg)$ to the region $\Re z_1\lsym<\Re z_n$ is given by
$$
\gather
I(W_\lg,w_\mg)\,=\int_{\II^{\lg_1}_1\lx\II^{\lg_n}_n}\!\!F(t)\>\dt
\\
\nn-3>
\Text{where}
\nn-3>
\II^{\lg_m}_m\,=\,\lb\)(t_{1+\lg^{m-1}}\lc t_{\lg^m})\in\C^{\,\lg_m}\vert
\Re t_a=\Re z_m\,,\ \,\lg^{m-1}<a\le\lg^m\)\rb
\endgather
$$
since the integrand has no poles at the hyperplanes $t_a=z_l-\La_l-sp$,
$s\in\Z$, for $\lg^{m-1}<a\le\lg^m$, $m>l$, has no poles at the hyperplanes
$t_a=z_l+\La_l+sp$, $s\in\Z$, for $\lg^{m-1}<a\le\lg^m$, $m<l$, and has
no poles at the hyperplanes $t_a=t_b+1+sp$, $s\in\Z$, for $a>b$.
\par
Let $z\to\8$ in $\AA^\nb$ so that $\Re(z_m-z_{m+1})\to-\8$ for all $\mn-1$.
Consider the case $\lg=\mg$. Transform the \hint/ $I(W_\lg,w_\lg)$ as above
and replace the integrand by its asymptotics as $z\to\8$ in $\AA^\nb$. Since
$$
w_\lg(\tell)\,=\>\pron w\"m\'{\lg_m}(t_{\lg^{m-1}+1}\lc t^{\lg_m})\,+\,o(1)
$$
as $z\to\8$ in $\AA^\nb$ and $t\in\II^{\lg_1}_1\lx\II^{\lg_n}_n$,
we obtain that
\ifMag
$$
\alignat2
I(W_\lg,w_\lg)\,&{}=\,\pi^{-\ell}\,\ell\)!\>\exp\bigl(\pii\tsun\lg_m z_m\bigr)
\plmn\!\bigl((z_l-z_m)/p\bigr)^{2(\lg_l\)\La_m+\lg_m\)\La_l-\lg_l\)\lg_m)/p}
\;\x &&
\\
& {}\>\x\pron\,\Bigl[\;\prod_{s=1}^{\lg_m}\,{\sin(\pi/p)\over\sin(\pi s/p)}
\;\int_{\II^{\lg_m}_m} w\"m\'{\lg_m}(t_{\lg^{m-1}+1}\lc t_{\lg_m})\>
\prod_{a\in\Gm_{\Rph lm}}\Bigl(\exp\bigl((\mu-\pii)\)t_a/p\bigr)\;\x &&
\\
\nn6>
&& \Llap{\x\;\Gmpp{t_a-z_m+\La_m}\>\Gmpb{1-(t_a-z_m-\La_m)/p}
\>\prod_{\tsize{b<a\atop b\in\Gm_{\Rph lm}}}
{\Gmpp{t_b-t_a-1}\over\Gmpp{t_b-t_a+1}}\;\Bigr)\>d^{\lg_m}t\>\Bigr]\;\x} &
\\
&& \Llap{\x\;\ono\,,} &
\endalignat
$$
\else
$$
\align
I(W_\lg,w_\lg)\, &{}=\,\pi^{-\ell}\,\ell\)!\>\exp\bigl(\pii\tsun\lg_m z_m\bigr)
\plmn\!\bigl((z_l-z_m)/p\bigr)^{2(\lg_l\)\La_m+\lg_m\)\La_l-\lg_l\)\lg_m)/p}
\;\x
\\
&\>{}\x\pron\,\Bigl[\;\prod_{s=1}^{\lg_m}\,{\sin(\pi/p)\over\sin(\pi s/p)}
\;\int_{\II^{\lg_m}_m} w\"m\'{\lg_m}(t_{\lg^{m-1}+1}\lc t_{\lg_m})\>
\prod_{a\in\Gm_{\Rph lm}}\Bigl(\exp\bigl((\mu-\pii)\)t_a/p\bigr)\;\x
\\
\nn6>
&\>{}\x\>\Gmpp{t_a-z_m+\La_m}\>\Gmpb{1-(t_a-z_m-\La_m)/p}\>
\prod_{\tsize{b<a\atop b\in\Gm_{\Rph lm}}}
{\Gmpp{t_b-t_a-1}\over\Gmpp{t_b-t_a+1}}\;\Bigr)\>d^{\lg_m}t\>\Bigr]\,\ono\,,
\endalign
$$
\fi
as $z\to\8$ in $\AA^\nb$. Here $|\arg\bigl((z_l-z_m)/p\bigr)|<\pi$.
Due to \(Phisym) the integrals over $\II^{\lg_m}_m$ are the \hint/s
$I(W\"m\'{\lg_m},w\"m\'{\lg_m})$ up to simple factors.
Hence, we finally obtain that
\ifMag
$$
\align
I(W_\lg,w_\lg)\,=\;{\ell\)!\over\lg_1!\ldots\lg_n!}\,
\plmn\!\bigl((z_l-z_m)/p\bigr)^{2(\lg_l\)\La_m+\lg_m\)\La_l-\lg_l\)\lg_m)/p}
\;\x &
\\
\x\;\bigl(\)\tpron I(W\"m\'{\lg_m},w\"m\'{\lg_m})+o(1)\bigr) &\,.
\endalign
$$
\else
$$
I(W_\lg,w_\lg)\,=\;{\ell\)!\over\lg_1!\ldots\lg_n!}
\plmn\!\bigl((z_l-z_m)/p\bigr)^{2(\lg_l\)\La_m+\lg_m\)\La_l-\lg_l\)\lg_m)/p}
\>\bigl(\)\tpron I(W\"m\'{\lg_m},w\"m\'{\lg_m})+o(1)\bigr)\,.
$$
\fi
The \hint/ $I(W_\lg,w_\mg)$ for $\lg\ne\mg$ can be treated similarly
to the \hint/ $I(W_\lg,w_\lg)$ considered above. The final answer is
$$
I(W_\lg,w_\mg)\,=\,I(W_\lg,w_\lg)\>o(1)\,,
$$
which completes the proof if the real parts of the \wt/s $\Lan$ are negative.
\Par
For general $\Lan$ the proof is similar. The \anco/ of $I(W_\lg,w_\mg)$
to the region $\Re z_1\lsym\ll\Re z_n$ is given by
\vv-.7>
$$
I(W_\lg,w_\mg)\,=\int_{\IIt^{\lg_1}_1\lx\IIt^{\lg_n}_n}\kern-.8em F(t)\>\dt
$$
\vvv-.7>
where $\IIt^{\lg_m}_m$ is the respective deformation of $\II^{\lg_m}_m$. On
every contour $\IIt^{\lg_m}_m$ the quantities $\Re(t_a-z_m)$ remain bounded as
$z\to\8$ in $\AA^\nb$ for all $a$ \st/ $\lg^{m-1}<a\le\lg^m$, and the rest of
the proof remains the same as before.
\par
Theorem~\[AsIWw] is proved.
\epf
Further in the proofs we will make use of the following identities:
\vvn-.5>
$$
\gather
\sum_{\si\in\S^l}\,\prod_{1\le j<k\le l}
{y_{\si_k}-\bt y_{\si_j}\over y_{\si_k}-y_{\si_j}}\;=\,
\prod_{s=1}^l\,{1-\Rph{\bt}{\bt^s}\over 1-\bt}\;,
\Tag{ident}
\\
\ald
\nn6>
\sum_{\si\in\S^l}\,\Bigl(\,\prod_{k=1}^l\,{1\over y_{\si_k}-\bt y_{\si_{k-1}}}
\prod_{1\le j<k\le l}{y_{\si_k}-\bt y_{\si_j}\over y_{\si_k}-y_{\si_j}}\,
\Bigr)\,=\,\prod_{k=1}^l\,{1\over y_k-\bt y_0}\;,\Rlap{\qqq\si_0=0\,.}
\Tag{identi}
\endgather
$$
\Pf of formula \{\(BarneS).
Consider the integral in \lhs/ of \(BarneS) as a \fn/ of $u$ and denote it
by $F(u)$. We will show that $F(u)$ satisfies a \difl/ \eq/
$$
(u+u\1)\,{d\over du}\,F(u)\,=\,
\ell\)\bigl(2a+(\ell-1)x\bigr)\>(1-u^{-2})\>F(u)\,.
\Tag{u-eq}
$$
\Rhs/ of formula \(BarneS) solves the same \difl/ \eq/. Therefore both sides
are proportional. The proportionality coefficient equals $1$ since,
as it is shown below, both sides have the same \as/s as $u\to +0$.
\par
Denote by $f(u;\sell)$ the integrand of integral \(BarneS):
$$
\gather
f(u;\sell)\,=\,\prod_{k=1}^\ell u^{2s_k}\>\Gm(a+s_k)\>\Gm(a-s_k)
\prod_{\tsize{j,k=1\atop j\ne k}}^\ell{\Gm(s_k-s_j+x)\over\Gm(s_k-s_j)}
\\
\nn-4>
\Text{so that}
\nn-4>
F(u)\,=\,\int_{-i\8}^{i\8}\!\cdots\!\int_{-i\8}^{i\8}\,f(u;\sell)\>d^\ell s\,.
\Tag{Fufu}
\endgather
$$
Differentiating the integral \wrt/ $u$ and using the identity
$$
\gather
\ell\)\bigl(a+(\ell-1)x/2\bigr)\,+\)\tsum_{k=1}^\ell s_k\,=\>
\sum_{k=1}^\ell\,(a+s_k)
\prod_{\tsize{j=1\atop j\ne k}}^\ell{s_k-s_j+x\over s_k-s_j}
\\
\Text{we obtain that}
\nn-4>
\bigl(\)\ell\)(2a+(\ell-1)x)+u\>{d\over du}\>\bigr)\>F(u)\,=\,
2\>\sum_{k=1}^\ell F_k(u)
\\
\Text{where}
\nn-4>
F_k(u)\,=\,\int_{-i\8}^{i\8}\!\cdots\!\int_{-i\8}^{i\8}\,
(a+s_k)\prod_{\tsize{j=1\atop j\ne k}}^\ell{s_k-s_j+x\over s_k-s_j}\,
f(u;\sell)\>d^\ell s\,.
\endgather
$$
The space
${\lb\)(\sell)\in\Cl\vert\Re s_k=1\,,\ \,\Re s_j=0\,,\ \,j\ne k\)\rb}$
is homologous to the imaginary space in the complement of the poles of the
integrand for $F_k(u)$. Therefore, changing the \^{$k$-th} integration \var/
$s_k\to s_k-1$ and using the \fn/al \eq/ for the gamma-\fn/ we obtain that
$$
\gather
F_k(u)\,=\,u^{-2}\int_{-i\8}^{i\8}\!\cdots\!\int_{-i\8}^{i\8}\,
(a-s_k)\prod_{\tsize{j=1\atop j\ne k}}^\ell{s_j-s_k+x\over s_j-s_k}\,
f(u;\sell)\>d^\ell s\,.
\\
\nn-4>
\Text{Now the identity}
\nn-4>
\ell\)\bigl(a+(\ell-1)x/2\bigr)\,-\)\tsum_{k=1}^\ell s_k\,=\>
\sum_{k=1}^\ell\,(a-s_k)
\prod_{\tsize{j=1\atop j\ne k}}^\ell{s_j-s_k+x\over s_j-s_k}
\\
\nn-4>
\Text{implies that}
\nn-4>
2\>\sum_{k=1}^\ell F_k(u)
\,=\,u^{-2}\>\bigl(\)\ell\)(2a+(\ell-1)x)-u\>{d\over du}\>\bigr)\>F(u)\,,
\endgather
$$
which complete the proof of equality \(u-eq).
\Par
To compute the \as/s of the integral \(BarneS) as $u\to +0$ we first suitably
transform the integrand. Taking identity \(identi) for $l=\ell$,
$\bt=\exp(2\pii x)$ and $y_k=\exp(2\pii s_k)$, $k=0\lc\ell$, we obtain that
\vv-.3>
$$
\align
1\, &{}=\,\exp\bigl(\pii\>\ell\)(\ell-1)\)x/2\bigr)\,
\prod_{k=1}^\ell\Sinpi{s_k-a}\;\x
\\
\nn3>
&\>{}\x\,\sum_{\si\in\Sl}\,\Bigl(\,\prod_{k=1}^\ell\,
{\exp\bigl(\pii\)(a-s_{\si_{k-1}}-x)\bigr)\>\over
\Sinpi{s_{\si_k}-s_{\si_{k-1}}-x}} \prod_{1\le j<k\le\ell}\!
{\Sinpi{s_{\si_k}-s_{\si_j}-x}\over\Sinpi{s_{\si_k}-s_{\si_j}}}\,\Bigr)
\endalign
$$
where $\si_0=0$, $s_0=a-x$. Substitute \rhs/ of the identity above into the
integral \(Fufu). Since $f(u;\sell)$ is a \sym/ \fn/ of the \var/s $\sell$ and
the imaginary space is invariant under \perm/s of $\sell$, we can keep in
the integral only one term of the sum, multiplying the result by ${\ell\)!}$.
Taking the term corresponding to the identity \perm/ we obtain that
$$
\gather
F(u)\,=\int_{\II^\ell_0} g(u;\sell)\>d^\ell s
\Tag{Fnew}
\\
\nn-6>
\Text{where}
\nn-2>
{\align
& g(u;\sell)\,=\,
\pi^\ell\,\ell\)!\,\exp\bigl(\pii\)\ell\)(\ell-1)\)x/2\bigr)\;\x
\\
\nn4>
&\>\x\;\prod_{k=1}^\ell\,\Bigl(u^{2s_k}\;{\Gm(a+s_k)\over\Gm(1+s_k-a)}\;
{\exp\bigl(\pii(a-s_{k-1}-x)\bigr)\over\Sinpi{s_{k-1}-s_k+x}}\,
\prod_{j=1}^{k-1}\,{(s_k-s_j)\>\Gm(s_k-s_j+x)\over\Gm(1+s_k-s_j-x)}\>\Bigr)\,.
\endalign}
\\
\nn6>
\Text{Here we use a notation}
\nn4>
\II^\ell_y\,=\,\lb\)(\sell)\in\Cl\vert\Re s_k=y\,,\ \,k=1\lc\ell\)\rb\,.
\endgather
$$
To compute the proportionality coefficient which we are interested in,
it suffices to study the \as/s of $F(u)$ as $u\to+0$ only for small real $x$
and real $a$ because both sides of formula \(BarneS) are analytic \fn/s of $x$
and $a$. Moreover, we can assume that $u$ is real. To find the \as/s of $F(u)$
we deform continuously the integration contour in integral \(Fnew).
Namely, we replace $\II^\ell_0$ by $\II^\ell_y$ and move $y$ from zero
to the positive direction.
\par
Since $\Re a>0$ and $\Re x>0$, there are no obstacles for
the deformation of the integration contour until $y$ becomes equal to $a$.
At this moment the integration contour touches the singularity hyperplane of
the integrand given by $s_1=a$. The next intersection of $\II^\ell_y$ with
a singularity hyperplane of the integrand appears at $y=a+1$ with
the hyperplane $s_1=a+1$. Therefore, we have
$$
F(u)\,=\,-2\pii\!\int_{\II^{\ell-1}_a} \Res_{s_1=a}g(u;\sell)\>d^{\ell-1}s\ +
\int_{\II^\ell_{a+\dl}}\! g(u;\sell)\>d^\ell s\,,
$$
where $0<\dl<1$ and $\II^{\ell-1}_y=\lb\)(s_2\lc s_\ell)\in\Cll\vert
\Re s_k=y\,,\ \,k=2\lc\ell\)\rb$.
We can estimate the second term above as follows,
$$
\Bigl|\,\int_{\II^\ell_{a+\dl}}\! g(u;\sell)\>d^\ell s\>\Bigr|\,\le\,
u^{a+\dl}\!\int_{\II^\ell_{a+\dl}}\! g(1;\sell)\>d^\ell s\,.
$$
In the first term we continue the deformation of the integration contour;
we replace $\II^{\ell-1}_a$ by $\II^{\ell-1}_y$ and move $y$ from $a$
to the positive direction. The first obstacle to the deformation appears
at $y=a+x$; at this moment $\II^{\ell-1}_y$ touches the hyperplane $s_2=a+x$.
Repeating the consideration $\ell$ times we finally obtain the following \as/s
for $F(u)$ as $u\to+0$:
$$
F(u)\,=\,(-2\pii)^\ell\,\Res_{a\>\smash\trileft\>x}\,g(u;\sell)\>\ono\,,
$$
where $\Res_{a\>\smash\trileft\>x}$ means the residue at the point
$s_k=a+(k-1)x$, $k=1\lc\ell$. Calculating the residue explicitly we find that
\vv-.5>
$$
F(u)\,=\,(2\pii)^\ell\>u^{\ell(2a+(\ell-1)x)}\,
\prod_{k=1}^\ell\,{\Gm(1+kx)\over\Gm(1+x)}\,\Gmpb{2a+(k-1)x}\,\ono\,,
$$
as $u\to+0$, which clearly coincides with the \as/s of \rhs/ of
formula \(BarneS). This means that the proportionality coefficient between
$F(u)$ and \rhs/ of formula \(BarneS) equals $1$. Formula \(BarneS) is proved.
\epf
\Pf of Theorem\{~\{\[mu<>0].
The proof is similar to the proof of Theorem~\[mu<>01]. Since both sides
of formula \(mu<>0) are analytic \fn/s of $\Lan$, it suffices to prove
the formula only for real negative $\Lan$.
\par
Denote by $F(z)$ the determinant $\DWwi_{\lg,\mg\in\Zln}$ and by $G(z)$ \rhs/
of formula \(mu<>0). Since for every $\lg\in\Zln$ the section $\Psi_{W_\lg}$
is a \^{$\Vl$-}valued \sol/ to the \qKZe/, $F(z)$ solves the following system
of \deq/s,
$$
F(\zmn)\,=\,\Det\ell K_m(\zn)\>F(\zn)\,.
$$
Here $\Det\ell K_m(z)$ stands for the determinant of the operator $K_m(z)$
\(Kmz) acting in the weight subspace $\Vl$. Using formula \(Rspec) we see that
$$
\align
& \Det\ell K_m(\zn)\,={}
\ifMag\,\ka^{\,{n+\ell-1\choose n}}_{\vp o}\;\x\fi
\\
\nn2>
& \ifMag\ {}\x{}\else\hp{\det}{}\,=\,\ka^{\,{n+\ell-1\choose n}}_{\vp o}\fi
\,\pros\,\Bigl[\prlm{z_m+\La_m-z_l+\La_l-s+p\over z_m-\La_m-z_l-\La_l+s+p}\>
\prod_{m<l\le n}{z_m+\La_m-z_l+\La_l-s\over z_m-\La_m-z_l-\La_l+s}
\,\Bigr]^{\tsize{n+\ell-s-2\choose n-1}}.
\endalign
$$
Therefore, the ratio $F(z)/G(z)$ is a \p-periodic \fn/ of each of the \var/s
$\zn$:
$$
{F\over G}\)(\zmn)\,=\,{F\over G}\)(\zn)\;.
$$
Theorem~\[asol] implies that the ratio $F(z)/G(z)$ tends to $1$ as $z$ tends
to infinity in the \azo/ $\Aid$. Hence, this ratio equals $1$ identically,
which completes the proof.
\epf
\Pf of Theorem\{~\{\[mu=0].
The idea of the proof is the same as in the case $\ell=1$. Since both sides of
formula \(mu=0) are analytic \fn/s of $\Lan$ and $p$, it suffices to prove
the formula for large negative $\Lan$ and large negative $p$. More precisely,
we assume that
$$
0\,<\,\tM\;-\,\ell(\ell-1)/p\,<\,1\,.
\Tag{0<<1}
$$
\par
We will construct certain bases in the \rat/ and \thgf/s of fibers \st/ in
these bases the \hpair/ has triangular \as/s as $\mu\to 0$. Using this fact
we will show that
$$
\align
\kern1em
\DWwi_{\lg,\mg\in\Zln}\, &{}=\,\Xi\;\mu^{-\M\cdot{\tsize{n+\ell-1\choose n}}
+2n/p\cdot\]{\tsize{n+\ell-1\choose n+1}}}\;\x
\kern-1em
\Tagg{asdet}
\\
\nn3>
&\>{}\x\,\prod_{k=1}^\ell
\det\bigl[I(\Wo_\ig,w_\jg)\bigr]_{\ig,\jg\in\Zc^{n-1}_k}\,\ono\,,
\kern-1em
\endalign
$$
where $0<\arg\mu<\pi$, $\jg\in\Zc^{n-1}_k\!$ is identified with
$(\jg,0)\in\Zc^n_k$, and $\Xi$ is a constant given by
\ifMag
\vvn-.7>
$$
\gather
\Rline{\Xi\,=\,(2i)^{\tsize{n+\ell-1\choose n}}\>
\exp\Bigl(\pii\>\Bigl(\,\tsun\La_m/p\cdot\!
{n+\ell-1\choose n}\)-\>n/p\cdot\!{n+\ell-1\choose n+1}\Bigr)\Bigr)\;\x}
\Tagg{Xic}
\\
\nn10>
\Rline{
\aligned
& \x\;\pros\Bigl(\,\Bigl[\,{p\>(\ell-s)\>
\Gmpb{-(s+1)/p}\over(2\La_n-s)\>\Gmpb{-1/p}}
\,\Bigr]^{\tsize{n+\ell-s-2\choose n-1}}\,\x
\\
\nn4>
&\quad\x\prod_{s<r\le\ell}\Gmpb{1+\tM+(s+2-2r)/p}\}
{\vp{\Big[}}^{\tsize{n+r-s-2\choose n-2}} \,\Bigr)\,.\!
\endaligned}
\endgather
$$
\else
$$
\gather
\Xi\,=\,(2i)^{\tsize{n+\ell-1\choose n}}\>
\exp\Bigl(\pii\>\Bigl(\,\tsun\La_m/p\cdot\!
{n+\ell-1\choose n}\)-\>n/p\cdot\!{n+\ell-1\choose n+1}\Bigr)\Bigr)\;\x
\Tag{Xic}
\\
\nn10>
\Rline{{}\>\x\,\pros\Bigl(\,\Bigl[\,
{p\>(\ell-s)\>\Gmpb{-(s+1)/p}\over(2\La_n-s)\>\Gmpb{-1/p}}
\,\Bigr]^{\tsize{n+\ell-s-2\choose n-1}}\!\!\prod_{s<r\le\ell}
\Gmpb{1+\tM+(s+2-2r)/p}\}{\vp{\Big[}}^{\tsize{n+r-s-2\choose n-2}}
\,\Bigr)\,.}
\endgather
$$
\fi
Formulae \(asdet) and \(Xic) imply Theorem~\[mu=0].
\Par
In the proof we use the \smod/ structure in the \rhgF/ which was defined
in Section~\SNo{2b}. We define \fn/s $w_\lg$, $W_\lg$ for an arbitrary vector
$\lg\in\Zp^n$ \resp/ by formulae \(wlga), \(Wlga), where we replace $\ell$ in
\lhs/s by the sum $\lg_1\lsym+\lg_n$. Similarly, we define \fn/s $\Wo_\lg$ for
arbitrary vector $\lg\in\Zp^{n-1}$ by formula \(Wlgo) replacing there $\ell$
in \lhs/ by the sum $\lg_1\lsym+\lg_{n-1}$.
\par
For any vector $\lg=(\lg_1\lc\lg_n)$ set $\lg'=(\lg_1\lc\lg_{n-1},0)$
and $\lge=(\lg_1\lc\lg_{n-1})$.
\Par
The required bases in the \rat/ and \thgf/s of a fiber are given by
\fn/s $w'_\lg$, $\lg\in\Zln$, and $W'_\lg$, $\lg\in\Zln$, \resp/.
The \fn/s $w'_\lg$ are defined by the rule:
$$
w'_\lg\,=\,F^{\lg_n}w_{\lg'}.
$$
where $F$ is the generator of $\gsl$ acting in the \rhgF/ $\Fg$
(\cf. \(HFE)\)). The \fn/s $W'_\lg$ are given by
\ifMag
$$
\align
\kern1.5em
& W'_\lg(\tell)\,=\,{1\over(\ell-\lg_n)!}\,
\prod_{s=1}^{\lg_n}\,{\sin(\pi/p)\over\sin(\pi s/p)}\,
\sum_{\si\in\Sl}\,\LBc\,\Wo_{\lge}(t_1\lc t_{\ell-\lg_n})\;\x{}
\Tagg{Wprime}
\\
\nn6>
& \x\prod_{\Lph{\ell-\lg}{\ell-\lg_n}<a\le\ell}\Bigl(\,
{\Expp{z_n-t_a}\over\Sinp{t_a-z_n-\La_n}}\,
\prm^{n-1}\,{\Sinp{t_a-z_m+\La_m}\over\Sinp{t_a-z_m-\La_m}}\,\Bigr)\>
\RBc_\si\,,
\kern-1.5em
\endalign
$$
\else
$$
\align
\kern1.2em
W'_\lg(\tell)\,&{}=\,{1\over(\ell-\lg_n)!}\,
\prod_{s=1}^{\lg_n}\,{\sin(\pi/p)\over\sin(\pi s/p)}\,
\sum_{\si\in\Sl}\,\LBc\,\Wo_{\lge}(t_1\lc t_{\ell-\lg_n})\;\x{}
\Tagg{Wprime}
\\
\nn6>
\x\kern.8em &\kern-.85em\prod_{\Lph{\ell-\lg}{\ell-\lg_n}<a\le\ell}\Bigl(\,
{\Expp{z_n-t_a}\over\Sinp{t_a-z_n-\La_n}}\,
\prm^{n-1}\,{\Sinp{t_a-z_m+\La_m}\over\Sinp{t_a-z_m-\La_m}}\,\Bigr)\>
\RBc_\si\,,
\kern-1.2em
\endalign
$$
\fi
By Lemmas~\[wbasis], \[Wbasis], for any $z\in\BB$ there are matrices $N(z)$,
$N^q(z)$ \st/
$$
\alignat2
w'_\lg(t,z)\, &{}=\!\sum_{\mg\in\Zln}N_{\lg\mg}(z)\>w_\mg(t,z)\,, &&
\lg\in\Zln,
\\
\nn4>
W'_\lg(t,z)\, &{}=\!\sum_{\mg\in\Zln}N^q_{\lg\mg}(z)\>W_\mg(t,z)\,,\qqq &&
\lg\in\Zln.
\endalignat
$$
\vsk>
\Lm{detN}
\vv-2>
$$
\gather
\det N(z)\,=\,
\pros\>\bigl((s+1)\>(2\La_n-s)\bigr)^{\tsize{n+\ell-s-2\choose n-1}}\,,
\\
\nn4>
\det N^q(z)\,=\,
\exp\Bigl(-\pii\tsmm^{n-1}\!\La_m/p\cdot\!{n+\ell-1\choose n}\)+
\>\pii\)(n-1)/p\cdot\!{n+\ell-1\choose n+1}\Bigr)\,.
\endgather
$$
\endpro
Now we study the \as/s of the \hpair/ as $\mu\to 0$. We will consider the total
family of the \hpair/s ${I:\Fql\ox\Fol\to\C}$, $l=0\lc\ell$, not indicating
explicitly the dependence on $l$. Recall that we assumed $\Lan$ to be
large negative.
\par
The first observation is that for any ${w\in\Fol}$, ${W\in\Fqsl}$ the \hint/
$I(W,w)$ has a finite limit as $\mu\to 0$, \cf. Lemma~\[Wwell]. Furthermore,
for any $w\in\Fo[\)l-1\)]$ and $W\in\Fqsl$ we have
$$
I(W,F\]w)\,=\,O(\mu)
\Tag{Omu}
$$
as $\mu\to 0$, because ${F\bigl(\Fo[\)l-1\)](z)\bigr)}$ is the \cosub/
${\Rc[\)l\>](z)\sub\Fol(z)}$, \cf. \(FfDf) and Lemma~\[IWw=0].
The \as/ behaviour of the \hint/ $I(W,w)$ for a general \fn/ $W$ is described
by the following lemma.
\Lm{asIW'w}
Let $\Lan$ be large negative. Assume that condition \(0<<1) holds.
Then for any $\lg\in\Zln$ and any $w\in\Fo$ the \hint/ $I(W'_\lg,w)$ has
the following \as/s as $\mu\to 0$, $0<\arg\mu<\pi${\rm:}
$$
I(W'_\lg,w)\,=\,\Xi'\;\mu^{\lg_n(2\ell-\lg_n-1-2\!\sun\!\La_m)/p}\>
I(\Wo_{\lge},E^{\lg_n}w)\>\ono\,.
$$
Here $E$ is the generator of $\gsl$ acting in the \rhgF/ $\Fg$
{\rm (\cf. \(HFE)\))}, and
\ifMag
$$
\align
\Xi'\, &{}=\;{(2i)^{\lg_n}\>\ell\)!\over(\ell-\lg_n)!}\;
\exp\bigl(\pii\>\lg_n\)(\La_n-(\lg_n-1)/2)/p\bigr)\,
\\
\nn4>
&\>{}\x\,\prod_{s=0}^{\lg_n-1}\,\Bigl[\,{\Gmpb{-(s+1)/p}\over\Gmpb{-1/p}}
\,\Gmpb{\tM+(2\lg_n-2\ell-s)/p}\>\Bigr]\,.
\endalign
$$
\else
$$
\;\Xi'\,=\;{(2i)^{\lg_n}\>\ell\)!\over(\ell-\lg_n)!}\;
\exp\bigl(\pii\)\lg_n\)(\La_n-(\lg_n-1)/2)/p\bigr)\,
\prod_{s=0}^{\lg_n-1}\>\Bigl[\,{\Gmpb{-(s+1)/p}\over\Gmpb{-1/p}}
\,\Gmpb{\tM+(2\lg_n-2\ell-s)/p}\>\Bigr]\,.
$$
\fi
\endpro
To obtain the required formulae \(asdet), \(Xic), we also need the following
lemma.
\Lm{shafl}
Let $\mu=0$. Then for any $k,l\in\Zp$ we have
\ifMag
$$
\align
& \det\bigr[I(\Wo_\lg,E^k\}F^k\]w_\mg]_{\lg,\mg\in\Zc^{n-1}_l}\;=
\\
\nn6>
& \hp{\det}\!\]\x\>\prod_{s=0}^{k-1}\>
\bigl((s+1)\>\bigl(\tLmn-2l-s\bigr)\bigr)_{\vp|}^{\tsize{n+l-2\choose n-2}}\,
\det\bigr[I(\Wo_\lg,w_\mg)]_{\lg,\mg\in\Zc^{n-1}_l}\,.
\endalign
$$
\else
$$
\det\bigr[I(\Wo_\lg,E^k\}F^k\]w_\mg]_{\lg,\mg\in\Zc^{n-1}_l}\,=\,
\prod_{s=0}^{k-1}\>\bigl((s+1)\>\bigl(\tLmn-2l-s\bigr)\bigr)_{\vp|}
^{\tsize{n+l-2\choose n-2}}\,
\det\bigr[I(\Wo_\lg,w_\mg)]_{\lg,\mg\in\Zc^{n-1}_l}\,.
$$
\fi
Here we identify $\mg\in\Zc^{n-1}_l\!$ and $(\mg,0)\in\Zc^n_l$.
\endpro
Now we will complete the proof of Theorem~\[mu=0] assuming that
Lemmas~\[detN], \[asIW'w] and \[shafl] hold, and then we will prove the lemmas.
\Par
Consider a matrix $U$ with the entries
$$
U_{\lg\)\mg}\,=\,I(W'_\lg,w_\mg)\,,\qqq\lg,\mg\in\Zln\,.
$$
Lemmas~\[FgVax], \[asIW'w] and formula \(Omu) imply that the \as/s of
the matrix $U$ as $\mu\to 0$ has a block-triangular structure,
\vv-.5>
$$
\alignat2
& U_{\lg\)\mg}\,=\,O(\mu^{\dl_{\lg_n}})\,, && \for\quad\lg_n\ge\mg_n\,,
\\
\nn-4>
\Text{and}
\nn-4>
& U_{\lg\)\mg}\,=\,O(\mu^{1+\dl_{\lg_n}})\,,\qqq && \for\quad\lg_n<\mg_n\,,
\endalignat
$$
where $\dl_l=l(2\ell-l-1-2\sun\La_m)/p$. Therefore, we see that
\vv-.5>
$$
\det U=O\bigl(\mu^{\;\sum_{l=0}^\ell d_l\dl_l}\bigr)\,,\qqq
d_l={n+\ell-l-2\choose n-2}\,,
$$
as $\mu\to 0$, where $d_0\lc d_\ell$ are the dimensions of the diagonal blocks.
Furthermore, the leading term of the \as/s of $\det U$ is given by the product
of the determinants of the diagonal blocks, which are described
by Lemma~\[shafl]. Finally, formula \(asdet) follows from the simple relation
$$
\DWwi_{\lg,\mg\in\Zln}\,=\;{\det U\over\det N\>\det N^q}
$$
and straightforward calculations.
\epf
\noPcd
\Pf of Lemma~\[detN].
The formula for $\det N$ is a corollary of Lemma~\[FgVax].
\par
To prove the formula for $\det N^q$, consider the points $y\"\lg\in\Cl$
defined below:
$$
y\"\lg_a\,=\,z_m-\La_m+\lg^m-a\,,\qqq \lg^{m-1}<a\le\lg^m\,,\qqq\aell\,.
$$
Recall that $\lg^m=\lg_1\lsym+\lg_m$, $\mn$. Let $L$ and $L'$ be matrices
with entries
$$
L_{\lg\)\mg}=W_\lg(y\"\mg)\,,\qqq L'_{\lg\)\mg}=W'_{\lg}(y\"\mg)\,,\qqq
\lg,\mg\in\Zln\,,
$$
\resp/. The matrices $L$ and $L'$ are triangular \wrt/ the following \lex/
order in $\Zln$: $\lg<\mg$ if $\lg_1<\mg_1$ or $\lg_1=\mg_1$, $\lg_2<\mg_2$
etc. Namely, $L_{\lg\)\mg}=0$ and $L'_{\lg\)\mg}=0$ for $\lg<\mg$.
Since $L'_{\lg\)\lg}=
\exp\bigl(\pii\>\lg_m\bigl((\lg_m-1)/2-\La_m\bigr)/p\bigr)\>L_{\lg\)\lg}$
and $N=L'\>L\1$, the formula for $\det N^q$ is proved.
\epf
\Pf of Lemma~\[shafl].
By Theorem~\[FgVax] the \rhgm/ is isomorphic to $\Vax$. We also have
${I(W,F\]w)=0}$ for any ${w\in\Fo[\)l-1\)]}$ and ${W\in\Fqsl}$. Therefore,
the coefficient of proportionality equals the determinant of the operator
${E^k\}F^k}$ acting in the quotient space $(\Vax)_l\:/F(\Vax)_{l-1}$.
This operator is isomorphic to the operator ${E^k\}F^k}$ acting
in the space of singular vectors $(\Vox)_l\sing\!$.
The last operator is simply the multiplication by
$\dsize\prod_{s=0}^{k-1}\>\bigl((s+1)\>\bigl(\tLmn-2l-s\bigr)\bigr)$
in the space of dimension $\dsize{n+l-2\choose n-2}$. The lemma is proved.
\epf
\Pf of Lemma~\[asIW'w].
First we will give another expression for the \fn/ $W'_\lg$ which will be
more convenient for our purpose.
\par
Taking identities \(ident), \(identi) for $\bt=\exp(-2\pii/p)$ and
$y_k=\exp(-2\pii t_k/p)$, $k=0\lc l$, we transform them \resp/
to the following form
$$
\gather
\sum_{\si\in\S^l}\lbc\)1\)\rbc\vpp\si\;=\,
\prod_{s=1}^l\,{\sin(\pi s/p)\over\sin(\pi/p)}\,
\prod_{1\le j<k\le l}\ {\Sinp{t_j-t_k}\over\Sinp{t_j-t_k+1}}\;,
\\
\nn8>
{\align
& \sum_{\si\in\S^l}\,\LBc\,
\prod_{k=1}^l{\Expp{t_{k-1}-t_k+1}\over\Sinp{t_{k-1}-t_k+1}}\,\RBc_\si\;=
\\
\nn6>
&\ \>=\,\exp\bigl(\pii\>l(l-1)/(2p)\bigr)\,\prod_{k=1}^l
{\Expp{z_n+\La_n-t_k}\over\Sinp{z_n+\La_n-t_k}}\,
\prod_{1\le j<k\le l}\ {\Sinp{t_j-t_k}\over\Sinp{t_j-t_k+1}}\;,
\endalign}
\endgather
$$
where $t_0=z_n+\La_n-1$. Subsequently using the identities above we replace
expression \(Wprime) for the \fn/ $W'_\lg$ by the following expression:
\ifMag
$$
\align
\kern1.7em
& W'_\lg(\tell)\,=\;
\Tagg{Wprimd}
\\
\nn6>
&=\;{(-1)^{\lg_n}\over(\ell-\lg_n)!}\,
\exp\bigl(\pii\>\lg_n\bigl((1-\lg_n)/2-\La_n\bigr)/p)\bigr)\,
\sum_{\si\in\Sl}\,\LBc\,\Wo_{\lge}(t_1\lc t_{\ell-\lg_n})\;\x
\\
\nn6>
&\>\x\prod_{\Lph{\ell-\lg}{\ell-\lg_n}<a\le\ell}\Bigl(\,
{\Expp{t'_{a-1}-t'_a+1}\over\Sinp{t'_{a-1}-t'_a+1}}\,\prm^{n-1}\,
{\Sinp{t_a-z_m+\La_m}\over\Sinp{t_a-z_m-\La_m}}\,\Bigr)\>\RBc_\si\,,
\kern-1.7em
\endalign
$$
\else
$$
\align
\kern1.6em
W'_\lg(\tell)\,&{}=\,{(-1)^{\lg_n}\over(\ell-\lg_n)!}\,
\exp\bigl(\pii\>\lg_n\bigl((1-\lg_n)/2-\La_n\bigr)/p)\bigr)\,
\sum_{\si\in\Sl}\,\LBc\,\Wo_{\lge}(t_1\lc t_{\ell-\lg_n})\;\x
\Tagg{Wprimd}
\\
\nn6>
&{}\>\x\!\prod_{\Lph{\ell-\lg}{\ell-\lg_n}<a\le\ell}\Bigl(\,
{\Expp{t'_{a-1}-t'_a+1}\over\Sinp{t'_{a-1}-t'_a+1}}\,\prm^{n-1}\,
{\Sinp{t_a-z_m+\La_m}\over\Sinp{t_a-z_m-\La_m}}\,\Bigr)\>\RBc_\si\,,
\kern-1.6em
\endalign
$$
\fi
where $t'_a=t_a$ for $\ell-\lg_n<a\le\ell$ and $t'_{\ell-\lg_n}=z_n+\La_n-1$.
\Par
Consider the \hint/ $I(W'_\lg,w)$,
$$
I(W'_\lg,w)\,=\int_{\){\II^\ell}\]}\Phi(t)\>w(t)\>W'_\lg(t)\>\dt\,.
$$
The imaginary space is invariant under \perm/s of the \var/s $\tell$.
By property \(Phisym) of \fn/ $\Phi(t)$, we can keep in the integral only
the term of the sum in \rhs/ of \(Wprimd) that corresponds to the identity
\perm/, multiplying the result by ${\ell\)!}$. Hence,
$$
\gather
I(W'_\lg,w)\,=\int_{\){\II^\ell}\]} F(t)\>\dt\,,
\Tag{I*}
\\
\nn-8>
\Text{where}
\nn-6>
\ifMag
{\align
& F(\tell)\,=\;{(-1)^{\lg_n}\>\ell\)!\over(\ell-\lg_n)!}\ \x
\\
\nn6>
&\x\;\exp\bigl(\pii\>\lg_n\bigl((1-\lg_n)/2-\La_n\bigr)/p)\bigr)\,
\Phi(\tell)\>w(\tell)\>\Wo_{\lge}(t_1\lc t_{\ell-\lg_n})\;\x
\\
\nn6>
& \x\prod_{\Lph{\ell-\lg}{\ell-\lg_n}<a\le\ell}\;
\prm^{n-1}\,{\Sinp{t_a-z_m+\La_m}\over\Sinp{t_a-z_m-\La_m}}
\ \ \prod_{\Lph{\ell-\lg}{\ell-\lg_n}<a\le\ell}
{\Expp{t'_{a-1}-t'_a+1}\over\Sinp{t'_{a-1}-t'_a+1}}\;\dt\,,
\endalign}
\else
\nn4>
{\align
F(\tell)\,=\;{(-1)^{\lg_n}\>\ell\)!\over(\ell-\lg_n)!}\;
\exp\bigl(\pii\>\lg_n\bigl((1-\lg_n)/2-\La_n\bigr)/p)\bigr)\,
\Phi(\tell)\>w(\tell)\>\Wo_{\lge}(t_1\lc t_{\ell-\lg_n})\;\x &
\\
\nn4>
\x\prod_{\Lph{\ell-\lg}{\ell-\lg_n}<a\le\ell}\;
\prm^{n-1}\,{\Sinp{t_a-z_m+\La_m}\over\Sinp{t_a-z_m-\La_m}}
\ \ \prod_{\Lph{\ell-\lg}{\ell-\lg_n}<a\le\ell}
{\Expp{t'_{a-1}-t'_a+1}\over\Sinp{t'_{a-1}-t'_a+1}}\;\dt\, &,
\endalign}
\fi
\endgather
$$
with the same convention for the \var/s $t'_{\ell-\lg_n}\lc t'_\ell$ as
in \(Wprimd).
\Par
Consider the \as/s of the integrand $F(t)$ for large $t$. Namely, assume that
$$
t_a=i\)(x_a-Au_a)\,,\qqq\Im u_a=0\qquad\Im x_a=0\,,\qqq\aell\,,
\Tag{t**}
$$
and $A\to+\8$. From the Stirling formula we find that $F(t)$ decays
exponentially for any $u=(\uell)$ which does not belong to the cone
$$
\lb\)u\in\R^\ell\vert
0=u_1\lsym=u_{\ell-\lg_n}\le u_{\ell-\lg_n+1}\lsym\le u_{\ell}\)\rb\,.
$$
The decay takes place for any $\mu$ including $\mu=0$. On the contrary,
if $u$ belongs to the cone
$$
\lb\)u\in\R^\ell\vert
0=u_1\lsym=u_{\ell-\lg_n}<u_{\ell-\lg_n+1}\lsym< u_{\ell}\)\rb\,,
\Tag{cone}
$$
then the \as/s of $F(t)$ depends essentially on whether $\mu$ equals zero
or not. If $\mu\ne 0$, then $F(t)$ decays exponentially due to the factor
$\exp\bigl(\mu\tsum_{a=1}^\ell t_a/p)$ and integral \(I*) converges.
If $\mu=0$, then $F(t)$ grows like a positive power of $A$
and integral \(I*) diverges.
\Par
So the leading term of the \as/s of the \hint/ $I(W'_\lg,w)$ as $\mu\to 0$
is given by the integral of the \as/s of the integrand $F(t)$ for large $t$
in the cone \(cone) (the justification of this fact is given below).
Explicitly computing the \as/s of $F(t)$ for large $t$, we obtain that
\ifMag
$$
\NN8>
\align
& F(\tell)\,=\;
{(-2i/p)^{\lg_n}\>\ell\)!\over(\ell-\lg_n)!}\;\exp\bigl(\pii\>\lg_n
\bigl((1-\lg_n)/2+2(\ell-\lg_n)+\La_n-2\tsun\La_m\bigr)/p\bigr)\;\x
\\
& {}\x\,\Phi_{[\>\ell-\lg_n]}(t_1\lc t_{\ell-\lg_n})\>
(E^{\lg_n}w)(t_1\lc t_{\ell-\lg_n})\>\Wo_{\lge}(t_1\lc t_{\ell-\lg_n})\;\x
\\
& \x\prod_{\Lph{\ell-\lg}{\ell-\lg_n}<a\le\ell}\!\exp(\mu t_a/p)\>
(t_a/p)^{\M-2(\ell-\lg_n)/p-1}
\prod_{\Lph{\ell-\lg}{\ell-\lg_n}<a,b\le\ell}\!\bigl((t_a-t_b)/p\bigr)^{-2/p}
\,\ono\,,
\endalign
$$
\else
$$
\NN8>
\align
F(\tell)\, &{}=\;
{(-2i/p)^{\lg_n}\>\ell\)!\over(\ell-\lg_n)!}\;\exp\bigl(\pii\>\lg_n
\bigl((1-\lg_n)/2+2(\ell-\lg_n)+\La_n-2\tsun\La_m\bigr)/p\bigr)\;\x
\\
&\>{}\x\,\Phi_{[\>\ell-\lg_n]}(t_1\lc t_{\ell-\lg_n})\>
(E^{\lg_n}w)(t_1\lc t_{\ell-\lg_n})\>\Wo_{\lge}(t_1\lc t_{\ell-\lg_n})\;\x
\\
&\>\x\prod_{\Lph{\ell-\lg}{\ell-\lg_n}<a\le\ell}\!\exp(\mu t_a/p)\>
(t_a/p)^{\M-2(\ell-\lg_n)/p-1}
\prod_{\Lph{\ell-\lg}{\ell-\lg_n}<a,b\le\ell}\!\bigl((t_a-t_b)/p\bigr)^{-2/p}
\,\ono\,,
\endalign
$$
\fi
as $A\to+\8$ and $\tell$ are described by \(t**).
Here the \fn/ $\Phi_{[\)l\>]}(t_1\lc t_l)$ is defined by formula \(Phi)
where $\ell$ is replaced by $l$.
\Par
Denote by $G(t)$ \rhs/ of the above formula without the factor $1+o(1)$.
Thus we have
$$
I(W'\lg,w)\,=\,\int_{\II^{\ell-\lg_n}\_{\vp|}\x\,\II^{\lg_n}_{\sss\top}}
\kern-.5em G(t)\>\dt\;\ono\,,
\Tag*
$$
as $\mu\to 0$, where
$\II^{\lg_n}_{\sss\top}=\lb\)(t_{\ell-\lg_n+1}\lc t_\ell\vert\Re t_a\le 0$,
$\Im t_a=0$, $\ell-\lg_n<a\le\ell\)\rb$.
\vv.2>
The integral \wrt/ the \var/s $t_1\lc t_{\ell-\lg_n}$ clearly gives
$I(\Wo_{\lge},E^{\lg_n}w)$. The integral \wrt/ the \var/s
$t_{\ell-\lg_n+1}\lc t_\ell$ can be calculated explicitly
via the Selberg integral and we obtain formula \(asIW'w).
\Par
The \as/s \(*) can be justified in a standard way. The main idea is the same
as in the \onedim/ case, \cf. Lemma~\[IWnwn]. We explain the details for
the example $n=1$, $\ell=2$. The general case is similar.
\par
To make formulae shorter we change notations and consider
the following integral
$$
\gather
J(\al)\;=\,\int_{\R^2}F(s_1,s_2)\exp\bigl(-\al(s_1+2s_2)\bigr)\,ds_1\>ds_2\,,
\\
\nn8>
\ifMag
{\align
F(s_1,s_2)\,=\,{\Gm(a+is_1)\over\Gm(1-a+is_1)}\,
{\Gm(b+is_2)\over\Gm(-b+is_2)}\,{\Gm(c+is_1+is_2)\over\Gm(1-c+is_1+is_2)}\ \x &
\\
\nn4>
\x\ {\exp\bigl(\pi(s_1+s_2+ia+ib)\bigr)\over
4\sinh\bigl(\pi(s_1+ia)\bigr)\>\sinh\bigl(\pi(s_2+ib)\bigr)} &
\endalign}
\else
F(s_1,s_2)\,=\,{\Gm(a+is_1)\over\Gm(1-a+is_1)}\,
{\Gm(b+is_2)\over\Gm(-b+is_2)}\,{\Gm(c+is_1+is_2)\over\Gm(1-c+is_1+is_2)}
\;{\exp\bigl(\pi(s_1+s_2+ia+ib)\bigr)\over
4\sinh\bigl(\pi(s_1+ia)\bigr)\>\sinh\bigl(\pi(s_2+ib)\bigr)}
\fi
\endgather
$$
\mgood
Our assumptions mean that parameters $a,b,c$ are small positive numbers \st/
$$
0\,<\,a+b+c\,<\,1/2\,.
$$
For $s_1=Au_1$, $s_2=Au_2$ and $A\to+\8$ the \fn/ $F(s_1,s_2)$ has
the following \as/s:
$$
F(s_1,s_2)\,=\,s_1^{2a-1}\>s_2^{2b}\>(s_1+s_2)^{2c-1}\,\ono
$$
if $u_1>0$ and $u_2>0$, and $F(s_1,s_2)$ decays exponentially if either
$u_1\le0$ or $u_2\le 0$. We have to show that
$$
J(\al)\;=\int_{\;\R^2_{\ge 0}}s_1^{2a-1}\>s_2^{2b}\>
(s_1+s_2)^{2c-1}\>\exp\bigl(-\al(s_1+2s_2)\bigr)\,ds_1\>ds_2\;\ono
\ifMag\kern-1.6em\fi
\Tag{asJal}
$$
as $\al\to 0$, $\Re \al>0$,
\Par
Fix a small positive number $\eps$ and decompose $\R^2$ into four parts:
$$
\NN3>
\align
Q_1(\eps)\, &{}=\,\lb\)(u_1,u_2)\in\R^2\vert u_1\ge 0\,,\ \,u_2\ge 0\,,\ \,
u_1\ge\eps u_2\,,\ \, u_2\ge\eps u_1\)\rb\,,
\\
Q_2(\eps)\, &{}=\,\lb\)(u_1,u_2)\in\R^2\vert u_1\ge |u_2|/\eps\)\rb\,,
\\
Q_3(\eps)\, &{}=\,\lb\)(u_1,u_2)\in\R^2\vert u_2\ge |u_1|/\eps\)\rb
\endalign
$$
and $Q_4(\eps)$ is the closure of
$\R^2\setminus\bigl(Q_1(\eps)\cup Q_2(\eps)\cup Q_3(\eps)\bigr)$.
The respective decomposition of the integral $J(\al)$ is
$$
\alignat2
J(\al)\; &{}=\int_{\;Q_1(\eps)}\!\!s_1^{2a-1}\>s_2^{2b}\>
(s_1+s_2)^{2c-1}\>\exp\bigl(-\al(s_1+2s_2)\bigr)\,ds_1\>ds_2\ +{}&&
\\
\nn6>
&\>{}+\int_{\;Q_1(\eps)}\!
\bigl(F(s_1,s_2)-s_1^{2a-1}\>s_2^{2b}\>(s_1+s_2)^{2c-1}
\bigr)\>\exp\bigl(-\al(s_1+2s_2)\bigr)\,ds_1\>ds_2\ +{}&&
\\
&& \Llap{+\ \Bigl(\!\int_{\;Q_2(\eps)} + \int_{\;Q_3(\eps)} +
\int_{\;Q_4(\eps)}\Bigr)\,F(s_1,s_2)\>\exp\bigl(-\al(s_1+2s_2)\bigr)
\,ds_1\>ds_2\,} & .
\endalignat
$$
The first integral equals
$$
\al^{-2(a+b+c)}\!\int_{\;Q_1(\eps)}\!\!s_1^{2a-1}\>s_2^{2b}\>
(s_1+s_2)^{2c-1}\>\exp\bigl(-s_1-2s_2)\bigr)\,ds_1\>ds_2\,,
$$
the second and fifth integrals have finite limits as $\al\to 0$,
the third and fourth integrals can be \resp/ estimated from above by
$$
\gather
A_1\,
\eps^{2b+1}\>\al^{-2(a+b+c)}\,\int_{0\;}^{+\8} s_1^{2(a+b+c)-1}\>
\exp(-s_1)\,ds_1
\\
\nn-2>
\Text{and}
\nn-2>
A_2\,
\eps^{2a}\>\al^{-2(a+b+c)}\,\int_{0\;}^{+\8} s_2^{2(a+b+c)-1}\>\exp(-2s_2)\,ds
\endgather
$$
as $\al\to 0$, where the constants $A_1,A_2$ do not depend on $\eps$.
The estimates can be obtained by means of Lemma~\[||<A]. Therefore,
\ifMag
$$
\align
\Bigl|\,\lim_{\al\to 0}\bigl(\al^{2(a+b+c)}\>J(\al)\bigr)\ -\!
\int_{\;Q_1(\eps)}\!\!s_1^{2a-1}\>s_2^{2b}\>
(s_1+s_2)^{2c-1}\>\exp\bigl(-s_1-2s_2)\bigr)\,ds_1\>ds_2\,\Bigr|\;<{}&
\\
\;
{}<\;\eps^{2b+1}+A_2\,\eps^{2a}\, & .
\endalign
$$
\else
$$
\Bigl|\,\lim_{\al\to 0}\bigl(\al^{2(a+b+c)}\>J(\al)\bigr)\ -\!
\int_{\;Q_1(\eps)}\!\!s_1^{2a-1}\>s_2^{2b}\>
(s_1+s_2)^{2c-1}\>\exp\bigl(-s_1-2s_2)\bigr)\,ds_1\>ds_2\,\Bigr|\;<\;
\eps^{2b+1}+A_2\,\eps^{2a}\,.
$$
\fi
Moving $\eps$ to zero we see that
$$
J(\al)\,=\,\al^{-2(a+b+c)}\int_{\;\R^2_{\ge 0}}s_1^{2a-1}\>s_2^{2b}\>
(s_1+s_2)^{2c-1}\>\exp\bigl(-s_1-2s_2)\bigr)\,ds_1\>ds_2\;\ono\,,
$$
which coincides with \(asJal).
\par
Lemma~\[asIW'w] is proved.
\epf
\makePcd
\Pf of Theorem\{~\{\[qcl0].
Let $h$ be a positive number. Consider the \hint/
$$
I(\Wo_\lg,w_\mg)\,=\int_{\){\II^\ell}\]}\Phi(t,z)\>w_\mg(t,z)\>
\Wo_\lg(t,z)\>\dt\,,
$$
where we substitute into formulae \(Phi), \(wlga), \(Wlgo), defining the \fn/s
$\Phi,w_\mg,\Wo_\lg$, a new parameter $\eta$ and new coordinates $\yn$:
$$
\mu=h\eta\,,\qqq z_m=y_m/h\,,\qquad\mn\,.
$$
We study the \qcl/ \as/s of the \hint/ $I(\Wo_\lg,w_\mg)$ as $h\to +0$
while the parameter $\eta$ and the coordinates $\yn$ remain fixed.
\par
For any $\lg\in\Zlm$ consider a region $\UUo_\lg$ and a domain
$\widehat\UU_\lg$ in the imaginary subspace $\II^\ell$ given by
$$
\alignat2
& \UUo_\lg\,=\,\lb\,u\in\II^\ell\,\vert\,
\Im y_m<\Im u_{1+\lg^{m-1}}\lsym<\Im u_{\lg^m}<\Im y_{m+1}\,, && \quad
\mn-1\)\rb\,,
\\
\nn6>
& \widehat\UU_\lg\,=\,
\lb\,u\in\II^\ell\,\vert\,\Im y_m\le\Im u_a\le\Im y_{m+1}\,,
& \Llap{\lg^{m-1}<a\le\lg^m\,,} & \quad \mn-1\)\rb\,.
\endalignat
$$
Recall that $\lg^m=\lg_1\lsym+\lg_m$.
\Par
Let $\S^\lg\sub\Sl$ be the following subgroup isomorphic to
$\S^{\lg_1}\lx\S^{\lg_n}$:
$$
\S^\lg\,=\,\lb\)\si\in\Sl\vert\lg^{m-1}<\si_a\le\lg^m\quad
\for\ \ \lg^{m-1}<a\le\lg^m\,,\quad\mn\)\rb\,.
$$
Set $\UUa_\lg=\lb\,u\in\Cl\vert(u_{\si_1}\lc u_{\si_\ell})\in\UUo_\lg$
for some $\si\in\S^\lg\)\rb$.
\Par
The imaginary subspace $\II^\ell$ is invariant under \perm/s
of the \var/s $\tell$. Using property \(Phisym) of \fn/ $\Phi(t)$ we see that
$$
\gather
I(\Wo_\lg,w_\mg)\,=\ell\)!\>\int_{\){\II^\ell}\]}\Phi(t,z)\>w_\mg(t,z)\>
\Wb_\lg(t,z)\,\dt\,,
\Tag{IWbw}
\\
\nn-2>
\Text{where}
\nn-2>
\alignedat2
\ifMag
&\>\Wb_\lg(\tell,\zn)\;= &&
\\
\nn6>
&\>{}\x\>\prmn\,\prod_{s=1}^{\lg_m}\,{\sin(\pi/p)\over\sin(\pi s/p)}\,
\Sinp{z_m-\La_m-z_{m+1}-\La_{m+1}+s-1}\;\x &&
\else
\Wb_\lg(\tell,\zn)\, &{}=
\,\prmn\,\prod_{s=1}^{\lg_m}\,{\sin(\pi/p)\over\sin(\pi s/p)}\,
\Sinp{z_m-\La_m-z_{m+1}-\La_{m+1}+s-1}\;\x &&
\fi
\\
\nn6>
&\>{}\x\>\prmn\,\prod_{a\in\Gm_{\Rph lm}}\Bigl(\,
{1\over\Sinp{t_a-z_m-\La_m}\>\Sinp{t_a-z_{m+1}-\La_{m+1}}}\ \x &&
\\
\nn6>
&& \Llap{\x\;\prlm\,{\Sinp{t_a-z_l+\La_l}\over\Sinp{t_a-z_l-\La_l}}\,\Bigr)} &
\endalignedat
\endgather
$$
\ifMag\else\vvv->\fi
and $\Gm_m=\lb\)\lg^{m-1}+1\lc\lg^m\)\rb$, $\mn$.
\par
The factors of the above integrand have the following \qcl/ \as/s as $h\to +0$
while the parameter $\eta$, the coordinates $\yn$ and the \var/s $u_a=ht_a$,
$\aell$, remain fixed:
$$
\alds
\gather
\Phi(u/h,y/h)\,=\,h^{\ell(\ell-1-\Lmn)/p}\,\Pht(u,y)\>\ono\,,
\\
\nn8>
w_\mg(u/h,y/h)\,=\,h^\ell\>\wti(u,h)\>\ono\,,
\\
\nn6>
\ifMag
\aligned
\Wb_\lg(u/h,y/h)\,=\,\prmn\Bigl(
\exp\bigl(\pii\bigl(2\tsun\La_m(\ell-\lg^{m-1})-\lg_m(\lg_m-1)/2\bigr)/p\bigr)
\>\prod_{s=1}^{\lg_m}\,{\sin(\pi/p)\over\sin(\pi s/p)}\,\Bigr)\;\x &
\\
\nn4>
\x\; \ono &
\endaligned
\else
\Wb_\lg(u/h,y/h)\,=\,\prmn\Bigl(
\exp\bigl(\pii\bigl(2\tsun\La_m(\ell-\lg^{m-1})-\lg_m(\lg_m-1)/2\bigr)/p\bigr)
\>\prod_{s=1}^{\lg_m}\,{\sin(\pi/p)\over\sin(\pi s/p)}\,\Bigr)\>\ono
\fi
\endgather
$$
for $u\in\UUo_\lg$ and $\Wb_\lg(u/h,y/h)=o(1)$ for $u\nin\UUa_\lg$.
Here the \fn/s $\Pti(u,y)$ and $\wti_\mg(u,y)$ are given by
formulae \(Phicl) and \(wticl), \resp/.
\par
The \qcl/ \as/s of the integral \(IWbw) is given by the integral
of the \qcl/ \as/s of the integrand, that is
\ifMag
$$
\align
\kern1.5em
I(\Wb_\lg,w_\mg)\,=\,\prmn\Bigl(
\exp\bigl(\pii\bigl(2\tsun\La_m(\ell-\lg^{m-1})-\lg_m(\lg_m-1)/2\bigr)/p\bigr)
\;\x &
\Tagg{asIWbw}
\\
\nn6>
\x\;\prod_{s=1}^{\lg_m}\,{\sin(\pi/p)\over\sin(\pi s/p)}\,\Bigr)\,
\int_{\UUa_\lg}\Pht(u,y)\>\wti(u,y)\>d^\ell u\,\ono\,. \kern-2.9em &
\endalign
$$
\else
$$
\align
\kern1.5em
I(\Wb_\lg,w_\mg)\,=\,\prmn\Bigl(
\exp\bigl(\pii\bigl(2\tsun\La_m(\ell-\lg^{m-1})-\lg_m(\lg_m-1)/2\bigr)/p\bigr)
\>\prod_{s=1}^{\lg_m}\,{\sin(\pi/p)\over\sin(\pi s/p)}\,\Bigr)\ \x &
\Tagg{asIWbw}
\\
\nn6>
\x\;\int_{\UUa_\lg}\Pht(u,y)\>\wti(u,y)\>d^\ell u\,\ono\, &.
\kern-1.5em
\endalign
$$
\fi
\vvv->
Taking into account that
$$
\int_{\UU_\lg}\Pht(u,y)\>\wti(u,y)\>d^\ell u\,=\,
\prmn\Bigl(\exp\bigl(\pii\>\lg_m(1-\lg_m)/(2p)\bigr)\>
\prod_{s=1}^{\lg_m}\,{\sin(\pi/p)\over\sin(\pi s/p)}\,\Bigr)\,
\int_{\UUa_\lg}\Pht(u,y)\>\wti(u,y)\>d^\ell u
$$
where $\UU_\lg$ is given by \(UUlg), we obtain formula \(qcl0).
\Par
Formula \(asIWbw) can be justified in a standard way,
similar to the proof of formula \(asJal).
\par
Theorem~\[qcl0] is proved.
\epf

\ifMag\else\fixedpage\fi
\vsk>

\myRefs
\widest{VW}

\ref\Key A
\by \Aomoto/
\paper $q$-analogue of de~Rham cohomology associated with Jackson integrals, I
\jour Proceedings of Japan Acad.{} \vol 66 {\rm Ser\&A} \yr 1990
\pages 161--164
\moreref \paper II \jour Proceedings of Japan Acad.{} \vol 66 {\rm Ser\&A}
\yr 1990 \pages 240--244
\endref

\ref\Key AK
\by \Aomoto/ and Y\]\&Kato
\paper Gauss decomposition of connection matrices for \sym/
{\sl A\/}-type Jackson integrals
\jour \SMNS/ \vol 1 \yr 1995 \issue 4 \pages 623--666
\endref

\ref\Key CP
\by V\]\&Chari and A\&Pressley
\book A guide to quantum groups
\yr 1994 \publ \CUPa/ 
\endref

\ref\Key D1
\by \Dri/
\paper Quantum groups
\inbook in Proceedings of the ICM, Berkley, 1986 \ed A\&M\&Gleason
\yr 1987 \pages 798--820 \publ \AMSa/
\endref

\ref\Key D2
\by \Dri/
\paper Quasi-Hopf algebras
\jour \LMJ/ \vol 1 \yr 1990 \pages 1419--1457
\endref

\ref\Key DF
\by J\&Ding and \Fre/
\paper Isomorphism of two realizations of \qaff/
$U_q\bigl(\widehat{\frak{gl}}(n)\bigr)$
\jour \CMP/ \vol 156 \yr 1993 \pages 277--300
\endref

\ref\Key F
\by \Feld/
\paper Conformal field theory and integrable systems associated to elliptic
curves
\inbook in Proceedings of the ICM, Z\"urich, 1994 \publ \Birk/ \yr 1994
\endref

\ref
\by\refin \Feld/
\paper Elliptic quantum groups
\inbook in Proceedings of the ICMP, Paris 1994 \ed D\&Iagolnitzer 
\yr 1995 \publ Intern.\ Press \publaddr Boston \pages 211--218
\endref

\ref\Key FR
\by \Fre/ and \Reshy/
\paper Quantum affine algebras and holonomic \dif/ \eq/s
\jour \CMP/ \vol 146 \yr 1992 \pages 1--60
\endref

\ref\Key FTT
\by \Fadd/, \Takh/ and \VoT/
\paper Local Hamiltonians for integrable models on a lattice
\jour \TMP/ \vol 57 \yr 1983 \pages 1059--1072
\endref

\ref\Key FTV1
\by \Feld/, \VT/ and \Varch/
\paper Solutions of elliptic \qKZBe/s and \Ba/ I
\jour \AMS/ Transl.,\ Ser\&\)2 \vol 180 \yr 1997
\endref

\ref\Key FTV2
\by \Feld/, \VT/ and \Varch/
\paper Monodromy of \sol/s of the elliptic quantum \KZvB/ \deq/s
\jour Preprint \yr 1997 \pages 1--20
\endref

\ref\Key FV
\by \Feld/ and \Varch/
\paper On \rep/s of the elliptic \qg/ $E_{\tau,\eta}(sl_2)$
\jour \CMP/ \vol 181 \yr 1996 \issue 3 \pages 741--761
\endref

\ref\Key G
\by\Gustaf/
\paper A generalization of Selbergs's beta integral
\jour Bull.\ \AMS/ \vol 22 \yr 1990 \pages 97--105
\endref

\ref\Key GR
\by G\&Gasper and M\&Rahman
\book Basic \hgeom/ series
\bookinfo Encycl.\ Math.\ Appl.{}
\yr 1990 \publ \CUPa/
\endref

\ref\Key IK
\by A\&Izergin and \Kor/
\paper The quantum scattering method approach to correlation \fn/s
\jour \CMP/ \vol 94 \yr 1984 \pages 67--92
\endref

\ref\Key J
\by M\&Jimbo
\paper Quantum \Rm/ for the generalized Toda system
\jour \CMP/ \vol 102 \yr 1986 \pages 537--547
\endref

\ref\Key JM
\by M\&Jimbo and T\]\&Miwa
\paper Algebraic analysis of solvable lattice models
\jour CBMS Regional Conf.\ Series in Math.{} \vol 85 \yr 1995
\endref

\ref\Key K
\by T\]\&Kohno
\paper Monodromy \rep/s of braid groups and \YB/s
\jour Ann.\ Inst.\ Fourier \vol 37 \yr 1987 \pages 139--160
\endref

\ref
\by\refin T\]\&Kohno
\paper Linear \rep/s of braid groups and classical \YB/s
\jour Contemp.\ Math.{} \vol 78 \yr 1988 \pages 339--363
\endref

\ref\Key Ka
\by K\&Kadell
\paper A proof of Askey's conjectured $q$-analogue of Selberg's integral and
a conjecture of Morris
\jour \SIAM/ \vol 19 \yr 1988 \pages 969--986
\endref

\goodbm
\ref\Key KL
\by \Kazh/ and \Lusz/
\paper Affine Lie algebras and \qg/s
\jour Intern.\ Math.\ Research Notices \vol 2 \yr 1991 \pages 21--29
\endref

\ref
\by\refin \Kazh/ and \Lusz/
\paper Tensor structures arising from affine Lie algebras, I
\jour J.\ \AMS/ \vol 6 \yr 1993 \pages 905--947
\moreref \paper II
\jour J.\ \AMS/ \vol 6 \yr 1993 \pages 949--1011
\endref

\ref\Key Ko
\by \Kor/
\paper Calculations of norms of Bethe wave \fn/s
\jour \CMP/ \vol 86 \yr 1982 \pages 391--418
\endref

\ref\Key KRS
\by P\]\&P\]\&Kulish, \Reshy/ and \Skl/
\paper \YB/ and \rep/ theory I
\jour \LMP/ \vol 5 \yr 1981 \pages 393--403
\endref

\ref\Key KS
\by \Kazh/ and Ya\&S\&\&Soibelman
\paper Representation theory of \qaff/s
\jour \SMNS/ \vol 3 \yr 1995
\endref

\ref\Key L
\by F\]\&Loeser
\paper Arrangements d'hyperplans et sommes de Gauss
\jour Ann.\ Scient.\ \'Ecole Normale Super.{} \vol {\rm 4-e serie, t\&24}
\yr 1991 \pages 379-400
\endref

\ref\Key Lu
\by S\&Lukyanov
\paper Free field \rep/ for massive integrable models
\jour \CMP/ \vol 167 \yr 1995 \pages 183--226
\endref

\ref\Key M
\by A\&Matsuo
\paper Jackson integrals of Jordan-Pockhammer type and quantum \KZv/ \eq/
\jour \CMP/ \vol 151 \yr 1993 \pages 263--274
\endref

\ref\Key R
\by \Reshy/
\paper Jackson-type integrals, \Bv/s, and \sol/s to a \dif/ analogue
of the \KZv/ system
\jour \LMP/ \vol 26 \yr 1992 \pages 153--165
\endref

\ref\Key RF
\by \Reshy/ and \Fadd/
\paper Hamiltonian structures for integrable field theory models
\jour \TMP/ \vol 56 \yr 1983 \pages 847--862
\endref

\ref\Key RS
\by \Reshy/ and M\&A\&Semenov-Tian-Shansky
\paper Central extensions of the quantum current groups
\jour \LMP/ \vol 19 \yr 1990 \pages 133--142
\endref

\ref\Key S
\by F\]\&A\&Smirnov
\book Form factors in completely integrable models of quantum field theory
\bookinfo Advanced Series in Math.\ Phys., vol\&14
\yr 1992 \publ \WSa/
\endref

\ref\Key Se
\by A\&Selberg
\paper Bemerkninger om et multipelt integral
\jour Norsk.\ Mat.\ Tidsskr.{} \vol 26 \yr 1944 \pages 71--78
\endref

\ref\Key SV1
\by \SchV/ and \Varn/
\paper Hypergeometric \sol/s of \KZv/ \eq/s
\jour \LMP/ \vol 20 \yr 1990 \pages 279--283
\endref

\ref\Key SV2
\by \SchV/ and \Varn/
\paper Arrangements of hyperplanes and Lie algebras homology
\jour \Inv/ \vol 106 \yr 1991 \pages 139--194
\endref

\ref
\by\refin \SchV/ and \Varn/
\paper Quantum groups and homology of local systems
\inbook in Algebraic Geometry and Analytic Geometry, Proceedings of the ICM
Satellite Conference, Tokyo, 1990
\yr 1991 \publ \Spria/ \pages 182--191
\endref

\ref\Key T
\by \VoT/
\paper Irreducible monodromy matrices for the \Rm/ of the {\sl XXZ}-model
and lattice local quantum Hamiltonians
\jour \TMP/ \vol 63 \yr 1985 \pages 440--454
\endref

\ref\Key TV1
\by \VoT/ and \Varn/
\paper Jackson integral \rep/s of \sol/s of the quantized \KZv/ \eq/
\jour \LpMJ/ \vol 6 \issue 2 \yr 1995 \pages 275--313
\endref

\ref\Key TV2
\by \VT/ and \Varch/
\paper Asymptotic \sol/ to the quantized \KZv/ \eq/ and \Bv/s
\jour \AMS/ Transl.,\ Ser\&\)2 \vol 174 \yr 1996 \pages 235--273
\endref

\ref\Key TV3
\by \VT/ and \Varch/
\paper Geometry of $q$-\hgeom/ \fn/s, \qaff/s and elliptic quantum groups
\jour Ast\'erisque \vol 246 \yr 1997 \pages 1--135
\endref

\ifMag\fixedpage\fi

\ref\Key V1
\by \Varn/
\paper The Euler beta-\fn/, the Vandermonde determinant, Legendre's equation,
and critical values of linear \fn/s on a configuration of hyperplanes, I
\jour Math.\ USSR, Izvestia \vol 35 \yr 1990 \pages 543--571
\moreref \paper II
\jour Math.\ USSR Izvestia \vol 36 \yr 1991 \pages 155-168
\endref

\ref
\by\refin \Varn/
\paper Determinant formula for Selberg type integrals
\jour Funct.\ Anal.\ Appl.{} \vol 4 \yr 1991 \pages 65--66
\endref

\ref\Key V2
\by \Varch/
\book Multidimensional \hgeom/ \fn/s and \rep/ theory of Lie algebras and
\qg/s \bookinfo Advanced Series in Math.\ Phys., vol\&21
\yr 1995 \publ \WS/ 
\endref

\ref\Key V3
\by \Varch/
\paper Quantized \KZv/ \eq/s, quantum \YB/, and \deq/s for $q$-\hgeom/ \fn/s
\jour \CMP/ \vol 162 \yr 1994 \pages 499--528
\endref

\ref\Key V4
\by \Varch/
\paper Asymptotic \sol/s to the \KZv/ \eq/ and crystal base
\jour \CMP/ \vol 171 \yr 1995 \pages 99-137
\endref

\ref\Key WW
\by E\&T\]\&Whittaker and G\&N\&Watson
\book A Course of Modern Analysis
\yr 1927 \publ \CUP/
\endref

\endRefs

\ifens
\else
\fixedpage
\nopagenumber

\ifx\Pcde\empty
\ifMag
\contents
\ContM
\endco
\fixedpage
\nopagenumber
\else
\contents
\ContS
\endco
\vsk2>
\fi
\else
\ifMag
\contents
\vbox{\ContMP\vsk-2>}
\endco
\else
\contents
\ContSP
\endco
\fi
\fixedpage
\nopagenumber
\fi
\abstext
\fi

\bye